%% file: Main.tex
\documentclass[12pt,twoside]{book} 
\usepackage{a4wide}
\usepackage{graphicx} 
\usepackage{times} 
\usepackage{latexsym} 
\usepackage{amsmath}
\usepackage{amsfonts}
\usepackage{amssymb}
\usepackage{amsthm} 
\usepackage{amscd} 	
\usepackage{pstricks} 
\usepackage{pst-text} 
\usepackage{pst-node} 
\usepackage{fancyvrb}
\RecustomVerbatimEnvironment%
  {Verbatim}{Verbatim}
  {numbers=left,numbersep=2mm,frame=lines,framerule=0.4pt,%
   firstnumber=last,numberblanklines=false,fontfamily=cmtt,%
   commandchars=\\\{\},framesep=4pt}

 \usepackage{mapleenv} 
\def\date{January 25, 2002}

\newtheorem{thrm}{Theorem}[chapter]
\newtheorem{dfn}[thrm]{Definiton}

\newtheorem{lmm}[thrm]{Lemma}

\def\nn{\nonumber \\}
\def\Cl{\mathcal{C}\kern -0.15em\ell}
\def\cl{\mathcal{C}\kern -0.15em\ell}
\def\JJ{\mathbin{\raisebox{0.25ex}{$\footnotesize
                  \rm\vphantom{I}%
                  \_\hskip -0.25em\_%
                  \vrule width 0.6pt$}}}           

\def\LL{\mathbin{\raisebox{0.25ex}{$\footnotesize
                  \rm\vphantom{I}%
                  \vrule width 0.6pt \hskip -0.5pt%
                  \_\hskip -0.25em\_$}}}          

\def\dwedge{\,\dot{\wedge}\,}

\def\JJBB{\mathop{\JJ \,}\displaylimits_B}

\def\trademark{$\,\bigcirc$\kern -1.75ex {\tiny R} \kern 0.5ex}
\def\openk{\Bbbk}

\def\openP{\mathbb{P}}

\def\openR{\mathbb{R}}
\def\openC{\mathbb{C}}
\def\openZ{\mathbb{Z}}

\def\vj{\u{$\kern -0.7ex\j$}}
\def\vjmath{\u{$\kern -0.8ex\jmath$}}
\def\XP{\phantom{X}}
\def\pprime{{\prime\prime}}
\def\sw{{\rm sw}}

\def\img{{\rm img\,}}
\def\ker{{\rm ker\,}}
\def\dim{{\rm dim\,}}
\def\End{{\rm End\,}}
\def\sign{{\rm sign\,}}
\def\BF{{\cal B\kern -0.5ex F}}

\def\linHom{\text{lin-Hom}\,}

\def\tr{{\rm tr}}
\def\Id{{\rm Id}}
\def\Hom{{\rm Hom}}
\def\Conv{{\rm Conv}}
\def\Sym{{\rm Sym}}
\def\Top{{\bf Top}}
\def\Cog{{\bf Cog}}

\def\pstlw{1.0pt}
\psset{nodesep=4pt}
\begin{document}
\pagestyle{empty}
\mbox{}
\vskip 0.5truecm
\begin{center}
\rule{\textwidth}{3pt}
\vskip 1.5truecm
\huge\bf\sf
A Treatise on\\[2ex] 
Quantum Clifford Algebras
\end{center}
\vspace{1truecm}
\rule{\textwidth}{3pt}
\vspace{2truecm}
\begin{center}
\Large\bf\sf
Habilitationsschrift\\
Dr. Bertfried Fauser
\end{center}
\vspace{3truecm}
\parbox{0.5\textwidth}{%
\hfil
\includegraphics[height=3.5cm]{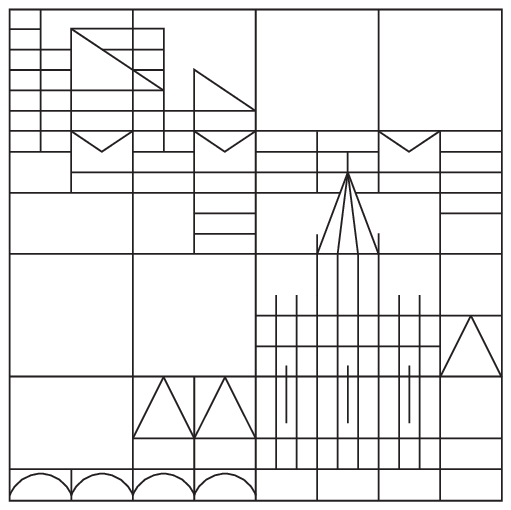}
\hfil}
\parbox{0.5\textwidth}{%
\begin{center}
\large\sf
Universit\"at Konstanz\\
Fachbereich Physik\\
Fach M 678\\
78457 Konstanz
\end{center}
}
\vspace{1truecm}
\begin{center}
\large\bf\sf
\date
\end{center}

\newpage
\mbox{}
\newpage
\mbox{}
\vskip0.3\textheight
\begin{center}
\Large\sf
To Dorothea Ida\\[2ex]
and Rudolf Eugen Fauser
\end{center}
\newpage
\mbox{}
\newpage

\setcounter{page}{1}
\def\thepage{\Roman{page}}
\pagestyle{myheadings}
\newpage
\addcontentsline{toc}{chapter}{Abstract}
\input{abstract.tex}

\newpage
\relax
\newpage

\addcontentsline{toc}{chapter}{Table of Contents}
\tableofcontents
%
\newpage
\newpage
\addcontentsline{toc}{chapter}{Preface}
\input{preface.tex}

%
\newpage
\setcounter{page}{1}
\def\thepage{\arabic{page}}
\input{peano_space.tex}

\input{cl_basics.tex}

\input{graphical_calc.tex}

\input{hopf_al_gebra.tex}

\input{hopf_gebra.tex}

\input{antipode.tex}

\input{gen_cli_prod.tex}

\input{qft.tex}

\begin{appendix}
\input{cliff_big.tex}

\end{appendix}
\newpage
\addcontentsline{toc}{chapter}{Bibliography}
\bibliographystyle{plain}
\bibliography{Main}

\end{document}

%% file: abstract.tex
\mbox{}
\vskip 3truecm
\hfill\parbox[t]{0.85\textwidth}{
{\sf ABSTRACT:} Quantum Clifford Algebras (QCA), i.e. Clifford
Hopf gebras based on bilinear forms of arbitrary symmetry, are treated 
in a broad sense. Five alternative constructions of QCAs are 
exhibited. Grade free Hopf gebraic product formulas are derived 
for meet and join of Gra{\ss}mann-Cayley algebras including co-meet 
and co-join for Gra{\ss}mann-Cayley co-gebras which are very efficient
and may be used in Robotics, left and right contractions, left and right 
co-contractions, Clifford and co-Clifford products, etc. The Chevalley 
deformation, using a Clifford map, arises as a special case.
We discuss Hopf {\it algebra} versus Hopf {\it gebra}, the latter emerging 
naturally from a bi-convolution. Antipode and crossing are consequences 
of the product and co-product structure tensors and not subjectable 
to a choice. A frequently used Kuperberg lemma is revisited 
necessitating the definition of non-local products and interacting 
Hopf gebras which are generically non-perturbative. A `spinorial' 
generalization of the antipode is given. The non-existence of non-trivial
integrals in low-dimensional Clifford co-gebras is shown. Generalized 
cliffordization is discussed which is based on non-exponentially
generated bilinear forms in general resulting in non unital, 
non-associative products. Reasonable assumptions lead to bilinear forms
based on 2-cocycles. Cliffordization is used to derive 
time- and normal-ordered generating functionals for the 
Schwinger-Dyson hierarchies of non-linear spinor field theory and 
spinor electrodynamics. The relation between the vacuum structure, the 
operator ordering, and the Hopf gebraic counit is discussed. QCAs are 
proposed as the natural language for (fermionic) quantum field theory. 
\vskip 2truecm

15-02 Research exposition (monographs, survey articles);\\
15A66 Clifford algebras, spinors;\\
15A75 Exterior algebra, Grassmann algebra;\\
81T15 Perturbative methods of renormalization
}

%% file: preface.tex
``Al-gebra and Co-gebra\\
are brother and sister``
\vskip 0.5truecm
\hskip 2truecm Zbigniew Oziewicz

\vskip 2truecm
{\small\sf\noindent%
Seht Ihr den Mond dort stehen\\
er ist nur halb zu sehen\\
und ist doch rund und sch\"on\\
so sind gar manche Sachen\\
die wir getrost belachen\\
weil unsre Augen sie nicht sehn.
\vskip 0.5truecm
\hskip 2truecm Matthias Claudius
}
\vskip 2truecm

\section*{Preface}

This `Habilitationsschrift' is the second incarnation of itself -- and still 
in a {\it status nascendi.} The original text was planned to contain Clifford
{\it algebras\/} of an arbitrary bilinear form, now called Quantum Clifford 
Algebras (QCA) and their beautiful application to quantum field theory 
(QFT). However, while proceeding this way, a major change in paradigm took 
place after the 5th Clifford conference held in Ixtapa 1999. As a
consequence the first incarnation of this work faded away without reaching
a properly typeset form, already in late 2000.

What had happened? During the 5th Clifford conference at Ixtapa a special 
session dedicated to Gian-Carlo Rota, who was assumed to attend the 
conference but died in Spring 1999, took place. Among other impressive
retrospectives delivered during this occasion about Rota and his work, 
Zbigniew Oziewicz explained the Rota-Stein cliffordization process and coined 
the term `Rota-sausage' for the corresponding tangle -- for obvious reason 
as you will see in the main text. This approach to the Clifford product 
turned out to be superior to all other previously achieved approaches in 
elegance, efficiency, naturalness and beauty -- for a discussion of 
`beautiness' in mathematics, see \cite{rota:1997a}, Chap. X, 
`{\it The Phenomenology of Mathematical Beauty\/}'.
So I had decided to revise the whole writing. During 2000, beside being 
very busy with editing \cite{ablamowicz:fauser:2000a}, it turned out, that 
not only a rewriting was necessary, but that taking a new starting point 
changes the whole tale!

A major help in entering the Hopf gebra business for Gra{\ss}mann 
and Clifford algebras and cliffordization was the CLIFFORD package 
\cite{CLIFFORD} developed by Rafa{\l} Ab{\l}amowicz. During a 
collaboration with him which took place in Konstanz in Summer 1999, major 
problems had been solved which led to the formation of the BIGEBRA 
package \cite{BIGEBRA} in December 1999. The package proved to be 
calculationable stable and useful for the first time in Autumn 2000 during
a joint work with Zbigniew Oziewicz, where many involved computations were 
successfully performed. The requirements of this lengthy computations 
completed the BIGEBRA package more or less. Its final form was produced 
jointly with Rafa{\l} Ab{\l}amowicz in Cookeville, September 2001.

The possibility of automated calculations and the knowledge of functional 
quantum field theory \cite{stumpf:borne:1994a,borne:lochak:stumpf:2001a} 
allowed to produce a first important result. The relation between time- and 
normal-ordered operator products and correlation functions was revealed to 
be a special kind of cliffordization which introduces an antisymmetric 
(symmetric for bosons) part in the bilinear form of the Clifford product 
\cite{fauser:2001b}. For short, QCAs deal with time-ordered monomials 
while regular Clifford algebras of a symmetric bilinear form deal with 
normal-ordered monomials.

It seemed to be an easy task to translate with benefits all of the work 
described in 
\cite{stumpf:fauser:pfister:1993a,fauser:1996b,fauser:stumpf:1997a,%
fauser:1998a,fauser:2000e,fauser:2001e} into the hopfish framework. But 
examining Ref. \cite{fauser:2001e} it showed up that the standard 
literature on Hopf algebras is set up in a too narrow manner so that some 
concepts had to be generalized first.

Much worse, Oziewicz showed that given an invertible scalar product $B$
the Clifford bi-convolution $\Cl(B,B^{-1})$, where the Clifford co-product 
depends on the co-scalar product $B^{-1}$, has {\it no antipode\/} and is 
therefore not a Hopf algebra at all. But the antipode played {\it the
central role\/} in Connes-Kreimer renormalization theory
\cite{kreimer:2000a,connes:kreimer:1999a,connes:kreimer:2000a,%
connes:kreimer:2001a}. Furthermore the topological meaning and the 
group-like structure are tied to Hopf algebras, not to convolution 
semigroups. This motivated Oziewicz to introduce a second 
{\it independent\/} bilinear form, the co-scalar product $C$ in the
Clifford bi-convolution $\Cl(B,C)$, $C\not= B^{-1}$ which is antipodal 
and therefore Hopf. A different solution was obtained jointly in 
\cite{fauser:oziewicz:2001a}.

Meanwhile QCAs made their way into differential geometry and showed up 
to be useful in Einstein-Cartan-K\"ahler theory with teleparallel 
connections developed by J. Vargas, see \cite{vargas:torr:2002a} and 
references therein. It was clear for some time that also differential forms,
the Cauchy-Riemann differential equations and cohomology have to be revisited
in this formalism. This belongs not to our main theme and will be 
published elsewhere \cite{fauser:dehnen:2002a}.

Another source supplied ideas -- geometry and robotics! -- the geometry of 
a Gra{\ss}mann-Cayley algebra, i.e. projective geometry is by 
the way the first {\it application\/} of Gra{\ss}mann's work by himself 
\cite{grassmann:1878a}. Nowadays these topics can be considered in their 
relation to Gra{\ss}mann Hopf gebras. The crucial `regressive product' 
of Gra{\ss}mann can easily be defined, again following Rota et al. 
\cite{doubilet:rota:stein:1974a,rota:stein:1976a,kung:rota:1984a,%
barnabei:brini:rota:1985a}, by Hopf algebra methods. A different route
also following Gra{\ss}mann's first attempt is discussed in Browne 
\cite{browne:2001a}. Rota et al., however, used a Peano space, a pair 
of a linear space $V$ and a volume to come up with invariant theoretic
methods. It turns out, and is in fact implemented in BIGEBRA this way 
\cite{ablamowicz:fauser:2002a,ablamowicz:fauser:2002b}, that meet and 
join operations of projective geometry are encoded most efficiently 
and mathematically sound using Gra{\ss}mann Hopf gebra. 
Gra{\ss}mannians, flag manifolds which are important in string theory, 
M-theory, robotics and various other objects from algebraic geometry 
can be reached in this framework with great formal and computational ease.

It turned out to be extremely useful to have geometrical ideas at hand 
which can be transformed into the QF theoretical framework. As a general 
rule, it is true that sane geometric concepts translate into sane concepts 
of QFT. However a complete treatment of the geometric background would 
have brought us too far off the road. Examples of such geometries would 
be M\"obius geometry, Laguerre geometry, projective and incidence 
geometries, Hijelmslev planes and groups etc. 
\cite{hilbert:1968a,benz:1973a,bachmann:1973a,bachmann:1989a,zaddach:1994a}.
I decided to come up with the algebraic part of
Peano space, Gra{\ss}mann-Cayley algebra, meet and join to have them 
available for later usage. Nevertheless, it will be possible for the 
interested reader to figure out to a large extend which geometric 
operations are behind many QF theoretical operations.  

In writing a treatise on QCAs, I assume that the reader is familiar with 
basic facts about Gra{\ss}mann and Clifford algebras. Reasonable 
introductions can be found in various text books, e.g.
\cite{riesz:1957a,porteous:1969a,benn:tucker:1987a,bourbaki:1989a,%
budinich:trautmann:1988a,crumeyrolle:1990a,lounesto:2001a}. A good source
is also provided by the conference volumes of the five international
Clifford conferences
\cite{canterbury:1985a,montpellier:1989a,deinze:1993a,aachen:1996a,%
ablamowicz:fauser:2000b,ryan:sproessig:2000a}. Nevertheless, the 
terminology needed later on is provided in the text. 

In this treatise we make to a large extend use of graphical calculi. 
These methods turn out to be efficient, inspiring and allow to memorize 
particular equations in an elegant way, e.g. the `Rota-sausage' of 
cliffordization which is explained in the text. Complicated calculations 
can be turned into easy manipulations of graphs. This is one key point 
which is already well established, another issue is to explore the 
topological and other properties of the involved graphs. This would 
lead us to graph theory itself, combinatorial topology, but also to the 
exciting topic of matroid theory. However, we have avoided graph theory,
topology and matroids in this work.

Mathematics provides several graphical calculi. We have decided to use 
three flavours of them. I: Kuperberg's translation of tensor algebra 
using a self-created very intuitive method because we require some of his 
important results. Many current papers are based on a couple of lemmas 
proved in his writings. II. Commutative diagrams constitute a sort of 
{\it lingua franca} in mathematics. III. Tangle diagrams turn out to be 
dual to commutative diagrams in a particular sense. From a physicist's 
point of view they constitute a much more natural way to display 
dynamical `processes'.

Of course, graphical calculi are present in physics too, especially in 
QFT and for the tensor or spinor algebra, e.g. \cite{penrose:rindler:1984a} 
appendix. The well known Feynman graphs are a particular case of a 
successful graphical calculus in QFT. Connes-Kreimer renormalization 
attacks QFT via this route. Following Cayley, rooted trees are taken to 
encode the complexity of differentiation which leads via the Butcher 
B-series \cite{butcher:1972a,butcher:1972b} and a `decoration' technique 
to the Zimmermann forest formulas of BPHZ 
(Bogoliubov-Parasiuk-Hepp-Zimmermann) renormalization in momentum space.

Our work makes contact to QFT on a different and very solid way not using
the mathematically peculiar path integral, but functional differential 
equations of functional quantum field theory, a method developed by Stumpf 
and coll. \cite{stumpf:borne:1994a,borne:lochak:stumpf:2001a}. This approach
takes its starting point in position space and proceeds by implementing 
an algebraic framework inspired by and closely related to $C^*$-algebraic 
methods without assuming positivity.

However, this method was not widely used in spite of reasonable and unique 
achievements, most likely due to its lengthy and cumbersome calculations.
When I became aware of Clifford algebras in 1993, as promoted by 
D. Hestenes \cite{hestenes:1966a,hestenes:sobczyk:1992a} for some decades 
now, it turns out that this algebraic structure is a key step to compactify 
notation and calculations of functional QFT \cite{fauser:1996c}. In the 
same time many {\it ad hoc\/} arguments have been turned into a 
mathematical sound formulation, see e.g. 
\cite{fauser:1996c,fauser:1996b,fauser:stumpf:1997a,fauser:1998a}. 
But renormalization was still not in the game, mostly since in Stumpf's 
group in T\"ubingen the main interest was laid on non-linear spinor 
field theory which has to be regularized since it is non-renormalizable.

While I was finishing this treatise Christian Brouder came up in 
January 2002 with an idea how to employ cliffordization in 
renormalization theory. He used the same transition as was employed in 
\cite{fauser:2001b} to pass from normal- to time-ordered operator products 
and correlation functions but implemented an additional bilinear form which
introduces the renormalization parameters into the theory but remains in 
the framework of cliffordization. This is the last part of a puzzle which 
is needed to formulate all of the algebraic aspects of (perturbative) QFT 
entirely using the cliffordization technique and therefore in the framework 
of a Clifford Hopf gebra (Brouder's term is `quantum field algebra', 
\cite{brouder:2002a}). This event caused a prolongation by a chapter on 
generalized cliffordization in the mathematical part in favour of some 
QFT which was removed and has to be rewritten along entirely hopfish lines. 
It does not make any sense to go with the {\it algebra only} description 
any longer. As a consequence, the discussion of QFT under the topic 
`QFT as Clifford Hopf gebra' will be a sort of second volume to this work. 
Nevertheless, we give a complete synopsis of QFT in terms of QCAs, i.e. 
in terms of Clifford Hopf gebras. Many results can, however, be found in 
a pre-Hopf status in our publications.
\bigskip

What is the content and what are the {\it main results}?
\begin{itemize}
\item The Peano space and the Gra{\ss}mann-Cayley algebra, also called
bracket algebra, are treated in its classical form as also in the Hopf 
algebraic context.
\item The bracket of invariant theory is related to a Hopf gebraic
integral.
\item Five methods are exhibited to construct (quantum) Clifford algebras, 
showing the outstanding beautiness of the Hopf gebraic method of 
cliffordization.
\item We give a detailed account on Quantum Clifford Algebras (QCA) based 
on an arbitrary bilinear form $B$ having no particular symmetry.
\item We compare Hopf {\it algebras} and Hopf {\it gebras}, the latter
providing a much more plain development of the theory.
\item Following Oziewicz, we present Hopf gebra theory. The crossing 
and the antipode are exhibited as dependent structures which have to be 
calculated from structure tensors of the product and co-product of
a bi-convolution and cannot be subjected to a choice.
\item We use Hopf algebraic methods to derive the basic formulas of
Clifford algebra theory (classical and QCA). One of them will be 
called {\it Pieri-formula of Clifford algebra\/}.
\item We discuss the Rota-Stein cliffordization and co-cliffordization,
which will be called, stressing an analogy, the
{\it Littlewood-Richardson rule of Clifford algebra\/}.
\item We derive {\it grade free\/} and very efficient product formulas 
for almost all products of Clifford and Gra{\ss}mann-Cayley algebras, e.g. 
Clifford product, Clifford co-product (time- and normal-ordered operator 
products and correlation functions based on dotted and undotted exterior 
wedge products), meet and join products, co-meet and co-join, 
left and right contraction by arbitrary elements, left and right 
co-contractions, etc.
\item We introduce non-interacting and interacting Hopf gebras which 
cures a drawback in an important lemma of Kuperberg which is frequently 
used in the theory of integrable systems, knots and even QFT as proposed 
by Witten. Their setting turns thereby out to be close to free theories.
\item We show in low dimensional examples that no non-trivial integrals 
do exist in Clifford co-gebras and conjecture this to be generally true.
\item A `spinorial' antipode, a convolutive unipotent, is given which 
symmetrizes the Kuperberg ladder.
\item We extend cliffordization to bilinear forms $\BF$ which are 
{\it not\/} derivable from the exponentiation of a bilinear form on the 
generating space $B$.
\item We discuss generalized cliffordization based on non-exponentially
generated bilinear forms. Assertions on the derived product show that
exponentially generated bilinear forms are related to 2-cocycles.  
\item An overview is presented on functional QFT. Generating functionals 
are derived for time- and normal-ordered non-linear spinor field theory 
and spinor electrodynamics.
\item A detailed account on the role of the counit as a `vacuum' state
is described. Two models with $U(1)$ and $U(2)$ symmetry are taken 
as examples.
\item It is shown how the quantization enters the cliffordization. 
Furthermore we explain in which way the vacuum is determined by the
propagator of the theory.
\item Quantum Clifford algebras are proposed as the algebras of QFT.
\end{itemize}

What is {\it not to be found} in this treatise? It was not intended to 
develop Clifford algebra theory from scratch, but to concentrate on the 
`quantum' part of this structure including the unavoidable hopfish methods. 
$q$-deformation, while possible and most likely natural in our framework 
is not explicitely addressed. However the reader should consult
our results presented in Refs. 
\cite{fauser:1999b,fauser:2000e,ablamowicz:fauser:2000b,%
fauser:2000d} where this topic is addressed. A detailed explanation why 
`quantum' has been used as prefix in QCA can be found in 
\cite{fauser:ablamowicz:2000c}. Geometry is reduced to algebra, which 
is a pity. A broader treatment, e.g. Clifford algebras over finite 
fields, higher geometries, incidence geometries, Hjielmslev planes etc. 
was not fitting coherently into this work and would have fatten it 
becoming thereby unhandsome. An algebro-synthetic approach to geometry 
would also constitute another volume which would be worth to be written.
This is not a work in mathematics, especially not a 
sort of `Bourbaki chapter' where a mathematical field is developed
straightforward to its highest extend providing all relevant definitions
and proving all important theorems. We had to concentrate on hot spots
for lack of time and space and to come to a status where the method
can be applied and prove its value. The symmetric group algebra and 
its deformation, the Hecke algebra, had to be postponed, as also a
discussion of Young tableaux and their relation to Specht modules
and Schubert varieties. And many more exciting topics \ldots 
\vskip 2truecm

\noindent
\addcontentsline{toc}{chapter}{Acknowledgement}
{\bf Acknowledgement:}
This work was created under the enjoyable support of many persons. 
I would like to thank a few of them personally, especially Prof. Stumpf 
for his outstanding way to teach and practise physics, Prof. Dehnen for the
patience with my hopfish exaggerations and his profound comments during 
discussions and seminars, Prof. Rafa{\l} Ab{\l}amowicz for helping me 
since 1996 with CLIFFORD, inviting me to be a co-author of this package
and most important becoming a friend in this turn. Prof. Zbigniew 
Oziewicz grew up most of my understanding about Hopf gebras.
Many thanks also to the theory groups in T\"ubingen and Konstanz which 
provided a inspiring working environment and took a heavy load of 
`discussion pressure'. Dr. Eva Ge{\ss}ner and Rafa{\l} Ab{\l}amowicz
helped with proof reading, however, the author is responsible for all
remaining errors.\\ 
My gratitude goes to my wife Mechthild for her support, to my children 
simply for being there, and especially to my parents to whom this work 
is dedicated.
\vskip 1truecm
\begin{flushright}
Konstanz, \date \\
Bertfried Fauser
\end{flushright}
\vfill
{\small\sf\noindent%
Wir armen Menschenkinder\\
sind eitel arme S\"under\\
und wissen garnicht viel\\
wir spinnen Luftgespinste\\
und suchen viele K\"unste\\
und kommen weiter von dem Ziel!
\vskip 0.5truecm
\hskip 2truecm Matthias Claudius
}

%% file: peano_space.tex
\chapter{Peano Space and Gra{\ss}mann-Cayley Algebra}

In this section we will turn our attention to the various possibilities 
which arise if additional structures are added to a linear space 
($\openk$-module or $\openk$-vector space). It will turn out that a second 
structure, such as a norm, a scalar product or a bracket lead to seemingly 
very different algebraic settings. To provide an overview, we review shortly 
normed spaces, Hilbert spaces, Weyl or symplectic spaces and concentrate on 
Peano or volume spaces which will guide us to projective geometry and the 
theory of determinants.

Let $\openk$ be a ring or a field. The elements of $\openk$ will be called
scalars, following Hamilton. Let $V$ be a linear space over $\openk$ having an
additively written group acting on it and a scalar multiplication. The elements
of $V$ are called vectors. Hamilton had a `vehend' also and his vectors
were subjected to a product and had thus an operative meaning, see e.g.
\cite{crowe:1994a}. We will also be interested mainly in the algebraic 
structure, but it is mathematical standard to disentangle the space 
underlying a `product' from the product structure. Scalar multiplication 
introduces `weights' on vectors sometimes also called `intensities'. As we 
will see later, the Gra{\ss}mann-Cayley algebra does not really need scalars 
and is strictly speaking not an algebra in the common sense. We agree that 
an algebra $A$ is a pair $A=(V,m)$ of a $\openk$-linear space $V$ and a 
product map $m : V\times V \rightarrow V$. Algebras are introduced more 
formally later. Products are mostly written in an {\it infix form\/}: 
$\;a \;m\; b \equiv m(a,b)$. Products are {\it defined\/} by Gra{\ss}mann 
\cite{grassmann:1878a} as those mappings which respect distributivity 
w.r.t. addition, $a,b,c \in V$:
\begin{align}
a \,m\, (b+c) &= a\,m\,b + a\,m\,c \nn
(a + b) \,m\, c &= a\,m\,c + b\,m\,c 
\end{align}
Hence the product is bilinear. Gra{\ss}mann does not assume associativity, 
which allows to drop parentheses
\begin{align}
 a \,m\, (b \,m\, c) &= (a\,m\,b) \,m\,c.
\end{align}
Usually the term algebra is used for `associative algebra' while 
`non-associative algebra' is used for the general case. We will be mostly
interested in associative algebras. 

\section{Normed space -- normed algebra}

Given only a linear space we own very few rules to manipulate its elements.
Usually one is interested in a reasonable extension, e.g. by a distance 
or length function acting on elements from $V$. In analytical applications 
it is very convenient to have a positive valued length function. A 
reasonable such structure is a norm $\|.\| : V \rightarrow \openk$, 
a linear map, defined as follows
\begin{align}
o)  &&\|\alpha a\| = \alpha\,\|a\| &&&\alpha\in \openk,\quad a\in V \nn
i)  &&\|a\| =0 &&&\text{if and only if $a\equiv 0$} \nn
ii) &&\|a\| \ge 0 &&&\forall a\in V \quad \text{positivity} \nn
iii)&&\|a+b\| \le \|a\| + \|b\| &&&\text{triangle relation} .
\end{align}
As we will see later this setting is to narrow for our purpose. Since 
it is a strong condition it implies lots of structure. Given an algebra 
$A=(V,m)$ over the linear space $V$, we can consider a {\it normed algebra} 
if $V$ is equipped additionally with a norm which fulfils
\begin{align}
\|ab\| &\le \|a\|\,\|b\|
\end{align} 
which is called {\it submultiplicativity}. Normed algebras provide a wealthy 
and well studied class of algebras \cite{gelfand:1964a}.

However, one can prove that on a finite dimensional vector space all norms
are equivalent. Hence we can deal with the prototype of a norm, the Euclidean 
length 
\begin{align}
\|x\|_2 &:=+\sqrt{\sum (x^i)^2}
\end{align}
where the $x^i \in \openk$ are the coefficients of $x\in V$ w.r.t. an
orthogonal generating set $\{e_i\}$ of $V$. We would need here the dual 
space $V^*$ of linear forms on $V$ for a proper description. From any norm
we can derive an {\it inner product} by polarization. We assume here that
$\openk$ has only trivial involutive automorphisms, otherwise the 
polarization is more complicated
\begin{align}
g(x,y) &: V \times V \rightarrow \openk \nn
g(x,y) &:= \| x-y\| .
\end{align}
A `distance' function also implies some kind of interpretation to the vectors
as `locations' in some space. 

Since the major part of the work will deal with algebras over finite vector
spaces or with formal power series of generating elements, i.e. without
a suitable topology, thus dropping convergence problems, we are not 
interested in normed algebras. The major playground for such a structure 
is over infinitely generated linear spaces of countable or continuous 
dimension. Banach and $C^*$-algebras are e.g. of such a type. The later is 
distinguished by a $C^*$-condition which provides a unique norm, 
the $C^*$-norm. These algebras are widely used in non-relativistic QFT and 
statistical physics, e.g. in integrable models, BCS superconductivity 
etc., see 
\cite{bratelli:robinson:1979a,bratelli:robinson:1981a,murphy:1990a}.

\section{Hilbert space, quadratic space -- classical Clifford algebra}

A slightly more general concept is to concentrate in the first place
on an inner product. Let 
\begin{align}
<.\mid.>  &\;:\; V\times V \rightarrow \openk \nn
<x\mid y> &\;=\; <y\mid x>
\end{align}
be a symmetric bilinear inner product. An inner product is called
positive semi definite if 
\begin{align}
<x\mid x> &\;\ge\; 0 
\end{align}
and positive definite if in the above equation equality holds if and only if
$x\equiv 0$. The pair of a finite or infinite linear space $V$ equipped with 
such a bilinear positive definite inner product $<.\mid .>$ is called a 
Hilbert space ${\cal H}=(V,<.\mid .>)$, if this space is closed in the
natural topology induced by the inner product. Hilbert spaces play a 
prominent role in the theory of integral equations, where they have been 
introduced by Hilbert, and in quantum mechanics. The statistical 
interpretation of quantum mechanics is directly connected to positivity. 
Representation theory of operator algebras benefits from positivity too, e.g. 
the important GNS construction \cite{murphy:1990a}. Of course one can add 
a multiplication to gain an algebra structure. This is a special case of 
a further generalization to quadratic spaces which we will consider now.

Let $Q$ be a quadratic form on $V$ defined as
\begin{align}
Q &\;:\; V \rightarrow \openk \nn
Q(\alpha x) &\;=\; \alpha^2 Q(x) &&\alpha\in \openk,\quad x\in V \nn
2\,B_p(x,y) &\;:=\; Q(x-y)-Q(x)-Q(y) &&\text{where $B_p$ is bilinear} .
\end{align}
The symmetric bilinear form $B_p$ is called polar bilinear form, the name
stems from the pol-polar relation of projective geometry, where the locus
of elements $x \in V$ satisfying $B_p(x,x)=0$ is called {\it quadric}. 
However, one should be careful and introduce dual spaces for the 
`polar elements', i.e. hyperplanes. It is clear that we have to assume that 
the characteristic of $\openk$ is not equal to 2.

We can ask what kind of algebras arise from adding this structure to 
and algebra having a product $m$. Such a structure $A=(V,m,Q)$ would e.g.
be an operator algebra where we have employed a non-canonical quantization,
as e.g. the Gupta-Bleuler quantization of electrodynamics.

However, it is more convenient to ask if the quadratic form can {\it imply} 
a product on $V$. In this case the product map $m$ is a consequence of the 
quadratic form $Q$ itself. As we will see later, classical Clifford algebras 
are of this type. From its construction, based on a quadratic form $Q$ having
a symmetric polar bilinear form $B_p$, it is clear that we can expect 
Clifford algebras to be related to orthogonal groups. Classical Clifford
algebras should thus be interpreted as a linearization of a quadratic
form. It was Dirac who used exactly this approach to postulate his 
equation. Furthermore, we can learn from the polarization process
that this type of algebra is related to anticommutation relations:
\begin{align}
Q(x) &= \sum_i x^i x^j e_i e_j \nn
2\,B_p(x,y) &= \sum_{i,j} x^i y^j (e_i e_j + e_j e_i)
\end{align}
which leads necessarily to
\begin{align}
e_i \,e_j + e_j\,e_i &= 2\,B_{p\,ij}.
\end{align}
Anticommutative such algebras are usually called (canonical) anticommutation 
algebras CAR and are related to fermions.

Classical Clifford algebras are naturally connected with the classical 
orthogonal groups and their double coverings, the {\bf pin} and {\bf spin}
groups, \cite{porteous:1969a,porteous:2000a,lounesto:2001a}.

Having generators $\{e_i\}$ linearly spanning $V$ it is necessary to pass 
over to the linear space $W=\bigwedge V$ which is the linear span of all 
linearly and algebraically independent products of the generators. 
Algebraically independent are such products of the $e_i$s which cannot 
be transformed into one another by using the (anti)commutation relations, 
which will be discussed later.

In the special case where the bilinear form on $W$, induced by this
construction, is positive definite we deal with a Hilbert space. That is, 
Clifford algebras with positive (or negative) definite bilinear forms 
{\it on the whole space $W$} are in fact $C^*$-algebras too, 
however of a special flavour.

\section{Weyl space -- symplectic Clifford algebras (Weyl algebras)}

While we have assumed symmetry in the previous section, it is equally
reasonable and possible to consider antisymmetric bilinear forms
\begin{align}
<.\mid.>  &: V\times V \rightarrow \openk \nn
<x\mid y> &= -<y\mid x> .
\end{align}
A linear space equipped with an antisymmetric bilinear inner product 
will be called Weyl space. The antisymmetry implies directly that
{\it all\/} vectors are null -- or synonymously isotrop: 
\begin{align}
<x\mid x> &= 0 &&\forall x\in V.
\end{align}
It is possible to define an algebra $A=(V,m,<. \mid .>)$, but once more 
we are interested in such products which are derived from the bilinear 
form. Using again the technique of polarization, one arrives this
time at a (canonical) commutator relation algebra CCR
\begin{align}
e_i \,e_j - e_j\,e_i &= 2\,A_{ij},
\end{align}
where $A_{ij}=-A_{ji}$. It should however be remarked, that this 
{\it symplectic Clifford algebras\/} are not related to classical 
groups in a such direct manner as the orthogonal Clifford algebras. 
The point is, that symplectic Clifford algebras do {\it not\/}
integrate to a group action if built over a field 
\cite{crumeyrolle:1990a,bourbaki:1989a}. In fact one awaits nevertheless
to deal with a sort of double cover of symplectic groups.

Such algebras are tied to bosons and occur frequently in quantum 
physics. Indeed, quantum physics was introduced for bosonic fields first
and studied these much more complicated algebras in the first place.

In literature one finds also the name Weyl algebra for this type of 
structure.

There is an odd relation between the scalars and the symmetry of the 
generators -- operators in quantum mechanics and quantum field theory. 
While for fermions the coefficients are commutative scalars forming a 
field and the generators are anticommutative we find in the case of bosons 
complicated scalars, at least a formal polynomial ring, or non-commutative 
coordinates. In combinatorics it is well known that such a vice-versa 
relation between coefficients and generators holds, see 
\cite{grosshans:rota:stein:1987a}.

Also looking at combinatorial aspects, symplectic Clifford algebras are 
much more complicated. This stems from two facts. One is that
one has to deal with multisets. The second is that the induced
bilinear forms on the space $W$ algebraically generated from $V$ have 
in the antisymmetric case the structure of minors and determinants which 
are related to Pfaffians and obey decomposition, while in the symmetric 
case one ends up with permanents and Hafnians. The combinatorics of 
permanents is much more complicated.

It was already noted by Caianiello \cite{caianiello:1973a} that such
structures are closely related to QFT calculations. We will however 
see below that his approach was not sufficient since he did not 
respect the symmetry of the operator product.

\section{Peano space -- Gra{\ss}mann-Cayley algebras}

In this section we recall the notion of a Peano space, as defined by Rota 
et al. \cite{doubilet:rota:stein:1974a,barnabei:brini:rota:1985a}, because
it provides the `classical' part of QFT as a good starting point. Furthermore
this notion is not well received. (In the older ref. 
\cite{doubilet:rota:stein:1974a} the term Cayley space 
was used). Peano space goes back to Giuseppe Peano's {\it Calcolo Geometrico}
\cite{peano:1888a}. In this important work, Peano managed to surmount the 
difficulties of Gra{\ss}mann's regressive product by setting up axioms in 
3-dimensional space. In later works this goes under the name of the 
{\it Regel des doppelten Faktors} [rule of the (double) common factor], 
see the discussion in \cite{browne:2001a} where this is taken as
an axiom to develop the regressive product. Gra{\ss}mann himself changed 
the way how he introduced the regressive product from the first {\bf A1}
({\bf Ai} is common for the $i$-th `lineale Ausdehnungslehre' [theory
 of extensions] from 1844 (A1) \cite{grassmann:1878a} and 1862 (A2) 
\cite{grassmann:1862a}) to the presentation in the {\bf A2} . Our goal is 
to derive the wealth of products accompanying the Gra{\ss}mann-Cayley algebra
of meet and join, emerging from a `bracket', which will later on be recast 
in Hopf algebraic terms. The bracket will show up as a Hopf algebraic 
integral of the exterior wedge products of its entries, see chapters below.
The Gra{\ss}mann-Cayley algebra is denoted bracket algebra in invariant 
theory.

\subsection{The bracket}

While we follow Rota et al. in their mathematical treatment, we 
separate explicitely from the comments about co-vectors and Hopf algebras
in their writing in the above cited references. It is less known that
also Rota changed his mind later. Unfortunately many scientists based their 
criticism of co-vectors or Hopf algebras on the above well received papers 
while the later change in the position of Rota was not appreciated, see 
\cite{grosshans:rota:stein:1987a,rota:stein:1994a} and many other joint
papers of Rota in the 90ies.

Let $V$ be a linear space of finite dimension $n$. Let lower case $x_i$ 
denote elements of $V$, which we will call also letters. We define a bracket 
as an alternating multilinear scalar valued function 
\begin{align}
{}[.,\ldots,.] &\;:\; V \times \ldots \times V \rightarrow \openk 
&&\text{$n$-factors} \nn
{}[x_1,\ldots,x_n] &= \sign(p)[x_{p(1)},\ldots,x_{p(n)}] \nn
{}[x_1,\ldots,\alpha x_r + \beta y_r,\ldots, x_n] &\;=\; 
\alpha [x_1,\ldots,x_r,\ldots,x_n] + 
\beta  [x_1,\ldots,y_r,\ldots,x_n] .
\end{align}
The sign is due to the permutation $p$ on the arguments of the bracket. 
The pair ${\cal P}=(V,[.,\ldots,.])$ is called a Peano space.

Of course, this structure is much weaker as e.g. a normed space or an
inner product space. It does not allow to introduce the concept of 
length, distance or angle. Therefore it is clear that a geometry based 
on this structure cannot be metric. However, the bracket can be addressed 
as a volume form. Volume measurements are used e.g. in the analysis of 
chaotic systems and strange attractors.

A {\it standard Peano space\/} is a Peano space over the linear space $V$ of
dimension $n$ whose bracket has the additional property that for every 
vector $x\in V$ there exist vectors $x_2,\ldots,x_n$ such that
\begin{align}
[x,x_2,\ldots,x_n] &\not= 0 .
\end{align}
In such a space the length of the bracket, i.e. the number of entries, 
equals the dimension of the space, and conversely. We will be concerned 
here with standard Peano spaces only. 

The notion of a bracket is able to encode linear independence. Let 
$x,y$ be elements of $V$ they are linearly independent if and only if
one is able to find $n-2$ vectors $x_3,...,x_n$ such that the bracket
\begin{align}
[x,y,x_3,\ldots,x_n] &\not= 0 .
\end{align}

A {\it basis} of $V$ is a set of $n$ vectors which have a non-vanishing
bracket. We call a basis {\it unimodular\/} or linearly ordered and 
normalized if for the ordered set $\{e_1,\ldots,e_n\}$, also called 
sequence in the following, we find the bracket
\begin{align}
[e_1,\ldots,e_n] &= 1 .
\end{align}
At this place we should note that an alternating linear form of rank
$n$ on a linear space of dimension $n$ is uniquely defined up to a 
constant. This constant is however important and has to be removed
for a fruitful usage, e.g. in projective geometry. This is done by 
introducing cross ratios. The group which maps two linearly ordered 
bases onto another is $gl_n$ and $sl_n$ for the mapping of unimodular
bases.

\subsection{The wedge product -- join}

To pass from a space to an algebra we need a product. For this reason we
introduce equivalence classes of ordered sequences of vectors using the 
bracket. We call two such sequences equivalent
\begin{align}
a_1,\ldots,a_k \cong b_1,\ldots,b_k
\end{align}
if for every choice of vectors $x_{k+1},\ldots,x_n$ the following 
equation holds
\begin{align}
[a_1,\ldots,a_k,x_{k+1},\ldots,x_n] &= [b_1,\ldots,b_k,x_{k+1},\ldots,x_n].
\end{align}
An equivalence class of this type will be called {\it extensor\/} or
{\it decomposable antisymmetric tensor\/} or {\it decomposable $k$-vector\/}.
The projection of the Cartesian product $\times$ (or the tensor product 
$\otimes$ if the $\openk$-linear structure is considered) under this 
equivalence class is called {\it exterior wedge product of points\/} or 
simply wedge product if the context is clear. Alternatively we use
the term {\it join\/} if geometrical applications are intended. 
In terms of formulas we find
\begin{align}
a\wedge b &:= \{a \;,\; b\} \mod\cong 
\end{align}
for the equivalence classes. The wedge product inherits antisymmetry
from the alternating bracket and associativity, since the bracket was 
`flat' (not using parentheses). Rota et al. write for the join the 
vee-product $\vee$ to stress the analogy to Boolean algebra, a connection 
which will become clear later. 
However, we will see that this identification is a matter of 
taste due to duality. For this reason we will stay with a wedge $\wedge$ 
for the `exterior wedge product of points'. Furthermore we will see later 
in this work that it is convenient to deal with {\it different exterior 
products\/} and to specify them in a particular context. In the course of
this work we even have occasion to use various exterior products at the 
same time which makes a distinction between them necessary. One finds 
$2^n$ linearly independent extensors. They span the linear space $W$ which
is denoted also as $\bigwedge V$. This space forms an algebra w.r.t. the 
wedge product, the {\it exterior algebra\/} or Gra{\ss}mann algebra. The 
exterior algebra is a graded algebra in the sense that the module 
$W=\bigwedge V$ is graded, i.e. decomposable into a direct sum of subspaces 
of words of the same length and the product respects this direct sum 
decomposition: 
\begin{align}
\wedge &\;:\; \bigwedge^r V\times \bigwedge^s V \rightarrow \bigwedge^{r+s} V.
\end{align}

The extensors of step $n$ form a one dimensional subspace. Gra{\ss}mann
tried to identify this space also with the scalars which is not convenient
\cite{zaddach:1994a}. Using an unimodular basis we can construct the element
\begin{align}
E &= e_1\wedge \ldots \wedge e_n 
\end{align} 
which is called {\it integral\/}, see \cite{sweedler:1969a}. Physicists
traditionally chose $\gamma^5$ for this element.

We allow extensors to be inserted into a bracket according to the 
following rule
\begin{align}
A &= a_1,\ldots,a_r, \quad B=b_1,\ldots,b_s,\quad C=c_1,\ldots,c_t \nn
{}[A,B,C] &=[a_1,\ldots,a_r,b_1,\ldots,b_s,c_1,\ldots,c_t] \nn
n &= r+s+t .
\end{align}
Since extensors are strictly speaking not generic elements, but 
representants of an equivalence class, it is clear that they are not 
unique. One can find quite obscure statements about this fact in literature,
especially at those places where an attempt is made to visualise extensors
as plane segments, even as circular or spherical objects etc. However an 
extensor $A$ defines uniquely a linear 
subspace $\bar{A}$ of the space $\bigwedge V$ underlying the Gra{\ss}mann 
algebra. The subspace $\bar{A}$ is called {\it support\/} of $A$. 

A geometrical meaning of the join can be derived from the following.
The wedge product of $A$ and $B$ is non-zero if and only if the supports of 
$A$ and $B$ fulfil $\bar{A}\cap\bar{B}=\emptyset$. In this case the 
support of $A\wedge B$ is the subspace $\bar{A}\cup\bar{B}$. Hence the join
is the union of $\bar{A}$ and $\bar{B}$ if they do {\it not\/} intersect
and otherwise zero -- i.e. disjoint union. The join is an incidence
relation.

If elements of the linear space $V$ are called `points', the join of
two points is a `line' and the join of three points is a `plane' etc. One 
has, however, to be careful since our construction is till now 
characteristic free and such lines, planes, etc. may behave very oddly. 

\subsection{The vee-product -- meet}

The wedge product with multiplicators of step greater or equal than $1$ 
raises the step of the multiplicand in any case. This is a quite 
asymmetric and geometrical unsatisfactory fact. It was already undertaken 
by Gra{\ss}mann in the {\bf A1} (`eingewandtes Produkt') to try to find 
a second product which lowers the step of the multiplicand extensor by the 
step of the multiplicator. Gra{\ss}mann changed his mind and based his step
lowering product in the {\bf A2} on another construction. He also changed 
the name to `regressives Produkt' [regressive product]. It might be
noted at this place, that Gra{\ss}mann denoted exterior products
as `combinatorisches Produkt' [combinatorial product] showing his 
knowledge about its link to this field.

Already in 1955 Alfred Lotze showed how the meet can be derived using
combinatorial methods only \cite{lotze:1955a}. Lotze considered this
formula superior to the `rule of the double factor' and called it
`Universalformel' [universal formula]. Lotze pointed clearly out that 
the method used by Gra{\ss}mann in the {\bf A2} needs a symmetric
correlation, i.e. a transformation in projective geometry which introduces
a quadric. However, Cayley and Klein showed that having a quadric is half 
the way done to pass over to metrical geometries. Mentioning this point 
seems to be important since in recent literature mostly the less general 
and less powerful method of the {\bf A2} is employed. Zaddach, who was
aware of Lotze's work \cite{zaddach:1994a}, seemed to have missed the 
importance of this approach. The reader should also consult the articles of 
Zaddach p. 285, Hestenes p. 243, and Brini et al. p. 231 in 
\cite{schubring:1996a} which exhibit tremendously different approaches.

We will shortly recall the second definition of the regressive
product as given in the {\bf A2} by Gra{\ss}mann. First of all 
we have to define the `Erg\"anzung' of an extensor $A$ denoted by a 
vertical bar $\vert A$. Let $A$ be an extensor, the Erg\"anzung $\vert A$ 
is defined using the bracket by
\begin{align}
[A,\vert A] &= 1.
\end{align} 
From this equation it is clear that the `Erg\"anzung' is a sort of
orthogonal (!) complement or negation. But due to the fact that we 
consider disjoint unions of linear spaces, the present notion is more 
involved. We find for the supports of $A$ and $\vert A$
\begin{align}
\bar{A} \cap \bar{\vert A} &= \emptyset \nn
\bar{A} \cup \overline{\vert A} &= \bar{E}
\end{align}
where $E$ is the integral. Furthermore one finds that the Erg\"anzung is
involutive up to a possible sign which depends on the dimension $n$ of $V$.
Gra{\ss}mann defined the regressive product, which we will call  {\it meet\/}
with Rota et al. and following geometrical tradition. The meet is derived 
from
\begin{align}
\vert( A \vee B) &:= (\vert A) \wedge (\vert B) \nn
\end{align}
which can be accompanied by a second formula
\begin{align}
\vert (A \wedge B) &= \pm (\vert A) \vee (\vert B) 
\end{align}
where the sign once more depends on the dimension $n$. The vee-product
$\vee$ is associative and anticommutative and thus another instance of
an exterior product. The di-algebra (double algebra by Rota et al.) having 
two associative multiplications, sometimes accompanied with a duality
map, is called {\it Gra{\ss}mann-Cayley algebra\/}. The two above displayed 
formulas could be addressed as {\it de Morgan laws\/} of Gra{\ss}mann-Cayley 
algebra. This implements a sort of logic on linear subspaces, a game which 
ships nowadays under the term quantum logic. It was Whitehead who emphasised 
this connection in his Universal Algebra. 

The geometric meaning of the meet, which we denote by a vee-product $\vee$,
is that of intersection. We give an example in $\dim V= 3$. Let 
$\{e_1,e_2,e_3\}$ be an unimodular basis, then we find
\begin{align}
 &\vert e_1 = e_2 \wedge e_3 && \vert e_2 = e_3 \wedge e_1 
&&\vert e_3 = e_1 \wedge e_2 .
\end{align} 
If we calculate the meet of the following two 2-vectors 
$e_1\wedge e_2$ and $e_2\wedge e_3$ we come up with
\begin{align}
& &&\vert( (e_1\wedge e_2) \vee (e_2 \wedge e_3) ) 
=( e3 \wedge e1) \,=\, \vert e_2 &&\nn
&\Rightarrow && (e_1\wedge e_2) \vee (e_2 \wedge e_3) = e_2 && 
\end{align}
which is the {\it common factor\/} of both extensors. The calculation
of the Erg\"anzung is one of the most time consuming operation in 
geometrical computations based on meet and join operations.
This renders the present definition of the meet as computational
inefficient. Moreover, it is unsatisfactory that the meet is a `derived'
product and not directly given as the join or wedge.

The rule of the double [middle / common] factor reads as follows.
Let $A,B,C$ be extensors of step $a+b+c=n$ one assumes
\begin{align}
(A \wedge C) \vee (B \wedge C) &=
(A \wedge B \wedge C) \vee C.
\end{align} 
Using this relation one can express all regressive products in
wedge products alone. Hence one is able to compute. However, also this
mechanism renders the meet to be a derived and not a generic product.

\paragraph{Splits and shuffles:}

We will not follow Lotze's presentation \cite{lotze:1955a} of his 
`universal formula' but for convenience the more recent presentation 
of Doubilet et al. \cite{doubilet:rota:stein:1974a}. First of all notation 
is much clearer there and secondly we will use their mechanism to derive 
a single wedge product of two factors, while Lotze computes a formula for 
the wedge product of $r$ factors, motivating his `universal' since it 
additionally does not need a symmetric correlation. Only thereafter the 
more general alternative laws could be derived which we have no occasion 
to consider in any depth here.

For convenience we {\it drop the wedge sign\/} for multiplication in 
the following. Note that the antisymmetry of elements allows to 
introduce a linear order in any sequence of vectors from $V$. We can 
e.g. use lexicographic ordering of letters or if we use indexed entities 
we can order by the value of the index. A word of $\bigwedge V$ (i.e. an 
extensor) is called {\it reduced\/} if it is ordered w.r.t. the 
chosen ordering. For instance 
\begin{align}
&A=abcde &&B=b_1b_2\ldots b_r &&C=c_1c_3c_6
\end{align}
are reduced words i.e. ordered, but 
\begin{align}
&A=acebd && B=b_4b_2b_3b_1b_5 && C=c_6c_1c_3
\end{align}
are not properly ordered w.r.t. the chosen ordering and need to be 
reordered. If one wants to come up with a basis for $\bigwedge V$ this 
is constituted by reduced words. Note that there are lots of orderings 
and it will be important to carefully distinguish them. In the following, 
we deal with reduced words (ordered basis extensors) only. A main problem 
in calculating the products is to expand the outcome into reduced basis 
elements. These are the straightening formulas of Rota et al. which could
be called Littlewood-Richardson rule for Gra{\ss}mann-Cayley algebra
equivalently.

It would be a nice sidestep to study Young-tableaux, symmetric group
representations and Specht modules, which we however resist to do in 
this work.

A {\it block} of an extensor is a subsequence (subword) extracted from 
the extensor (word). A {\it $(\lambda_{i_1},\ldots,\lambda_{i_k})$-split\/}
of an extensor $A$ is the decomposition of the reduced word representing
$A$ into $k$ blocks of length $\lambda_{i_j}$ where $\sum \lambda_{i_s}
= $ step $A$. E.g. $A =a\ldots bc\ldots de\ldots f$ is decomposed into
$B_1=(a\ldots b)$, $B_2=(c\ldots d), \ldots$, $B_k=(e\ldots f)$. A
{\it shuffle\/} of the $(\lambda_{i_1},\ldots,\lambda_{i_k})$-split
of $A$ is a permutation $p\in S_{\text{step} A+1}$ of $A$ such that 
every block $B_s$ remains to be reduced. In other words, the blocks
$B_s$ consist of ordered subsequences of letters from the word 
representing $A$.

The meet of $k$ factors can be defined along the lines of Lotze
using these shuffles and splits into $k$ blocks. Rota et al. call 
these products bracket products. We will restrict ourselves to 
consider only ($s,t$)-splits into two blocks. Let 
$A=a_1\ldots a_k$ and $B=b_1\ldots b_s$ with $ \dim V=n$ and $k+s\ge n$. 
We define the {\it meet\/} $\vee$ as
\begin{align}
A \vee B &= \sum_{\text{shuffles}} \sign(p)
[a_{p(1)},\ldots,a_{p(n-s)},b_{1},\ldots,b_{s}]\;
a_{p(n-s+1)}\wedge \ldots \wedge a_{p(k)}
\end{align} 
where the permutations $p$ range over all ($n-s,k-n+s$)-shuffles of
$a_1\ldots a_k$. Note the order of factors inside the bracket, which
is given sometimes differently.

We introduce a {\it co-product\/} $\Delta \;:\; W \rightarrow 
W\otimes W$, which we will discuss later in detail, as the mapping of 
extensors $A$ into a sum of tensor products of its ($n-s,k-n+s$)-shuffles 
of subsequences
\begin{align}
\Delta(A) &=  \sum_{\text{shuffles}} \sign(p)\,
a_{p(1)}\ldots a_{p(n-s)} \otimes a_{p(n-s+1)}\ldots a_{p(k)} \nn
&= a_{(1)} \otimes a_{(2)}
\end{align}
where we have introduced a shorthand known as Sweedler notation which
implies the sum and the signs of the split as a sort of summation convention. 
Using this shorthand notation, the meet can be written as
\begin{align}
A \vee B &= [A_{(1)},B]\, A_{(2)} \;=\; 
B_{(1)}[A, B_{(2)}].
\end{align} 
The second identity holds if and only if the particular order
of factors is employed, otherwise a difference in sign may occur.

From this construction of the meet it is clear that no symmetric
correlation is needed and consequently no Erg\"anzungs operator 
has to be employed. The BIGEBRA package \cite{BIGEBRA} has both
versions implemented as {\tt meet} and {\tt \&v} products. There 
one can check the above identity on examples. Furthermore it turns 
out that the combinatorial implementation, which is ultimately based on 
Hopf algebra methods, is far more efficient than the above given
and widely utilized method using the Erg\"anzung. Especially in robotics,
where meet and join operations are frequently needed, this should
speed up calculations dramatically \cite{ablamowicz:fauser:2002b}. 
For benchmarks see the online help-page of {\tt meet} or {\tt \&v} from 
the BIGEBRA package.

\subsection{Meet and join for hyperplanes and co-vectors}

In projective geometry one observes a remarkable duality. If we consider
a 3-dimensional projective space a correlation maps points into planes 
and planes into points. It is hence possible to consider planes as
elementary objects and to construct lines and points by `joining' planes.
Projective duality shows that this geometry is equivalent to the geometry 
which considers points as basic objects and constructs lines and planes 
as joins of points. Recently projective duality was studied in terms
of Clifford algebras \cite{conradt:2000a,conradt:2000b,conradt:2000c}.
Clifford algebras have been employed for projective geometry in e.g.  
\cite{hestenes:ziegler:1991a}. However, the Clifford structure
is essentially not needed, but was only introduced to compute the 
Erg\"anzung. Ziegler has described the history of
classical mechanics in the 19th century \cite{ziegler:1985a} and showed
there, that screw theory and projective methods have influenced
the development of algebraic systems too. Gra{\ss}mann considered
(projective) geometry to be the first field to employ and exemplify 
his `new brach of mathematics', see {\bf A1}. Projective methods
are widely used in image processing, camera calibration, robotics etc.
\cite{bayrocorrochano:sobczyk:2001a}.

However in these fields, engineers and applied mathematicians do not like
co-vectors or tensor products, not to mention Hopf algebras. Rota et al.
tried to cure the case by introducing co-vectors using the bracket,
see \cite{barnabei:brini:rota:1985a}, p. 122. They black-listed Bourbaki's
treatment \cite{bourbaki:1989a} of co-vectors as follows: ``Unfortunately, 
with the rise of functional analysis, another dogma was making headway at 
the time, namely, the distinction between a vector space $V$ and its dual 
$V^*$, and the pairing of the two viewed as a bilinear form.'' A few lines
later, Hopf algebras are ruled out by stating that the ``common presentation
of both [interior and exterior products, BF] in the language of Hopf
algebras, further obscures the basic fact that the exterior algebra is 
a bird of a different feather. ... If one insists in keeping interior
products, one is sooner or later faced with the symmetry of exterior 
algebra as a Hopf algebra''. They develop a sort of co-vectors
inside the bracket or Gra{\ss}mann-Cayley algebra. We will see later, 
and Rota changed his mind also 
\cite{grosshans:rota:stein:1987a,rota:stein:1994a}, that this is not 
the proper way to deal with the subject. Indeed we have to reject even 
the term co-vector for this construction. We will call dual vectors
introduced by the bracket as {\it reciprocal vectors}. It will turn out 
that reciprocal vectors need implicitly the Erg\"anzung and imply
therefore the usage of a symmetric correlation. This introduces a 
distinguished quadric and spoils invariance under general projective
transformations. Our criticism applies for the now frequently used
homogenous models of hyperbolic spaces in terms of Clifford algebras
\cite{bayrocorrochano:sobczyk:2001a}.

If we identify vectors of the space $V$ of dimension $n$ with points, a 
hyperplane is represented by an extensor of step $n-1$. In other
words, $n-1$ linearly independent points span a hyperplane. If 
hyperplanes are identified with reciprocal vectors, one can define
an action of reciprocal vectors on vectors which yields a scaler. 
This motivated the misnaming of reciprocal vectors as co-vectors. We find
using summation convention and an unimodular basis $\{e_i\}$ and 
the Erg\"anzung
\begin{align}
& x\in V && x= x^i \, e_i \nn
& u\in \bigwedge^{n-1} V && u= u^{i_1,\ldots,i_{n-1}} 
e_{i_1}\wedge \ldots\wedge e_{i_{n-1}} \nn
& && u=u_{k} \vert e^{k}
\end{align}
where the Erg\"anzung yields the vector 
$e_{i_1}\wedge \ldots\wedge e_{i_{n-1}}\in \bigwedge^{n-1} V$, 
which we identify with the reciprocal vector $e^k$ and the coefficients 
$ u^{i_1,\ldots,i_{n-1}}$ are identified with $u_k$ accordingly. Using 
the bracket one gets $[e_i,e^k] = \delta_i^k$. This reads for a vector 
$x$ and a reciprocal vector $u$ 
\begin{align}
&u=u_k e^k\, .
\end{align}
We are able to use the bracket to write for the action $\bullet$ of a 
vector $x$ on a reciprocal vector $u$
\begin{align}
x\bullet u &= x^i u_k [e_i, e^k] \nn
&= x^i u_k [e_i, \vert e_k] = x^i u_k \delta_i^k \nn
&= x^i u_i\; \in \openk.
\end{align}
This mechanism can be generalized to an action of $\bigwedge^{r} V$ on
$\bigwedge^{n-r} V$. The usage of the Erg\"anzung implying a quadric 
is pretty clear. This construction is used in \cite{hestenes:sobczyk:1992a}
to derive `co-vectors'. Hence all their formulas are not applicable
in projective geometry which does not single out the Erg\"anzung or a 
symmetric correlation which implies a quadric.

However, we can follow Lotze, \cite{lotze:1955a} note added in prove,
to do the same construction but starting this time from the space of planes.
Let $\vartheta\in V^*$ be a co-vector and $\{\vartheta^a\}$ be a set of 
canonical co-vectors dual to a basis $\{x_i\}$ of vectors spanning $V$, 
i.e. $\vartheta^a x_i = \delta^a_i$. One can form a Gra{\ss}mann algebra 
on $V^*$ along the same lines as given above by introducing a bracket 
on $\times^n V^*$. We denote the exterior product of this particular 
Gra{\ss}mann algebra by vee $\vee$ that is the meet (join of hyperplanes). 
This reflects the fact that if it is allowed in a special case that 
co-vectors and reciprocal vectors are identified, their product is the meet. 
We can derive along the same lines as above a dual product called join.
This join plays the same role to the above meet as the meet played beforehand
to the join. It is denoted as `join' (meet of hyperplanes) and uses the 
wedge $\wedge$ symbol, using splits and shuffles. It turns out, as our 
notation has anticipated, that this operation is the join of points, if 
points are identified as $n-1$-reciprocal vectors of co-vectors. We have
experienced an instance of product co-product duality here, which will
be a major topic in the later development of this treatise.

This consideration, which is exemplified to some detail in the online 
help-page of the {\tt meet} and {\tt\&v} products of the BIGEBRA
package \cite{BIGEBRA}, shows that it is a matter of choice which 
exterior product is used as meet and which as join by dualizing. This is 
the reason why we did not follow Rota et al. to use the vee-product $\vee$
for the join of points to make the analogy to Boolean algebra perfect.

However, we can learn an important thing. It is possible and may be
necessary to implement an exterior algebra on the vector space
$V$ and the co-vector space $V^*$ {\it independently}. This will
give us a great freedom in the Hopf gebraic structure studied below.
Moreover, it will turn out to be of utmost importance in QFT. 
Reordering and renormalization problems are hidden at this place.
After our remarks it might not surprise that also classical differential
geometry can make good use of such a general structure 
\cite{vargas:torr:2002a,fauser:dehnen:2002a}.

%
%

%% file: cl_basics.tex
\chapter{Basics on Clifford algebras}

\section{Algebras recalled}

In this section we recall some definitions and facts from module and ring 
theory. In the sense we use the terms `algebra' and `ring', they are synonyms. 
We want to address the structure of the scalars as ring and the 
additive and distributive multiplicative structure on a module as algebra.
The following statements about rings hold also for algebras.

From any book on module theory, e.g \cite{wisbauer:1988a}, one can take 
the following definitions:

\begin{dfn} A {\em ring} is a non-empty set $R$ with two morphisms
$+,\cdot\;:\; R\times R \rightarrow R$ fulfilling 
\begin{align}
i)& &&(R,+) &&\text{is an abelian group, $0$ its neutral element} \nn
ii)& &&(R,\cdot) &&\text{is a semigroup} \nn
iii)& && (a+b)c = ab+bc&&\forall\;a,b,c \in R \nn
     &&& a(b+c) = ab + ac
\end{align}
\end{dfn}  

A ring $R$ (same symbol for the underlying set and the ring) is called 
commutative, if $(R,\cdot)$ is commutative. If the multiplication map 
$\,\cdot\,$ enjoys associativity, the ring is called associative. We will 
assume associativity for rings.

An element $e\in R$ is called {\it left (right) unit\/} if $ea=a$ ($ae=a$) 
for all $a\in R$. A {\it unit\/} is a left and a right unit. A ring
with unit is denoted {\it unital ring\/}.

The {\it opposite ring $R^{op}$} of $R$ is defined to be the additive 
group $(R,+)$ with the opposite multiplication
\begin{align}
a\circ^{op} b &= b\cdot a .
\end{align}

A subgroup $I$ of $(R,+)$ is called {\it left ideal\/} if 
$R\cdot I \subset I$ holds and {\it right ideal\/} if 
$I\cdot R \subset I$ holds. An {\it ideal\/} (also bilateral ideal)
is at the same time a left and right ideal. If $(R,\cdot)$ is commutative
then every ideal is a bilateral ideal. The intersection of left (right) 
ideals is again a left (right) ideal.

A {\it morphism of (unital) rings\/} is a mapping $f : (R,+,\cdot)
\rightarrow (S,+,\circ)$ satisfying
\begin{align}
f(a+b) &= f(a) + f(b) \nn
f(a\cdot b) &= f(a) \circ f(b) \nn
f(e_R) &= e_S &&\text{if $e_R,e_S$ do exist} .
\end{align}
The kernel of a ring homomorphism $f : R\rightarrow S$ is an ideal
\begin{align}
I_f &= \ker f \,=\, \{ a\in R \mid f(a)=0 \}.
\end{align}
The converse is true, every ideal is the kernel of an appropriate
homomorphism. The {\it canonical projection\/} is given as 
$\pi_I : R\rightarrow R/I$ where $R/I$ is the {\it residue class ring}. 
The ring structure in $R/I$ is given as ($a,b\in R$)
\begin{align}
(a+I)+(b+I) &= (a+b+I) \nn
(a+I)(b+I)  &= (ab +I).
\end{align}
$R/I$ is also called a {\it factor ring\/}.

Let $A$ be a subset of $R$. An (left/right) ideal $I_A$ is called 
{\it generated by $A$\/} if it is the smallest (left/right) ideal $I_A$
with $A\subset I_A$. If $A$ has finite cardinality we call $I_A$ 
{\it finitely generated\/}. $I_A$ is the intersection of all ideals 
which contain $A$.

The {\it direct sum\/} $A\oplus B$ of two ideals $A,B$ is defined to be
their Cartesian product $A\times B$ under the condition $A\cap B =\emptyset$.
The ring $R$ is called {\it decomposable\/} if it is a direct sum of 
(left/right) ideals $R=A\oplus B \oplus\ldots$,
$A\cap B = \emptyset$, etc. In such rings every element $r$ can be 
{\it uniquely\/} decomposed as
\begin{align}
R \ni r&= a+b+\ldots \nn
a &\subset A,\;\; b\subset B,\ldots
\end{align}
A ring is called (left/right) {\it indecomposable\/} if it cannot
be written as a direct sum of (left/right) ideals. An analogous definition
applies for ideals.

We define some special elements which will be needed later. An element
$a$ of the ring $R$ is denoted as
\begin{itemize}
\item {\it left divisor of zero\/} if it exists a $b\not=0$ such that
$ab=0$.
\item {\it right divisor of zero\/} if it exists a $b\not=0$ such that
$ba=0$.
\item {\it divisor of zero\/} if it is a left and right divisor of zero.
\item {\it idempotent\/} if $a^2=a$.
\item {\it nilpotent\/} (of order $k$) if $a^k=0$.
\item {\it unipotent\/} if R is unital and $a^2=e$.
\item {\it regular\/} if it exists an element $b\in R$ with $aba=a$.
\item {\it left (right) invertible\/} if $R$ is unital and it 
exists an element $b\in R$ such that $ab=e$ ($ba=e$).
\item {\it invertible\/} if it is left and right invertible.
\item {\it central\/} if for all $b\in R$ holds $ab-ba=0$.
\end{itemize}

Two idempotents $f_1,f_2$ are called orthogonal if $f_1f_2=0=f_2f_1$.
An idempotent is called primitive if it cannot be written as the 
orthogonal sum of idempotents.

A subset $A$ of $R$ is called 
{\it left annulator\/} if $An^l(A):=\{b\in R\mid ba=0,\quad\forall a\in R\}$,
{\it right annulator\/} if $An^r(A):=\{b\in R\mid ab=0,\quad\forall a\in R\}$,
or {\it annulator\/} if $An^l(A):=An^l(A) \cap An^r(A)$.

\begin{thrm}[Direct decomposition]
Let $R$ be a ring, it holds
\begin{itemize}
\item[1)] If the left ideal $I\subset R$ is generated by an idempotent
$f\in R$, $I=Rf$, then $R$ is decomposable into left ideals
$R = A\oplus An^l(f)$.
\item[2)] If the ideal $J$ is generated by a central idempotent $f$
then $R$ is decomposable into $R=J+An(f)$.
\item[3)] Let $R$ be an unital ring. Every (left/right) ideal $I$
which is a direct summand is generated by an idempotent element $f$.
If $I$ is an ideal then $f$ is central. The decomposition is
$R=Rf+An^l(F)$, where $An^l(f)=R(1-f)$.
\end{itemize}
\end{thrm}

\noindent 
{\bf Proof:} see \cite{wisbauer:1988a}.

\section{Tensor algebra, Gra{\ss}mann algebra, Quadratic forms}

Our starting point to construct Clifford algebras and later on Clifford 
Hopf gebras will be the Gra{\ss}mann algebra. We have already used the
language of an alphabet having letters which do form words to introduce 
this mathematical structure in the chapter on the Peano bracket. Hence we
will introduce here the same structure by factoring out an ideal from 
tensor algebra. We will have occasion to use this technique later on.

Let $\openk$ be an unital commutative ring and let $V$ be a $\openk$-linear 
space. The tensor algebra $T(V)$ is formed by the direct sum of tensor 
products of $V$
\begin{align}
T(V) &= \openk\, \oplus\;V\; \oplus\; (V\otimes V)\; \oplus \ldots \nn
     &= \oplus_r T^r(V) \;=\; \oplus_r\;\otimes^r V.
\end{align}
We identify $\openk$ with $V^0$ in a canonical way. The unit of $\openk$
in $T(V)$ is denoted as $\Id$. The injection $\eta \;:\; \openk\rightarrow 
T(V)$ into the tensor algebra will be needed below and is called 
{\it unit map}, also denoted $\Id_V$. The elements of the set $\{e_i\}$ 
of linearly independent elements which span $V$ are called set of
generators. The words obtained from these generators by concatenation 
yields a basis of $T(V)$. All elements of $V$ are called letters, 
decomposable elements of $\otimes^r V$, i.e. 
$a_1\otimes\ldots\otimes a_r \in T^r(V)\,\equiv\,\otimes ^r V$. The 
number of factors is called {\it length} of the word or {\it rank} of the 
tensor. One can add words of the same length which will in general lead to 
an indecomposable tensor, but still a tensor of the same rank. Sums of words 
of same or arbitrary length might be called sentences. The tensors of a 
particular rank form a linear subspace of $T(V)$. Products of tensors are 
formed by concatenation of words,
\begin{align}
(a_1\otimes\ldots\otimes a_r) \otimes (b_1\otimes\ldots\otimes b_s)
&=
a_1\otimes\ldots\otimes a_r \otimes b_1\otimes\ldots\otimes b_s .
\end{align}
Concatenation is by definition associative. $T(V)$ is naturally 
graded by the length or rank of the tensors. i.e. products of $r$-tensors
and $s$-tensors are $r+s$-tensors. For a precise definition of the 
tensor product look up any algebra book, 
\cite{scheja:storch:1980a,scheja:storch:1988a}.

The Gra{\ss}mann algebra is obtained by projecting the tensor product
onto the antisymmetric wedge product $\pi(\otimes) \rightarrow \wedge$.
In the case of the Gra{\ss}mann algebra, we can either describe the 
equivalence class or deliver relations among some generators. If relations 
hold, not all words of the tensor algebra which can be formed by 
concatenation remain to be independent. The problem to identify two words 
w.r.t. given relations is called the {\it word problem}. It can in general 
not be solved, however, we will deal with solvable cases here. To be able 
to pick a representant from an equivalence class, we have to define 
{\it reduced words}. A reduced word is semi {\it ordered} in a certain 
sense. One has to use the relations to establish such a semi ordering, 
sometimes called {\it term ordering\/} in the theory of Gr\"obner bases.

We define the following ideal which identifies all but antisymmetric
tensors
\begin{align}
{\cal I}_{\wedge} &= \{ a\otimes x\otimes x\otimes b \mid 
a,b\in T(V),\quad x\in V\} .
\end{align} 
The Gra{\ss}mann algebra $\bigwedge V$ is the factor algebra of 
$T(V)$ where the elements of the above given ideal are identified 
to zero.
\begin{align}
\bigwedge V &= \frac{T(V)}{{\cal I}_{\wedge}} \nn
            &= \pi_\wedge\,( T(V))
\end{align}
where $\pi_\wedge$ is the canonical projection from $T(V)$ onto $\bigwedge V$.
From this construction it is easy to show by means of categorial
methods that a Gra{\ss}mann algebra over a space $V$ is a {\it universal
object} and is defined uniquely up to isomorphy.

The relations which are equivalent to the above factorization read
\begin{align}
e_i \otimes e_i &= 0 \mod {\cal I}_{\wedge}\nn
\pi_\wedge(e_i \otimes e_i) &= e_i \wedge e_i \,=\,0 \\
\pi_\wedge(e_i \otimes e_j) &= e_i \wedge e_j \,=\, -e_j\wedge e_i .
\end{align}
While the tensor algebra had essentially no calculational rules to 
manipulate words or sentences, beside multilinearity, one has to respect 
such relations after factorization. We can introduce reduced words by
asserting that words of generators are ordered by ascending (descending)
indices. A basis of $\bigwedge V$ is given as 
\begin{align}
GB &= \{\Id;e_1,\ldots,e_n;e_1\wedge e_2,\ldots e_{n-1}\wedge e_n;\ldots
;e_1\wedge\ldots\wedge e_n\}
\end{align}
where we have separated words of a different length by a semi-colon. Due 
to the relations we find for a finitely generated space $V$ of dimension 
$n$ a {\it finite number} of reduced words only. Their
number is $\sum {n\choose r} = 2^n$. The space spanned by these generators 
will be called $W=\bigwedge V$. In analogy to the group theory
\cite{johnson:1997a} we can define a {\it presentation of an algebra} 
over $V$ spanned by the set of generators $X$ as follows:
\begin{align}
{\rm Alg}(V) &= \; \langle X,{\bf R}\rangle \nn
             &= \{ \,(\{e_i\},\{{\bf R_i}\}) \mid V=\,\text{span}\{e_i\}, 
                \text{${\bf R_i}$ relations}\,\}. 
\end{align}
We will freely pass from one picture to the other as it is convenient. 
The techniques from group presentations and terminology, e.g. word problem,
generator, etc. can be applied to algebras by analogy. E.g. a free algebra
is an algebra generated by a set $X$ of generators $e_i$ which span $V$
having no relations at all. A free Lie algebra has of course relations
which renders it to be a Lie algebra, but no further constraints among
its `Lie words'.

We had already occasion to define quadratic forms previously, so we
recall here only the basis free definition 
\begin{align}
Q(\alpha\,x) &=\alpha^2\,Q(x) \nn
2\;B_p(x,y)  &=Q(x-y)-Q(x)-Q(y) &&\text{$B_p$ bilinear.}  
\end{align}
As we pointed out, the addition of a quadratic form to a linear space
yields a quadratic space. The main idea of a Clifford algebra is to form
an algebra in a natural way from this building blocks.
One can show that there is a functorial relation between quadratic spaces
and associative unital algebras. This functor is injective and denoted as 
$\Cl$. It is clear from this observation that the classification of Clifford
algebras is essentially given by the classification of the quadratic
forms used in their construction. If $\openk$ is $\openR$ or $\openC$,
this can be readily done by signature and dimension in the case of 
$\openR$ or dimension only in the case of $\openC$.

In the following sections we will provide some possible methods to 
establish this functorial relation. Each method has its advantages 
in certain circumstances, so none has to be abandoned, however, we will 
spend lots of efforts to provide a universal, computationally efficient, 
and sound approach to Clifford algebras, which will turn out to be Rota-Stein
cliffordization. We will in the same time generalize the term 
Clifford algebra to {\it Quantum Clifford algebra\/} (QCA) if we consider 
algebras built from spaces having a bilinear form of arbitrary symmetry.
It will turn out during our treatment of the subject that we will need 
necessarily the co-algebra and Hopf algebra structure which is hidden or 
implicit in the more basic approaches. Hopf techniques will be extremely 
helpful in applications, speeding up actual computations, e.g. of meet and 
join, used in robotics. The same holds true for Clifford products,
\cite{ablamowicz:fauser:2002a,ablamowicz:fauser:2002b}. Cliffordization
turns out to be a neat device to describe normal-, time-, and even
renormalized time-ordered operator products and correlation functions 
in QFT.

\section{Clifford algebras by generators and relations}

The generator and relation method is the historical root of several
algebraic systems. Hamilton's quaternion units ${\bf i,j,k}$ are still
used in vector analysis, Gra{\ss}mann used basis vectors $e_i$ to 
generate his `Hauptgebiet', our linear space $V$. A basis independent
method was in general not available during these times, hence, also 
Clifford introduced and studied algebras in terms of generators and 
relations. The presentation of a Clifford algebra is as follows:
\begin{align}
\Cl(V,Q) &= <X,R> \nn
&= <\,\{e_i\}, e_ie_j + e_je_i = 2\,g_{ij}>
\end{align}
where the $e_i\in X$ span $V$ and $g_{ij}$ is the symmetric polar
bilinear form which represents $Q$ in the basis of the generators.
These relations are usually called (anti)commutation relations. In physics
only the commutation relations of the generators are usually given to 
define algebras, hence one writes
\begin{align}
e_ie_j + e_je_i &= 2\,g_{ij}.
\end{align}
Synonymous notations are $\Cl(V,Q)$, $\Cl(Q)$ if $V$ is clear,
$\Cl_{p,q,r}$ if $V$ is an $\openR$-linear space of dimension $p+q+r$,
while the quadratic form has $p$ positive, $q$ negative eigenvalues and 
a radical of dimension $r$, and $\Cl_n$ if $V$ is a $\openC$-linear space 
of dimension $n$. The Clifford product is denoted by juxtaposition or
if we want to make it explicite by a circle $\circ$, sometimes called 
circle product
\cite{rota:stein:1994a}.
A natural basis for this algebra would be the Clifford basis, ordered
by ascending indices
\begin{align}
CB &= \{\Id;e_1,\ldots,e_n;e_1e_2,\ldots e_{n-1}e_n;\ldots
;e_1\ldots e_n\}
\end{align}
which does not resort to the Gra{\ss}mann exterior product. But most
applications actually use a Gra{\ss}mann basis. Such a basis is obtained
by antisymmetrization of the Clifford basis elements, e.g.
\begin{align}
e_i \wedge e_j &= \frac{1}{2}(e_i e_j - e_j e_i).
\end{align}
It was shown by Marcel Riesz \cite{riesz:1957a} that a wedge product
can be consistently developed in a Clifford algebra. This basis is 
isomorphic to a basis of a Gra{\ss}mann algebra. Hence it is clear, that 
Clifford and Gra{\ss}mann algebras have the same dimension. We will see 
below, that one can construct Clifford algebras as a subalgebra of the 
endomorphism algebra of an underlying Gra{\ss}mann algebra.

The most remarkable changes between a Gra{\ss}mann and a Clifford algebra 
are, that the latter has a richer representation theory. This stems from the 
fact that in a Gra{\ss}mann algebra $\Id$ and $0$ are the only idempotent
elements. That is $\bigwedge V$ is an indecomposable algebra.
One finds beside nilpotent ideals only trivial ideals. 
Clifford algebras have idempotent elements which generate various spinor
representations. This fact follows directly from the quadratic form 
introduced in the Clifford algebra.

We had noted that the Gra{\ss}mann basis $GB$ spans a $\openZ$-graded
linear space. The exterior wedge product was graded too. Since the Clifford
algebra can be described using a Gra{\ss}mann basis, it seems to be possible
to introduce a $\openZ$-grading here also. However, a short calculation
shows that the Clifford product does {\it not} respect this grading, but
only a weaker filtration, see later chapters. Let $u,v$ be extensors
of step $r$ and $s$ one obtains
\begin{align}
u \circ v & \in \oplus_{n=\vert r-s\vert}^{r+s}\;
\bigwedge^{n}V .
\end{align}
This is not an accident of the foreign basis, but remains to be true
in a Clifford basis also. The terms of lower step emerge from the necessary
commutation of some generators to the proper place in a reduced word.
For instance 
\begin{align}
(e_1e_2e_5)\circ(e_3e_4e_6) &= (e_1e_2e_3e_4e_5e_6) +
g_{35}(e_1e_2e_4e_6) - g_{45}(e_1e_2e_3e_6) .
\end{align}
As a matter of fact a Clifford algebra is only $Z_2$-graded since
even- and oddness of the length of words is preserved. The commutation 
relation contracts two generators for each commutation.

The usually defined {\it grade projection operators} 
$<\ldots>_r \;:\; \bigwedge V\rightarrow \bigwedge^r V$ are foreign to
the concept of a Clifford algebra and belongs to the underlying 
Gra{\ss}mann algebra. We will see later, that one is able to employ 
various $\openZ$-gradings at the same time. It will be of 
great importance to keep track of the grading which is inherited from the
Gra{\ss}mann algebra. However, the mere {\it choice} of a set of
generators $\{e_i\}$ induces a $\openZ_n$-grading w.r.t. an underlying 
Gra{\ss}mann algebra. The question if such representations are equivalent 
is known as {\it isomorphy problem}
in the theory of group presentations \cite{johnson:1997a}. In fact it
is easy to find, e.g. using CLIFFORD \cite{CLIFFORD}, non grade
preserving transformations of generators. This is well known from the 
group theory. E.g. the braid group on three strands has presentations
\begin{align}
B_3 &= <\{x,y\},xyx=yxy>\qquad\text{or} \nn
B_3 &= <\{a,b\},a^3=b^2>
\end{align}
where one sets with $xy=a$ and $x=a^{-1}b$ 
\begin{align}
&y=x^{-1}a,&&ax=x^{-1}a^2
\end{align}
and finds that the length function w.r.t. the generators $x,y$ is different
to that w.r.t. $a,b$.

This observation is crucial for any attempt to identify algebraic expressions
with geometric objects. The same will hold in QFT when identifying operator
products.

\section{Clifford algebras by factorization}

Clifford algebras can be approached in a basis free manner which for 
obvious reasons avoids the problems discussed in the previous section. 
While generators can be used very conveniently in actual calculations,
the strength of the basis free method is to achieve general statements
about the structure of Clifford algebras.

Following the procedure which led to the Gra{\ss}mann algebra, we can 
introduce an ideal $I_{\Cl}$ and factor out the Clifford algebra from the 
tensor algebra $T(V)$. This ideal has to introduce the quadratic form and 
reads
\begin{align}
{\cal I}_{\Cl} &= \{ a\otimes(x\otimes y + y\otimes x)\otimes b-
2g(x,y) a\otimes b \vert a,b \in T(V),\quad x,y\in V\}
\end{align}
where $g(x,y)$ is the basis free symmetric polar bilinear form corresponding
to $Q$. Inspection of the elements in this ideal shows that they are not 
homogeneous and identify elements of different rank. This ideal is not 
$\openZ_n$-graded. Since even- and oddness is preserved by the ideal, it 
remains to be $\openZ_2$-graded.

We arrive at the Clifford algebra via the following factorization
\begin{align}
\Cl(V,Q) &= \frac{T(V)}{{\cal I}_{\Cl}}\,.
\end{align}
Following Chevalley \cite{chevalley:1997a} (see ``The construction and study 
of certain important algebras'') one is able to show that Clifford algebras
are universal, which allows to speak about {\it the} Clifford algebra 
(up to isomorphy). Existence is also proved in this approach. 

The most important and structural interesting observation may be however
the identification of $\Cl$ as a functor. We call
a space {\it reflexive\/} if its dual has a set of generators of the same 
cardinality. All finite dimensional spaces are reflexive in this sense. 
Infinite dimensional spaces are usually not, but if generators are used, 
we want to have an isomorphism between generators for the spaces $V$ and 
$V^*$. Let ${\cal H}$ be a reflexive quadratic space, i.e. a pair of a 
linear space $V$ and a symmetric quadratic form $Q$. We find that $\Cl$ 
is an injective functor from the category (see Chapter 4) of quadratic 
spaces {\bf Quad} into the category of associative unital algebras {\bf Alg}.
\begin{align} 
\Rnode{A}{\text{\bf Quad}}\hskip 3truecm 
\Rnode{B}{\text{\bf Alg}}
\ncline{->}{A}{B}
\Aput{\Cl}
\end{align}
In the same manner we could have introduced a Gra{\ss}mann functor
${\bf\bigwedge}$. Functorial investigations would lead us also to the 
cohomology of these algebras. In fact, we will need the functorial
approach later to define the concept of a co-algebra, co-products
etc. by a simple duality argument.

\section{Clifford algebras by deformation -- Quantum Clifford algebras}

The previous section is to some extend unsatisfactory since it does not
allow to {\it compute} in a plain way. Even the generator and relation 
method suffers from computational difficulties. It is quite not easy, 
to Clifford multiply e.g. two extensors $u,v$. As an example we
compute
\begin{align}
e_1\circ(e_2\wedge e_3)
&=\frac{1}{2}e_1\circ(e_2\circ e_3 - e_3\circ e_2) \nn
&=\frac{1}{6}(
 e_1\circ e_2\circ e_3 + e_2\circ e_3\circ e_1 + e_3\circ e_1\circ e_2  \nn
&\phantom{=\frac{1}{6}(}
-e_1\circ e_2\circ e_3 - e_2\circ e_3\circ e_1 - e_3\circ e_1\circ e_2)  \nn
&\phantom{=\frac{1}{6}(}
+ \frac{4}{6}g_{12}e_3 - \frac{2}{6}g_{13}e_2
- \frac{4}{6}g_{13}e_2 + \frac{2}{6}g_{12}e_3 \nn
&= e_1 \wedge e_2 \wedge e_3 + g_{12}e_3 - g_{13}e_2
\end{align}
which is cumbersome due to the fact that we have to recast exterior products
into Clifford products where we can use the (commutator) relations. Finally
one has to transform back at the end into the wedge basis of reduced words
of the Gra{\ss}mann basis. Furthermore, that factor $2$ occurring in the 
(anti)commutation relations prevents an application of this mechanism to
rings of characteristic $2$. Claude Chevalley developed a method which
is applicable to this case and which provides an efficient method to
evaluate the Clifford product in a Gra{\ss}mann basis \cite{chevalley:1997a}.
An emphatic article of Oziewicz \cite{oziewicz:1986a} generalized Chevalley's
method from quadratic forms to bilinear forms. This will be a key
point in later applications to QFT.

Chevalley's observation was that it is possible to implement the Clifford 
algebra as an endomorphism algebra of the Gra{\ss}mann algebra
\begin{align}
\Cl &\subset \End \bigwedge V.
\end{align}
This inclusion is strict. To be able to define an endomorphism on
$\bigwedge V$, we have to introduce a dual basis and a dual Gra{\ss}mann
algebra $\bigvee V$. Let $\varepsilon^i(e_j) = \delta^i_j$, where 
$\delta^i_j$ is the Kronecker symbol, and let $\{\Id;\varepsilon^i;
\varepsilon^i\vee\varepsilon^j (i<j);\ldots \}$ be a Gra{\ss}mann co-basis
w.r.t. the vee-product. An endomorphism on $\bigwedge V$ can be written
as
\begin{align}
R &: \bigwedge V\rightarrow \bigwedge V \nn
R &= \sum_{IK} R^I_{\phantom{I}K}e_I\otimes \varepsilon^K
\end{align}
where $I,K$ are multi-indices of ordered basis words (basis monomials).

\subsection{The Clifford map} 

Let $B$ be a scalar product $B : V \times V \rightarrow \openk$. $B$ is 
at the same time a map $B : V \rightarrow V^*$. The action of the 
co-vectors $\varepsilon^i$ on vectors $e_j$ does form a pairing 
$<.\mid .> : V^* \times V \rightarrow \openk$.
\begin{dfn}[contraction]
Using the pairing $<.\mid .>_B$, where the scalar product $B$ is used to mediate 
the adjoint map,  
a {\em left (right) contraction} $\JJ$ (\,$\LL$) is defined as
\begin{align}
<\varepsilon^i\mid e_j >_B &= < \Id \mid B(\varepsilon^i) \JJ_{\delta}
e_j> \,=\, <\Id, e_i \JJBB e_j> \nn
<\varepsilon^i\mid e_j >_B &= < \varepsilon^i\LL_{\delta} 
B^{-1}(e_j) \mid_{\delta} \Id> \,=\, 
<\varepsilon^i \LL_{B^{-1}} \varepsilon^j, \Id>
\end{align}
\end{dfn}

\begin{dfn}[Clifford map] 
A {\em Clifford map} $\gamma_x : \bigwedge V \rightarrow \bigwedge V$ 
is an endomorphism parameterized by a 1-vector $x\in V$ of the following 
form
\begin{align}
\gamma_x &= x\JJBB\;+\; x\wedge
\end{align}
obeying the following calculational rules ($x,y\in V$, $u,v,w\in \bigwedge V$):
\begin{align}
i)& &&&x\JJBB y &= B(x,y) &&&\nn
ii) & &&& x \JJBB (u\wedge v) &= (x\JJBB u) \wedge v 
   + \hat{u} \wedge (x\JJBB v) &&&\nn
iii)& &&& (u\wedge v) \JJBB w &= u\JJBB (v \JJBB w) &&&  
\end{align}
where $\hat{~}$ is the main involution $\hat{~} : V \rightarrow -V$,
extended to $\bigwedge V$, also called grade involution. One obtains 
$\hat{u} = (-1)^{\text{length}(u)}u$.
\end{dfn}

We decompose $B=g+F$ into a symmetric part $g^T=g$ and an antisymmetric
part $F^T=-F$. The Clifford maps $\{\gamma_{e_i}\}$ of the generators 
$\{e_i\}$ of $V$ generate the Clifford algebra $\Cl(V,B)$. Let $\Id\;$ be 
the identity morphism, we find in a basis free notation
\begin{align}
\gamma_x\gamma_y + \gamma_y\gamma_x &= 2\, g(x,y) \Id\,.
\end{align}
It is remarkable, that in the anticommutation relation only the 
symmetric part of $B$ occurs. However, the anticommutators are 
altered
\begin{align}
\gamma_x\gamma_y -\gamma_y\gamma_x &= 2\, x\wedge y + 2\, F(x,y)\Id\,.
\end{align}
This shows that the $Z_n$-grading depends directly on the presence of the
antisymmetric part. If we compute a Clifford basis with or without 
an antisymmetric part $F$ we get ($\gamma^g\in \Cl(V,g)$,
$\gamma^B\in \Cl(V,B)$)
\begin{align}
\Id & &&&\Id& &&\nn
\gamma^g_{e_i}\Id &= e_i &&& \gamma^B_{e_i}\Id& = e_i &&\nn
\gamma^g_{e_i}\gamma^g_{e_j}\Id &= e_i\wedge e_j+g_{i,j} 
&&& \gamma^B_{e_i}\gamma^B_{e_j}\Id& = e_i\wedge e_j +B_{ij}&& \nn
\text{etc.} &&&& &&
\end{align}
If $g$ is identical zero $g\equiv 0$ we find two different 
{\it Gra{\ss}mann algebras}! One is $Z_n$-graded w.r.t. the exterior 
wedge products $\wedge$ while the other is not! It is however possible
to introduce a second {\it dotted wedge\/} $\dwedge$, also an exterior 
product, which is the $Z_n$-graded product under the presence of the 
antisymmetric part $F$.
\begin{align}
x\dwedge y &= x\wedge y +F(x,y)\Id \nn
x\dwedge y \dwedge z &= x\wedge y\wedge z
  +F(x,y)z + F(y,z)x + F(z,x)y \nn
\text{etc.}&
\end{align}
This structure was employed to obtain Hecke algebra representations
\cite{fauser:1999b,ablamowicz:fauser:2000b} and is crucial to the
compact formulation of Wick's theorem in QFT
\cite{fauser:1996c,fauser:1998a,fauser:2001b}.

\subsection{Relation of $\Cl(V,g)$ and $\Cl(V,B)$}

\begin{thrm}[Wick theorem]
The Clifford algebras $\Cl(V,g)$ and $\Cl(V,B)$ are isomorphic
as Clifford algebras. The isomorphisms in $Z_2$-graded.
\end{thrm}

\noindent
{\bf Proof:} see 
\cite{fauser:1996c,fauser:1998a,fauser:ablamowicz:2000c,fauser:2001b}.

\begin{thrm}[Chevalley \cite{chevalley:1997a}]
The opposite Clifford algebra $\Cl^{op}(V,g)$ of $\Cl(V,g)$ is isomorphic
to $\Cl(V,-Q)$.
\end{thrm}

This can be generalized to

\begin{thrm}
The opposite Clifford algebra $\Cl^{op}(V,B)$ of $\Cl(V,B)$ is isomorphic
to $\Cl(V,-B^T)$.
\end{thrm}

\noindent
{\bf Proof:} see 
\cite{fauser:stumpf:1997a,fauser:1998a}.

\vspace{0.2cm}

One obtains that 

\begin{align}
\End \bigwedge V 
&= \bigwedge V \otimes \bigvee V^*
&= \Cl(V,B)\,\hat{\otimes}\, \Cl(V,-B^T) 
&= \Cl(V\oplus V,B\oplus -B^T)
\end{align}
where $\hat{\otimes}$ is a $\openZ_2$-graded tensor product. In terms of 
commutation relations this reads 
\begin{align}
\gamma_x\gamma_y + \gamma_y\gamma_x &= 2\, g(x,y) \nn
\gamma_x\gamma^{op}_y + \gamma^{op}_y\gamma_x &= 0 \nn
\gamma^{op}_x\gamma^{op}_y + \gamma^{op}_y\gamma^{op}_x &= -2\, g(x,y) .
\end{align}

\section{Clifford algebras of multivectors}

An intriguing approach to Clifford algebras was developed by Oziewicz
and will be called {\it Clifford algebra of multivectors\/}. This method 
originated out of a discussion of QF theoretic composite particle 
calculations \cite{fauser:stumpf:1997a} which was elevated in 
\cite{oziewicz:1997b}
to a mathematical setting. We recall this approach here for completeness
and because of its extraordinary character and generality.

It was Woronowicz \cite{woronowicz:1987a,woronowicz:1989a} who studied
systematically the theory of deformed Gra{\ss}mann algebras. As we 
discussed above, Gra{\ss}mann algebras are obtained by factorization
w.r.t. an antisymmetrizer, which projects out all symmetric tensors
from tensor algebra. The canonical projection $\pi_\wedge$ maps the
tensor product $\otimes$ onto the exterior wedge product 
$\pi_\wedge(\otimes)\rightarrow \wedge$. If one proceeds to deformed
symmetries, e.g. Hecke algebras, one obtains deformed Gra{\ss}mann
algebras $\bigwedge_q V$. The presentation of the symmetric algebra
reads
\begin{align}
S_n &= \;<X,\{{\bf R_1},{\bf R_2},{\bf R_3}\}> \nn
{\bf R_1} &:  s_1^2=1 \nn
{\bf R_2} &:  s_i s_j s_i = s_j s_i s_j \nn
{\bf R_3} &:  s_i s_j = s_j s_i \qquad \text{if~~} \vert s_i-s_j\vert \ge 2 .
\end{align}
$X$ contains $n-1$ generators, $s_i$. This is a restriction of the Artin 
braid group, resp. its group algebra, by asserting additionally the 
relation ${\bf R_1}$. The projection operator onto the alternating part 
reads 
\begin{align}
\pi_\wedge &= \frac{1}{n!}\sum_{\text{red. words}} (-1)^{\text{length}(w)} w
\end{align} 
where $w$ runs in the set of $n!$ reduced words. For $S_3$ we find
\begin{align}
\pi_\wedge &= \frac{1}{3!}(1-s_1-s_2+s_1s_2+s_2s_1-s_1s_2s_1) .
\end{align}
A slight generalization of this setting is to allow a quadratic relation 
for the transposition which leads to the Hecke algebra
\begin{align}
H_n &=\; <X,\{{\bf R_1},{\bf R_2},{\bf R_3}\}> \nn
{\bf R_1} &:  \tau_1^2=a\,\tau + b \nn
\end{align}
where ${\bf R_2}$ and ${\bf R_3}$ are still the braid relations. Since the
cardinality of the set $Y=\{\text{red. words}\}$ which is generated from 
the generators $\tau_i\in X$ does not change, one proceeds as above, but 
has to take care of the additional parameters. Let $a=(1-q)$ and $b=q$, 
one ends up with a projection operator \cite{ablamowicz:fauser:2000b}
\begin{align}
\pi_{\wedge_q} &= \frac{1-\tau_1-\tau_2+\tau_1\tau_2
+\tau_2\tau_1-\tau_1\tau_2\tau_2}{(1+q+q^2)(1+q)} .
\end{align}
It is a remarkable fact, that these generators can be found also in 
an undeformed Clifford algebra if it has a carefully chosen 
non-symmetric bilinear form 
\cite{fauser:1999a,fauser:1999b,fauser:2000d,ablamowicz:fauser:2000b}.

Woronowicz showed that factoring the tensor algebra by such deformed
switch generators yields in a functorial way a $q$-deformed exterior 
algebra
\begin{align}
{\bigwedge}_q V &= \frac{T(V)}{{\cal I}_{\pi_q}}\,=\, T(V) \mod \pi_q .
\end{align}
It should be noted that the relations for such algebras look quite 
different, involving $q$s. Moreover, the parameter $q$ has to be treated
as a formal variable and deformed exterior algebras have to be built 
over $\openk[[q]]$.

Oziewicz`s idea was to study non-grade preserving isomorphisms $j$ of 
$T(V)$ and their projection under an {\it ungraded\/} switch onto exterior
algebra. This can be displayed by the following diagram
\begin{align}
\begin{array}{c@{\hskip 2truecm}c}
\Rnode{a}{T(V)} 	& \Rnode{b}{T_j(V)} \\[8ex]
\Rnode{c}{\bigwedge V}	& \Rnode{d}{\bigwedge_j V}
\end{array}
\ncline{->}{a}{b}
\Aput{j}
\ncline{->}{b}{d}
\Aput{\pi_\wedge}
\ncline{->}{a}{c}
\Aput{\pi_\wedge}
\ncline{->}{c}{d}
\Aput{\gamma}              
\end{align}
The aim is to define the map $\gamma$ by this diagram and to study
the properties of the algebra $\bigwedge_j V$. The main and astonishing
outcome is, that if $j^2$ is a $\openZ_2$-graded mapping which respects 
a filtration
\begin{align}
j^2 &: T(V) \rightarrow T_j(V) \nn
j^2 &: T^k(V) \rightarrow T^k_j(V) \oplus T^{k-2}_J(V) \oplus
\ldots \oplus \left\{
\begin{array}{cl}
T^1_j(V) & \qquad\text{if $k$ is odd} \\
T^0_j(V) & \qquad\text{if $k$ is even} \\
\end{array}\right.
\end{align}
Oziewicz proved that $\bigwedge_j V$ is a Clifford algebra w.r.t. an 
arbitrary bilinear form induced by $j^2$. Since we have no occasion 
to follow this interesting path, the reader is invited to consult
the original work \cite{oziewicz:1997b}. We will deliver an example
which provides some evidence that the above described mechanism
works.

\noindent
{\bf Example:}
Let $a,b\in T^1(V)$ and $j^2 : V \otimes V \rightarrow T(V)$ be defined as
$j^2(a\otimes b) = a\otimes b + B_{ab}$ where $B$ is an arbitrary bilinear 
form. We compute the above given commutative diagram on these elements
\begin{align}
\begin{array}{c@{\hskip 2truecm}c}
\Rnode{a}{(a\otimes b)} 	& \Rnode{b}{a\otimes b + B_{ab}} \\[8ex]
\Rnode{c}{a\wedge b = \frac{1}{2}(a\otimes b - b \otimes a)}	
& \Rnode{d}{a\wedge b +B_{ab} = \frac{1}{2}(a\otimes b - b \otimes a)+B_{ab}}
\end{array}
\ncline{->}{a}{b}
\Aput{j^2}
\ncline{->}{b}{d}
\Aput{\pi_\wedge}
\ncline{->}{a}{c}
\Aput{\pi_\wedge}
\ncline{->}{c}{d}
\Aput{\gamma^B}              
\end{align}
If $\gamma$ is interpreted as the action of $a$ on $b$ it constitutes
a Clifford map $\gamma_a b = a\wedge b + B_{ab}$. The general case is given 
in Oziewicz \cite{oziewicz:1997b}.

Relevant to our consideration is that this construction can be interpreted 
as a product mutation or the other way around a homomorphism of 
algebras. Let $\Gamma$ be the map $\gamma$ extended to $\bigwedge V$,
we find
\begin{align}
\Gamma &: \bigwedge V \rightarrow \bigwedge_j V \cong \Cl(V,B) \nn 
\Gamma(a\wedge b) &= \Gamma(a) \circ \Gamma(b)
\end{align} 
where $\circ$ is the product of the new algebra, in our case a Clifford 
product. An analogous mechanism was used by Brouder to introduce renormalized
time-ordered products in QFT. 

A further remarkable fact is that one can discuss deformation {\it versus\/}
quantization. It might be even surprising that a Clifford algebra can be 
considered as an exterior algebra w.r.t. a different $\openZ$-grading. This
is obtained from the identification $\bigwedge_j V \cong \Cl(V,B)$. Such an 
outcome depends strongly on the properties of $j^2$. Oziewicz's method is
much more general and various algebras may be generated along this lines.
It is obvious that such a construction holds for the symmetric algebras
and Weyl algebras also.

\section{Clifford algebras by cliffordization}

Studying cliffordization is a major aspect of this treatise. We postpone
its precise elaboration to later chapters. In this section we discuss 
cliffordization in a non-technical way and try to highlight the advantages
of cliffordization and to make contact to some notions from the group theory.
This will help to recognize the fundamental nature of cliffordization not 
only in our case.

The Clifford map $\gamma_x$ introduced by Chevalley is a mapping
\begin{align}
\gamma_x &: V \times \bigwedge V \rightarrow \bigwedge V
\end{align}
and thus quite asymmetric in the structure of its factors. Stressing an 
analogy, we will call the process induced by the Clifford map as 
{\it Pieri formula\/} of Clifford algebra. In the theory of the symmetric 
group (alternating groups included) a Pieri formula allows to add a single 
box to a standard Young  tableaux and gives the result expanded into such 
standard tableaux, see e.g. \cite{fulton:1999a}. Denote a partition of 
the natural number $n$ into $k$ parts as 
$\lambda=(\lambda_1\ge\ldots\ge\lambda_k\ge0)$, with $\sum \lambda_i=n$.
Young operators can be constructed which are projection operators allowing
a decomposition of the representation space. The formulas
which allow to add one box (possibly in each row) to a Young tableau 
is a Pieri formula
\begin{align}
Y^1 \circ Y^{(\lambda_1\ldots\lambda_k)}
&=\sum_K a_K Y^{\lambda_K}
\end{align}
where $K$ runs over all partitions of the standard Young tableaux obtained
by adding the box. The crucial point is to have an explicite rule to 
calculate the coefficients $a_K$ in this expansion, the {\it branching rule\/} 
and {\it branching coefficients\/}. In the case of a Clifford map these 
coefficients emerge from the contractions and the involved bilinear form.

Recursive application of the Pieri formula allows to calculate products of 
arbitrary Young tableaux. A {\it closed\/} formula for such a product is 
called a {\it Littlewood-Richardson rule\/}. The question is, if such a 
formula can be given for Clifford algebras too. The affirmative answer was 
given by Rota and Stein \cite{rota:stein:1994a,rota:stein:1994b}.

Using the co-product which we introduced by employing shuffles of 
($r,s-r$)-splits one can give the following formula for a Clifford product 
of two reduced words
\begin{align}
u \circ v &= B^\wedge(u_{(2)},v_{(1)})\, u_{(1)}\wedge v_{(2)}
\end{align}
where $B^\wedge : \bigwedge V \times \bigwedge V \rightarrow \openk$
is the extension of $B : V \times V \rightarrow\openk$ by exponentiation.
The product is extended to $\bigwedge V$ by bilinearity. Hence we identify 
$B^\wedge(u_{(2)},v_{(1)})$ with the branching coefficients.

It would be misleading to recognize `cliffordization' as closely tied 
to `Clifford' algebras. Rota and Stein showed in Refs.
\cite{rota:stein:1994a,rota:stein:1994b} that this is in fact a general 
mechanism and that e.g. the Littlewood-Richardson rule emerges as a special 
case. Cliffordization provides a direct and computational very efficient
approach to various product formulas of deformed structures. The language of
cliffordization is that of Hopf gebras.

The above example using Young tableaux is not far away from our topic.
Representation theory of $gl_n(\openC)$ is closely related to this topic.
Therewith related irreducible representations of the symmetric group are 
called Specht modules or Schur modules if finite representations over 
$\openC$ are considered. If one fills Young diagrams not by numbers but 
by vectors, the resulting spaces are the Schubert varieties, which are 
extremely useful in algebraic geometry. Gra{\ss}mannians, flag manifolds and 
cohomological aspects can be treated along this route.

Given the variety of approaches to `Quantum Clifford Algebras' it is
clear, that we have to study cliffordization as the most general and 
promising tool for a great bunch of interesting mathematical and 
physical problems. Especially quantum field theory will benefit 
extraordinarily from cliffordization.

\section{Dotted and un-dotted bases}

It is a triviality that one can choose various bases to span the
linear space underlying an algebra. In our case, it is convenient
to use reduced words w.r.t. the wedge product $\wedge$, the Clifford 
product $\circ$ or the dotted wedge product $\dwedge$ which leads
to bases of the following type
\begin{align}
i)&&& GB = \{\Id;e_i;{e_i\wedge e_j}_{i<j};\ldots \}
&&\text{Gra{\ss}mann basis}&& \nn
ii)&&& CB = \{\Id;e_i;{e_i\circ e_j}_{i<j};\ldots \} 
&&\text{Clifford basis}&& \nn
iii)&&& dGB = \{\Id;e_i;{e_i\dwedge e_j}_{i<j};\ldots \} 
&&\text{dotted Gra{\ss}mann basis.}&&
\end{align}
We will investigate a few cases where a choice of the basis
leads to a different outcome.

\subsection{Linear forms}

A Clifford algebra comes with a unique linear form. We
have identified the scalars by the unit map $\eta : \openk \rightarrow 
\Id = V^0$. It is now convenient to introduce the inverse mapping
$\epsilon$ such that $\epsilon\circ\eta = \Id_{\openk}$
\begin{align}
&&\Rnode{a}{\openk}&& 
\Rnode{b}{\Cl(V,B)}&& 
\Rnode{c}{\openk}&& 
\ncline{->}{a}{b}
\Aput{\eta}
\ncline{->}{b}{c}
\Aput{\epsilon}
\end{align} 
Technically speaking, the linear form $\epsilon$ comes up with the 
coefficient of $\Id$ in the expansion of an element $u$ in terms
of the basis
\begin{align}
u &= u^0\Id + u^i e_i + u^{ij} e_i\wedge e_j + \ldots \nn
\epsilon(u) &= u^0
\end{align}
But, this outcome depends strongly on the basis chosen to expand
an element. Converting an expansion from one product to another
will change the value of the linear form. Let us take 
$u=u^0\Id+u^ie_i+u^{ij}e_i\wedge e_j$
and compute
\begin{align}
GB&&& u=u^0\Id+u^ie_i+u^{ij} e_i\wedge e_j &&\nn
  &&& \epsilon(u) = u^0 &&\nn
CB&&& u=u^0\Id+u^ie_i+u^{ij} e_i\circ e_j &&\nn
  &&& \phantom{u}=u^0\Id +u^ie_i+ u^{ij}(e_i\wedge e_j + B_{ij}) &&\nn
  &&& \epsilon(u) = u^0+u^{ij}B_{ij} &&\nn
dGB&&& u=u^0\Id+u^ie_i+u^{ij} e_i{\dwedge} e_j &&\nn
  &&& \phantom{u}=u^0\Id +u^ie_i+ u^{ij}(e_i\wedge e_j + F_{ij}) &&\nn
  &&& \epsilon(u) = u^0+u^{ij}F_{ij} . && 
\end{align}
It is thus convenient to introduce a unique linear form $\epsilon^\wedge$,
$\epsilon^\circ$, and $\epsilon^{\dwedge}$ for every basis. These are
different linear forms and their appearance is important in quantum 
physics. Moreover, we saw that $\Cl(V,g)$ and $\Cl(V,B)$ are isomorphic
as $\openZ_2$-graded algebras, but they possess different canonical 
linear structures. The isomorphy, mediated by the Wick theorem, is nothing
but a change in the product from the wedge to the dotted wedge and vice versa.
Of course, one could introduce a new unit map $\eta^{\dwedge}$ or $\eta^\circ$
to avoid these problems, however, we will see that there are other 
obstructions which prevent this.

\subsection{Conjugation}

The main involution of a Gra{\ss}mann or Clifford algebra was given 
as the the map $\hat{~} : V\rightarrow -V$ extended to $\bigwedge V$.
In terms of reduced words, this reads
\begin{align}
\hat{u} &= (-1)^{\text{length}(u)}\,u .
\end{align}
It is now obvious, that the length function depends also on the 
chosen expansion to a basis, if the basis is not of the same 
$\openZ_2$-grade (parity). Since we do not investigate such supersymmetric
transformations here, all bases behave similar under this involution.

\subsection{Reversion}

This changes if we come to the reversion $\tilde{~} : \bigwedge V 
\rightarrow\bigwedge V$. Reversion is an anti algebra homomorphisms
\begin{align}
(u \wedge v)\tilde{~} &= \tilde{v} \wedge \tilde{u} \nn
\tilde{~}\vert_{\Id\oplus V^1} &= \Id_{\Id\oplus V^1} .
\end{align}
This operation is quite sensible to the chosen basis, as we will exemplify
once more on the element $u=u^0\Id + u^ie_i + u^{ij}e_i\wedge e_j$.
\begin{align}
GB&&& u=u^0\Id+u^i e_i+ u^{ij} e_i\wedge e_j &&\nn
  &&& \tilde{u} = u^0\Id +u^ie_i - u^{ij} e_i\wedge e_j &&\nn
CB&&& u=u^0\Id+u^i e_i+u^{ij} e_i\circ e_j &&\nn
  &&& \phantom{u}=u^0\Id +u^i e_i+ u^{ij}(e_i\wedge e_j + B_{ij}) &&\nn
  &&& \tilde{u} = u^0\Id +u^i e_i- u^{ij}(e_i\wedge e_j - B_{ij}) &&\nn
  &&& \phantom{u} = (u^0-2u^{ij}B_{ij})\Id +u^i e_i- u^{ij} e_i\circ e_j &&\nn
dGB&&& u=u^0\Id+u^ie_i+u^{ij} e_i\dwedge e_j &&\nn
  &&& \tilde{u}= u^0\Id +u^ie_i- u^{ij}(e_i\wedge e_j - F_{ij}) &&\nn
  &&& \phantom{u}= (u^0+2u^{ij}F_{ij})\Id +u^ie_i- u^{ij}e_i\dwedge e_j .&& 
\end{align}
Since the reversion is needed to form spinor inner products,
this is an outcome of major importance. Also Clifford-Lipschitz,
{\bf pin} and {\bf spin} groups will be altered by this mechanism.
Indeed we have been able to employ this type of transformation
to study $q$-spin groups, and Hecke algebras
\cite{fauser:1999a,fauser:1999b,ablamowicz:fauser:2000b}.

Note, that one is once more able to define a reversion w.r.t. any product
which was chosen to build the reduced words or even a different one. In the
above calculation we used reversion w.r.t. the wedge product. Regarding 
the group structures coming with Clifford algebras, it might be convenient 
to use the reversion w.r.t. the Clifford product.
Reversion may be called {\it Gra{\ss}mann reversion\/}, {\it dotted
Gra{\ss}mann reversion\/} or {\it Clifford reversion\/} to indicate w.r.t.
which product it acts as antihomomorphism.

Given a reversion, say w.r.t. the Clifford product, one finds an
exterior product which is stable under reversion in that sense that
is does not pick up additional terms of a different grade. This is the
dotted wedge product. This will justify the identification of dotted wedge
products with normal-ordered operator products while the undotted wedge
will be related to time-ordered operator product. This relations was 
established in 
\cite{fauser:1996c,fauser:stumpf:1997a,fauser:1998b,fauser:1998a,fauser:2001b}.
Clifford reversion acts as:
\begin{align}
(e_i \wedge e_j)\tilde{~} &= - e_i \wedge e_j + 2F_{ij} \nn
(e_i \dwedge e_j)\tilde{~} &= - e_i \dwedge e_j .
\end{align}

%% file: graphical_calc.tex
\chapter{Graphical calculi}
 
\section{The Kuperberg graphical method}

\subsection{Origin of the method}

In 1991, Kuperberg introduced a graphical method to visualize tensorial 
equations \cite{kuperberg:1991a}. His method received some recognition, 
e.g. \cite{kauffman:1999a,kauffman:radford:1999a}, since he derived a valuable
set of lemmas and theorems in the course of calculating an invariant 
for 3-manifolds. While the first paper mentioned above deals with 
involutory Hopf algebras, the second paper \cite{kuperberg:1996a} generalized
the method to the non-involutory case.
\index{Kuperberg, Greg}
\index{Kuperberg, graphical method}
\index{algebra!tensor}

Tensors appear quite naturally at nearly any place in physics. Kuperberg's 
starting point is the theory of {\it state models\/} $M$. Such models
\index{state model}
consist of a commutative ring {\bf R\/} (usually the field $\openC$)
and a bi-partite graph $G$, the {\it connectivity graph\/}, whose vertices
are labelled as {\it atoms\/} and {\it interactions\/}; a set $S_A$ for
each atom $A$, called the {\it state set\/} of $A$; and a function
$\omega_I : A_1 \times A_2 \times \ldots \times A_n \rightarrow {\bf R}$
for each interaction $I$ (where $A_1,\ldots,A_n$ are {\it neighbours\/}
of $I$), called {\it weight function\/} or the {\it Boltzmann weights\/}
of $I$. A {\it state\/} of $M$ is a function $s$ on the atoms of $M$ such
that $s(A) \in S_A$. The {\it weight\/} $\omega(s)$ of a state $s$ is defined
as the product of the $\omega_I$'s evaluated at the state $s$ when this
product converges, and in particular when $G$ is finite. Finally, the
{\it partition function\/} $Z(M)$ is defined to be the sum of $\omega(s)$
over all states $s$ when this sum converges, and in particular when all 
state sets are finite.

There are more types of models like face type models or ice type models,
which however can be handled along the same lines \cite{baxter:1984a}.

In topology, a state model is connected to {\it knots\/} and {\it links\/},
if a projection $P$ of a link is given, one declares the arcs between 
{\it crossings\/} to be atoms and the crossing (which contains the 
information which arc crosses over and which one under the other one)
themselves to be interactions. Crossings may also be called 
{\it scatterings\/} which is derived from a particle interpretation 
of these models.
\index{scattering}
\index{crossing}
\index{knots}
\index{links}
A {\it weight function\/} on states $a,b,c,d$ may be defined as
\begin{align}
\omega(a,b,c,d) &=
t \delta(a,b)\delta(c,d) + t^{-1}\delta(a,c)\delta(b,d)
\end{align}
where $t$ is chosen so that $n=-(t^2+t^{-2})$, and $\delta(a,b)=1$ when
$a=b$ and $0$ otherwise. This state model is a 'link covariant' called
the Kauffman bracket \cite{kauffman:1991a}, which is essentially the 
Jones polynomial up to normalization. Pictographically the crossing is
written as
\index{Kauffman bracket}
\index{Jones polynomial}
\begin{align}
\pspicture[0.5](0,0)(2,2)
\psset{nodesep=5pt}
\psline{-}(0,0)(2,2)
\psline[border=4pt,bordercolor=white]{-}(2,0)(0,2)
\rput*(0,2){a}
\rput*(0,0){c}
\rput*(2,0){d}
\rput*(2,2){b}
\endpspicture
\quad&=\quad
t\quad
\pspicture[0.5](0,0)(2,2)
\psset{nodesep=5pt}
\psbezier(0,2)(0,1)(2,1)(2,2)
\psbezier(0,0)(0,1)(2,1)(2,0)
\rput*(0,2){a}
\rput*(0,0){c}
\rput*(2,0){d}
\rput*(2,2){b}
\endpspicture
\quad+t^{-1}\,\,
\pspicture[0.5](0,0)(2,2)
\psset{nodesep=5pt}
\psbezier(0,2)(1,2)(1,0)(0,0)
\psbezier(2,2)(1,2)(1,0)(2,0)
\rput*(0,2){a}
\rput*(0,0){c}
\rput*(2,0){d}
\rput*(2,2){b}
\endpspicture
\end{align}
This diagram is also an instance of a {\it skein relation\/} which
\index{skein relation} 
allows to cut knots and links into smaller and more elementary objects. 
It can easily be checked, that if dim$\,V = n$ one finds
the inverse scattering from the substitution $t\rightarrow 
t^{-1}$. Moreover this scattering is a braid, i.e. it satisfies 
the braid relations of the Artin braid group.
\index{scattering!inverse}

For our purpose, we review Kuperberg's graphical method, which is then 
afterwards compared with the method using tangles 
\cite{yetter:1990a,lyubashenko:1995b,lyubashenko:1995a}.
Moreover, we are interested in some basic results derived by Kuperberg,
e.g. Lemma 3.1. \cite{kuperberg:1991a}. A further interesting result, 
which will be discussed in a subsequent chapter, is the fact that
quantum Clifford Hopf gebras and quantum Gra{\ss}mann Hopf gebras
provide counterexamples to Kuperberg's Lemma 3.2., which has to be reformulated
to hold on {\it non-interacting\/} Hopf gebras. This analysis will allow
us to describe a distinction of {\it interacting\/} and 
\index{Hopf gebra!interacting}
\index{Hopf gebra!non-interacting}
{\it non-interacting\/} products, co-products and Hopf gebras.

\subsection{Tensor algebra}

Tensorial equations use an index notation which is common in physics 
and is mainly used in hydrodynamics, electrodynamics, special and general 
relativity. The invariant objects like {\it vectors\/} and 
{\it tensors\/} are displayed via their components w.r.t. a (commonly
not written down) basis,
\index{algebra!tensor!index notation}
\index{tensor algebra}
\begin{align}
T &= T_{\phantom{kl}ij}^{kl} e_k\otimes e_l \otimes 
\epsilon^i \otimes \epsilon^j
\quad\mapsto\quad T^{kl}_{\phantom{kl}ij}.
\end{align}
A basis is assumed, but it needs not even to be a holonomic basis on 
a manifold. A thoughtful introduction to abstract vectors and the usage
of indices can be found in Penrose Rindler 
\cite{penrose:rindler:1984a,penrose:rindler:1986a}.

Let now $V$ be a module (vector space), elements of $V$ are written
as $v^a$. The dual space $V^\ast$ of linear forms gives rise to elements
\index{dual space}
\index{linear forms}
$\omega_b$ acting on the $v^a$. In abstract index notation, $a$ and
$b$ are `placeholders of names' of a vector and a linear form 
(co-vector) respectively, we will also write the index $a,b,c,\ldots$
directly to denote a vector or co-vector there using Greek letters.
On the other hand using letters as indices of kernel symbols, e.g. $v^a$,
one assumes a vector to be a $n$-tuple of $\openk$-numbers, then the
$v^a$ are simply the components of the vector. With Kuperberg we will from 
now on assume a canonical basis in $V$ and $V^\ast$ denoted by $\{e_i\}$ 
and $\{\epsilon^j\}$. All spaces are assumed to have finite dimension, 
i.e. the index set ${\cal I}$ is of finite order $\#{\cal I} < \infty$.

Vectors are then described as tuples of coefficients, indexed by the
'names' of the basis vectors:
\begin{align}
v &\cong (v^1,\ldots,v^n).
\end{align}
The components of $v$ can be obtained by applying the canonical co-vectors
$\epsilon^i$ on $v$,
assuming the relation
\begin{align}
\epsilon^i(e_j) &= \delta^i_j ,
\end{align}
which fixes $\#{\cal I}$ linear forms $\epsilon^i \in V$. Obviously we
have
\begin{align}
\epsilon^i(v) &= \epsilon^i(v^j e_j) = v^j \delta^i_j = v^i,
\end{align}
where the Einstein summation convention is in force between upper and 
lower indices.

An endomorphism $S : V \rightarrow V$ is an element of $\linHom(V,V)
\cong V\otimes V^\ast$ and it has therefore the index structure
\begin{align}
S^a_{\phantom{a}b} &\cong S = S^a_{\phantom{a}b} e_a \otimes \epsilon^b.
\end{align}
The action of an endomorphism is translated (via summation convention) 
in this method into {\it matrix multiplication\/} $\bullet$ of the 
coefficients.
\begin{align}
S\bullet v &= (S^a_{\phantom{a}b} e_a \otimes \epsilon^b) (v^c e_c) \nn
&= S^a_{\phantom{a}b} v^c e_a \epsilon^b(e_c) \nn
&= S^a_{\phantom{a}b} v^c e_a \delta^b_c \nn
&= S^a_{\phantom{a}b} v^b e_a,
\end{align}
which reads after dropping the basis vectors as usual
\begin{align}
v^{\prime a} &= S^a_{\phantom{a}b} v^b .
\end{align}
We have to distinguish four type of maps, which are different in their
index structure:
\index{map!linear!type of}

\begin{align}
S^a_{\phantom{a}b} & &\Rnode{A}{V^{\phantom{*}}} &&& \Rnode{B}{V} 
	\ncline{->}{A}{B}	
	\Aput{S} 
	&& \text{endomorphism on $V$} \nn
T_a^{\phantom{a}b} & &\Rnode{A}{V^*} &&& \Rnode{B}{V^*}    
	\ncline{->}{A}{B} 
	\Aput{T}
	&& \text{endomorphism on $V^*$} \nn
B_{ab} & &\Rnode{A}{V^{\phantom{*}}} &&& \Rnode{B}{V^*} 
	\ncline{->}{A}{B}
	\Aput{B}
	&& \text{scalar product} \nn
C^{ab} & &\Rnode{A}{V^*} &&& \Rnode{B}{V}   
	\ncline{->}{A}{B}
	\Aput{C}
	&& \text{co-scalar product.} 
\end{align}
Note that the usual symmetry types can be established on type-changing
operations having two indices of the same type, i.e. $B$ or $C$:
\begin{align}
B_{(ab)} &= g_{ab} = \frac{1}{2}(B_{ab} + B_{ab}^T), \quad
B_{ab}^T= B_{ba} \nn
B_{[ab]} &= A_{ab} = \frac{1}{2}(B_{ab} - B_{ab}^T),
\end{align}
where we have introduced the Bach-brackets common in tensor calculus of 
general relativity. Vector (component) indices are called 
{\it contra-variant\/} and co-vector (component) indices are called 
{\it co-variant}. This notion
reflects the fact that under a (linear) change of the basis the vector 
coefficients transform in an inverse (contra $\cong$ against) way as
the basis itself. Co-vector indices transform covariant (co $\cong$ with).
Hence this notion resorts to the invariance of the vectors (tensors)
\index{indices!covariant}
\index{indices!contravariant}
\begin{align}
v &= v^a e_a = (v^a {T^{-1}}_{a}^{\phantom{a}b})
	       (T_{b}^{\phantom{b}c} e^\prime_c) \nn
&= v^{\prime b} e^\prime_b.
\end{align} 
The attentive reader will have noticed that $T^{-1}$ acts from the 
right, since it has to mimic a map
\begin{align}
&\Rnode{A}{V^\ast} \hskip 2truecm \Rnode{B}{V^\ast}&
\ncline{->}{A}{B}
\Aput{T}
\end{align}
which would have the index structure $T^{b}_{\phantom{b}a}$. This shows
directly the pitfall to look at coefficients (and tuples of coefficients)
as constituting `vectors', but see e.g. Hilbert \cite{hilbert:1968a}.
The $v^b$ are simply elements of the number field (or a commutative ring
${\bf R}$) and obey {\it no\/} vectorial transformation law at all. Tensor 
calculus, by omitting the basis, shifts in a peculiar manner the vector 
character of the object $v= v^a e_a$ to the index (position) of the component.

We introduce some more notations. A tensor is said to have {\it step $n$}
\index{tensor!step of a}
\index{tensor!type of a}
if it has $n$ indices. The terms {\it rank\/}, {\it degree} or 
{\it grade} are sometimes used also. It is said to have 
{\it type $(p,q)$} if it has $p$ contravariant and $q$ covariant indices.

Basic actions with tensors are:
\index{tensor!calculational rules}

\begin{itemize}
\item[i)]
Tensors may be added if and only if they have the same index structure
(type) including the names of these indices
\begin{align}
 A^{i_1,\ldots,i_r}_{j_1,\ldots,j_s}
+B^{i_1,\ldots,i_r}_{j_1,\ldots,j_s}
=C^{i_1,\ldots,i_r}_{j_1,\ldots,j_s}.
\end{align}
\item[ii)]
Tensors may be multiplied. The product tensor picks up all indices
in their mutual order and gets a new kernel symbol
\begin{align}
v^a w^b &= U^{ab} \nn
A^{i_1,\ldots,i_r}_{j_1,\ldots,j_s}
B^{k_1,\ldots,k_n}_{l_1,\ldots,l_n}
&=
C^{i_1,\ldots,i_r,k_1,\ldots,k_n}_{j_1,\ldots,j_s,l_1,\ldots,l_n}.
\end{align}
\item[iii)] Factor switching (transposition) is given by the map 
$\Rnode{A}{a\otimes b} \hskip 1truecm \Rnode{B}{b\otimes a}
\ncline{->}{A}{B}\Aput{T}$ of adjacent or non-adjacent indices of 
the same type. This device allows to speak about symmetry, i.e.
\begin{align}
g_{(ab)} &= g_{ba} = g^T_{ab}, 
\quad T : V^*\otimes V^* \rightarrow V^*\otimes V^*
\end{align} 
is a symmetric tensor. This map will be called {\it switch\/} if adjacent
indices (elements) are interchanged.
\item[iv)]
The canonical map from $V^*\otimes V$ into $\openk$ which is called 
{\it trace map\/} or {\it evaluation map\/} is denoted by repeated indices.
Note that this map implicitly uses the isomorphism
\begin{align}
&\Rnode{A}{V} \hskip 2truecm \Rnode{B}{V^*,}&
\ncline{->}{A}{B}
\Aput{\star} 
\end{align}
which is called Euclidean dual isomorphism by Saller \cite{saller:brevier}.
This map is needed to establish the correspondence
\begin{align}
&\Rnode{A}{e_i} \hskip 2truecm \Rnode{B}{\epsilon^i}&
\ncline{->}{A}{B}
\Aput{\star} 
\end{align}
and is assumed to be bijective allowing to establish $\star^{-1}$
\begin{align}
&\Rnode{A}{\epsilon^i} \hskip 2truecm \Rnode{B}{e_i}&
\ncline{->}{A}{B}
\Aput{\star^{-1}} 
\end{align}
{\bf Examples} are $\omega_a v^a$ which is the value of $\omega$ at 
the point $v$, $S^b_{\phantom{b}a} v^a$ is the action of the 
endomorphism $S$ on $v$, $S^a_{\phantom{a}a}$ is the {\it trace\/}
of $S$.\\
The trace or evaluation is commonly called {\it contraction\/} in 
tensor calculus.
\end{itemize}
Free (open) i.e. uncontracted indices will be called {\it boundary\/}
indices, while contracted indices are called {\it inner}.
\index{indices!boundary}
\index{indices!inner}

\subsection{Pictographical notation of tensor algebra}

Kuperberg's translation of tensor equations into a graphical language
is now as follows:
\index{Kuperberg!translation rules}
\begin{itemize}
\item[i)]
Every tensor is represented by its kernel symbol. Scalars are usually
not drawn at all, since they have no connectivity, i.e. no open or
boundary indices. In graphical terms this corresponds to in- or out-going 
arrows.
\item[ii)]
Every contravariant index is represented as an arrow pointing towards
the kernel symbol.
\item[iii)]
Every covariant index is represented as an arrow pointing away from
the kernel symbol.\\
{\bf Example:} A tensor of step 4 and type $(2,2)$ i.e. $T^{ab}_{cd}$
is iconographically represented as:
\begin{align}
\pspicture[0.5](0,0)(2,2)
\psset{nodesep=5pt}
\psline{->}(0,2)(0.75,1.25)
\psline{->}(0,0)(0.75,0.75)
\psline{->}(1.25,1.25)(1.75,1.75)
\psline{->}(1.25,0.75)(1.75,0.25)
\rput(1,1){T}
\rput*(0,2){a}
\rput*(0,0){b}
\rput*(2,0){c}
\rput*(2,2){d}
\endpspicture
\quad &\cong \quad
\pspicture[0.5](0,0)(2,2)
\psline{->}(0,2)(0.75,1.25)
\psline{->}(0,0)(0.75,0.75)
\psline{->}(1.25,1.25)(2,2)
\psline{->}(1.25,0.75)(2,0)
\rput(1,1){T}
\endpspicture
\end{align}
An endomorphism is given as $S^a_{\phantom{a}b}$ 
\begin{align}
& \rightarrow S \rightarrow .
\end{align}
A vector $v$ or a co-vector $\omega$ appears as source or sink of an arrow
\begin{align}
&v \rightarrow, \quad\quad  \rightarrow \omega .
\end{align}
\item[iv)]
Contraction of tensors translates into connecting diagrams.\\
{\bf Example:} $S^a_{\phantom{a}b}v^b$ the endomorphic product as
\begin{align}
&v \rightarrow S \rightarrow .
\end{align}
The trace $S^a_{\phantom{a}a}$ and $\Id^a_{\phantom{a}a} = \text{tr}(\Id)$
are depicted as
\begin{align}
\pspicture[0.5](0,0)(2,1)
\psset{xunit=0.5,yunit=0.5,runit=0.5}
\psbezier{->}(1,2)(2,2)(2,0)(1,0)
\psbezier{-}(1,0)(0,0)(0,2)(1,2)
\rput*(1,2){$S$}
\endpspicture
\quad &\quad
\pspicture[0.5](0,0)(2,1)
\psset{xunit=0.5,yunit=0.5,runit=0.5}
\psbezier{->}(1,2)(2,2)(2,0)(1,0)
\psbezier{-}(1,0)(0,0)(0,2)(1,2)
\endpspicture
\end{align}
where these diagrams having no in- or outgoing arrows represent scalars,
i.e. the trace of $S$ and the trace of $\Id$ which is $\dim V$.
\item[v)]
Arrows are allowed to and will cross.
\end{itemize}
If boundary arrows are not labelled there is an ambiguity in their
re-labelling. However, if we adopt the rule that external lines will be named
counter clockwise starting at the top-left arrow, and that arrows of
diagrams which will be subjected to an equality have to end at equal places,
this ambiguity is removed.\\
{\bf Examples:} A linear equation results in (l.h.s. vector equation, 
r.h.s. equation for a scalar coefficient) 
\begin{align}
v \rightarrow S \rightarrow  \;=\; v^\prime \rightarrow 
\quad \Rightarrow \quad
v \rightarrow S \rightarrow \epsilon &= v^\prime \rightarrow \epsilon
\end{align}
and the symmetry of a bilinear form translates into
\begin{align}
\pspicture[0.5](0,0)(1,2)
\psline{->}(0,2)(0.75,1.25)
\psline{->}(0,0)(0.75,0.75)
\rput(1,1){B}
\endpspicture
\quad=\quad
\pspicture[0.5](0,0)(2,2)
\psbezier{->}(0,0)(0.5,2)(1.5,2)(1.75,1.25)
\psbezier{->}(0,2)(0.5,0)(1.5,0)(1.75,0.75)
\rput(2,1){B}
\endpspicture
\quad&~\quad\Rightarrow
\quad B_{ab}=B_{ba}.
\end{align}

\subsection{Some particular tensors and tensor equations}

The multiplication of an algebra can be described as tensor of rank 3
with valence $1,2$. The components $M^i_{jk}$ are called multiplication
coefficients and $\{M^i_{jk}\}$ is the multiplication table which
uniquely defines the product structure of the algebra.
\index{map!multiplication}
\index{multiplication!table}
\begin{align}
\pspicture[0.5](0,0)(2,2)
\psline{->}(0,2)(0.75,1.25)
\psline{->}(0,0)(0.75,0.75)
\psline{->}(1.25,1)(1.75,1)
\rput*(1,1){M}
\rput*(0,2){a}
\rput*(0,0){b}
\rput*(2,1){$\rho$}
\endpspicture
\quad &\cong \quad
\rho_l M^l_{ij}a^i b^j.
\end{align}
Note that $\rho$ is a co-vector while $a,b$ are vectors and the equation
holds between scalars.\\
Associativity is represented by
\index{associativity}
\begin{align}
\begin{array}{c@{\hskip 0.75truecm}c@{\hskip 0.75truecm}c@{\hskip 0.75truecm}%
c@{\hskip 0.75truecm}c@{\hskip 0.75truecm}c@{\hskip 0.75truecm}c%
c@{\hskip 0.75truecm}c}
\Rnode{a}{\XP} & & & & & & \Rnode{o}{\XP} & &\\[2ex]
& \Rnode{c}{M} & & &= & \Rnode{n}{\XP} & & \Rnode{q}{M}& \Rnode{r}{\XP}\\[2ex]
\Rnode{b}{\XP} & & \Rnode{e}{M} & \Rnode{f}{\XP}& & & \Rnode{p}{M} & &\\[2ex]
& \Rnode{d}{\XP} & & & & \Rnode{m}{\XP} & & &
\ncline{->}{a}{c}
\ncline{->}{b}{c}
\ncline{->}{c}{e}
\ncline{->}{d}{e}
\ncline{->}{e}{f}
\ncline{->}{n}{p}
\ncline{->}{m}{p}
\ncline{->}{o}{q}
\ncline{->}{p}{q}
\ncline{->}{q}{r}
\end{array}
\end{align}
which permits one to drop braces and to condense the diagrams as follows
\begin{align}
\pspicture[0.5](0,0)(5,2)
\psline{->}(0,1)(0.75,1)
\psline{->}(0,0)(0.75,0.75)
\psline{->}(1,1)(1.75,1)
\psline{->}(1,0)(1.75,0.75)
\psline[linestyle=dotted]{-}(2,1)(2.75,1)
\psline{->}(3,1)(3.75,1)
\psline{->}(3,0)(3.75,0.75)
\psline{->}(4.25,1)(5,1)
\rput*(1,1){M}
\rput*(2,1){M}
\rput*(4,1){M}
\endpspicture
\quad &= \quad
\pspicture[0.5](0,0)(3,2)
\psline{->}(0,2)(0.7,1.1)
\psline{->}(0,1.5)(0.7,1.0)
\psarc[linestyle=dotted](1,1){1.0}{170}{235}
\psline{->}(0.5,0)(0.7,0.9)
\psline{->}(1.25,1)(2,1)
\rput*(1,1){M}
\endpspicture
\end{align}

Co-algebras and co-gebras can be defined by a certain categorial 
duality described in some detail in the next chapter. 
In the Kuperberg graphical calculus, this results in the reversion 
of the arrows and obviously in renaming of the structure elements. 
As an example, we can define a co-product $\Delta$ which splits up a 
single line (module) into a tensor product of two lines (modules)
\index{categorial duality}
\index{co-product}
\index{co-product!Sweedler notation}
\begin{align}
\begin{array}{c@{\hskip 0.75truecm}c@{\hskip 0.75truecm}c}
 & & \Rnode{c}{\XP} \\[2ex]
\Rnode{a}{a} & \Rnode{b}{\Delta} & \\[2ex]
 & & \Rnode{d}{\XP}
\ncline{->}{a}{b}
\ncline{->}{b}{c}
\ncline{->}{b}{d}
\end{array}
\quad&\Leftrightarrow\quad
\Delta(a) \,=\, \sum_{(a)} a_{(1)} \otimes a_{(2)}.
\end{align}
We have employed the so called Sweedler notation \cite{sweedler:1969a} 
to indicate the elements of the co-product in the first and second 
tensor slot. As a rule for translation, one might associate the terms 
$a_{(i)}$ in counter clockwise order of arrows from left to right
in the tensor products. In our case $a_{(1)}$ would be the lower
outgoing arrow etc. Internal names (indices) have no particular meaning 
at all. The co-product is also called diagonalization, since some well 
recognized co-products, e.g. that of groups, simply double the entry 
element $\Delta(a) = a \otimes a$. The associativity of co-products,
i.e. co-associativity, is derived from the associativity if we
replace 
\index{diagonalization}
$\rightarrow\,\Leftrightarrow\,\leftarrow$ and $M\,\Leftrightarrow\,\Delta$
and it reads:
\begin{align}
\label{eqn:4-31}
\begin{array}{%
c@{\hskip 0.75truecm}c@{\hskip 0.75truecm}c@{\hskip 0.75truecm}%
c@{\hskip 0.75truecm}c@{\hskip 0.75truecm}c@{\hskip 0.75truecm}%
c@{\hskip 0.75truecm}c@{\hskip 0.75truecm}c}
 & & & \Rnode{e}{\XP}& & & & \Rnode{j}{\XP} & \\[2ex]
 & & \Rnode{c}{\Delta} & & = & \Rnode{h}{\XP}& \Rnode{i}{\Delta}& & 
\Rnode{l}{\XP} \\[2ex]
\Rnode{a}{\XP} & \Rnode{b}{\Delta}& & \Rnode{f}{\XP}& & & & 
\Rnode{k}{\Delta}& \\[2ex]
 & & \Rnode{d}{\XP}& & & & & & \Rnode{m}{\XP}
\ncline{->}{a}{b}
\ncline{->}{b}{c}
\ncline{->}{b}{d}
\ncline{->}{c}{e}
\ncline{->}{c}{f}
\ncline{->}{h}{i}
\ncline{->}{i}{j}
\ncline{->}{i}{k}
\ncline{->}{k}{l}
\ncline{->}{k}{m}
\end{array}
\end{align}
If a product has an unit, this is pictographically represented as
\index{unit}
\begin{align}
\label{eqn:4-32}
\begin{array}{%
c@{\hskip 0.75truecm}c@{\hskip 0.75truecm}c@{\hskip 0.75truecm}%
c@{\hskip 0.75truecm}c@{\hskip 0.75truecm}c@{\hskip 0.75truecm}%
c@{\hskip 0.75truecm}c@{\hskip 0.75truecm}c@{\hskip 1.5truecm}c}
\Rnode{a}{\eta} & & & & \Rnode{f}{\XP}& & & & & \\[2ex]
 &\Rnode{c}{M} & \Rnode{d}{\XP}& = & & \Rnode{h}{M}& \Rnode{i}{\XP}&= & \Rnode{l}{\XP} & \Rnode{m}{\XP} \\[2ex]
\Rnode{b}{\XP}    & & & & \Rnode{g}{\eta}& & & & &
\ncline{->}{a}{c}
\ncline{->}{b}{c}
\ncline{->}{c}{d}
\ncline{->}{f}{h}
\ncline{->}{g}{h}
\ncline{->}{h}{i}
\ncline{->}{l}{m}
\end{array}
\end{align}
and dualizing yields the definition of a counit $\epsilon$
\index{counit}
\begin{align}
\label{eqn:4-33}
\begin{array}{%
c@{\hskip 0.75truecm}c@{\hskip 0.75truecm}c@{\hskip 0.75truecm}%
c@{\hskip 0.75truecm}c@{\hskip 0.75truecm}c@{\hskip 0.75truecm}%
c@{\hskip 0.75truecm}c@{\hskip 0.75truecm}c@{\hskip 1.5truecm}c}
 & & \Rnode{c}{\epsilon} & & & & \Rnode{h}{\XP} & & & \\[2ex]
\Rnode{a}{\XP} & \Rnode{b}{\Delta}& & = & \Rnode{f}{\XP}& \Rnode{g}{\Delta}& & = & \Rnode{l}{\XP} & \Rnode{m}{\XP}\\[2ex]
 & & \Rnode{d}{\XP}        & & & & \Rnode{i}{\epsilon} & & &
\ncline{->}{a}{b}
\ncline{->}{b}{c}
\ncline{->}{b}{d}
\ncline{->}{f}{g}
\ncline{->}{g}{h}
\ncline{->}{g}{i}
\ncline{->}{l}{m}
\end{array}
\end{align}
A further prominent structure element is the antipode $S$, an endomorphism,
which, if it exists, fulfils the following defining relations
\index{antipode}
\begin{align}
\label{eqn:4-34}
\begin{array}{%
c@{\hskip 0.75truecm}c@{\hskip 0.75truecm}c@{\hskip 0.75truecm}%
c@{\hskip 0.75truecm}c@{\hskip 0.25truecm}c@{\hskip 0.25truecm}%
c@{\hskip 0.75truecm}c@{\hskip 0.25truecm}c@{\hskip 0.75truecm}%
c@{\hskip 0.25truecm}c@{\hskip 0.25truecm}c@{\hskip 0.75truecm}%
c@{\hskip 0.75truecm}c@{\hskip 0.75truecm}c@{\hskip 0.75truecm}c}
 & & \Rnode{c}{S}& & & & & & & & & & & & & \\[2ex]
\Rnode{a}{\XP} & \Rnode{b}{\Delta} & & \Rnode{d}{M} & \Rnode{e}{\XP}
& = & 
\Rnode{g}{\XP}& \Rnode{h}{\epsilon}& \Rnode{i}{\eta} & \Rnode{j}{\XP}
& = & 
\Rnode{l}{\XP} & \Rnode{m}{\Delta} & & \Rnode{o}{M} & \Rnode{p}{\XP} \\[2ex]
 & & & & & & & & & & & & &\Rnode{n}{S} & & 
\ncline{->}{a}{b}
\ncline{->}{b}{c}
\ncline{->}{c}{d}
\ncline{->}{d}{e}
\nccurve[angleA=-45,angleB=225]{->}{b}{d}
\ncline{->}{g}{h}
\ncline{->}{i}{j}
\ncline{->}{l}{m}
\ncline{->}{m}{n}
\ncline{->}{n}{o}
\ncline{->}{o}{p}
\nccurve[angleA=45,angleB=135]{->}{m}{o}
\end{array}
\end{align}
In other words, $S$ is the left and right {\it convolutive\/} inverse of the
linear identity $\Id \in \End V$. Having the Hopf algebra as our goal in mind, 
it is convenient to introduce a further relation, which is the graphical 
counterpart of the fact that $M$ is asserted to be a co-algebra morphism 
and $\Delta$ is an algebra morphism. This can be postulated as an axiom,
which will result in a constraint to choose $M$ and $\Delta$ or can be 
checked to be true or false for an arbitrary given pair of structure
tensors $(M,\Delta)$. 
\index{multiplication!as co-algebra morphism}
\index{co-multiplication!as algebra morphism}
\begin{align}
\label{eqn:4-35}
\begin{array}{%
c@{\hskip 0.75truecm}c@{\hskip 0.75truecm}c@{\hskip 0.75truecm}%
c@{\hskip 0.75truecm}c@{\hskip 0.75truecm}c@{\hskip 0.75truecm}%
c@{\hskip 0.75truecm}c@{\hskip 0.75truecm}c}
\Rnode{a}{\XP} & & &\Rnode{e}{\XP} & &\Rnode{h}{\XP} & \Rnode{j}{\Delta}& \Rnode{l}{M}&\Rnode{n}{\XP} \\[2ex]
 & \Rnode{c}{M}& \Rnode{d}{\Delta} & & = & & & & \\[2ex]
\Rnode{b}{\XP} & & &\Rnode{f}{\XP} & & \Rnode{i}{\XP}& \Rnode{k}{\Delta}& \Rnode{m}{M}& \Rnode{o}{\XP}
\ncline{->}{a}{c}
\ncline{->}{b}{c}
\ncline{->}{c}{d}
\ncline{->}{d}{e}
\ncline{->}{d}{f}
\ncline{->}{h}{j}
\ncline{->}{j}{l}
\ncline{->}{l}{n}
\ncline{->}{i}{k}
\ncline{->}{k}{m}
\ncline{->}{m}{o}
\ncline{->}{k}{l}
\ncline{->}{j}{m}
\end{array}
\end{align}
Note that a crossing of arrows occurs in the right hand side of this 
equation. This will be allowed and it is described in more detail below.

\index{Hopf algebra}
\index{algebra!Hopf}
The notion of a Hopf algebra is equivalent to the assertion of the 
relations in Eqns. \ref{eqn:4-31} to \ref{eqn:4-35}. Indeed, some rules 
to manipulate the Kuperberg graphs have to be given since the mere 
notion of {\it calculating\/} with them means that there are rules to 
manipulate them.

Hopf algebras have been derived from some properties of groups and 
group manifolds \cite{hopf:1941a,milnor:moore:1965a}. For groups
one obtains the structure relations
\index{Hopf algebra!of groups}
\begin{align}
\begin{array}{%
c@{\hskip 0.75truecm}c@{\hskip 0.75truecm}c@{\hskip 0.25truecm}%
c@{\hskip 0.25truecm}c@{\hskip 0.75truecm}c@{\hskip 0.75truecm}%
c@{\hskip 0.75truecm}c@{\hskip 0.75truecm}c@{\hskip 0.75truecm}%
c@{\hskip 0.25truecm}c@{\hskip 0.25truecm}c@{\hskip 0.75truecm}%
c}
\Rnode{a}{g} & & & & & & & & &\Rnode{k}{\XP} & & \Rnode{n}{g}& \Rnode{p}{\XP}\\[2ex]
 & \Rnode{c}{M}& \Rnode{d}{\XP}&=& \Rnode{f}{gh} & \Rnode{g}{\XP}& & \Rnode{i}{g}& \Rnode{j}{\Delta}& & = & & \\[2ex]
\Rnode{b}{h} & & & & & & & & &\Rnode{l}{\XP} & &\Rnode{o}{g} &\Rnode{q}{\XP}\\[3ex] 
\ncline{->}{a}{c}
\ncline{->}{b}{c}
\ncline{->}{c}{d}
\ncline{->}{f}{g}
\ncline{->}{i}{j}
\ncline{->}{j}{k}
\ncline{->}{j}{l}
\ncline{->}{n}{p}
\ncline{->}{o}{q}
\Rnode{a}{~g} & \Rnode{b}{S}& \Rnode{c}{\XP}&=& \Rnode{e}{g^{-1}} & \Rnode{f}{\XP} & & & & & & &
\ncline{->}{a}{b}
\ncline{->}{b}{c}
\ncline{->}{e}{f}
\end{array}
\end{align}
having the nice property that Hopf algebra morphisms induce as restrictions
morphisms of the underlying groups.

\subsection{Duality}

There are different ways to introduce new Hopf algebra structure
tensors related to the given one. Indeed, we have already used the fact that 
we can exchange by categorial duality $\Delta \Leftrightarrow M$, 
$\epsilon \Leftrightarrow \eta$, $\rightarrow\,\Leftrightarrow\,\leftarrow$ 
etc. This Hopf algebra is denoted by $H^*$ and $H$ may be denoted as 
$H_*$ for topological reasons \cite{hopf:1941a,milnor:moore:1965a}.
\index{categorial duality}

A further possibility to introduce new structure elements related with the
old ones is to introduce opposite algebras and opposite co-gebras, which
are given via
\index{algebra!opposite}
\index{multiplication!opposite}
\begin{align}
&\begin{array}{%
c@{\hskip 0.75truecm}c@{\hskip 0.75truecm}c@{\hskip 0.75truecm}%
c@{\hskip 0.75truecm}c@{\hskip 0.75truecm}c@{\hskip 0.75truecm}%
c@{\hskip 0.25truecm}c@{\hskip 0.25truecm}c@{\hskip 1.25truecm}%
c@{\hskip 1.25truecm}c}
\Rnode{a}{\XP} & & & &\Rnode{f}{\XP} & & & &\Rnode{k}{\XP} & & \\[2ex]
 & \Rnode{c}{M}& \Rnode{d}{\XP}& \Leftrightarrow& & \Rnode{h}{M^{op}} &
\Rnode{i}{\XP} & = & & \Rnode{m}{M}& \Rnode{n}{\XP}\\[2ex]
\Rnode{b}{\XP} & & & &\Rnode{g}{\XP} & & & &\Rnode{l}{\XP} & & 
\ncline{->}{a}{c}
\ncline{->}{b}{c}
\ncline{->}{c}{d}
\ncline{->}{f}{h}
\ncline{->}{g}{h}
\ncline{->}{h}{i}
\nccurve[angleA=-45,angleB=225]{->}{k}{m}
\nccurve[angleA=45,angleB=135]{->}{l}{m}
\ncline{->}{m}{n}
\end{array}
\\
&\begin{array}{%
c@{\hskip 0.75truecm}c@{\hskip 0.75truecm}c@{\hskip 0.75truecm}%
c@{\hskip 0.75truecm}c@{\hskip 0.75truecm}c@{\hskip 0.75truecm}%
c@{\hskip 0.25truecm}c@{\hskip 0.25truecm}c@{\hskip 1.25truecm}%
c@{\hskip 1.25truecm}c}
 & & \Rnode{a}{\XP} & & & &\Rnode{f}{\XP} & & & &\Rnode{k}{\XP} \\[2ex]
\Rnode{d}{\XP} & \Rnode{c}{\Delta}& & \Leftrightarrow& \Rnode{i}{\XP}& \Rnode{h}{\Delta^{op}} & & = &\Rnode{n}{\XP} & \Rnode{m}{M}& \\[2ex]
 & & \Rnode{b}{\XP} & & & &\Rnode{g}{\XP} & & & &\Rnode{l}{\XP} 
\ncline{->}{c}{a}
\ncline{->}{c}{b}
\ncline{->}{d}{c}
\ncline{->}{h}{f}
\ncline{->}{h}{g}
\ncline{->}{i}{h}
\nccurve[angleA=-45,angleB=225]{->}{m}{k}
\nccurve[angleA=45,angleB=135]{->}{m}{l}
\ncline{->}{n}{m}
\end{array}
\\
&\begin{array}{%
c@{\hskip 0.75truecm}c@{\hskip 0.75truecm}c@{\hskip 0.75truecm}%
c@{\hskip 0.75truecm}c@{\hskip 0.75truecm}c@{\hskip 0.75truecm}c}
\Rnode{a}{\XP} &\Rnode{b}{S} &\Rnode{c}{\XP} & \Leftrightarrow& 
\Rnode{d}{\XP}& \Rnode{e}{S^{op}}&\Rnode{f}{\XP}
\ncline{->}{a}{b}
\ncline{->}{b}{c}
\ncline{->}{d}{e}
\ncline{->}{e}{f}
\end{array}
\end{align}
If the opposite algebra and co-algebra structures have units $\eta^{op}$
and counits $\epsilon^{op}$ depends on the crossing of arrows,
but are taken usually to be the same units $\eta$ and counits $\epsilon$
as in the untwisted case, which is an assumption about the crossing.

\subsection{Kuperberg's Lemma 3.1.}
 
In 1991 \cite{kuperberg:1991a} Kuperberg introduced so called 
{\it ladder\/} diagrams, according to their graphical representation, 
which are extremely useful for proving further identities.
A further Lemma from this important paper will be considered below.
\index{Kuperberg!Lemma 3.1}
\index{diagram!ladder}
\index{ladder}

\begin{lmm}[Kuperberg]
The tensors of (bi-associative, bi-unital) Hopf objects
\begin{align}
\begin{array}{%
c@{\hskip 0.75truecm}c@{\hskip 0.75truecm}c@{\hskip 1.5truecm}%
c@{\hskip 0.75truecm}c@{\hskip 0.75truecm}c@{\hskip 0.75truecm}c}
\Rnode{a}{\XP} &\Rnode{c}{M} &\Rnode{e}{\XP} & & \Rnode{g}{\XP}& \Rnode{i}{M}& \Rnode{k}{\XP}\\[2ex]
\Rnode{b}{\XP} &\Rnode{d}{\Delta} &\Rnode{f}{\XP} & & \Rnode{h}{\XP}& \Rnode{j}{\Delta}& \Rnode{l}{\XP}
\ncline{->}{a}{c}
\ncline{->}{c}{e}
\ncline{->}{b}{d}
\ncline{->}{d}{c}
\ncline{->}{d}{f}
\ncline{->}{g}{i}
\ncline{->}{i}{k}
\ncline{->}{j}{h}
\ncline{->}{j}{i}
\ncline{->}{l}{j}
\end{array}
\end{align}
when viewed as vector space endomorphisms of $H\otimes H$ and
$H\otimes H^*$, are invertible.
\end{lmm}
\noindent
{\bf Proof:}
\begin{align}
\begin{array}{%
c@{\hskip 0.5truecm}c@{\hskip 0.5truecm}c@{\hskip 0.5truecm}%
c@{\hskip 0.25truecm}c@{\hskip 0.25truecm}c@{\hskip 0.5truecm}%
c@{\hskip 0.5truecm}c@{\hskip 0.5truecm}c@{\hskip 0.5truecm}%
c@{\hskip 0.25truecm}c@{\hskip 0.25truecm}c@{\hskip 0.5truecm}%
c@{\hskip 0.5truecm}c@{\hskip 0.5truecm}c@{\hskip 0.25truecm}%
c@{\hskip 0.25truecm}c@{\hskip 1.0truecm}c}
\Rnode{a}{\XP} &\Rnode{c}{M} &\Rnode{f}{M} &\Rnode{h}{\XP} & & & &\Rnode{X}{\XP} &\Rnode{o}{M} &\Rnode{q}{\XP} & & &\Rnode{s}{\XP} &\Rnode{v}{M} &\Rnode{y}{\XP} & &\Rnode{z}{\XP} & \Rnode{B}{\XP}\\[2ex]
 &\Rnode{d}{S} & & &=& &\Rnode{k}{\Delta} &\Rnode{m}{S} &\Rnode{p}{M} & &=& &\Rnode{t}{\epsilon} & \Rnode{w}{\eta}& &=& & \\[2ex]
\Rnode{b}{\XP} &\Rnode{e}{\Delta} &\Rnode{g}{\Delta} &\Rnode{i}{\XP} & &\Rnode{j}{\XP} &\Rnode{l}{\Delta} &\Rnode{n}{\XP} & & & &\Rnode{r}{\XP} &\Rnode{u}{\Delta} &\Rnode{x}{\XP} & & & \Rnode{A}{\XP}& \Rnode{C}{\XP}
\ncline{->}{a}{c}
\ncline{->}{c}{f}
\ncline{->}{f}{h}
\ncline{->}{b}{e}
\ncline{->}{e}{d}
\ncline{->}{d}{c}
\ncline{->}{e}{g}
\ncline{->}{g}{f}
\ncline{->}{g}{i}
\ncline{->}{j}{l}
\ncline{->}{l}{n}
\ncline{->}{l}{k}
\ncline{->}{k}{m}
\nccurve[angleA=-45,angleB=225]{->}{k}{p}
\ncline{->}{m}{p}
\ncline{->}{p}{o}
\ncline{->}{X}{o}
\ncline{->}{o}{q}
\ncline{->}{r}{u}
\ncline{->}{u}{x}
\ncline{->}{u}{t}
\ncline{->}{s}{v}
\ncline{->}{w}{v}
\ncline{->}{v}{y}
\ncline{->}{z}{B}
\ncline{->}{A}{C}
\end{array}
\\[2ex]
\begin{array}{%
c@{\hskip 0.5truecm}c@{\hskip 0.5truecm}c@{\hskip 0.5truecm}%
c@{\hskip 0.25truecm}c@{\hskip 0.25truecm}c@{\hskip 0.5truecm}%
c@{\hskip 0.5truecm}c@{\hskip 0.5truecm}c@{\hskip 0.5truecm}%
c@{\hskip 0.25truecm}c@{\hskip 0.25truecm}c@{\hskip 0.5truecm}%
c@{\hskip 0.5truecm}c@{\hskip 0.5truecm}c@{\hskip 0.25truecm}%
c@{\hskip 0.25truecm}c@{\hskip 1.0truecm}c}
\Rnode{a}{\XP} &\Rnode{c}{M} &\Rnode{f}{M} &\Rnode{h}{\XP} & & & &\Rnode{X}{\XP} &\Rnode{o}{M} &\Rnode{q}{\XP} & & &\Rnode{s}{\XP} &\Rnode{v}{M} &\Rnode{y}{\XP} & &\Rnode{z}{\XP} & \Rnode{B}{\XP}\\[2ex]
 &\Rnode{d}{S} & & &=& &\Rnode{k}{\Delta} &\Rnode{m}{S} &\Rnode{p}{M} & &=& &\Rnode{t}{\epsilon} & \Rnode{w}{\eta}& &=& & \\[2ex]
\Rnode{b}{\XP} &\Rnode{e}{\Delta} &\Rnode{g}{\Delta} &\Rnode{i}{\XP} & &\Rnode{j}{\XP} &\Rnode{l}{\Delta} &\Rnode{n}{\XP} & & & &\Rnode{r}{\XP} &\Rnode{u}{\Delta} &\Rnode{x}{\XP} & & & \Rnode{A}{\XP}& \Rnode{C}{\XP}
\ncline{->}{a}{c}
\ncline{->}{c}{f}
\ncline{->}{f}{h}
\ncline{->}{e}{b}
\ncline{->}{e}{d}
\ncline{->}{d}{c}
\ncline{->}{g}{e}
\ncline{->}{g}{f}
\ncline{->}{i}{g}
\ncline{->}{l}{j}
\ncline{->}{n}{l}
\ncline{->}{l}{k}
\ncline{->}{k}{m}
\nccurve[angleA=-45,angleB=225]{->}{k}{p}
\ncline{->}{m}{p}
\ncline{->}{p}{o}
\ncline{->}{X}{o}
\ncline{->}{o}{q}
\ncline{->}{u}{r}
\ncline{->}{x}{u}
\ncline{->}{u}{t}
\ncline{->}{s}{v}
\ncline{->}{w}{v}
\ncline{->}{v}{y}
\ncline{->}{z}{B}
\ncline{->}{C}{A}
\end{array}
\end{align}
The first equality is due to associativity, the second holds because
of the antipode axiom while the third reflects the unit and counit
properties. 

\section{Commutative diagrams versus tangles}

\subsection{Definitions}

We do not intend to go into details of category theory, hence the interested 
reader may consult e.g. Mac~Lane \cite{maclane:1971a}, whom we follow
\index{Mac Lane, Saunders}
in our presentation. Since some notions of category theory are, however, 
frequently used in physics, we want to give definitions for the most 
frequently used terms which are also freely used in this work. Especially 
the literature defining Hopf algebras is full of {\it commutative diagrams\/},
and uses the notion of {\it categorial duality\/}, i.e. reversing of arrows 
in diagrams. At the same time, we introduce {\it tangles\/}. Tangles will 
be seen to be opposite (might be also called dual, but should not be 
confused with categorial duality) to commutative diagrams in a certain sense.
\index{commutative!diagram}
\index{diagram!commutative}
\index{tangles}

Categories and functors are most often described by graphical
methods. It is therefore appropriate to define a {\it metagraph\/}
which consists of {\it objects\/} $a,b,c,\ldots$ and {\it arrows\/}
$f,g,h,\ldots$. The arrows depict morphisms of some structure of
the objects. Every arrow has a {\it domain\/} or {\it source\/}
and a {\it codomain\/} or {\it target}. We have thus
\index{arrows}
\index{objects}

\begin{align}
&a = \text{dom}\, f, \quad\quad b = \text{codom}\,f.
\end{align}
A morphism $f$ can be graphically represented in two ways:
\begin{align}
&f~:~a\rightarrow b &&\text{or}&&\Rnode{a}{a}\hskip 2truecm \Rnode{b}{b.}
\ncline{->}{a}{b}
\Aput{f}
\end{align}

\noindent
A finite graph may be composed from such objects and arrows. A 
{\it metacategory\/} is a metagraph with two additional operations.
The {\it identity\/} which assigns to each object the morphism
$\Id_a ~:~ a \rightarrow a$ and the composition, which assigns to each 
pair of morphisms $f,g$ having $\text{dom}\,f = \text{codom}\,g$
the {\it composite\/} morphism $h = g\circ f$. One finds that
$\text{dom}\,h = \text{dom}\,f$ and $\text{codom}\,h = \text{codom}\,g$.
This operation can be most clearly displayed by a
{\it commutative diagram\/} which relates the arrows as follows
\begin{align}
\label{eqn-CD}
\begin{array}{%
c@{\hskip 2truecm}c@{\hskip 2truecm}c}
 & \Rnode{b}{b} & \\[8ex]
\Rnode{a}{a} && \Rnode{c}{c}
\ncline{->}{a}{b}
\Aput{f}
\ncline{->}{b}{c}
\Aput{g}
\ncline{->}{a}{c}
\Bput{h=g\circ f}
\end{array}
\end{align}
\medskip

\noindent
The diagrammatic description contains full information about
all arrows, their domains and co-domains, and objects
involved in the transformations. The diagram is called {\it commutative\/}
since we end up with the result $c$ if we take either route
$g\circ f$ or $h$. The composition of morphisms, i.e. arrows,
is taken to be {\it associative\/}, i.e.

\begin{align}
\Rnode{a}{a}\hskip 1truecm \Rnode{b}{b}\hskip 1truecm
\Rnode{c}{c}\hskip 1truecm \Rnode{d}{d} 
\ncline{->}{a}{b}
\Aput{f}
\ncline{->}{b}{c}
\Aput{g}
\ncline{->}{c}{d}
\Aput{h}
\nn
h\circ(g\circ f) = (h\circ g)\circ f.
\end{align} 
Of course this definition is restricted to the case where composition
can be performed, that is when co-domains and domains are compatible.

A {\it category\/} is the restriction of a metacategory to the case
where the objects are sets. A {\it graph\/}, also called diagram scheme,
is a set $O$ of objects and a set $A$ of arrows (morphisms) and two
functions
\index{category}
\begin{align}
\Rnode{A}{A}\hskip 1.5truecm \Rnode{O}{O.}
\ncline[offset=4pt]{->}{A}{O}
\Aput{\text{dom}}
\ncline[offset=-4pt]{->}{A}{O}
\Bput{\text{codom}}
\end{align} 

\noindent
The arrows which can be composed are elements of the set of ordered pairs
\begin{align}
A\times_O A &= \{ (g,f) \mid g,f \in A,\text{~and~dom}\,g=\text{codom}\,f\}. 
\end{align}
This is the product over $O$. A {\it category\/} is a graph having the 
two additional functions
\begin{align}
&\Rnode{O}{O}\hskip 1.5truecm \Rnode{A}{A}, &&\Rnode{AA}{A\times_O A}
\hskip 1.5truecm\Rnode{A2}{A}, 
\ncline{->}{O}{A}\Aput{\Id}
\ncline{->}{AA}{A2}\Aput{\circ}
\nn
&\Rnode{O}{c}\hskip 1.5truecm \Rnode{A}{\Id_c}, &&\Rnode{AA}{g\times f}
\hskip 1.5truecm\Rnode{A2}{h = g\circ f}, 
\ncline{->}{O}{A}
\ncline{->}{AA}{A2}
\end{align}
with compatible domains and co-domains and assuming associativity. A 
further notion in the category $C$ is
\begin{align}
&\Hom(b,c) \,=\, \{ f\mid f \text{~in~} C, \text{~dom}\, f=b,
\text{~codom}\,f=c\}
\end{align}
which is equivalent to the set of arrows of $C$. A generalization 
of these notions is possible, see Oziewicz \cite{oziewicz:2001a}, 
where one finds $n$-categories, sketches and operads related to 
ideas we are using here. The idea is to study graphs like
\begin{align}
\ldots &\Rnode{C}{\XP}\hskip 1.5truecm\Rnode{B}{B}
  \hskip 1.5truecm\Rnode{A}{A}\hskip 1.5truecm \Rnode{O}{O}    
\ncline[offset=4pt]{->}{A}{O}
\Aput{\text{dom}}
\ncline[offset=-4pt]{->}{A}{O}
\Bput{\text{codom}}
\ncline[offset=4pt]{->}{C}{B}
\ncline[offset=-4pt]{->}{C}{B}
\ncline[offset=4pt]{->}{B}{A}
\ncline[offset=-4pt]{->}{B}{A}
\end{align} 
where one examines morphisms of arrows and morphisms of morphisms of $\ldots$
of arrows. In fact, it is possible to define a category consisting only
of arrows and morphisms of arrows.

We have to remark, that the sets used in categories here will be so called
{\it small sets\/} which are tame sets in the sense of set theory, i.e.
one disallows pathological sets as the 'set of all sets' etc. to avoid 
antinomies. Such categories are called {\it small categories}.

Examples of categories are among the following:
\begin{itemize}
\item[{\bf 0}] the empty category, no objects no arrows.
\item[{\bf 1}] the identity category, one object, one (identity) arrow.
\item[{\bf Set}] {\it Set:\/} Objects are small sets, arrows are 
functions between them.
\item[{\bf Mon}] {\it Monoid:\/} A monoid can be addressed itself as a
category with one object and arrows, among them the identity arrow. 
The category {\bf Mon} is that where the objects are small monoids and 
the arrows are morphisms of monoids.
\item[$\openk$-{\bf Mod}] {\it $\openk$-Modules:} small modules over the 
commutative ring $\openk$.
\item[$A$-{\bf Mod}] {\it $A$-modules:} small left $A$-modules and 
morphisms of modules.
\item[{\bf Grp}] {\it Groups:\/} Objects are small groups, arrows are 
morphisms of groups. A group itself is a category with one object where 
all morphisms are isomorphisms.
\item[{\bf Top}] {\it Topological spaces:\/} small topological spaces 
and continuous maps.
\end{itemize}

\index{functor}
A {\it functor\/} is a morphism of categories. A functor consists therefore
of two morphisms, a morphism of objects and a morphism of arrows.
This opens two possibilities, either the directions of the arrows 
before and after the transformation are 'parallel', i.e. the direction
is not changed, such functors are called covariant. Otherwise
the functor reverses the direction of arrows and is called 
contravariant. We will not have much opportunity to enter this topic, but
the notion of opposite algebras, products etc. touches this fact.

Some more notation is appropriate to make contact with the current literature. 
If an arrow $f~:~a\rightarrow b$ in a category {\bf C} is invertible,
i.e. there exists $f^{-1}~:~b\rightarrow a$ with $f\circ f^{-1} = \Id =
f^{-1}\circ f$ including the domain and codomain restriction, then
$a$ and $b$ are {\it isomorphic\/} written as $a\cong b$. An arrow 
$m~:~a\rightarrow b$ is called {\it monic\/}, and the corresponding morphism 
is a {\it monomorphism\/} in {\bf C}, if for any two parallel arrows 
$f,g~:~d\rightarrow a$ the equality $f\circ m = g\circ m$ implies $f = g$. 
That is there holds a right cancellation law. An arrow $h$ is called 
{\it epi}, and its corresponding morphism is an {\it epimorphism\/} if there 
holds analogously a left cancellation law. In {\bf Set} and in {\bf Grp} 
monic arrows, i.e. monomorphisms, are injections and epi arrows, i.e. 
epimorphisms, are surjections.
\index{morphism!epi}
\index{morphism!mono}

We will have below occasion to see that graphical calculus is dangerous
in that sense that it is difficult to keep track of the epi- and mono 
morphism properties of arrows. The discussion of Kuperberg's Lemma 3.2
will show how graphical calculus can be misleading. One pays a price for 
a nice representation with a certain peculiarity in calculating with them.

Note, that in a commutative diagram as given in Eqn. \ref{eqn-CD}
the objects are localized and the morphisms are given as arrows. This
goes directly with intuition where one expects objects to be 'solid'
or 'material' as indicated by the name and arrows are seen to by 
'dynamical', 'operating' or 'transforming' entities. That this is not
a necessity can be shown directly by categories where objects are
themselves morphisms. One is therefore able to develop a graphical
notation which reverses this assumption and depicts the objects as
lines and the morphisms as points on this lines or as other localized
operations as connecting or splitting lines. Such diagrams will
be called {\it tangles\/} 
\cite{yetter:1990a,lyubashenko:1995a,lyubashenko:1995b}. 
As an example we give the notion of a product in both representations
\begin{align}
&\begin{array}{c}
\Rnode{AA}{A\otimes A}\\[8ex]
\Rnode{A}{A}
\ncline{->}{AA}{A}
\Aput{m}
\end{array}
&&
\pspicture[0.5](0,0)(1,2.5)
\psset{linewidth=\pstlw,xunit=0.5,yunit=0.5,runit=0.5}
\psset{arrowsize=2pt 2,arrowinset=0.2} 
\psline{-}(0,5)(0,3)
\psline{-}(2,5)(2,3)
\psarc(1,3){1.0}{180}{360}
\psline{-}(1,2)(1,0)
\pscircle[linewidth=0.4pt,fillstyle=solid,fillcolor=white](1,2){0.2}   
\rput(1,2.75){m}
\endpspicture
\end{align}
The tangle is read downwards unless arrows are used to indicate which
lines have to be traversed in which direction. As commutative diagrams 
can be read for elements or for sets of objects, or even for 
categories themselves, i.e. diagrams which contain functors, this is 
true also for tangle relations. Actually a tangle can be seen as 
constituting a dual type of graph. We think that tangles can be 
intuitively understood as a {\it process\/}. A set of objects 
can move along the line of a tangle suffering certain operations 
(morphisms). This is much more a dynamical picture closely related to
physical processes where also objects like point particles or quanta
'move' around subjected to forces or quantum processes such as creation or
annihilation. The multiplication shown  above can be seen as the
annihilation of two factors into a newly created entity, the product,
of possibly different type.   

The notion of a category such allows to speak about Hopf objects, which are
elements of the objects of the category of Hopf algebras etc.

\subsection{Tangles for knot theory}

As an example we examine tangles of a special kind arising from knot 
theory and link invariants, see e.g. \cite{oziewicz:2001a}. The projection 
of 3-dimensional knots into a plane, yields a planar graph which contains 
information about the knot in three space. If additionally the information 
is conserved which string of the knot crosses over and which under the other 
w.r.t. the particular projection, the planar graph contains all topological 
information about the knot. A classification of such planar graphs would 
then classify the knots. We can state an alphabet, i.e. a basic set, 
for (graphs of) knots and links in the following way:
\index{alphabet!for knots and links}
\begin{align}
\label{eqn-TGL}
\left\{\hskip 0.5truecm
\pspicture[0.5](0,0)(1,1.5)
\psset{linewidth=\pstlw,xunit=0.5,yunit=0.5,runit=0.5}
\psset{arrowsize=2pt 2,arrowinset=0.2} 
\psline{-}(0,3)(2,1)
\psline[border=4pt,bordercolor=white]{-}(2,3)(0,1)
\rput(0,0){\Rnode{1}{2}}
\rput(2,0){\Rnode{2}{2}}
\ncline{->}{1}{2}
\endpspicture
,\quad
\pspicture[0.5](0,0)(1,1.5)
\psset{linewidth=\pstlw,xunit=0.5,yunit=0.5,runit=0.5}
\psset{arrowsize=2pt 2,arrowinset=0.2} 
\psline{-}(2,3)(0,1)
\psline[border=4pt,bordercolor=white]{-}(0,3)(2,1)
\rput(0,0){\Rnode{1}{2}}
\rput(2,0){\Rnode{2}{2}}
\ncline{->}{1}{2}
\endpspicture
,\quad
\pspicture[0.5](0,0)(1,1.5)
\psset{linewidth=\pstlw,xunit=0.5,yunit=0.5,runit=0.5}
\psset{arrowsize=2pt 2,arrowinset=0.2} 
\psbezier(0,3)(0,1)(2,1)(2,3)
\rput(0,0){\Rnode{1}{2}}
\rput(2,0){\Rnode{2}{0}}
\ncline{->}{1}{2}
\endpspicture
,\quad
\pspicture[0.5](0,0)(1,1.5)
\psset{linewidth=\pstlw,xunit=0.5,yunit=0.5,runit=0.5}
\psset{arrowsize=2pt 2,arrowinset=0.2} 
\psbezier(0,1)(0,3)(2,3)(2,1)
\rput(0,0){\Rnode{1}{0}}
\rput(2,0){\Rnode{2}{2}}
\ncline{->}{1}{2}
\endpspicture
\hskip 0.5truecm\right\}
\end{align}
Knots and links are composed from these basic words so that no
loose ends are open. However, since a knot can have different planar 
projections, one has to introduce rules which allow to transform
one such representation into another. In this case, these are the
Reidemeister moves. The first assertion is:
\index{moves!Reidemeister}
\begin{align}
\label{eqn:4-54}
\pspicture[0.5](0,0)(2,2)
\psset{linewidth=\pstlw,xunit=0.5,yunit=0.5,runit=0.5}
\psset{arrowsize=2pt 2,arrowinset=0.2} 
\psline(0,4)(0,2)
\psline(4,2)(4,0)
\psarc(1,2){1.0}{180}{360}
\psarc(3,2){1.0}{0}{180}
\endpspicture
\quad=\quad
\pspicture[0.5](0,0)(1,2)
\psset{linewidth=\pstlw,xunit=0.5,yunit=0.5,runit=0.5}
\psset{arrowsize=2pt 2,arrowinset=0.2} 
\psline(1,4)(1,0)
\endpspicture
\quad=\quad
\pspicture[0.5](0,0)(2,2)
\psset{linewidth=\pstlw,xunit=0.5,yunit=0.5,runit=0.5}
\psset{arrowsize=2pt 2,arrowinset=0.2} 
\psline(4,4)(4,2)
\psline(0,2)(0,0)
\psarc(1,2){1.0}{0}{180}
\psarc(3,2){1.0}{180}{360}
\endpspicture
\end{align} 
which states that lines can be bent at will like rubber strings.
The second move is
\begin{align}
\label{eqn:4-55}
\pspicture[0.5](0,0)(2,2)
\psset{linewidth=\pstlw,xunit=0.5,yunit=0.5,runit=0.5}
\psset{arrowsize=2pt 2,arrowinset=0.2} 
\psline(0,4)(0,3)
\psarc(2,2){1.0}{270}{90}
\psbezier(2,3)(1,3)(0,2)(0,1)
\psbezier[border=4pt,bordercolor=white](0,3)(0,2)(1,1)(2,1)
\psline(0,1)(0,0)
\endpspicture
\quad=\quad
\pspicture[0.5](0,0)(1,2)
\psset{linewidth=\pstlw,xunit=0.5,yunit=0.5,runit=0.5}
\psset{arrowsize=2pt 2,arrowinset=0.2} 
\psline(1,4)(1,0)
\endpspicture
\quad=\quad
\pspicture[0.5](0,0)(2,2)
\psset{linewidth=\pstlw,xunit=0.5,yunit=0.5,runit=0.5}
\psset{arrowsize=2pt 2,arrowinset=0.2} 
\psline(5,4)(5,3)
\psarc(3,2){1.0}{90}{270}
\psbezier(3,3)(4,3)(5,2)(5,1)
\psbezier[border=4pt,bordercolor=white](5,3)(5,2)(4,1)(3,1)
\psline(5,1)(5,0)
\endpspicture
\end{align}
and an analogous tangle equation for the second crossing, stating
that a single loop can be untwisted. The last move is a braid
relation and shows how to move a single line in the middle of two
other lines through a crossing:
\begin{align}
\label{eqn:4-56}
\pspicture[0.5](0,0)(2,3)
\psset{linewidth=\pstlw,xunit=0.5,yunit=0.5,runit=0.5}
\psset{arrowsize=2pt 2,arrowinset=0.2} 
\psbezier(2,6)(2,5)(0,5)(0,4)
\psbezier[border=4pt,bordercolor=white](0,6)(0,5)(2,5)(2,4)
\psline(4,6)(4,4)
\psbezier(4,4)(4,3)(2,3)(2,2)
\psbezier[border=4pt,bordercolor=white](2,4)(2,3)(4,3)(4,2)
\psline(0,4)(0,2)
\psbezier(2,2)(2,1)(0,1)(0,0)
\psbezier[border=4pt,bordercolor=white](0,2)(0,1)(2,1)(2,0)
\psline(4,2)(4,0)
\endpspicture
\quad&=\quad
\pspicture[0.5](0,0)(2,3)
\psset{linewidth=\pstlw,xunit=0.5,yunit=0.5,runit=0.5}
\psset{arrowsize=2pt 2,arrowinset=0.2} 
\psbezier(4,6)(4,5)(2,5)(2,4)
\psbezier[border=4pt,bordercolor=white](2,6)(2,5)(4,5)(4,4)
\psline(0,6)(0,4)
\psbezier(2,4)(2,3)(0,3)(0,2)
\psbezier[border=4pt,bordercolor=white](0,4)(0,3)(2,3)(2,2)
\psline(4,4)(4,2)
\psbezier(4,2)(4,1)(2,1)(2,0)
\psbezier[border=4pt,bordercolor=white](2,2)(2,1)(4,1)(4,0)
\psline(0,2)(0,0)
\endpspicture
\end{align}
This is in fact the braid relation of an Artin braid group
which is generated by the crossings at the position $i\in 
\{1,\ldots,n-1\}$, if $n$-strings are given.
\index{braid relation}
If one adds tangles, projections of knots or links
and to multiplies them by numbers, the set of basic tangles
as given in Eqn. \ref{eqn-TGL} and composed so that no
outgoing lines occur, then constitutes a basis of an infinite dimensional
free module. The relations asserted by the Reidemeister moves 
lead to equivalence classes of tangles which also constitute
a basis of a module. If relations between the letters of the basic
alphabet hold, such relations are called {\it skein relations\/.}
A module having skein relations on its constituting alphabet is a 
{\it skein module}.
\index{skein module}

\subsection{Tangles for convolution}

The convolution alphabet \cite{oziewicz:2001a} is given by a multiplication 
map which is of type $2 \rightarrow 1$, and a dual structure called 
{\it co-product\/}. The co-product arrises from the product by 
categorial duality. Product and co-product form together the 
convolution alphabet.
\index{alphabet!for convolution}
\begin{align}
\label{eqn-PCP}
\begin{array}{c}
\Rnode{AA}{A\otimes A}\\[8ex]
\Rnode{A}{A}
\ncline{->}{AA}{A}
\Aput{M}
\end{array}
\hskip 0.75truecm\cong\hskip 0.75truecm
\pspicture[0.5](0,0)(1,1.5)
\psset{linewidth=\pstlw,xunit=0.5,yunit=0.5,runit=0.5}
\psset{arrowsize=2pt 2,arrowinset=0.2} 
\psline{-}(0,3)(0,2)
\psline{-}(2,3)(2,2)
\psarc(1,2){1.0}{180}{360}
\psline{-}(1,1)(1,0)
\pscircle[linewidth=0.4pt,fillstyle=solid,fillcolor=white](1,1){0.2}   
\rput(1,1.75){M}
\endpspicture
\quad\quad\quad\quad
\begin{array}{c}
\Rnode{AA}{A}\\[8ex]
\Rnode{A}{A\otimes A}
\ncline{->}{AA}{A}
\Aput{\Delta}
\end{array}
\hskip 1truecm\cong\hskip 0.75truecm
\pspicture[0.5](0,0)(1,1.5)
\psset{linewidth=\pstlw,xunit=0.5,yunit=0.5,runit=0.5}
\psset{arrowsize=2pt 2,arrowinset=0.2} 
\psline{-}(0,0)(0,1)
\psline{-}(2,0)(2,1)
\psarc(1,1){1.0}{0}{180}
\psline{-}(1,2)(1,3)
\pscircle[linewidth=0.4pt,fillstyle=solid,fillcolor=white](1,2){0.2}   
\rput(1,1){$\Delta$}
\endpspicture
\end{align}
The convolution product is defined in the following way using either
commutative diagrams or tangles:
\begin{align}
\label{eqn-CONV}
\begin{array}{%
c@{\hskip 2truecm}c}
\Rnode{1}{A} & \Rnode{2}{A\otimes A} \\[8ex]
\Rnode{3}{A} & \Rnode{4}{A\otimes A}
\ncline{->}{1}{2}
\Aput{\Delta}
\ncline{->}{2}{4}
\Aput{f\otimes g}
\ncline{->}{4}{3}
\Bput{m}
\ncline{->}{1}{3}
\Bput{\star}
\end{array}
\quad\hskip 2truecm\quad
\pspicture[0.5](0,0)(1,3)
\psset{linewidth=\pstlw,xunit=0.5,yunit=0.5,runit=0.5}
\psset{arrowsize=2pt 2,arrowinset=0.2} 
\psline{-}(1,6)(1,4)
\psline{-}(1,2)(1,0)
\psarc(1,3){1.0}{0}{360}
\pscircle[linewidth=0.4pt,fillstyle=solid,fillcolor=black](0,3){0.2}   
\pscircle[linewidth=0.4pt,fillstyle=solid,fillcolor=black](2,3){0.2}   
\pscircle[linewidth=0.4pt,fillstyle=solid,fillcolor=white](1,4){0.2}   
\pscircle[linewidth=0.4pt,fillstyle=solid,fillcolor=white](1,2){0.2}   
\rput(-0.75,3){$f$}
\rput(2.75,3){$g$}
\endpspicture
\end{align}
The convolution product turns the module of endomorphisms into an 
algebra, the {\it convolution algebra\/}. In terms of algebraic 
symbolism we can write down this process as
\begin{align}
(f\star g)(x) &= m\circ(f\otimes g)\circ \Delta(x)
\,=\, \sum_{(x)} f(x_{(1)}) g(x_{(2)}).
\end{align}
\index{convolution product}
The element $x$ is from the objects of a category, $f,g$ are 
morphisms. The convolution product $\star$ is therefore from 
the morphisms of morphisms and we would enter here the case
of a 2-category, but see Oziewicz for the general case and
proper definitions \cite{oziewicz:2001a}.

The commutative diagrams and tangles of Hopf algebras and Hopf gebras 
will be given below where these objects will be defined. However,
we have seen some of them in the Kuperberg notation already, see
Eqns. \ref{eqn:4-31} to \ref{eqn:4-35}.
Products and co-products have been displayed already in Eqn. \ref{eqn-PCP}, 
while the convolution was given in Eqn. \ref{eqn-CONV}.

We should finally remark that the Kuperberg graphical notation is
a variant of the tangle notation since it localizes also the morphisms
and depicts the objects by arrows. However, there are differences,
and as we will see later, the notion of duality etc. is not
so clearly expressed in Kuperberg's representation.

%% file: hopf_al_gebra.tex
\chapter{Hopf algebras}

Hopf algebras were introduced by Heinz Hopf \cite{hopf:1941a}
\index{Hopf, Heinz}
to study topological aspects of group manifolds and their generalizations.
From this origin it is clear that Hopf algebras are closely related to 
{\it groups\/}, as we will see more clearly below, and {\it topology}.
Indeed Hopf algebras play an important role in the theory of
link invariants, knot theory and lattice models of various types.
For our purpose this fact does not directly come into play, but should
kept in mind.

The name {\it Hopf algebra\/}\index{Hopf algebra}\index{algebra!Hopf}
was, for obvious reasons, not coined by Hopf, but by Milnor and Moore 
\cite{milnor:moore:1965a}.
The usage of this term by Borel \cite{borel:1953a} for a structure
without the diagonalization (complication, diagonal map or co-product)
\index{complication}
\index{diagonalization}
\index{map!diagonal}
\index{co-product}
has to be rejected. Moreover, Milnor and Moore's Hopf algebra comes nowadays
under the name of a {\it bialgebra\/}\index{bialgebra}
\index{bialgebra!antipodal} and only antipodal, i.e. group like see below, 
such structures are named {\it Hopf algebras\/} today. Detailed 
descriptions of Hopf algebras may be found e.g. in Milnor and Moore 
\cite{milnor:moore:1965a}, the standard reference by Sweedler 
\cite{sweedler:1969a} and Abe \cite{abe:1980a}. Our presentation
is along the trail of Milnor and Moore.

We will later on contrast the notion of a {\it Hopf al-gebra\/} with
that of a {\it Hopf gebra\/}\index{Hopf gebra}\index{gebra!Hopf}
which we will derive from a {\it convolution algebra\/}
\index{convolution algebra}
\index{algebra!convolution}
following Oziewicz \cite{oziewicz:1997a,oziewicz:2001b}. This approach
seems to be mathematically more sound and provides us a
better classification and understanding of various {\it al-} 
and {\it co-}gebraic structures including Hopf gebras. Moreover,
following Bourbaki \cite{bourbaki:1989a} the naming {\it gebra}
is linguistically correct and does not abuse its Arabic origin\footnote{%
The name algebra goes back to a book of Abu Ja'far Mohammad ibn Mussa 
al Khw\=arizmi, (780--850), see
\cite{vanderwaerden:1985a,historyofmath:web}, 
named {\it Al-{\vjmath}abr w'al muq\=abala\/}, Baghdad, Iraq, where for
\index{Khwarizmi@Al Khw\=arizmi, Mohammad (780--850)} 
the first time algebraic methods have been systematically developed.
From 'al-{\vj}abr' , transcribed as 'the completion' the word 
algebra is derived. In German transcriptions one finds even a
\u{g} in place of the English \vj. Unfortunately the article {\it al\/} 
was incorporated. The Bourbaki group had already suggested to use
therefore al-gebra, but co-gebra, bi-gebra, Hopf gebra etc. which comes
to meet the usage of the Arabic language without abusing it.
From a further book of Kw\=arizmi available only in a Latin translation, 
{\it Algorithmi de numero indorum\/}, i.e. {\it Al Khw\=arizmi on the 
Hindu Art of Reckoning}, the term {\it algorithm\/} was taken.
\index{algorithm!origin of}
}. 
 
Before we start to resume the mathematical details needed for our 
treatise, we give a genealogy of al-gebras and co-algebras to display 
the further development in a diagram:
\medskip

\begin{center}
\begin{tabular}{c@{\hskip 1truecm}c@{\hskip 1truecm}c}
\rnode{1}{Al-Gebra} & & \rnode{2}{Co-Algebra} \\[6ex]
 & \rnode{3}{Bi-Algebra} & \\[8ex]
 & \rnode{4}{Hopf Algebra} &
\end{tabular}
\ncline[nodesep=6pt]{->}{1}{3}
\ncline[nodesep=6pt]{->}{2}{3}
\ncline[nodesep=6pt]{<->}{1}{2}
\Bput{Mutual Compatibility}
\Aput{duality}
\ncline[nodesep=6pt]{->}{3}{4}
\Aput*{$\exists$~Antipode}
\end{center}
\medskip

{\it Algebras\/} and {\it Co-algebras\/}\index{co-algebra} are mutually
dual structures in a certain sense as we will explain below. If certain 
compatibility laws hold, the selfdual structure is called a bi-algebra. 
If in addition an antipode exists, one obtains a Hopf algebra. 

\section{Algebras}

\subsection{Definitions}

Let $\openk$ be a commutative ring or even a field which is chosen
once and fixed thereafter. We consider finitely generated $\openk$-modules,
denoted by capital letters $A,B,C,\ldots,V,W,\ldots$. Tensor products
of two $\openk$-modules are taken over $\openk$ and will be denoted
as $A\otimes_{\openk} B \equiv A\otimes B$. These $\openk$-modules
constitute a category. $\hom(A,B)$ denotes the morphisms of $A$ into $B$
in this category, while $A,B,C,\ldots$ are the objects.

\index{kmodule@$\openk$-module}
\index{module!over $\openk$}
\index{module!graded}
A {\it graded\/} $\openk$-module is a (finite) family of $\openk$-modules
$\{ A_n \}$ where $n$ runs through the non-negative integers. $n$ is 
called the {\it degree}. If $A,B$ are graded $\openk$-modules, a graded 
morphism of graded $\openk$-modules $f~:~A\rightarrow B$ is a family of 
morphisms $\{ f_n \}$ such that $f_n~:~ A_n \rightarrow B_n$ is a morphism 
of $\openk$-modules.

If $A$ and $B$ are graded $\openk$-modules, then $A \hat\otimes B$ is
a graded $\openk$-module such that $(A\hat\otimes B)_n = \oplus_{i+j=n}
A_i \hat\otimes B_j$. Commonly the {\it hat} in $\hat{\otimes}$ is
dropped if the grading is clear. The graded tensor product is a 
particular case of a crossed product, see below. If $f~:~A \rightarrow 
A^\prime$ and $g~:~ B \rightarrow B^\prime$ are graded morphisms of 
graded $\openk$-modules then $(f \hat\otimes g) ~:~ A\,\hat\otimes\, B
\rightarrow A^\prime \,\hat\otimes\, B^\prime$ is the graded morphism
of graded $\openk$-modules such that $(f \hat\otimes g)_n \,=\,
\oplus_{i+j=n} f_i \hat \otimes g_j$. 
\index{module!morphism}

If $A$ is a graded $\openk$-module, we denote by $A^*$ the graded 
$\openk$-module such that $A^*_n = \hom(A_n,\openk)$. $A^*$ is called the
dual $\openk$-module, elements of $A^*$ are forms. If these forms
\index{module!dual}
are linear and $\openk\;$ a field, the elements of $A^*$ are called 
co-vectors. Morphisms of graded dual $\openk$-modules are defined as in 
the case of graded $\openk$-modules. The identity morphism will be denoted 
by $\Id$ or eventually by the same symbol as for the $\openk$-module:
\begin{align}
\Rnode{A}{A}\hskip 2truecm \Rnode{B}{A}
\quad&\equiv\quad 
\Rnode{C}{A}\hskip 2truecm \Rnode{D}{A.}
\ncline{->}{A}{B}
\Aput{A}
\ncline{->}{C}{D}
\Aput{\Id}
\end{align}
If necessary, $\openk$ will be considered as graded $\openk$-module
which is defined to be the $0$-module in all degrees except $0$, and
the ring (field) $\openk$ in degree 0. This definition is equivalent
to $A\,\hat\otimes\,\openk = A = \openk\,\hat\otimes\,A$, where $A$ and 
$\openk$ are (graded) $\openk$-modules.

\index{multiplication!as def. by Gra{\ss}mann}
The notion of an algebra emerges directly from the writings of
H. Gra{\ss}mann \cite{grassmann:1878a,grassmann:1862a}. Gra{\ss}mann
denotes any {\it bilinear\/}, that is left and right distributive, map
a {\it multiplication\/}, hence allowing non-associative
multiplications, but keeping linearity. In our treatise, we
will be interested in linear and associative multiplications
most of the time and will explicitly state when these assumptions
are not met. Hopf algebras are usually defined, however, as
associative, but linearity of the multiplication map may not be
assumed.

\begin{dfn}
A {\em graded unital algebra} over $\openk$ is a graded $\openk$-module
$A$ together with morphisms of graded $\openk$-modules
$m~:~A\hat\otimes A \rightarrow A$ (multiplication) and 
$\eta~:~ \openk \rightarrow A$ (unit) such that the diagrams
\index{algebra!definition of}
\index{algebra!graded}
\begin{align}
\begin{array}{c@{\hskip 2truecm}c}
\Rnode{a}{A\otimes A \otimes A} & \Rnode{b}{A \otimes A} \\[8ex]
\Rnode{c}{A\otimes A}           & \Rnode{d}{A}
\end{array}
\ncline{->}{a}{b}
\Aput{\Id \otimes m}
\ncline{->}{b}{d}
\Aput{m}
\ncline{->}{a}{c}
\Aput{m \otimes \Id}
\ncline{->}{c}{d}
\Aput{m}
\quad&&~\quad
\pspicture[0.5](0,0)(1.5,2.0)
\psset{linewidth=\pstlw,xunit=0.5,yunit=0.5,runit=0.5}
\psset{arrowsize=2pt 2,arrowinset=0.2}
\psline(0,5)(0,4)
\psline(2,5)(2,4)
\psline(3,5)(3,2)
\psarc(1,4){1}{180}{360}
\psline(1,3)(1,2)
\psarc(2,2){1}{180}{360}
\psline(2,1)(2,0)
\pscircle[linewidth=0.4pt,fillstyle=solid,fillcolor=white](1,3){0.2}
\pscircle[linewidth=0.4pt,fillstyle=solid,fillcolor=white](2,1){0.2}
\rput(1,3.75){m}
\rput(2,1.75){m}
\endpspicture
\quad=\quad
\pspicture[0.5](0,0)(1.5,2.0)
\psset{linewidth=\pstlw,xunit=0.5,yunit=0.5,runit=0.5}
\psset{arrowsize=2pt 2,arrowinset=0.2}
\psline(0,5)(0,2)
\psline(1,5)(1,4)
\psline(3,5)(3,4)
\psarc(2,4){1}{180}{360}
\psline(2,3)(2,2)
\psarc(1,2){1}{180}{360}
\psline(1,1)(1,0)
\pscircle[linewidth=0.4pt,fillstyle=solid,fillcolor=white](1,1){0.2}
\pscircle[linewidth=0.4pt,fillstyle=solid,fillcolor=white](2,3){0.2}
\rput(1,1.75){m}
\rput(2,3.75){m}
\endpspicture
\end{align}
and
\begin{align}
\label{eqn:5-3}
\begin{array}{c@{\hskip 2truecm}c@{\hskip 2truecm}c}
 & \Rnode{a}{A\otimes A} & \\[8ex]
\Rnode{b}{A\otimes\openk} & \Rnode{c}{A} & \Rnode{d}{\openk\otimes A}
\ncline{->}{a}{c}
\Aput{m}
\ncline{->}{b}{a}
\Aput{A \otimes\eta}
\ncline{->}{d}{a}
\Bput{\eta\otimes A}
\ncline{->}{c}{b}
\Bput{\approx}
\ncline{->}{c}{d}
\Aput{\approx}
\end{array}
\quad&&~\quad
\pspicture[0.5](0,0)(1,1.0)
\psset{linewidth=\pstlw,xunit=0.5,yunit=0.5,runit=0.5}
\psset{arrowsize=2pt 2,arrowinset=0.2}
\psline(0,3)(0,2)
\psarc(1,2){1.0}{180}{360}
\psline(1,1)(1,0)
\pscircle[linewidth=0.4pt,fillstyle=solid,fillcolor=white](1,1){0.2}
\pscircle[linewidth=0.4pt,fillstyle=solid,fillcolor=white](2,2){0.2}
\rput(1,1.75){m}
\rput(2,2.75){$\eta$}
\endpspicture
=
\pspicture[0.5](0,0)(1,1.0)
\psset{linewidth=\pstlw,xunit=0.5,yunit=0.5,runit=0.5}
\psset{arrowsize=2pt 2,arrowinset=0.2}
\psline(1,3)(1,0)
\endpspicture
=
\pspicture[0.5](0,0)(1,1.0)
\psset{linewidth=\pstlw,xunit=0.5,yunit=0.5,runit=0.5}
\psset{arrowsize=2pt 2,arrowinset=0.2}
\psline(2,3)(2,2)
\psarc(1,2){1.0}{180}{360}
\psline(1,1)(1,0)
\pscircle[linewidth=0.4pt,fillstyle=solid,fillcolor=white](1,1){0.2}
\pscircle[linewidth=0.4pt,fillstyle=solid,fillcolor=white](0,2){0.2}
\rput(1,1.75){m}
\rput(0,2.75){$\eta$}
\endpspicture
\end{align}
are commutative.
\end{dfn}

If the multiplication $m$ has no unit, we speak about a non-unital algebra.
A prominent example of such algebras are the Lie algebras, which however 
fail to be associative hence they are not algebras in the above defined sense.
\index{algebra!non-unital}

\begin{dfn}
The {\em graded switch\/} (crossing, scattering) 
$\hat{\tau}~:~ A\,\hat\otimes\, B \rightarrow A\,\hat\otimes\, B$
is a morphism of graded $\openk$-modules $A$ and $B$ such that 
$\hat{\tau}_n(a\,\hat\otimes\, b) = (-1)^{\partial a\partial b}\;
b\,\hat\otimes\, a$ for $a\in A_p$, $b\in B_q$ and $p+q=n$. The degree 
of a homogenous element $c$ is denoted by $\partial c$. 
\index{switch!graded}
\index{crossing}
\index{scattering}
\end{dfn}
As a tangle the crossing reads:
\begin{align}
\pspicture[0.5](0,0)(1,1)
\psset{linewidth=\pstlw,xunit=0.5,yunit=0.5,runit=0.5}
\psset{arrowsize=2pt 2,arrowinset=0.2}
\psbezier(0,2)(0,1)(2,1)(2,0)
\psbezier[border=4pt,bordercolor=white](2,2)(2,1)(0,1)(0,0)
\rput(2.5,1.0){$\hat{\tau}$}
\endpspicture
\end{align}

\noindent
If the grading is trivial, i.e. $A = A_0$ and $B=B_0$, one obtains
through this definition the usual switch, as used e.g. by Sweedler
\cite{sweedler:1969a}

\begin{dfn}
A graded unital algebra is called {\em commutative} if the diagram
\index{algebra!commutative}
\begin{align}
\begin{array}{c@{\hskip 2truecm}c}
\Rnode{a}{A\otimes A} & \\[4ex]
& \Rnode{b}{A} \\[4ex]
\Rnode{c}{A\otimes A} & 
\end{array}
\ncline{->}{a}{b}
\Aput{m}
\ncline{->}{a}{c}
\Aput{\hat{\tau}}
\ncline{->}{c}{b}
\Bput{m}
\quad&&~\quad
\pspicture[0.5](0,0)(1,2)
\psset{linewidth=\pstlw,xunit=0.5,yunit=0.5,runit=0.5}
\psset{arrowsize=2pt 2,arrowinset=0.2}
\psline(0,2.5)(0,4)
\psline(2,2.5)(2,4)
\psarc(1,2.5){1}{180}{360}
\psline(1,0)(1,1.5)
\pscircle[linewidth=0.4pt,fillstyle=solid,fillcolor=white](1,1.5){0.2}
\rput(1,2.75){m}
\endpspicture
\quad=\quad
\pspicture[0.5](0,0)(1,2)
\psset{linewidth=\pstlw,xunit=0.5,yunit=0.5,runit=0.5}
\psset{arrowsize=2pt 2,arrowinset=0.2}
\psbezier(0,4)(0,3)(2,3)(2,2)
\psbezier[border=4pt,bordercolor=white](2,4)(2,3)(0,3)(0,2)
\psarc(1,2){1}{180}{360}
\psline(1,0)(1,1)
\pscircle[linewidth=0.4pt,fillstyle=solid,fillcolor=white](1,1){0.2}
\rput(2.5,3.0){$\hat{\tau}$}
\rput(1.0,1.75){m}
\endpspicture
\end{align}
is commutative (the tangle identity holds).
\end{dfn}
We follow the convention to talk about commutativity if the standard
crossing is employed in the above tangle, the precise term would however
be `graded commutative'.

\noindent
{\bf Note:} If the graded $\openk$-module $A$ is of the form
$A=\openk\oplus A_1$ such algebras are called usually anticommutative.
If one finds $A=A_0 \oplus A_1 = (\openk\oplus A^\prime_0)
\oplus A_1$ one deals with a supersymmetric algebra, which obeys
a $\openZ_2$-grading. In the tangle diagram one sees that in a commutative 
algebra the crossing is absorbed by the product morphism $m$
\index{algebra!anticommutative}
\index{algebra!supersymmetric}
\index{algebra!graded}
\begin{align}
m &= m\circ \hat{\tau}.
\end{align}

\begin{dfn}[Crossed product]
If $A$ and $B$ are graded unital algebras over $\openk$, $A\,\hat\otimes\,B$
is turned into an algebra with multiplication $M$ and unit $\eta$
via
\index{product!crossed}
\begin{align}
M &: A\hat\otimes B \times A \hat\otimes B \rightarrow A \hat\otimes B \nn
M &= (m_A\otimes m_B)\circ(\Id \otimes \hat\tau \otimes \Id) \\
\eta &: \openk = \openk\otimes \openk \rightarrow A\,\hat\otimes\, B\nn
\eta &= \eta_A\otimes \eta_B.   
\end{align} 
\end{dfn}
\noindent
In terms of tangles this reads:
\begin{align}
\pspicture[0.5](0,0)(3,2.5)
\psset{linewidth=\pstlw,xunit=0.5,yunit=0.5,runit=0.5}
\psline(2,5)(2,4)
\psline(4,5)(4,4)
\psline(0,5)(0,2)
\psline(6,5)(6,2)
\psline{c-c}(4,4)(2,2)
\psline[border=4pt,bordercolor=white]{c-c}(2,4)(4,2)
\psarc(1,2){1.0}{180}{360}
\psarc(5,2){1.0}{180}{360}
\psline(1,1)(1,0)
\psline(5,1)(5,0)
\pscircle[linewidth=0.4pt,fillstyle=solid,fillcolor=white](1,1){0.2}
\pscircle[linewidth=0.4pt,fillstyle=solid,fillcolor=white](5,1){0.2}
\rput(1.0,1.75){m}
\rput(3.0,3.75){$\hat\tau$}
\rput(5.0,1.75){m}
\endpspicture
&&\text{~and~}&&
\pspicture[0.5](0,0)(1,1.0)
\psset{linewidth=\pstlw,xunit=0.5,yunit=0.5,runit=0.5}
\psset{arrowsize=2pt 2,arrowinset=0.2}
\psline(0,3)(0,0)
\psline(2,3)(2,0)
\pscircle[linewidth=0.4pt,fillstyle=solid,fillcolor=white](0,3){0.2}
\pscircle[linewidth=0.4pt,fillstyle=solid,fillcolor=white](2,3){0.2}
\rput(0,3.75){$\eta_A$}
\rput(2,3.75){$\eta_B$}
\endpspicture
\end{align}
A morphism of algebras $f~:~ A\rightarrow B$ is a morphism of graded
$\openk$-modules which fulfils
\index{morphism!of algebras}
\begin{align}
\label{eqn:5-9}
\begin{array}{c@{\hskip 2truecm}c}
\Rnode{a}{A\otimes A} & \Rnode{b}{A} \\[8ex]
\Rnode{c}{B\otimes B} & \Rnode{d}{B} 
\ncline{->}{a}{b}
\Aput{m_A}
\ncline{->}{c}{d}
\Aput{m_B}
\ncline{->}{a}{c}
\Aput{f\otimes f}
\ncline{->}{b}{d}
\Aput{f}
\end{array}
&&
\pspicture[0.5](0,0)(1,2.0)
\psset{linewidth=\pstlw,xunit=0.5,yunit=0.5,runit=0.5}
\psset{arrowsize=2pt 2,arrowinset=0.2}
\psline(0,4)(0,3)
\psline(2,4)(2,3)
\psarc(1,3){1}{180}{360}
\psline(1,0)(1,2)
\pscircle[linewidth=0.4pt,fillstyle=solid,fillcolor=white](1,2){0.2}
\pscircle[linewidth=0.4pt,fillstyle=solid,fillcolor=black](1,1){0.2}
\rput(1,2.75){$m_A$}
\rput(0.5,1){$f$}
\endpspicture
\quad=\quad
\pspicture[0.5](0,0)(1,2.0)
\psset{linewidth=\pstlw,xunit=0.5,yunit=0.5,runit=0.5}
\psset{arrowsize=2pt 2,arrowinset=0.2}
\psline(0,4)(0,2)
\psline(2,4)(2,2)
\psarc(1,2){1}{180}{360}
\psline(1,0)(1,1)
\pscircle[linewidth=0.4pt,fillstyle=solid,fillcolor=white](1,1){0.2}
\pscircle[linewidth=0.4pt,fillstyle=solid,fillcolor=black](0,3){0.2}
\pscircle[linewidth=0.4pt,fillstyle=solid,fillcolor=black](2,3){0.2}
\rput(1,1.75){$m_B$}
\rput(-0.5,3){$f$}
\rput( 1.5,3){$f$}
\endpspicture
\end{align}
In other words, the multiplications in $A$ and $B$ are related via
\begin{align}
f\circ m_A &= m_B\circ(f\otimes f).
\end{align}
Furthermore, an algebra $A$ is commutative --in the sense defined above--
if and only if the product $m_A ~:~ A\otimes A \rightarrow A$ is a 
morphism of algebras
\index{algebra!commutativity of}
\begin{align}
\begin{array}{c@{\hskip 3truecm}c}
\Rnode{a}{(A\otimes A)\otimes (A\otimes A)} & 
\Rnode{b}{A\otimes A} \\[8ex]
\Rnode{c}{A\otimes A} &
\Rnode{d}{A} 
\ncline{->}{a}{b}
\Aput{m_A\otimes m_A}
\ncline{->}{a}{c}
\Aput{m_A\otimes m_A}
\ncline{->}{b}{d}
\Aput{m_A}
\ncline{->}{c}{d}
\Aput{m_A}
\end{array}
&&\Rightarrow&&
\pspicture[0.5](0,0)(1,2)
\psset{linewidth=\pstlw,xunit=0.5,yunit=0.5,runit=0.5}
\psset{arrowsize=2pt 2,arrowinset=0.2}
\psline(0,2.5)(0,4)
\psline(2,2.5)(2,4)
\psarc(1,2.5){1}{180}{360}
\psline(1,0)(1,1.5)
\pscircle[linewidth=0.4pt,fillstyle=solid,fillcolor=white](1,1.5){0.2}
\rput(1,2.75){m}
\endpspicture
\quad=\quad
\pspicture[0.5](0,0)(1,2)
\psset{linewidth=\pstlw,xunit=0.5,yunit=0.5,runit=0.5}
\psset{arrowsize=2pt 2,arrowinset=0.2}
\psbezier(0,4)(0,3)(2,3)(2,2)
\psbezier[border=4pt,bordercolor=white](2,4)(2,3)(0,3)(0,2)
\psarc(1,2){1}{180}{360}
\psline(1,0)(1,1)
\pscircle[linewidth=0.4pt,fillstyle=solid,fillcolor=white](1,1){0.2}
\rput(2.5,3.0){$\hat{\tau}$}
\rput(1.0,1.75){m}
\endpspicture
\end{align}
where the conclusion follows by a simple calculation using associativity.
The crossed product is in a certain sense a canonical generalization 
of the product $m_A$ on the algebra $A\otimes A$.
\index{product!crossed}
\index{crossed product}

The {\it augmentation\/} of an algebra $A$ is a morphism of algebras
$\epsilon ~:~ A \rightarrow \openk$
\index{augmentation}
\begin{align}
\begin{array}{c@{\hskip 2truecm}c}
\Rnode{a}{A\otimes A} & \Rnode{b}{A} \\[8ex]
\Rnode{c}{\openk\otimes \openk} & \Rnode{d}{\openk} 
\ncline{->}{a}{b}
\Aput{m_A}
\ncline{->}{a}{c}
\Aput{\epsilon\otimes \epsilon}
\ncline{->}{b}{d}
\Aput{\epsilon}
\ncline{->}{c}{d}
\Aput{m_{\openk}}
\end{array}
&&
\pspicture[0.5](0,-0.5)(1,1.0)
\psset{linewidth=\pstlw,xunit=0.5,yunit=0.5,runit=0.5}
\psset{arrowsize=2pt 2,arrowinset=0.2}
\psline(0,3)(0,2)
\psline(2,3)(2,2)
\psarc(1,2){1}{180}{360}
\psline(1,1)(1,0)
\pscircle[linewidth=0.4pt,fillstyle=solid,fillcolor=white](1,0){0.2}
\pscircle[linewidth=0.4pt,fillstyle=solid,fillcolor=white](1,1){0.2}
\rput(1,1.75){m}
\rput(1,-0.75){$\epsilon_A$}
\endpspicture
\quad=\quad
\pspicture[0.5](0,-0.5)(1,1.0)
\psset{linewidth=\pstlw,xunit=0.5,yunit=0.5,runit=0.5}
\psset{arrowsize=2pt 2,arrowinset=0.2}
\psline(0,3)(0,0)
\psline(2,3)(2,0)
\pscircle[linewidth=0.4pt,fillstyle=solid,fillcolor=white](0,0){0.2}
\pscircle[linewidth=0.4pt,fillstyle=solid,fillcolor=white](2,0){0.2}
\rput(0,-0.75){$\epsilon_A$}
\rput(2,-0.75){$\epsilon_A$}
\endpspicture
\end{align}
Obviously $\epsilon_A$ is a (linear) form on the algebra. An algebra $A$ 
together with an augmentation $\epsilon_A$ is called an augmented
(or supplemented) algebra. If $A$ is an augmented algebra, we denote
\index{algebra!augmented}
\index{algebra!supplemented}
\index{augmentation!ideal}
\index{ideal!augmentation}
by $I(A)$ the kernel of $\epsilon_A~:~ A\rightarrow \openk$.
Clearly one has $I(A)_q = A_q$ for $q>0$ and $I(A)_0$ is the kernel of
$\epsilon_0 ~:~ A_0\rightarrow \openk$. The ideal $I(A)$ is called 
the augmentation ideal of $A$. As a graded $\openk$-module $A$ may
be decomposed as the direct sum
\begin{align}
A &= \img \eta \oplus \ker \epsilon
\end{align} 
or identifying $\openk = \img \eta$
\begin{align}
A &= \openk \oplus \ker \epsilon = \openk \oplus I(A),
\end{align}
where the fact was used that $\epsilon\circ \eta ~:~ \openk\rightarrow\openk$
is the identity map on $\openk$.

\noindent
{\bf Note:} The fact that $I(A)$ is an ideal in $A$, is directly
related to the property that the augmentation form $\epsilon$ is
an algebra homomorphism. One may note the similarity between the 
real and imaginary parts of complex, quaternionic and octonionic division 
algebras. An augmented algebra generalizes such a structure. We can 
therefore loosely speak about the real part $\img \eta = \openk$
and the imaginary part $\ker \epsilon = I(A)$ for connected such algebras.
\begin{dfn}
An algebra over $\openk$ is {\em connected\/} if
$\eta ~:~ \openk \rightarrow A_0$ is an isomorphism.
\end{dfn}
\index{algebra!connected}\index{augmentation} 
In that case the algebra has a unique augmentation 
$\epsilon ~:~A \rightarrow \openk$ such that $\openk \cong A_0$ and 
$\openk\overset{\eta}{\rightarrow} A_0\overset{\epsilon}{\rightarrow}\openk$
where $\epsilon_0\eta_0 = \Id_\openk$.

\subsection{$A$-modules}

$A$-modules are needed to study the representation theory of algebras
also group algebras, and of various groups connected with them. From a 
mathematical point of view every module is well suited for this purpose 
as long as an action can be defined on it. We will, however, be mostly 
interested in representations of the multiplicative (semi) group on the 
module the algebra is built over. Concerning physics, there is no reason to 
introduce a new entity (remember that 'particles' are defined
as elements of irreducible modules), i.e. a new type of module, which
is foreign to the problem at hand. Moreover, for our cases we find all 
fundamental irreducible representations as (sub)spaces of the module $A$ 
and can build up any representation from them. 

To distinguish $A$-modules and algebras, we use letters $N,M,\ldots$ for 
$A$-modules and do not assume here that $N$ is a submodule of the
algebra $N\subseteq A$.

\begin{dfn}
If $A$ is a graded unital algebra over $\openk$, a graded left $A$-module
is a graded $\openk$-module $N$ together with a morphism 
$m_N ~:~ A\otimes N \rightarrow N$ such that the diagrams
\index{amodule@$A$-module}
\index{module@module!$A$}
\begin{align}
\begin{array}{c@{\hskip 2truecm}c}
\Rnode{a}{A\otimes A \otimes N} & \Rnode{b}{A \otimes N} \\[8ex]
\Rnode{c}{A\otimes N}           & \Rnode{d}{N}
\end{array}
\ncline{->}{a}{b}
\Aput{\Id \otimes m_N}
\ncline{->}{b}{d}
\Aput{m_N}
\ncline{->}{a}{c}
\Aput{m_A \otimes \Id}
\ncline{->}{c}{d}
\Aput{m_N}
\quad&&~\quad
\pspicture[0.5](0,0)(1.5,2.0)
\psset{linewidth=\pstlw,xunit=0.5,yunit=0.5,runit=0.5}
\psset{arrowsize=2pt 2,arrowinset=0.2}
\psline(0,5)(0,4)
\psline(2,5)(2,4)
\psline(3,5)(3,2)
\psarc(1,4){1}{180}{360}
\psline(1,3)(1,2)
\psarc(2,2){1}{180}{360}
\psline(2,1)(2,0)
\pscircle[linewidth=0.4pt,fillstyle=solid,fillcolor=white](1,3){0.2}
\pscircle[linewidth=0.4pt,fillstyle=solid,fillcolor=white](2,1){0.2}
\rput(1,3.75){$m_A$}
\rput(2,1.75){$m_N$}
\endpspicture
\quad=\quad
\pspicture[0.5](0,0)(1.5,2.0)
\psset{linewidth=\pstlw,xunit=0.5,yunit=0.5,runit=0.5}
\psset{arrowsize=2pt 2,arrowinset=0.2}
\psline(0,5)(0,2)
\psline(1,5)(1,4)
\psline(3,5)(3,4)
\psarc(2,4){1}{180}{360}
\psline(2,3)(2,2)
\psarc(1,2){1}{180}{360}
\psline(1,1)(1,0)
\pscircle[linewidth=0.4pt,fillstyle=solid,fillcolor=white](1,1){0.2}
\pscircle[linewidth=0.4pt,fillstyle=solid,fillcolor=white](2,3){0.2}
\rput(1,1.75){$m_N$}
\rput(2,3.75){$m_N$}
\endpspicture
\end{align}
and
\begin{align}
\begin{array}{c@{\hskip 2truecm}c}
 & \Rnode{a}{A\otimes N} \\[8ex]
\Rnode{b}{\openk\otimes N} & \Rnode{c}{N} 
\ncline{->}{a}{c}
\Aput{m_N}
\ncline{->}{c}{b}
\Bput{\approx}
\ncline{->}{b}{a}
\Aput{\eta\otimes \Id}
\end{array}
\quad&&~\quad
\pspicture[0.5](0,0)(1,1.0)
\psset{linewidth=\pstlw,xunit=0.5,yunit=0.5,runit=0.5}
\psset{arrowsize=2pt 2,arrowinset=0.2}
\psline(2,3)(2,2)
\psarc(1,2){1.0}{180}{360}
\psline(1,1)(1,0)
\pscircle[linewidth=0.4pt,fillstyle=solid,fillcolor=white](1,1){0.2}
\pscircle[linewidth=0.4pt,fillstyle=solid,fillcolor=white](0,2){0.2}
\rput(1,1.75){$m_N$}
\rput(0,2.75){$\eta$}
\endpspicture
=
\pspicture[0.5](0,0)(1,1.0)
\psset{linewidth=\pstlw,xunit=0.5,yunit=0.5,runit=0.5}
\psset{arrowsize=2pt 2,arrowinset=0.2}
\psline(1,3)(1,0)
\endpspicture
\end{align}
are commutative.
\end{dfn}

\noindent
A morphism of graded $A$-modules $f ~:~ N\rightarrow N^\prime$ is a 
morphism of graded $\openk$-modules $N,N^\prime$ such that the diagrams
\index{module!morphism}
\begin{align}
\begin{array}{c@{\hskip 2truecm}c}
\Rnode{a}{A\otimes N} & \Rnode{b}{N} \\[8ex]
\Rnode{c}{A\otimes N^\prime} & \Rnode{d}{N^\prime}
\ncline{->}{a}{b}
\Aput{m_N}
\ncline{->}{a}{c}
\Aput{\Id\otimes f}
\ncline{->}{b}{d}
\Aput{f}
\ncline{->}{c}{d}
\Aput{m_{N^\prime}}
\end{array}
\quad&&~\quad
\pspicture[0.5](0,0)(1,2.0)
\psset{linewidth=\pstlw,xunit=0.5,yunit=0.5,runit=0.5}
\psset{arrowsize=2pt 2,arrowinset=0.2}
\psline(0,4)(0,3)
\psline(2,4)(2,3)
\psarc(1,3){1}{180}{360}
\psline(1,0)(1,2)
\pscircle[linewidth=0.4pt,fillstyle=solid,fillcolor=white](1,2){0.2}
\pscircle[linewidth=0.4pt,fillstyle=solid,fillcolor=black](1,1){0.2}
\rput(1,2.75){$m_N$}
\rput(0.5,1){$f$}
\endpspicture
\quad=\quad
\pspicture[0.5](0,0)(1,2.0)
\psset{linewidth=\pstlw,xunit=0.5,yunit=0.5,runit=0.5}
\psset{arrowsize=2pt 2,arrowinset=0.2}
\psline(0,4)(0,2)
\psline(2,4)(2,2)
\psarc(1,2){1}{180}{360}
\psline(1,0)(1,1)
\pscircle[linewidth=0.4pt,fillstyle=solid,fillcolor=white](1,1){0.2}
\pscircle[linewidth=0.4pt,fillstyle=solid,fillcolor=black](2,3){0.2}
\rput(1,1.75){$m_{N^\prime}$}
\rput( 1.5,3){$f$}
\endpspicture
\end{align}
and
\begin{align}
\begin{array}{c@{\hskip 2truecm}c}
\Rnode{a}{\openk} & \Rnode{b}{N} \\[8ex]
\Rnode{c}{\openk} & \Rnode{d}{N^\prime}
\ncline{->}{a}{b}
\Aput{\eta_N}
\ncline{->}{c}{d}
\Aput{\eta_{N^\prime}}
\ncline{->}{a}{c}
\Aput{\Id_\openk}
\ncline{->}{b}{d}
\Aput{f}
\end{array}         
&&
\pspicture[0.5](0,0)(1,2)
\psset{linewidth=\pstlw,xunit=0.5,yunit=0.5,runit=0.5}
\psset{arrowsize=2pt 2,arrowinset=0.2}
\psline(1,0)(1,2)
\pscircle[linewidth=0.4pt,fillstyle=solid,fillcolor=white](1,2){0.2}
\rput(1,3){$\eta_{N^\prime}$}
\endpspicture 
\quad=\quad
\pspicture[0.5](0,0)(1,2)
\psset{linewidth=\pstlw,xunit=0.5,yunit=0.5,runit=0.5}
\psset{arrowsize=2pt 2,arrowinset=0.2}
\psline(1,3)(1,0)
\pscircle[linewidth=0.4pt,fillstyle=solid,fillcolor=white](1,3){0.2}
\pscircle[linewidth=0.4pt,fillstyle=solid,fillcolor=black](1,2){0.2}
\rput(1,3.5){$\eta_N$}
\rput(0.5,2){$f$}
\endpspicture     
\end{align}
are commutative. Compare this relation with that of Eqn. \ref{eqn:5-9}.
Morphisms of $A$-modules can be added $(f+g) ~:~ N\rightarrow N^\prime$ 
via adding the particular degrees of $f,g$. Kernels and cokernels are
defined as $(\ker f)_q = \ker f_q$ etc.

Right $A$-modules are defined via the right action in a similar way.

Analogously to algebras we define connected modules as a $\openk$-module 
if $N_0 \cong \openk$. We will see later on that Gra{\ss}mann algebras 
are connected and have connected modules, while Clifford algebras are 
not connected.

\section{Co-algebras}

\subsection{Definitions}

We use now categorial duality, i.e. the reversion of arrows in 
commutative diagrams or the horizontal mirroring of tangles, to
define a dual structure called co-algebra.
\index{duality!categorial}

\begin{dfn}
An unital {\em co-algebra\/} over $\openk$ is a graded $\openk$-module
$C$ together with morphisms $\Delta$, the co-multiplication, and 
$\epsilon$, the counit, of graded $C$-modules
\index{co-algebra}
\index{counit}
\index{co-multiplication}\index{diagonalization}
\begin{align}
\Delta ~&:~ C \rightarrow C\otimes C \nn
\epsilon ~&:~ C \rightarrow \openk 
\end{align}
such that the diagrams
\begin{align}
\begin{array}{c@{\hskip 2truecm}c}
\Rnode{a}{C} & \Rnode{b}{C \otimes C} \\[8ex]
\Rnode{c}{C\otimes C} & \Rnode{d}{C\otimes C\otimes C}
\end{array}
\ncline{->}{a}{b}
\Aput{\Delta}
\ncline{->}{b}{d}
\Aput{\Id \otimes \Delta}
\ncline{->}{a}{c}
\Bput{\Delta}
\ncline{->}{c}{d}
\Aput{\Delta \otimes \Id}
\quad&&~\quad\quad
\pspicture[0.5](0,0)(1.5,2.5)
\psset{linewidth=\pstlw,xunit=0.5,yunit=0.5,runit=0.5}
\psset{arrowsize=2pt 2,arrowinset=0.2}
\psline(0,0)(0,1)
\psline(2,0)(2,1)
\psline(3,0)(3,3)
\psarc(1,1){1}{360}{180}
\psline(1,2)(1,3)
\psarc(2,3){1}{0}{180}
\psline(2,4)(2,5)
\pscircle[linewidth=0.4pt,fillstyle=solid,fillcolor=white](1,2){0.2}
\pscircle[linewidth=0.4pt,fillstyle=solid,fillcolor=white](2,4){0.2}
\rput(1,1.25){$\Delta$}
\rput(2,3.25){$\Delta$}
\endpspicture
\quad&=\quad
\pspicture[0.5](0,0)(1.5,2.5)
\psset{linewidth=\pstlw,xunit=0.5,yunit=0.5,runit=0.5}
\psset{arrowsize=2pt 2,arrowinset=0.2}
\psline(0,0)(0,3)
\psline(1,0)(1,1)
\psline(3,0)(3,1)
\psarc(2,1){1}{0}{180}
\psline(2,2)(2,3)
\psarc(1,3){1}{0}{180}
\psline(1,4)(1,5)
\pscircle[linewidth=0.4pt,fillstyle=solid,fillcolor=white](1,4){0.2}
\pscircle[linewidth=0.4pt,fillstyle=solid,fillcolor=white](2,2){0.2}
\rput(1,3.25){$\Delta$}
\rput(2,1.25){$\Delta$}
\endpspicture       
\end{align}
and
\begin{align}
\label{eqn:5-21}
\begin{array}{c@{\hskip 2truecm}c@{\hskip 2truecm}c}
\Rnode{a}{\openk\otimes C} & \Rnode{b}{C\otimes C} & 
\Rnode{c}{C\otimes\openk} \\[8ex]
 & \Rnode{d}{C} &
\ncline{->}{b}{a}
\Bput{\epsilon\otimes \Id}
\ncline{->}{b}{c}
\Aput{\Id\otimes \epsilon}
\ncline{->}{d}{b}
\Aput{\Delta}
\ncline{->}{a}{d}
\Bput{\approx}
\ncline{->}{c}{d}
\Aput{\approx}
\end{array}
\quad&&~
%
\pspicture[0.5](0,0)(1,1.5)
\psset{linewidth=\pstlw,xunit=0.5,yunit=0.5,runit=0.5}
\psset{arrowsize=2pt 2,arrowinset=0.2}
\psline(2,0)(2,1)
\psarc(1,1){1.0}{0}{180}
\psline(1,2)(1,3)
\pscircle[linewidth=0.4pt,fillstyle=solid,fillcolor=white](0,1){0.2}
\pscircle[linewidth=0.4pt,fillstyle=solid,fillcolor=white](1,2){0.2}
\rput(1,1.25){$\Delta$}
\rput(0,0.25){$\epsilon$}
\endpspicture
=
\pspicture[0.5](0,0)(1,1.5)
\psset{linewidth=\pstlw,xunit=0.5,yunit=0.5,runit=0.5}
\psset{arrowsize=2pt 2,arrowinset=0.2}
\psline(1,3)(1,0)
\endpspicture
=
\pspicture[0.5](0,0)(1,1.5)
\psset{linewidth=\pstlw,xunit=0.5,yunit=0.5,runit=0.5}
\psset{arrowsize=2pt 2,arrowinset=0.2}
\psline(0,0)(0,1)
\psarc(1,1){1.0}{0}{180}
\psline(1,2)(1,3)
\pscircle[linewidth=0.4pt,fillstyle=solid,fillcolor=white](2,1){0.2}
\pscircle[linewidth=0.4pt,fillstyle=solid,fillcolor=white](1,2){0.2}
\rput(1,1.25){$\Delta$}
\rput(2,0.25){$\epsilon$}
\endpspicture
\end{align}
are commutative.
\end{dfn}

\noindent
The co-algebra is called co-commutative (or simply commutative in a
unique context) if the diagram
\index{co-commutativity}
\begin{align}
\begin{array}{c@{\hskip 2truecm}c}
\Rnode{a}{C\otimes C} & \\[4ex]
& \Rnode{b}{C} \\[4ex]
\Rnode{c}{C\otimes C} &
\end{array}
\ncline{->}{b}{a}
\Bput{\Delta}
\ncline{->}{a}{c}
\Bput{\hat{\tau}_C}
\ncline{->}{b}{c}
\Bput{\Delta}
\quad&&~\quad
%
\pspicture[0.5](0,0)(1,2)
\psset{linewidth=\pstlw,xunit=0.5,yunit=0.5,runit=0.5}
\psset{arrowsize=2pt 2,arrowinset=0.2}
\psbezier(2,0)(2,2)(0,0)(0,2)
\psbezier[border=4pt](0,0)(0,2)(2,0)(2,2)
\psarc(1,2){1}{0}{180}
\psline(1,3)(1,4)
\pscircle[linewidth=0.4pt,fillstyle=solid,fillcolor=white](1,3){0.2}
\rput(1,0.25){$\hat{\tau}_C$}
\rput(1,2.25){$\Delta$}
\endpspicture
\quad=\quad
\pspicture[0.5](0,0)(1,2)
\psset{linewidth=\pstlw,xunit=0.5,yunit=0.5,runit=0.5}
\psset{arrowsize=2pt 2,arrowinset=0.2}
\psline(2,0)(2,2)
\psline(0,0)(0,2)
\psarc(1,2){1}{0}{180}
\psline(1,3)(1,4)
\pscircle[linewidth=0.4pt,fillstyle=solid,fillcolor=white](1,3){0.2}
\rput(1,2.25){$\Delta$}
\endpspicture    
\end{align}
is commutative.

A {\it crossed co-product\/} of co-algebras $C$ and $D$ over $\openk$
is defined as
\index{crossed co-product}
\index{co-product!crossed}
\begin{align}
\Rnode{1}{C\otimes D}\hskip 2truecm
\Rnode{2}{C\otimes C\otimes D\otimes D}\hskip 2truecm
\Rnode{3}{(C\otimes D)\otimes (C\otimes D)}
\ncline{->}{1}{2}
\Aput{\Delta_C\otimes\Delta_D}
\ncline{->}{2}{3}
\Aput{\Id\otimes\hat{\tau}\otimes\Id}
\end{align}
where the co unit is defined as
\begin{align}
\Rnode{1}{C\otimes D}\hskip 2truecm
\Rnode{2}{\openk\otimes\openk \approx \openk.}
\ncline{->}{1}{2}
\Aput{\epsilon_C\otimes\epsilon_D}
&&
\pspicture[0.5](0,0)(1,2.0)
\psset{linewidth=\pstlw,xunit=0.5,yunit=0.5,runit=0.5}
\psset{arrowsize=2pt 2,arrowinset=0.2}
\psline(0,4)(0,1)
\psline(2,4)(2,1)
\pscircle[linewidth=0.4pt,fillstyle=solid,fillcolor=white](0,1){0.2}
\pscircle[linewidth=0.4pt,fillstyle=solid,fillcolor=white](2,1){0.2}
\rput(0,0.25){$\epsilon$}
\rput(2,0.25){$\epsilon$}
\endpspicture
\end{align}
In terms of tangle one displays the crossed co-product as:
\begin{align}
\pspicture[0.4](0,0)(3,2.5)
\psset{linewidth=\pstlw,xunit=0.5,yunit=0.5,runit=0.5}
\psline(2,0)(2,1)
\psline(4,0)(4,1)
\psline(0,0)(0,3)
\psline(6,0)(6,3)
\psline{c-c}(4,1)(2,3)
\psline[border=4pt,bordercolor=white]{c-c}(2,1)(4,3)
\psarc(1,3){1.0}{0}{180}
\psarc(5,3){1.0}{0}{180}
\psline(1,4)(1,5)
\psline(5,4)(5,5)
\pscircle[linewidth=0.4pt,fillstyle=solid,fillcolor=white](1,4){0.2}
\pscircle[linewidth=0.4pt,fillstyle=solid,fillcolor=white](5,4){0.2}
\rput(1.0,3.0){$\Delta$}
\rput(3.0,1.0){$\hat\tau$}
\rput(5.0,3.0){$\Delta$}
\endpspicture
\end{align}

\noindent
{\bf Note:}
Observe that in older literature e.g. Sweedler \cite{sweedler:1969a}
and Milnor and Moore \cite{milnor:moore:1965a} it is assumed that the
crossing $\hat{\tau}$ is an involution, i.e. that 
$\hat{\tau}^2=\Id\otimes\Id$ hence $\hat{\tau}^{-1}=\hat{\tau}$.
If we do not impose this restriction, it could be convenient to introduce
the crossed co-product using the inverse crossing $\hat{\tau}^{-1}$
in spite of the fact that mirroring the tangle
\index{crossing!involutive}
\index{crossing!inverse}
\begin{align}
\pspicture[0.5](0,0)(1.5,1)
\psset{linewidth=\pstlw,xunit=0.5,yunit=0.5,runit=0.5}
\psset{arrowsize=2pt 2,arrowinset=0.2}
\psline(2,2)(0,0)
\psline[border=4pt,bordercolor=white](0,2)(2,0)
\psline{<->}(3,2)(3,0)
\endpspicture
\quad&\Rightarrow\quad
\pspicture[0.5](0,0)(1,1)
\psset{linewidth=\pstlw,xunit=0.5,yunit=0.5,runit=0.5}
\psset{arrowsize=2pt 2,arrowinset=0.2}
\psline(2,2)(0,0)
\psline[border=4pt,bordercolor=white](2,0)(0,2)
\endpspicture
\quad=\quad
\left(\quad
\pspicture[0.5](0,0)(1,1)
\psset{linewidth=\pstlw,xunit=0.5,yunit=0.5,runit=0.5}
\psset{arrowsize=2pt 2,arrowinset=0.2}
\psline(2,0)(0,2)
\psline[border=4pt,bordercolor=white](2,2)(0,0)
\endpspicture
\quad\right)^{-1}
\end{align}
yields the crossing back, assuming that the over and under information
is corerctly encoded.

A morphism of co-algebras $f~:~ C\rightarrow D$ is a morphism of
graded $\openk$-modules such that the diagrams
\index{morphism!of co-algebras}
\begin{align}
\begin{array}{c@{\hskip 2truecm}c}
\Rnode{a}{C} & \Rnode{b}{C\otimes C} \\[8ex]
\Rnode{c}{D} & \Rnode{d}{D\otimes D}
\ncline{->}{a}{b}
\Aput{\Delta_C}
\ncline{->}{c}{d}
\Aput{\Delta_D}
\ncline{->}{a}{c}
\Aput{f}
\ncline{->}{b}{d}
\Aput{f\otimes f}
\end{array}         
&&
\pspicture[0.5](0,0)(1,2)
\psset{linewidth=\pstlw,xunit=0.5,yunit=0.5,runit=0.5}
\psset{arrowsize=2pt 2,arrowinset=0.2}
\psline(0,0)(0,1)
\psline(2,0)(2,1)
\psarc(1,1){1}{0}{180}
\psline(1,4)(1,2)
\pscircle[linewidth=0.4pt,fillstyle=solid,fillcolor=white](1,2){0.2}
\pscircle[linewidth=0.4pt,fillstyle=solid,fillcolor=black](1,3){0.2}
\rput(1,1){$\Delta_D$}
\rput(0.5,3){$f$}
\endpspicture
\quad=\quad
\pspicture[0.5](0,0)(1,2)
\psset{linewidth=\pstlw,xunit=0.5,yunit=0.5,runit=0.5}
\psset{arrowsize=2pt 2,arrowinset=0.2}
\psline(0,0)(0,2)
\psline(2,0)(2,2)
\psarc(1,2){1}{0}{180}
\psline(1,4)(1,3)
\pscircle[linewidth=0.4pt,fillstyle=solid,fillcolor=white](1,3){0.2}
\pscircle[linewidth=0.4pt,fillstyle=solid,fillcolor=black](0,1){0.2}
\pscircle[linewidth=0.4pt,fillstyle=solid,fillcolor=black](2,1){0.2}
\rput(1,2.0){$\Delta_C$}
\rput(-0.5,1){$f$}
\rput( 1.5,1){$f$}
\endpspicture              
\end{align}
and
\begin{align}
\begin{array}{c@{\hskip 2truecm}c}
\Rnode{a}{C} & \Rnode{b}{\openk} \\[8ex]
\Rnode{c}{D} & \Rnode{d}{\openk}
\ncline{->}{a}{b}
\Aput{\epsilon_C}
\ncline{->}{c}{d}
\Aput{\epsilon_D}
\ncline{->}{a}{c}
\Aput{f}
\ncline{->}{b}{d}
\Aput{\Id_\openk}
\end{array}         
&&
\pspicture[0.5](0,0)(1,2)
\psset{linewidth=\pstlw,xunit=0.5,yunit=0.5,runit=0.5}
\psset{arrowsize=2pt 2,arrowinset=0.2}
\psline(1,4)(1,1)
\pscircle[linewidth=0.4pt,fillstyle=solid,fillcolor=white](1,1){0.2}
\rput(1,0.15){$\epsilon_C$}
\endpspicture 
\quad=\quad
\pspicture[0.5](0,0)(1,2)
\psset{linewidth=\pstlw,xunit=0.5,yunit=0.5,runit=0.5}
\psset{arrowsize=2pt 2,arrowinset=0.2}
\psline(1,4)(1,1)
\pscircle[linewidth=0.4pt,fillstyle=solid,fillcolor=white](1,1){0.2}
\pscircle[linewidth=0.4pt,fillstyle=solid,fillcolor=black](1,3){0.2}
\rput(1,0.15){$\epsilon_D$}
\rput(0.5,3){$f$}
\endpspicture     
\end{align}
are commutative. A co-algebra is co-commutative if and only if 
$\Delta ~:~ C\rightarrow C\otimes C$ is a morphism of co-algebras.
The proof follows by duality from the proof for the case of algebras.
\index{commutativity!of co-algebras}
\index{co-commutativity}

One can consider $\openk$ to be a co-algebra in a canonical way. This 
allows to introduce an {\it augmentation\/} of co-algebras as a morphism
of co-algebras $\eta ~:~ \openk \rightarrow C$. If $C$ is an augmented
co-algebra, i.e. a co-algebra with augmentation $\eta$, we denote by
$J(A)$ the cokernel of $\eta$. Considering $C$ as a $\openk$-module we
find 
\index{co-algebra!augmented}
\index{augmentation!cokernel}
\begin{align}
C &= \openk \oplus J(A).
\end{align}

\subsection{$C$-comodules}

Let $C$ be a co-algebra over $\openk$. A left $C$-comodule is a graded 
$\openk$-module $N$ together with morphisms 
$\Delta_N ~:~ N \rightarrow C\otimes N$ such that the diagrams
\index{comodule}
\begin{align} 
\begin{array}{c@{\hskip 2truecm}c}
\Rnode{a}{N} & \Rnode{b}{C \otimes N} \\[8ex]
\Rnode{c}{C\otimes N} & \Rnode{d}{C\otimes C\otimes N}
\end{array}
\ncline{->}{a}{b}
\Aput{\Delta_N}
\ncline{->}{b}{d}
\Aput{\Id \otimes \Delta_N}
\ncline{->}{a}{c}
\Bput{\Delta_N}
\ncline{->}{c}{d}
\Aput{\Delta_C \otimes \Id}
\quad&&~\quad\quad
\pspicture[0.5](0,0)(1.5,2.5)
\psset{linewidth=\pstlw,xunit=0.5,yunit=0.5,runit=0.5}
\psset{arrowsize=2pt 2,arrowinset=0.2}
\psline(0,0)(0,1)
\psline(2,0)(2,1)
\psline(3,0)(3,3)
\psarc(1,1){1}{360}{180}
\psline(1,2)(1,3)
\psarc(2,3){1}{0}{180}
\psline(2,4)(2,5)
\pscircle[linewidth=0.4pt,fillstyle=solid,fillcolor=white](1,2){0.2}
\pscircle[linewidth=0.4pt,fillstyle=solid,fillcolor=white](2,4){0.2}
\rput(1,1.25){$\Delta_C$}
\rput(2,3.25){$\Delta_N$}
\endpspicture
\quad&=\quad
\pspicture[0.5](0,0)(1.5,2.5)
\psset{linewidth=\pstlw,xunit=0.5,yunit=0.5,runit=0.5}
\psset{arrowsize=2pt 2,arrowinset=0.2}
\psline(0,0)(0,3)
\psline(1,0)(1,1)
\psline(3,0)(3,1)
\psarc(2,1){1}{0}{180}
\psline(2,2)(2,3)
\psarc(1,3){1}{0}{180}
\psline(1,4)(1,5)
\pscircle[linewidth=0.4pt,fillstyle=solid,fillcolor=white](1,4){0.2}
\pscircle[linewidth=0.4pt,fillstyle=solid,fillcolor=white](2,2){0.2}
\rput(1,3.25){$\Delta_N$}
\rput(2,1.25){$\Delta_N$}
\endpspicture           
\end{align}
and
\begin{align}
\begin{array}{c@{\hskip 2truecm}c}
\Rnode{a}{\openk\otimes C} & \Rnode{b}{C\otimes C} \\[8ex]
 & \Rnode{d}{C} 
\ncline{->}{b}{a}
\Bput{\epsilon\otimes\Id}
\ncline{->}{d}{b}
\Bput{\Delta}
\ncline{->}{a}{d}
\Bput{\approx}
\end{array}
\hskip 2truecm&&~        
\pspicture[0.5](0,0)(1,1.5)
\psset{linewidth=\pstlw,xunit=0.5,yunit=0.5,runit=0.5}
\psset{arrowsize=2pt 2,arrowinset=0.2}
\psline(2,0)(2,1)
\psarc(1,1){1.0}{0}{180}
\psline(1,2)(1,3)
\pscircle[linewidth=0.4pt,fillstyle=solid,fillcolor=white](0,1){0.2}
\pscircle[linewidth=0.4pt,fillstyle=solid,fillcolor=white](1,2){0.2}
\rput(1,1.25){$\Delta$}
\rput(0,0.25){$\epsilon$}
\endpspicture
\quad=
\pspicture[0.5](0,0)(1,1.5)
\psset{linewidth=\pstlw,xunit=0.5,yunit=0.5,runit=0.5}
\psset{arrowsize=2pt 2,arrowinset=0.2}
\psline(1,3)(1,0)
\endpspicture 
\end{align} 
are commutative.

Let $N,N^\prime$ are left $C$-modules. A morphism $f ~:~ N \rightarrow
N^\prime$ of left $C$-modules is a morphism of graded $\openk$-modules
such that the diagram
\index{morphism!of comodules}
\begin{align}
\begin{array}{c@{\hskip 2truecm}c}
\Rnode{a}{N} & \Rnode{b}{C\otimes N} \\[8ex]
\Rnode{c}{N^\prime} & \Rnode{d}{C\otimes N^\prime}
\ncline{->}{a}{b}
\Aput{\Delta_N}
\ncline{->}{c}{d}
\Aput{\Delta_{N^\prime}}
\ncline{->}{a}{c}
\Aput{f}
\ncline{->}{b}{d}
\Aput{\Id\otimes f}
\end{array}         
&&
\pspicture[0.5](0,0)(1,2)
\psset{linewidth=\pstlw,xunit=0.5,yunit=0.5,runit=0.5}
\psset{arrowsize=2pt 2,arrowinset=0.2}
\psline(0,0)(0,1)
\psline(2,0)(2,1)
\psarc(1,1){1}{0}{180}
\psline(1,4)(1,2)
\pscircle[linewidth=0.4pt,fillstyle=solid,fillcolor=white](1,2){0.2}
\pscircle[linewidth=0.4pt,fillstyle=solid,fillcolor=black](1,3){0.2}
\rput(1,1){$\Delta_{N^\prime}$}
\rput(0.5,3){$f$}
\endpspicture
\quad=\quad
\pspicture[0.5](0,0)(1,2)
\psset{linewidth=\pstlw,xunit=0.5,yunit=0.5,runit=0.5}
\psset{arrowsize=2pt 2,arrowinset=0.2}
\psline(0,0)(0,2)
\psline(2,0)(2,2)
\psarc(1,2){1}{0}{180}
\psline(1,4)(1,3)
\pscircle[linewidth=0.4pt,fillstyle=solid,fillcolor=white](1,3){0.2}
\pscircle[linewidth=0.4pt,fillstyle=solid,fillcolor=black](2,1){0.2}
\rput(1,2.0){$\Delta_N$}
\rput( 1.5,1){$f$}
\endpspicture 
\end{align}
is commutative.

\noindent
{\bf Note:} Beside the completely parallel developments of algebras,
co-algebras and modules, comodules due to categorial duality, they exhibit, 
after a close look, some different features. As an example, one should
note, that the left $A$-modules constitute an abelian category, while 
the left $C$-comodules do not unless further conditions are assumed
\cite{milnor:moore:1965a}. Since we do not need such sophisticated
facts here, the interested reader should consult the original literature.
\index{module!abelian category of}

\section{Bialgebras}

\subsection{Definitions}

A graded $\openk$-module $A$ is of {\it finite type\/} if each $A_n$
is finitely generated as $\openk$-module and only a finite number of 
$A_n$'s are not zero. It is called projective if each $A_n$ is projective.
\index{module!finite}
\index{module!projective}

Under these conditions one obtains that the morphism of $\openk$-modules
\begin{align}
\lambda ~&:~ A \rightarrow A^{**}
\end{align}
defined by $\lambda(x)a^* = a^*(x)$ for $x\in A_n$, $a^* \in A^*_n$
is an isomorphism, and the morphism of graded $\openk$-modules
\begin{align}
\alpha ~&:~ A^*\otimes B^* \rightarrow (A\otimes B)^*
\end{align}
defined as $\alpha(a^*\otimes b^*) (x\otimes y) = a^*(y) b^*(x)$
for $a^*\in A^*_p$, $b^*\in B^*_q$, $y\in A_p$, $x\in B_q$ is an
isomorphism. That is the tensor product of linear forms is a linear 
form itself acting on a tensor product in a canonical way. This allows 
us to write $A^{**}=A$ and $A^* \otimes B^* = (A\otimes B)^*$.
One could prove, that $A^*$ is of projective finite type if $A$ is
of projective finite type and vice versa.

The product $m_{A}$ from $A$ induces now a {\it co-product\/} $m^*_{A^*}$
by categorial duality also abusing our notation. In other words, a product 
of vectors implies a co-product on co-vectors and a products on co-vectors 
implies a co-product on vectors by so called product co-product duality. 
Having a pair of a product and co-product on $A$ and $A^*$, we find four 
products which might be independently chosen.

\index{duality!of algebra and co-gebra}
\begin{thrm}[Duality]
Suppose that $A$ is a graded $\openk$-module of projective and finite 
type, then:
\begin{itemize}
\item[(1)] $m_A~:~ A\otimes A \rightarrow A$ is a multiplication in $A$ 
if and only if $m^*_{A^*} ~:~ A^* \rightarrow A^* \otimes A^*$ is a 
co-multiplication in $A^*$.
\item[(2)] $m_A$ is associative if and only if $m^*_{A^*}$ is 
co-associative.
\item[(3)] $\eta ~:~ \openk\rightarrow A$ is a unit for the 
multiplication $m_A$ if and only if $\eta^* ~:~ A^*\rightarrow \openk^*
\cong \openk$ is a counit for the co-multiplication $m^*_{A^*}$.
\item[(4)] $(A,m_A,\eta)$ is an associative unital algebra if and only if 
$(A^*,m^*_{A^*},\eta^*)$ is a co-associative counital co-algebra.
\item[(5)] $\epsilon ~:~ A\rightarrow \openk$ is an augmentation
of the algebra $(A,m_A,\eta)$ if and only if $\epsilon^* ~:~
\openk \rightarrow A^*$ is an augmentation of the co-algebra 
$(A^*,m^*_{A^*},\eta^*)$.
\item[(6)] The algebra $(A,m_A,\eta)$ is commutative if and only
if the co-algebra  $(A^*,m^*_{A^*},\eta^*)$ is cocommutative.
\end{itemize}
\end{thrm}
The preceeding theorem exhibits a beautiful symmetry. Moreover, one 
notes, that if an algebra structure is employed on the module $A^*$
of linear forms, this induces automatically an co-algebra structure
on the (double) dual module $A^{**}=A$.

At this point one should reconsider the algebra of meet and join,
i.e. the Gra{\ss}mann-Cayley algebra. Since the join acts on 
{\it points\/} $\in A$ and the meet acts on {\it planes\/} $\in A^*$,
the join induces a co-algebra structure $\Delta_\vee$ on the planes, 
i.e. on $A^*$, and the meet induces a co-algebra structure $\Delta_\wedge$
on points, i.e. on $A$. The duality between meet and join is mediated by 
Gra{\ss}mann's Erg\"anzung, i.e. the orthogonal complement and the duality
between products and co-products by categorial duality. This can be 
displayed in the diagram
\index{algebra@algebra!Gra{\ss}mann-Cayley}
\index{meet}\index{join}\index{Ergaenzung@Erg\"{a}nzung}
\begin{align}
\begin{array}{c@{\hskip 3truecm}c}
\Rnode{1}{GC(\wedge,\vee,*)} & \Rnode{2}{GC(\Delta_\vee,\Delta_\wedge,*)} 
\\[8ex]
\Rnode{3}{H_\wedge(\wedge,\Delta_\vee,*)} &
\Rnode{4}{H_\vee(\vee,\Delta_\wedge,*)}
\ncline{<->}{1}{2}
\Aput{\star}
\ncline{<->}{3}{4}
\Aput{\star}
\ncline{<->}{1}{3}
\ncline{<->}{2}{4}
\end{array}
\end{align}
Hence the Gra{\ss}mann-Cayley algebra is built over a pair of spaces 
$A\oplus A^*$, both seen as algebras, and this structure is equivalent
to a Gra{\ss}mann-Hopf algebra, see below.
\index{algebra@algebra!Gra{\ss}mann-Hopf}

Now, it is clear that one is tempted to complete the structure to go
over to a four-fold algebra, see \cite{doubilet:rota:stein:1974a}. 
In fact this is also done using the quantum double 
\cite{drinfeld:1987a,majid:1995a}
\begin{align}
H\oplus H^* &\cong GC(\wedge,\Delta_\vee,\vee,\Delta_\wedge,*)
\end{align}

\noindent
{\bf Note} that the introduction of two independent products
on $A$ and $A^*$ results in a completely {\it independent\/} 
structure $(\wedge,\Delta_\vee)$ on $A$ and $(\vee,\Delta_\wedge)$
on $A^*$. Hence product and co-product on {\it one\/} space are not related
at all. This will motivate later on the study of convolution algebras
having deformed products and {\it independently\/} deformed co-products.
Unfortunately, because of the canonical identification of $A^*$ with $A$
via the co-vector basis $\varepsilon^i$ fulfilling 
$\varepsilon^i(e_j)=\delta^i_j$, hides the fact that one deals with 
two independent structures.

\index{bialgebra}
\begin{dfn}[Bialgebra]
A {\em bialgebra} over $\openk$ is a graded $\openk$-module $B$
together with morphisms of graded $\openk$-modules
\begin{align}
m_B~&:~ B\otimes B \rightarrow B,\qquad \eta_B~:~\openk\rightarrow B \nn
\Delta_B~&:~ B\rightarrow B\otimes B,\qquad \epsilon_B~:~B\rightarrow \openk
\end{align} 
such that
\begin{itemize}
\item[(1)] $(B,m_B,\eta_B)$ is an augmented $\openk$-algebra,
\item[(2)] $(B,\Delta_B,\epsilon_B)$ is an augmented $\openk$-co-algebra,
\item[(3)] the diagram
\begin{align}
\label{eqn:5-38}
\hskip -2truecm
\begin{array}{c@{\hskip 1.0truecm}c@{\hskip 1.0truecm}c}
\Rnode{1}{B\otimes B} & \Rnode{2}{B} & \Rnode{3}{B\otimes B} \\[8ex]
\Rnode{4}{B\otimes B\otimes B\otimes B} & & 
\Rnode{5}{B\otimes B\otimes B\otimes B}
\ncline{->}{1}{2}
\Aput{m_B}
\ncline{->}{2}{3}
\Aput{\Delta_B}
\ncline{->}{1}{4}
\Aput{\Delta_B \otimes \Delta_B}
\ncline{->}{4}{5}
\Aput{\Id \otimes \hat{\tau}\otimes \Id}
\ncline{->}{5}{3}
\Bput{m_B \otimes m_B}
\end{array}
\quad&&~
\pspicture[0.5](0,0)(1,3)
\psset{linewidth=\pstlw,xunit=0.5,yunit=0.5,runit=0.5}
\psline(0,6)(0,5)
\psline(2,6)(2,5)
\psarc(1,5){1.0}{180}{360}
\psline(1,4)(1,2)
\psarc(1,1){1.0}{0}{180}
\psline(0,1)(0,0)
\psline(2,1)(2,0)
\pscircle[linewidth=0.4pt,fillstyle=solid,fillcolor=white](1,4){0.2}
\pscircle[linewidth=0.4pt,fillstyle=solid,fillcolor=white](1,2){0.2}
\rput(1.0,4.75){$m_B$}
\rput(1.0,1.25){$\Delta_B$}
\endpspicture
=\quad 
\pspicture[0.5](0,0)(3,3)
\psset{linewidth=\pstlw,xunit=0.5,yunit=0.5,runit=0.5}
\psline(1,6)(1,5)
\psline(5,6)(5,5)
\psarc(1,4){1.0}{0}{180}
\psarc(5,4){1.0}{0}{180}
\psline(0,4)(0,2)
\psline(6,4)(6,2)
\psline{c-c}(4,4)(2,2)
\psline[border=4pt]{c-c}(2,4)(4,2)
\psarc(1,2){1.0}{180}{360}
\psarc(5,2){1.0}{180}{360}
\psline(1,1)(1,0)
\psline(5,1)(5,0)
\pscircle[linewidth=0.4pt,fillstyle=solid,fillcolor=white](1,1){0.2}
\pscircle[linewidth=0.4pt,fillstyle=solid,fillcolor=white](1,5){0.2}
\pscircle[linewidth=0.4pt,fillstyle=solid,fillcolor=white](5,1){0.2}
\pscircle[linewidth=0.4pt,fillstyle=solid,fillcolor=white](5,5){0.2}
\rput(1.0,4.25){$\Delta_B$}
\rput(1.0,1.75){$m_B$}
\rput(3.0,3.75){$\hat{\tau}$}
\rput(5.0,4.25){$\Delta_B$}
\rput(5.0,1.75){$m_B$}
\endpspicture          
\end{align}
is commutative.
\end{itemize}
\end{dfn}

\noindent
Condition (3) states that $\Delta_B$ is an algebra homomorphism since
it preserves also units and $m_B$ is a co-algebra homomorphism preserving 
also counits. This follows from the fact that augmented algebras and
co-algebras have been considered.

\noindent
{\bf Note:}
This definition, even adopted in the vast majority of literature will
be too narrow for our purpose. We will see that we have to drop the fact
that $B$ is an augmented algebra or co-algebra. In this case, $\Delta_B$ 
and $m_B$ do not preserve the counit and unit respectively.
Furthermore in older literature this structure is already called Hopf 
algebra, while we reserve this term for a still more restrictive setting.

Having defined the notion of a bialgebra, unnaturally emphasizing the 
algebra part, we can speak in a more symmetric fashion about 
bi-associativity if $B$ is associative and co-associative, about $B$ 
being bi-unital if it is unital and counital etc.

Using crossed products and crossed co-products one can establish a 
bialgebra action  and coaction on the tensor product $N\otimes M$ 
of two left $B$-modules over $\openk$. In tangle notation this reads:
\begin{align}
\pspicture[0.5](0,0)(3,3)
\psset{linewidth=\pstlw,xunit=0.5,yunit=0.5,runit=0.5}
\psline(1,6)(1,5)
\psarc(1,4){1.0}{0}{180}
\psline(4,6)(4,4)
\psline(6,6)(6,4)
\psline(0,4)(0,2)
\psline(6,4)(6,2)
\psline{c-c}(4,4)(2,2)
\psline[border=4pt]{c-c}(2,4)(4,2)
\psarc(1,2){1.0}{180}{360}
\psarc(5,2){1.0}{180}{360}
\psline(1,1)(1,0)
\psline(5,1)(5,0)
\pscircle[linewidth=0.4pt,fillstyle=solid,fillcolor=white](1,1){0.2}
\pscircle[linewidth=0.4pt,fillstyle=solid,fillcolor=white](1,5){0.2}
\pscircle[linewidth=0.4pt,fillstyle=solid,fillcolor=white](5,1){0.2}
\rput(1.5,5.5){$B$}
\rput(4.5,5.5){$N$}
\rput(6.5,5.5){$M$}
\rput(1.5,0.25){$N$}
\rput(5.5,0.05){$M$}
\endpspicture          
\quad&&~\quad
\pspicture[0.5](0,0)(3,3)
\psset{linewidth=\pstlw,xunit=0.5,yunit=0.5,runit=0.5}
\psline(1,6)(1,5)
\psline(5,6)(5,5)
\psarc(1,4){1.0}{0}{180}
\psarc(5,4){1.0}{0}{180}
\psline(0,4)(0,2)
\psline(6,4)(6,2)
\psline{c-c}(4,4)(2,2)
\psline[border=4pt]{c-c}(2,4)(4,2)
\psarc(1,2){1.0}{180}{360}
\psline(4,2)(4,0)
\psline(6,2)(6,0)
\psline(1,1)(1,0)
\pscircle[linewidth=0.4pt,fillstyle=solid,fillcolor=white](1,5){0.2}
\pscircle[linewidth=0.4pt,fillstyle=solid,fillcolor=white](1,1){0.2}
\pscircle[linewidth=0.4pt,fillstyle=solid,fillcolor=white](5,5){0.2}
\rput(1.5,5.5){$N$}
\rput(1.5,0.25){$B$}
\rput(5.5,5.5){$M$}
\rput(4.5,0.25){$N$}
\rput(6.5,0.25){$M$}
\endpspicture          
\end{align}
where one has to use $\Delta_N$ and $\Delta_M$ in the r.h.s. tangle.

\index{bialgebra!quasi}
Milnor and Moore introduce a {\it quasi bialgebra\/} (quasi
Hopf algebra in their notation) which does not assume associativity
of multiplication and co-multiplication and where the augmentation
is replaced by the condition
\begin{align}
(4)~&~\epsilon\circ\eta ~=~ \Id_\openk,\qquad
\Rnode{1}{\openk}\hskip 2truecm\Rnode{2}{\openk .}
\ncline{->}{1}{2}
\Aput{\epsilon\circ\eta}
\end{align}
However, we are interested in associative multiplications exclusively
and we will not follow this track.

\index{filtration}
If $A$ is a graded $\openk$-module, a {\it filtration\/} of $A$ is a 
family $\{F_p A\}$ of sub-graded $\openk$-modules of $A$, indexed by
the integers such that $F_p A \subset F_{p+1}A$. The filtration 
$\{F_p A\}$ of the graded module $A$ is {\it complete\/} if
\begin{align}
&&&(1)\quad A = \lim_{p\rightarrow  \infty} F_p A &&&&\nn
&\text{or}&&(2)\quad A = \lim_{p\rightarrow -\infty} A/F_p A. &&&&
\end{align}
\index{algebra!filtered}\index{module!filtered}%
A {\it filtered algebra\/} is an algebra over a graded $\openk$-module
such that the multiplication map $m_A~:~A\otimes A\rightarrow A$
is a morphism of filtered graded modules. A filtered left $A$-module
$M$ is a graded left $A$-module with filtration on its underlying
graded $\openk$-module such that the action $m_A~:~ A\otimes M 
\rightarrow M$ is a morphism of filtered graded $A$-modules. One may
easily generalize this notion to a filtered bialgebra.

If $A$ is an augmented algebra over $\openk$, let $Q(A) = \openk
\otimes_A I(A)$. The elements of the graded $\openk$-module $Q(A)$
are called the {\it indecomposable elements\/} of $A$. If $C$ is an
\index{element!indecomposable}
augmented co-algebra over $\openk$, let $P(C) = \openk \Box_A J(A)$, 
where $\Box_A$ is the co-tensor product, see \cite{milnor:moore:1965a}.
The elements of the graded $\openk$-module $P(C)$ are called 
{\it primitive elements\/} of $C$. A bialgebra $B$ is said to be
{\it primitively generated\/} if the smallest sub-bialgebra of $B$ 
containing $P(B)$ is $B$ itself.
\index{element!primitive}

The notion of a filtration can now be used to generalize that of 
grading, connectedness and primitivity. However, while the augmentation
was sufficient to prove the facts about the kernel $I(B)$ and the
cokernel $J(B)$, i.e. determines the structure of $Q(B)$ and $P(B)$,
this is no longer true for filtered (quasi) bialgebras. There one has to
impose further conditions, that is to assume that certain exact sequences
of the $Q(A)$ and $P(Q)$ modules are split, to be able to draw the 
conclusions. This fact will face us below, but we do not develop the 
corresponding theory since we have not yet applied it to the examples
in physics.

\section{Hopf algebras i.e. antipodal bialgebras}

\subsection{Morphisms of connected co-algebras and connected algebras : 
group like convolution}

Let $C$ be a connected co-algebra and $A$ be a connected algebra,
i.e. $\eta :\openk\rightarrow A_0$ is an isomorphism, let 
$\Conv(C,A)$ be the set of morphisms $f~:~C\rightarrow A$ such that
$f_0$ is the identity morphism of $\openk$. If $f,g \in \Conv(C,A)$, the 
{\it convolution product\/} $f\star g$ is defined as the composition
\index{convolution}
\index{product!convolution}
\begin{align}
\Rnode{1}{C}\hskip 1truecm
\Rnode{2}{C\otimes C}\hskip 1truecm
\Rnode{3}{A\otimes A}\hskip 1truecm
\Rnode{4}{A}
\ncline{->}{1}{2}
\Aput{\Delta_C}
\ncline{->}{2}{3}
\Aput{f\otimes g}
\ncline{->}{3}{4}
\Aput{m_A}
&&~
\pspicture[0.5](0,0)(1,3)
\psset{linewidth=\pstlw,xunit=0.5,yunit=0.5,runit=0.5}
\psset{arrowsize=2pt 2,arrowinset=0.2}
\psline{-}(1,6)(1,4)
\psline{-}(1,2)(1,0)
\psarc(1,3){1.0}{0}{360}
\pscircle[linewidth=0.4pt,fillstyle=solid,fillcolor=black](0,3){0.2}
\pscircle[linewidth=0.4pt,fillstyle=solid,fillcolor=black](2,3){0.2}
\pscircle[linewidth=0.4pt,fillstyle=solid,fillcolor=white](1,4){0.2}
\pscircle[linewidth=0.4pt,fillstyle=solid,fillcolor=white](1,2){0.2}
\rput(-0.75,3){$f$}
\rput(2.75,3){$g$}
\rput(1.5,5.5){$C$}
\rput(1.5,0.5){$A$}
\endpspicture            
\end{align}

\begin{thrm}[Milnor \& Moore]
\label{thrm:5-9}
If $C$ is a connected co-algebra and $A$ is a connected algebra, 
then $\Conv(C,A)$ is a group under the convolution product $\star$
with identity (convolution unit)
\index{group!convolution}
\begin{align}
\Rnode{1}{C}\hskip 2truecm
\Rnode{2}{\openk}\hskip 2truecm
\Rnode{3}{A}
\ncline{->}{1}{2}
\Aput{\epsilon}
\ncline{->}{2}{3}
\Aput{\eta}
\end{align}
\end{thrm}
\label{proof-antipode}
{\bf Proof:} $\Conv(C,A)$ is a monoid regarding its definition, and one 
needs to prove the existence of $f^{-1}$ only. Suppose now that the action
of $f^{-1}$ is defined on degrees less than $n$, $x\in A_n$ and $\Delta(x) = 
x\otimes 1 + 1\otimes x +\sum^\prime_{(x)}x_{(1)}\otimes x_{(2)}$,
where the prime indicates a sum over proper cuts, i.e. $x_{(1)}\not=1$
and $x_{(2)}\not=1$. Assume $n>0$ and recall that 
$0 < \text{degree~} x_{(2)}< n$ for proper cuts and all $i$. 
Let $f^{-1}(x) = -x - \sum^\prime_{(x)} x_{(1)}f^{-1}(x_{(2)})$, which is 
the recursive definition of $f^{-1}$, since $f_0\cong \Id_\openk$.
I.e. we have $f\star f^{-1}=\eta\circ\epsilon$.

\noindent
{\bf Note} that this recursive definition of $f^{-1}$ is used, and was 
reinvented, by various authors. The most complete treatment might be
found in Schmitt \cite{schmitt:1987a}, where the antipode is constructed
in this way, but also a non-recursive formula is given. The Connes-Kreimer
antipode is, up to the renormalization scheme, calculated this way
\cite{connes:kreimer:1999a,kreimer:2000a,brouder:1999a,brouder:1999b%
,brouder:2000a}. In fact $f^{-1}$ is the inverse of $\,\Id$ and hence
the antipode, but any other inverse can be obtained in this way.
\index{inverse!convolutive}
\index{antipode!recursive def.}
\index{renormalization}

If one considers morphisms $f~:~ C\rightarrow C^\prime$ of connected 
co-algebras and $g~:~A\rightarrow A^\prime$ of connected algebras
there is an induced morphism of groups $\Conv(f,g)~:~ \Conv(C,A)
\rightarrow \Conv(C^\prime,A^\prime)$
\begin{align}
\Conv(f,g)h = ghf &&\quad 
\pspicture[0.5](0,0)(1,2)
\psset{linewidth=\pstlw,xunit=0.5,yunit=0.5,runit=0.5}
\psset{arrowsize=2pt 2,arrowinset=0.2}
\psline{-}(1,4)(1,0)
\pscircle[linewidth=0.4pt,fillstyle=solid,fillcolor=black](1,1){0.2}
\pscircle[linewidth=0.4pt,fillstyle=solid,fillcolor=black](1,2){0.2}
\pscircle[linewidth=0.4pt,fillstyle=solid,fillcolor=black](1,3){0.2}
\rput(1.5,1){$g$}
\rput(1.5,2){$h$}
\rput(1.5,3){$f$}
\endpspicture 
\end{align} 

\begin{dfn}[Antipode]
\index{antipode}
Let $B$ be a biconnected bialgebra, the {\em antipode\/} or 
{\em conjugation\/} of $B$ is the (unique) inverse in $\Conv(B,B)$ of 
the identity morphism of $B$. The antipode is denoted as $S_B$.
\end{dfn}

Let $U=\eta\circ\epsilon$ be the convolution unit as defined above. 
The defining relations of the antipode reads in tangle notation:
\begin{align}
\pspicture[0.5](0,0)(1,2.5)
\psset{linewidth=\pstlw,xunit=0.5,yunit=0.5,runit=0.5}
\psset{arrowsize=2pt 2,arrowinset=0.2}
\psline(1,5)(1,4)
\psarc(1,3){1}{0}{180}
\psline(2,3)(2,2)
\psarc(1,2){1}{180}{360}
\psline(1,1)(1,0)
\pscircle[linewidth=0.4pt,fillstyle=solid,fillcolor=white](1,1){0.2}
\pscircle[linewidth=0.4pt,fillstyle=solid,fillcolor=white](1,4){0.2}
\rput(1,1.75){m}
\rput(1,3.25){$\Delta$}
\rput(0,2.5){S}
\endpspicture
\quad&=\quad
\pspicture[0.5](0,0)(1,2.5)
\psset{linewidth=\pstlw,xunit=0.5,yunit=0.5,runit=0.5}
\psset{arrowsize=2pt 2,arrowinset=0.2}
\psline(1,5)(1,3)
\psline(1,0)(1,2)
\pscircle[linewidth=0.4pt,fillstyle=solid,fillcolor=white](1,3){0.2}
\pscircle[linewidth=0.4pt,fillstyle=solid,fillcolor=white](1,2){0.2}
\rput(1.5,3){$\epsilon$}
\rput(1.5,2){$\eta$}
\endpspicture
\quad=\quad
\pspicture[0.5](0,0)(1,2.5)
\psset{linewidth=\pstlw,xunit=0.5,yunit=0.5,runit=0.5}
\psset{arrowsize=2pt 2,arrowinset=0.2}
\psline(1,5)(1,4)
\psarc(1,3){1}{0}{180}
\psline(0,3)(0,2)
\psarc(1,2){1}{180}{360}
\psline(1,1)(1,0)
\pscircle[linewidth=0.4pt,fillstyle=solid,fillcolor=white](1,1){0.2}
\pscircle[linewidth=0.4pt,fillstyle=solid,fillcolor=white](1,4){0.2}
\rput(1,1.75){m}
\rput(1,3.25){$\Delta$}
\rput(2,2.5){S}
\endpspicture
\qquad\qquad
\text{U}=\eta\circ\epsilon\,=\,\Id_\openk 
\end{align}

The following two theorems show that the antipode is intimately related
with the notion of opposite products and opposite co-products. In fact,
this establishes a further {\it duality\/} connecting the four
spaces $H,H^{op},H^{*},H^{*op}$.

\begin{thrm}
If $B$ is a {\em biconnected bialgebra}, then the co-product diagram
\begin{align}
\begin{array}{c@{\hskip 2truecm}c}
\Rnode{1}{B} & \Rnode{2}{B\otimes B} \\[4ex]
& \Rnode{3}{B\otimes B} \\[4ex]
\Rnode{4}{B} & \Rnode{5}{B\otimes B}
\ncline{->}{1}{2}
\Aput{\Delta}
\ncline{->}{1}{4}
\Aput{S}
\ncline{->}{4}{5}
\Aput{\Delta}
\ncline{->}{2}{3}
\Aput{\hat{\tau}}
\ncline{->}{3}{5}
\Aput{S\otimes S}
\end{array}
&&\quad
\pspicture[0.5](0,0)(1,3)
\psset{linewidth=\pstlw,xunit=0.5,yunit=0.5,runit=0.5}
\psset{arrowsize=2pt 2,arrowinset=0.2}
\psline(0,0)(0,2)
\psline(2,0)(2,2)
\psline(1,3)(1,6)
\psarc(1,2){1}{0}{180}
\pscircle[linewidth=0.4pt,fillstyle=solid,fillcolor=white](1,3){0.2}
\rput*(1,5){S}
\endpspicture
\quad=\quad
\pspicture[0.5](0,0)(1,3)
\psset{linewidth=\pstlw,xunit=0.5,yunit=0.5,runit=0.5}
\psset{arrowsize=2pt 2,arrowinset=0.2}
\psline(0,0)(0,2)
\psline(2,0)(2,2)
\psbezier(2,2)(2,3)(0,3)(0,4)
\psbezier[border=4pt,bordercolor=white](0,2)(0,3)(2,3)(2,4)
\psarc(1,4){1}{0}{180}
\psline(1,5)(1,6)
\pscircle[linewidth=0.4pt,fillstyle=solid,fillcolor=white](1,5){0.2}
\rput*(0,1){S}
\rput*(2,1){S}
\endpspicture
\end{align}
is commutative.
\end{thrm}

\begin{thrm}
If $B$ is a {\em biconnected bialgebra}, then the product diagram
\begin{align}
\begin{array}{c@{\hskip 2truecm}c}
\Rnode{1}{B\otimes B} & \Rnode{2}{B} \\[4ex]
\Rnode{3}{B\otimes B} & \\[4ex]
\Rnode{4}{B\otimes B} & \Rnode{5}{B}
\ncline{->}{1}{2}
\Aput{m}
\ncline{->}{2}{5}
\Aput{S}
\ncline{->}{1}{3}
\Aput{S\otimes S}
\ncline{->}{3}{4}
\Aput{\hat{\tau}}
\ncline{->}{4}{5}
\Aput{m}
\end{array}
&&\quad
\pspicture[0.5](0,0)(1,3)
\psset{linewidth=\pstlw,xunit=0.5,yunit=0.5,runit=0.5}
\psset{arrowsize=2pt 2,arrowinset=0.2}
\psline(0,6)(0,4)
\psline(2,6)(2,4)
\psline(1,3)(1,0)
\psarc(1,4){1}{180}{360}
\pscircle[linewidth=0.4pt,fillstyle=solid,fillcolor=white](1,3){0.2}
\rput*(1,1){S}
\endpspicture
\quad=\quad
\pspicture[0.5](0,0)(1,3)
\psset{linewidth=\pstlw,xunit=0.5,yunit=0.5,runit=0.5}
\psset{arrowsize=2pt 2,arrowinset=0.2}
\psline(0,6)(0,4)
\psline(2,6)(2,4)
\psbezier(2,4)(2,3)(0,3)(0,2)
\psbezier[border=4pt,bordercolor=white](0,4)(0,3)(2,3)(2,2)
\psarc(1,2){1}{180}{360}
\psline(1,1)(1,0)
\pscircle[linewidth=0.4pt,fillstyle=solid,fillcolor=white](1,1){0.2}
\rput*(0,5){S}
\rput*(2,5){S}
\endpspicture
\end{align}
is commutative.
\end{thrm}

\begin{thrm}
If $B$ is a biconnected bialgebra where either $m_B$ or $\Delta_B$ is 
commutative, then $S\circ S ~:~ B\rightarrow B$ is the identity 
morphism of $B$, i.e. $S\circ S = \Id_B$ is an involution.
\end{thrm}

The proofs will be discussed together with Kuperberg's Lemma 3.2, see 
below.

\subsection{Hopf algebra definition}

\begin{dfn}[Hopf algebra, Milnor Moore]
\index{Hopf algebra}\index{algebra!Hopf}
A {\em Hopf algebra\/} is an {\it antipodal biconnected bialgebra\/}, 
i.e. a bialgebra which possesses an antipode.
\end{dfn}

\noindent
In fact, this raises a question if any antipodal bialgebra is a
biconnected bialgebra, see also the discussion in section 6 of 
Fauser and Oziewicz \cite{fauser:oziewicz:2001a}. We will later coin
the term `Hopf gebra' as used by Oziewicz which will not imply
connectedness. In fact, the antipode definition etc. does not depend
on connectedness.

The importance of these definitions for topology comes from the
following: Consider the category of augmented co-algebras
with (graded) commutative co-multiplication $\Cog(\eta,\hat{\tau})$.
This category carries a product just by taking the tensor product 
$A\otimes B$ which needs essentially the commutativity of the 
co-multiplication. If $\openk$ is a point in this category, one has the 
morphisms $\eta~:~\openk\rightarrow C$, $\epsilon~:~C \rightarrow
\openk$ which turn the category $\Cog(\eta,\hat{\tau},
\epsilon\circ\eta,\otimes)$ with product into a monoid. We have thus 
seen above, that connected bialgebras with commutative co-multiplication
are groups in the category  $\Cog(\eta,\hat{\tau},\epsilon\circ
\eta,\otimes) = {\bf Hopf}$, i.e. Hopf algebras with involutive antipode.
\index{Hopf algebra!in topology}

If now $\Top_*$ is the category of topological spaces with base point, 
then if $\openk$ is a field, there is a natural functor $H_*(~~,\openk)
~:~\Top_* \rightarrow {\bf Hopf}$ which to every space $X$ assigns its 
singular homology with coefficients in $\openk$. The co-multiplication
$H_*(X,\openk)\rightarrow H_*(X,\openk)\otimes H_*(X,\openk)$ is
the morphism induced by the diagonal map $\Delta~:~X\rightarrow
X\otimes X$. This was the starting point of Hopf \cite{hopf:1941a}
and motivated the works of Milnor and Moore \cite{milnor:moore:1965a},
Kuperberg \cite{kuperberg:1991a,kuperberg:1996a} and others. 

Further notions like integrals will be defined below where they are 
explored in some examples.

%% file: hopf_gebra.tex
\chapter{Hopf gebras}

In this chapter we will develop the theory of {\it Hopf gebras\/} as 
opposed to that of {\it Hopf al-gebras\/} which was developed in the preceeding 
chapter. To some extend, Hopf gebras and Hopf al-gebras are 
equivalent, but it will turn out, that the notion of Hopf gebra
allows a much clearer genealogy of al- and co-gebraic structures.

As we saw at various places, the Hopf algebras over a graded $A$-module
are defined by the following structure (tensors): the associative 
multiplication $m_A~:~A\otimes A\rightarrow A$, the 
unit $\eta~:~\openk\rightarrow A$ 
with $m(x\otimes \eta(1)) = x = m(\eta(1)\otimes x)$ $\forall x\in A$,
the associative co-multiplication $\Delta_A~:~A\rightarrow A\otimes A$,
the counit $\epsilon~:~A\rightarrow \openk$
with $(\Id\otimes \epsilon)\circ\Delta_A = \Id = 
      (\epsilon\otimes \Id)\circ\Delta_A$,
the antipode $S~:~A\rightarrow A$, an antihomomorphism and finally
the crossing $\hat{\tau}$. One can summarize this as 
$H(A,m_A,\Delta_A,\eta,\epsilon;S;\hat{\tau})$. 

The question which will be the guiding principle in this section is:
Are the structure tensors independent? Already the presentation
of Hopf algebras in the last chapter showed, that topological
requirements as connectedness or the splitting of certain exact
sequences played an important role to be able to show the
existence of inverses which turned the convolution into a group.
\index{tensors!structure!independence of}

The idea is not to start from algebras and co-algebras, but to 
take possibly non-unital and non-associative products and
co-products to form a convolution algebra. Then it is a naturally
given way to add structures unless one arrives at a Hopf gebra.
We will restrict this general setting by assuming bi-associativity,
i.e. associative products and co-associative co-products. Furthermore,
we assume here that the product and co-product are endomorphic, so that
source and target are the same $\openk$-module $A$.

The fact that an antipodal convolution is already a Hopf gebra, follows
from the theorem on the crossing derived by Oziewicz 
\cite{oziewicz:1997a,oziewicz:guzman:2001a}, see below. One finds
that an antipodal convolution  has a unique crossing derived from 
the antipode since the antipode is unique. This idea will be 
generalized in the next chapter using convolutive idempotents. 

In \cite{fauser:oziewicz:2001a} co-convolutions have been
introduced. It is clear that a convolution algebra is turned by
categorial duality into a co-convolution co-algebra, however, we will 
not develop a theory of co-convolutions here.
\index{co-convolution co-gebra}

A map of our further development might help to see how we proceed to
obtain Hopf gebras.
\index{convolution!genealogy of}

\begin{align}
\begin{array}{c@{\hskip 0.15truecm}c@{\hskip 0.15truecm}c}
\Rnode{1}{\text{al-gebra, $(A,m_A)$}} & & 
\Rnode{2}{\text{co-gebra, $(C,\Delta_C)$}} \\[6ex]
& \Rnode{3}{\text{convolution algebra, $\Conv(C,A)$}} & \\[9ex]
\Rnode{4}{\text{unital convolution}} &
\Rnode{5}{\text{non-unital bigebra}}  &
\Rnode{6}{\text{non-unital convolution}} \\[15ex]
\Rnode{7}{\parbox[t]{4truecm}{antipodal Convolution\\i.e. Hopf gebra}} &
\Rnode{8}{\text{unital bigebra}} &
\Rnode{9}{\text{non-antipodal convolution}} 
\ncline{<->}{1}{2}
\Aput{\text{categorial duality}}
\ncline{->}{1}{3}
\ncline{->}{2}{3}
\ncline{->}{3}{4}
\mput*{\exists U}
\ncline{->}{3}{5}
\mput*{\not\exists U,\,\exists\hat{\tau}}
\ncline{->}{3}{6}
\mput*{\not\exists U,\,\not\exists\hat{\tau}}
\ncline{->}{4}{7}
\Bput{\exists S}
\ncline{->}{4}{8}
\mput*{\not\exists S,\,\exists\hat{\tau}}
\ncline{->}{4}{9}
\mput*{\not\exists S,\,\not\exists\hat{\tau}}
\end{array}
\end{align}

We will derive some facts about the leftmost trail down
to the Hopf gebras. There will not be an opportunity to go into the
details of the other structures which occur in the above displayed
diagram. The classification proposed here is different to that of
the previous chapter. A bi-associative bi-connected Hopf gebra
would fulfil the axioms of a Hopf algebra. We use the term
Hopf gebra as synonym for bi-associative, not necessarily
bi-connected antipodeal convolution. This is a bi-algebra as we will 
show below.

{\bf Note} that at every point the arrows describe {\it yes\/} or 
{\it no\/} questions which renders the structures in a single line
as being disjoint. Hence a Hopf gebra is {\it not\/} a unital bigebra,
etc. This is a major difference to the usual treatment of Hopf 
al-gebras, where a Hopf algebra is in the same time a bialgebra, and the
notions there are inclusive and not exclusive.

\section{Cup and cap tangles}

\subsection{Evaluation and co-evaluation}

While in the alphabet for knots and links, Eqn. \ref{eqn-TGL}, cup and 
cap tangles already occurred, the convolution alphabet, Eqn. \ref{eqn-PCP},
does not contain such tangles.
\index{tangle!cup}
\index{tangle!cap}

However, dealing with endomorphisms, we have implicitly assumed an 
action of a dual space since
$f \in \End V \cong \Hom(V,V) \cong V\otimes V^*$. Then, the action of $V^*$
on $V$ is described by a cup tangle, the {\it evaluation map\/} of type
\index{map!evaluation}
\index{evaluation map}
$2\rightarrow 0$ denoted as $eval~:~ V^*\otimes V \rightarrow \openk$. $V$ 
or elements $x$ of $V$ are represented by downward pointing arrows, 
while $V^*$ or elements $\omega$ from $V^*$ are drawn as upward pointing 
arrows, i.e. we use oriented lines which already occurred in the Kuperberg 
graphical method.
\begin{align}
\pspicture[0.5](0,0)(1,1.5)
\psset{linewidth=\pstlw,xunit=0.5,yunit=0.5,runit=0.5}
\psline(0,3)(0,2)
\psline(2,3)(2,2)
\psline{<-}(0,2.5)(0,2.25)
\psline{->}(2,2.5)(2,2.25)
\psarc(1,2){1}{180}{360}
\pscircle[linewidth=0.4pt,fillstyle=solid,fillcolor=black](1,1){0.2}
\rput(1,0.25){eval}
\endpspicture
\quad\quad
\pspicture[0.5](0,0)(1,1.5)
\psset{linewidth=\pstlw,xunit=0.5,yunit=0.5,runit=0.5}
\psline(0,3)(0,2)
\psline(2,3)(2,2)
\psline{->}(0,2.5)(0,2.25)
\psline{<-}(2,2.5)(2,2.25)
\psarc(1,2){1}{180}{360}
\pscircle[linewidth=0.4pt,fillstyle=solid,fillcolor=black](1,1){0.2}
\rput(1,0.25){eval}
\endpspicture
\qquad&~
\pspicture[0.5](0,0)(1,2.5)
\psset{linewidth=\pstlw,xunit=0.5,yunit=0.5,runit=0.5}
\psline(1,5)(1,0)
\psline{->}(1,3.6)(1,3.5)
\psline{->}(1,1.6)(1,1.5)
\pscircle[linewidth=0.4pt,fillstyle=solid,fillcolor=black](1,2.5){0.2}
\rput(0.25,2.5){$f$}
\endpspicture
\quad\cong\quad
\pspicture[0.5](0,0)(2,2.5)
\psset{linewidth=\pstlw,xunit=0.5,yunit=0.5,runit=0.5}
\psline(4,5)(4,3)
\psline{->}(4,4.25)(4,4)
\psarc(3,3){1}{270}{360}
\psbezier(1,3)(2.1,3)(1.8,2)(3,2)
\psline{->}(2.05,2.45)(1.95,2.55)
\psarc(1,2){1}{90}{180}
\psline(0,2)(0,0)
\psline{->}(0,1.25)(0,1)
\pscircle[linewidth=0.4pt,fillstyle=solid,fillcolor=black](1,3){0.2}
\pscircle[linewidth=0.4pt,fillstyle=solid,fillcolor=black](3,2){0.2}
\rput(3,1.25){eval}
\rput(1,3.75){$f$}
\endpspicture  
\qquad
\parbox[t]{3truecm}{left action by\\ evaluation}                     
\end{align}
Co-evaluation is displayed by a cap tangle. This can be seen
considering the identity map in the above given description:
\index{tangle!cap}
\index{co-evaluation}
\begin{align}
\pspicture[0.5](0,0)(1,1.5)
\psset{linewidth=\pstlw,xunit=0.5,yunit=0.5,runit=0.5}
\psarc(1,1){1}{0}{180}
\psline(0,1)(0,0)
\psline(2,1)(2,0)
\psline{<-}(0,0.5)(0,0.25)
\psline{->}(2,0.5)(2,0.25)
\pscircle[linewidth=0.4pt,fillstyle=solid,fillcolor=black](1,2){0.2}
\rput(1,2.75){coeval}
\endpspicture
\quad~\quad
\pspicture[0.5](0,0)(1,1.5)
\psset{linewidth=\pstlw,xunit=0.5,yunit=0.5,runit=0.5}
\psarc(1,1){1}{0}{180}
\psline(0,1)(0,0)
\psline(2,1)(2,0)
\psline{->}(0,0.5)(0,0.25)
\psline{<-}(2,0.5)(2,0.25)
\pscircle[linewidth=0.4pt,fillstyle=solid,fillcolor=black](1,2){0.2}
\rput(1,2.75){coeval}
\endpspicture       
\qquad&~
\pspicture[0.5](0,0)(1,2.5)
\psset{linewidth=\pstlw,xunit=0.5,yunit=0.5,runit=0.5}
\psline(1,5)(1,0)
\psline{->}(1,2.6)(1,2.5)
\endpspicture
\quad=\quad
\pspicture[0.5](0,0)(2,2.5)
\psset{linewidth=\pstlw,xunit=0.5,yunit=0.5,runit=0.5}
\psline(4,5)(4,3)
\psline{->}(4,4.25)(4,4)
\psarc(3,3){1}{270}{360}
\psbezier(1,3)(2.1,3)(1.8,2)(3,2)
\psline{->}(2.05,2.45)(1.95,2.55)
\psarc(1,2){1}{90}{180}
\psline(0,2)(0,0)
\psline{->}(0,1.25)(0,1)
\pscircle[linewidth=0.4pt,fillstyle=solid,fillcolor=black](1,3){0.2}
\pscircle[linewidth=0.4pt,fillstyle=solid,fillcolor=black](3,2){0.2}
\rput(3,1.25){eval}
\rput(1,3.75){coeval}
\endpspicture          
\end{align}
Cup and cap tangles constitute a so called {\it closed structure\/}
\cite{kelly:laplaza:1980a,lyubashenko:1995a}.
\index{closed structure}

Having introduced an irrelevant basis $\{e_i\}$ in $V$ and a canonical
dual basis  $\{\epsilon^j\}$ in $V^*$, i.e. $\epsilon^i(e_j)=\delta^i_j$,
one can easily compute the action of $\omega \in V^*$ on $v \in V$ using 
$(\omega_i\epsilon^i)(v^j e_j)= \omega_i v^j \epsilon^i(e_j)
= \omega_i v^j \delta^i_j =  \omega_i v^i$. However, this relation does
{\it not\/} introduce a duality operation $*(e_i) = \epsilon^i$. This
isomorphism, simply keeping the coefficients, is called Euclidean dual 
isomorphism \cite{saller:1993a,saller:1994a}. If we introduce Gra{\ss}mann 
exterior\index{dual isomorphism!Euclidean}
algebras one has to define the action of a multi-co-vector $\in \bigvee V^*$
on multivectors $\in \bigwedge V$, where we used a {\it vee\/} $\vee$
to denote the exterior product of co-vectors. Following standard conventions
\cite{bourbaki:1989a,sweedler:1969a,%
doubilet:rota:stein:1974a,rota:stein:1994a}
\index{pairing!by evaluation}
one introduces the following {\it pairing\/} on homogenous elements 
(extensors) and extends it by bilinearity
\begin{align}
&\langle\,\mid\,\rangle ~:~ 
\bigvee V^* \times \bigwedge V \rightarrow \openk \nn
&\langle V^*\mid V\rangle~\cong \text{eval} \nn
&\langle\omega_1\vee\ldots\vee\omega_n\mid x_1\wedge\ldots\wedge x_m\rangle =
\left\{ 
\begin{array}{ll}
\pm \det(\langle\omega_i\mid x_j\rangle) & \text{ if $n=m$}\\
0 & \text{otherwise}
\end{array}
\right.
\end{align}
The $\pm$ sign has to be arranged due to the involved permutations.
We use sometime a slightly different setting, where the indices in the
first argument are in reversed order and no sign occurs in front of det. 
Using this construction the space underlying the Gra{\ss}mann algebra 
can be turned into a Hilbert space \cite{depills:1968a}.
In fact this is nothing but the Laplace expansion. In Hopf algebraic 
terms, see \cite{rota:stein:1994a}, the pairing can be expanded as
\index{Laplace expansion}
\begin{align}
\langle \omega^\prime\vee\omega^\pprime\mid x\wedge y\rangle
&= \langle \omega^\prime\otimes\omega^\pprime\mid 
\Delta_\vee(x\wedge y)\rangle \nn
&= \langle \omega^\prime\mid (x\wedge y)_{(1)}\rangle 
   \langle \omega^\pprime\mid (x\wedge y)_{(2)}\rangle \nn
&= \langle \omega^\prime\mid x_{(1)}\rangle
   \langle \omega^\pprime\mid y_{(2)}\rangle
+(-)^{\partial x\partial y}
   \langle \omega^\prime\mid y_{(1)}\rangle
   \langle \omega^\pprime\mid x_{(2)}\rangle
\end{align}
Note that since wedge and vee are independent, $(\Delta_\vee,\wedge)$ is 
also a pair of an independent Gra{\ss}mann co-product and product. If
$\omega^\prime, \omega^\pprime \in V^*$ and $x,y\in V$, this is 
the particular case of a $2\times2$-determinant.

\subsection{Scalar and co-scalar products}

To be able to introduce Clifford algebras and Clifford co-gebras, we
need to introduce scalar and co-scalar products. Let $B\in V^*\otimes V^*$ 
be a scalar product and $D\in V\otimes V$ be a co-scalar product:
\index{scalar product}
\index{co-scalar product}
\index{product!scalar}
\index{product!co-scalar}
\begin{align}
\Rnode{1}{V}\hskip 2truecm
\Rnode{2}{V^*}
\ncline[offset=1ex]{->}{1}{2}
\Aput{B}
\ncline[offset=1ex]{->}{2}{1}
\Aput{D}
\end{align}
Equivalently using the action by evaluation this can be denoted as
\begin{align}
\Rnode{1}{V\otimes V}\hskip 2truecm
\Rnode{2}{\openk}\hskip 2truecm
\Rnode{3}{V^*\otimes V^*}
\ncline{->}{1}{2}
\Aput{B}
\ncline{->}{3}{2}
\Bput{D}
\end{align}
These actions are depicted also by cup and cap tangles, but with two
ingoing and two outgoing lines. The tangles are also decorated by the 
map in use
\index{tangle!cup}
\begin{align}
\pspicture[0.5](0,0)(1,2)
\psset{linewidth=\pstlw,xunit=0.5,yunit=0.5,runit=0.5}
\psline(0,4)(0,2)
\psline(2,4)(2,2)
\psline{->}(0,2.5)(0,2.25)
\psline{->}(2,2.5)(2,2.25)
\psarc(1,2){1}{180}{360}
\pscircle[linewidth=0.4pt,fillstyle=solid,fillcolor=black](1,1){0.2}
\rput(1,0.25){B}
\endpspicture
\quad\cong\quad
\pspicture[0.5](0,0)(1,2)
\psset{linewidth=\pstlw,xunit=0.5,yunit=0.5,runit=0.5}
\psline(0,4)(0,2)
\psline(2,4)(2,2)
\psline{->}(0,2.5)(0,2.25)
\psline{->}(2,3.7)(2,3.5)
\psline{<-}(2,2.5)(2,2.25)
\psarc(1,2){1}{180}{360}
\pscircle[linewidth=0.4pt,fillstyle=solid,fillcolor=black](1,1){0.2}
\pscircle[linewidth=0.4pt,fillstyle=solid,fillcolor=black](2,3){0.2}
\rput(1.25,3){B}
\rput(1,0.25){eval}
\endpspicture  
\qquad&&~
\pspicture[0.5](0,0)(1,2)
\psset{linewidth=\pstlw,xunit=0.5,yunit=0.5,runit=0.5}
\psline(0,4)(0,2)
\psline(2,4)(2,2)
\psline{<-}(0,2.5)(0,2.25)
\psline{<-}(2,2.5)(2,2.25)
\psarc(1,2){1}{180}{360}
\pscircle[linewidth=0.4pt,fillstyle=solid,fillcolor=black](1,1){0.2}
\rput(1,0.25){D}
\endpspicture
\quad\cong\quad
\pspicture[0.5](0,0)(1,2)
\psset{linewidth=\pstlw,xunit=0.5,yunit=0.5,runit=0.5}
\psline(0,4)(0,2)
\psline(2,4)(2,2)
\psline{<-}(0,2.5)(0,2.25)
\psline{<-}(2,3.75)(2,3.5)
\psline{->}(2,2.5)(2,2.25)
\psarc(1,2){1}{180}{360}
\pscircle[linewidth=0.4pt,fillstyle=solid,fillcolor=black](1,1){0.2}
\pscircle[linewidth=0.4pt,fillstyle=solid,fillcolor=black](2,3){0.2}
\rput(1.25,3){D}
\rput(1,0.25){eval}
\endpspicture             
\end{align}
Using categorial duality one introduces the corresponding cap tangles
\index{tangle!cap}
\begin{align}
\Rnode{1}{V\otimes V}\hskip 2truecm
\Rnode{2}{\openk}\hskip 2truecm
\Rnode{3}{V^*\otimes V^*}
\ncline{<-}{1}{2}
\Aput{C}
\ncline{<-}{3}{2}
\Bput{E}
\end{align}
\begin{align}
\pspicture[0.5](0,-0.5)(1,2)
\psset{linewidth=\pstlw,xunit=0.5,yunit=0.5,runit=0.5}
\psarc(1,1){1}{0}{180}
\psline(0,1)(0,0)
\psline(2,1)(2,0)
\psline{->}(0,0.5)(0,0.25)
\psline{->}(2,0.5)(2,0.25)
\pscircle[linewidth=0.4pt,fillstyle=solid,fillcolor=black](1,2){0.2}
\rput(1,2.75){C}
\endpspicture
\qquad&&~
\pspicture[0.5](0,-0.5)(1,2)
\psset{linewidth=\pstlw,xunit=0.5,yunit=0.5,runit=0.5}
\psline(0,1)(0,0)
\psline(2,1)(2,0)
\psline{<-}(0,0.5)(0,0.25)
\psline{<-}(2,0.5)(2,0.25)
\psarc(1,1){1}{0}{180}
\pscircle[linewidth=0.4pt,fillstyle=solid,fillcolor=black](1,2){0.2}
\rput(1,2.75){E}
\endpspicture                     
\end{align}
Observe, that $C^{-1} \not= D$ and $E^{-1}\not= B$ in general. That
means, that also the Reidemeister moves are not in general valid
and the present tangles are not 'knottish', e.g.:
\begin{align}
\pspicture[0.5](0,0)(2,2.5)
\psset{linewidth=\pstlw,xunit=0.5,yunit=0.5,runit=0.5}
\psline(0,5)(0,3)
\psline{->}(0,4.25)(0,4)
\psarc(1,3){1}{180}{270}
\psbezier(1,2)(2.2,2)(1.7,3)(3,3)
\psline{<-}(1.95,2.45)(2.05,2.55)
\psarc(3,2){1}{0}{90}
\psline(4,2)(4,0)
\psline{->}(4,1.25)(4,1)
\pscircle[linewidth=0.4pt,fillstyle=solid,fillcolor=black](1,2){0.2}
\pscircle[linewidth=0.4pt,fillstyle=solid,fillcolor=black](3,3){0.2}
\rput(1,1.25){B}
\rput(3,3.75){C}
\endpspicture
\quad&\not=\quad
\pspicture[0.5](0,0)(1,2.5)
\psset{linewidth=\pstlw,xunit=0.5,yunit=0.5,runit=0.5}
\psline(1,5)(1,0)
\psline{->}(1,2.6)(1,2.5)
\endpspicture
\quad\not=\quad  
\pspicture[0.5](0,0)(2,2.5)
\psset{linewidth=\pstlw,xunit=0.5,yunit=0.5,runit=0.5}
\psline(4,5)(4,3)
\psline{->}(4,4.25)(4,4)
\psarc(3,3){1}{270}{360}
\psbezier(1,3)(2.1,3)(1.8,2)(3,2)
\psline{<-}(2.05,2.45)(1.95,2.55)
\psarc(1,2){1}{90}{180}
\psline(0,2)(0,0)
\psline{->}(0,1.25)(0,1)
\pscircle[linewidth=0.4pt,fillstyle=solid,fillcolor=black](1,3){0.2}
\pscircle[linewidth=0.4pt,fillstyle=solid,fillcolor=black](3,2){0.2}
\rput(3,1.25){B}
\rput(1,3.75){C}
\endpspicture
\end{align}        
The condition for Reidemeister moves to hold is hence $C\circ B = Id_V$
and $D\circ E = \Id_{V^*}$. In the case of Clifford products, we will 
learn, that this condition prevents the existence of an antipode.

\subsection{Induced graded scalar and co-scalar products}

Till now, we have not made any assumptions about the scalar and co-scalar
products. However, since we will deal mainly with Clifford algebras,
our particular scalar and co-scalar products will have a quite special
structure.
\index{scalar product!induced}
\index{co-scalar product!induced}

It is convenient to introduce the following rules for a {\it pairing\/}
of $n$-co-tensors on $n$-tensors
\begin{align}
&\langle \rnode{A}{V^*} \mid \rnode{B}{V}\rangle \rightarrow\openk
\ncbar[linewidth=\pstlw,nodesep=5pt,angle=-90,arm=0.3,linearc=0.2]{-}{A}{B}
\nonumber\\[2ex]
&\langle \rnode{C}{V^*}\otimes \ldots \otimes \rnode{D}{V^*} \mid
        \rnode{E}{V}  \otimes \ldots \otimes \rnode{F}V \rangle 
\rightarrow \openk
\ncbar[linewidth=\pstlw,nodesep=5pt,angle=-90,arm=0.3,linearc=0.2]{-}{D}{E}
\ncbar[linewidth=\pstlw,nodesep=5pt,angle=-90,arm=0.5,linearc=0.2]{-}{C}{F}   
\end{align}
\vskip 3ex
\noindent
and to agree that the pairing of $m$-co-tensors on $n$-tensors is
zero for $m\not=n$. This notion suggests that we index $n$-co-tensors
in a reverse way as $n$-tensors, i.e.
\begin{align}
\langle \omega_n\otimes\ldots\otimes \omega_1\mid 
         x_1 \otimes\ldots\otimes x_n \rangle &=
\langle \omega_1\mid x_1\rangle \ldots \langle \omega_n\mid x_n\rangle\,.
\end{align}
This setting reflects the practice in the theory of knots and links,
where the tangles are closed by cup and cap tangles of adjacent open
ends of a braid to form a knot or link. The above pairing between
antisymmetric exterior products is the Gra{\ss}mann Hopf version of 
this definition. 
\index{pairing!of tensors}

Now let us turn to the case of Gra{\ss}mann Hopf algebras. 
Let $V$ be a (finitely generated) $\openk$-module and $\bigwedge V$
the Gra{\ss}mann algebra built over this space. Furthermore let
$V^*$ be the dual space and $\bigvee V^*$ the Gra{\ss}mann algebra
over that space. Usually one introduces there a bilinear form,
i.e. a scalar product, $B:~V\otimes V\rightarrow \openk$, or a
bilinear form, i.e. a co-scalar product, $C:~V^*\otimes V^*\rightarrow 
\openk$. The question arrises, in which way the bilinear forms 
are lifted to the whole space $\bigwedge V$ or $\bigvee V^*$.
Let us denote this lifted scalar and co-scalar products by
$B^{\wedge}$ and $C^{\vee}$

Remember that $B \in V^*\otimes V^*$, and that $B^{\wedge}$ will
live in $\bigvee V^* \otimes \bigvee V^*$. We require that this 
extension is a graded morphism $B^{\wedge} \in \Hom(\bigwedge V
\otimes \bigwedge V,\openk)$. The required extension can be given, 
see Oziewicz \cite{oziewicz:1997a,oziewicz:2001a} and 
\cite{stumpf:borne:1994a,fauser:stumpf:1997a,borne:lochak:stumpf:2001a}, as:
\index{tangle!induced scalar product}
\begin{align}
B^{\wedge} &= \exp(B) =
\epsilon\otimes\epsilon + B_{ij}\epsilon^i\otimes \epsilon^j
+ B_{[i_1i_2],[j_1j_2]} 
\epsilon^{i_1}\wedge\epsilon^{i_2} \otimes \epsilon^{j_1}\wedge\epsilon^{j_2}
+ \ldots 
\nonumber \\[3ex]
\pspicture[0.5](0,0)(1,2)
\psset{linewidth=\pstlw,xunit=0.5,yunit=0.5,runit=0.5}
\psset{arrowsize=2pt 2,arrowinset=0.2}
\psline(0,4)(0,1)
\psline(2,4)(2,1)
\psline{->}(0,4)(0,3.5)
\psline{->}(2,4)(2,3.5)
\psarc(1,1){1}{180}{360}
\pscircle[linewidth=0.4pt,fillstyle=solid,fillcolor=black](1,0){0.2}
\rput(1,0.75){$B^{\wedge}$}
\endpspicture
\quad&=\quad
\pspicture[0.5](0,0)(1,2)
\psset{linewidth=\pstlw,xunit=0.5,yunit=0.5,runit=0.5}
\psset{arrowsize=2pt 2,arrowinset=0.2}
\psline(0,4)(0,0)
\psline(2,4)(2,0)
\psline{->}(0,4)(0,3.5)
\psline{->}(2,4)(2,3.5)
\pscircle[linewidth=0.4pt,fillstyle=solid,fillcolor=white](0,0){0.2}
\pscircle[linewidth=0.4pt,fillstyle=solid,fillcolor=white](2,0){0.2}
\endpspicture
\,\oplus\,
\pspicture[0.5](0,0)(1,2)
\psset{linewidth=\pstlw,xunit=0.5,yunit=0.5,runit=0.5}
\psset{arrowsize=2pt 2,arrowinset=0.2}
\psline(0,4)(0,1)
\psline(2,4)(2,1)
\psline{->}(0,4)(0,3.5)
\psline{->}(2,4)(2,3.5)
\psarc(1,1){1}{180}{360}
\pscircle[linewidth=0.4pt,fillstyle=solid,fillcolor=black](1,0){0.2}
\rput(1,0.75){$B$}
\endpspicture
\,\oplus\,
\frac{1}{2!}\,
\pspicture[0.5](0,0)(3,2.5)
\psset{linewidth=\pstlw,xunit=0.5,yunit=0.5,runit=0.5}
\psset{arrowsize=2pt 2,arrowinset=0.2}
\psline(1,5)(1,4)
\psline(5,5)(5,4)
\psline{->}(1,5)(1,4.5)
\psline{->}(5,5)(5,4.5)
\psarc(1,3){1}{0}{180}
\psarc(5,3){1}{0}{180}
\psbezier(2,3)(2,1.5)(2,1.5)(3,1.5)
\psbezier(4,3)(4,1.5)(4,1.5)(3,1.5)
\psbezier(0,3)(0,0)(0,0)(3,0)
\psbezier(6,3)(6,0)(6,0)(3,0)
\pscircle[linewidth=0.4pt,fillstyle=solid,fillcolor=white](1,4){0.2}
\pscircle[linewidth=0.4pt,fillstyle=solid,fillcolor=white](5,4){0.2}
\pscircle[linewidth=0.4pt,fillstyle=solid,fillcolor=black](3,1.5){0.2}
\pscircle[linewidth=0.4pt,fillstyle=solid,fillcolor=black](3,0){0.2}
\rput(3,0.75){$B$}
\rput(3,2.25){$B$}
\endpspicture
\,\oplus\,
\frac{1}{3!}\,
\pspicture[0.5](0,0)(4,3)
\psset{linewidth=\pstlw,xunit=0.5,yunit=0.5,runit=0.5}
\psset{arrowsize=2pt 2,arrowinset=0.2}
\psline(2,6)(2,5)
\psline(6,6)(6,5)
\psline{->}(2,6)(2,5.5)
\psline{->}(6,6)(6,5.5)
\psarc(2,4){1}{0}{180}
\psarc(1,3){1}{0}{180}
\psarc(6,4){1}{0}{180}
\psarc(7,3){1}{0}{180}
\psbezier(3,4)(3,3)(3,3)(4,3)
\psbezier(5,4)(5,3)(5,3)(4,3)
\psbezier(2,3)(2,1.5)(2,1.5)(4,1.5)
\psbezier(6,3)(6,1.5)(6,1.5)(4,1.5)
\psbezier(0,3)(0,0)(0,0)(4,0)
\psbezier(8,3)(8,0)(8,0)(4,0)
\pscircle[linewidth=0.4pt,fillstyle=solid,fillcolor=white](1,4){0.2}
\pscircle[linewidth=0.4pt,fillstyle=solid,fillcolor=white](2,5){0.2}
\pscircle[linewidth=0.4pt,fillstyle=solid,fillcolor=white](6,5){0.2}
\pscircle[linewidth=0.4pt,fillstyle=solid,fillcolor=white](7,4){0.2}
\pscircle[linewidth=0.4pt,fillstyle=solid,fillcolor=black](4,3){0.2}
\pscircle[linewidth=0.4pt,fillstyle=solid,fillcolor=black](4,1.5){0.2}
\pscircle[linewidth=0.4pt,fillstyle=solid,fillcolor=black](4,0){0.2}
\rput(4,0.75){$B$}
\rput(4,2.25){$B$}
\rput(4,3.75){$B$}
\endpspicture
\quad\ldots
\end{align}
And the same holds true for co-scalar products:
\index{tangle!induced co-scalar products}
\begin{align}
C^{\vee} &= \exp(C) =
\eta\otimes\eta + C^{ij}e_i\otimes e_j
+ C^{[i_1i_2],[j_1j_2]} 
e_{i_1}\wedge e_{i_2} \otimes e_{j_1}\wedge e_{j_2}
+ \ldots
\nonumber
\end{align}
\begin{align}
\pspicture[0.5](0,0)(1,2)
\psset{linewidth=\pstlw,xunit=0.5,yunit=0.5,runit=0.5}
\psset{arrowsize=2pt 2,arrowinset=0.2}
\psline(0,0)(0,3)
\psline(2,0)(2,3)
\psline{->}(0,0.5)(0,0)
\psline{->}(2,0.5)(2,0)
\psarc(1,3){1}{0}{180}
\pscircle[linewidth=0.4pt,fillstyle=solid,fillcolor=black](1,4){0.2}
\rput(1,3.25){$C^{\vee}$}
\endpspicture
\quad&=\quad
\pspicture[0.5](0,0)(1,2)
\psset{linewidth=\pstlw,xunit=0.5,yunit=0.5,runit=0.5}
\psset{arrowsize=2pt 2,arrowinset=0.2}
\psline(0,0)(0,4)
\psline(2,0)(2,4)
\psline{->}(0,0.5)(0,0)
\psline{->}(2,0.5)(2,0)
\pscircle[linewidth=0.4pt,fillstyle=solid,fillcolor=white](0,4){0.2}
\pscircle[linewidth=0.4pt,fillstyle=solid,fillcolor=white](2,4){0.2}
\endpspicture
\,\oplus\,
\pspicture[0.5](0,0)(1,2)
\psset{linewidth=\pstlw,xunit=0.5,yunit=0.5,runit=0.5}
\psset{arrowsize=2pt 2,arrowinset=0.2}
\psline(0,0)(0,3)
\psline(2,0)(2,3)
\psline{->}(0,0.5)(0,0)
\psline{->}(2,0.5)(2,0)
\psarc(1,3){1}{0}{180}
\pscircle[linewidth=0.4pt,fillstyle=solid,fillcolor=black](1,4){0.2}
\rput(1,3.25){$C$}
\endpspicture
\,\oplus\,
\frac{1}{2!}\,
\pspicture[0.5](0,0)(3,2.5)
\psset{linewidth=\pstlw,xunit=0.5,yunit=0.5,runit=0.5}
\psset{arrowsize=2pt 2,arrowinset=0.2}
\psline(1,0)(1,1)
\psline(5,0)(5,1)
\psline{->}(1,0.5)(1,0)
\psline{->}(5,0.5)(5,0)
\psarc(1,2){1}{180}{360}
\psarc(5,2){1}{180}{360}
\psbezier(2,2)(2,3.5)(2,3.5)(3,3.5)
\psbezier(4,2)(4,3.5)(4,3.5)(3,3.5)
\psbezier(0,2)(0,5)(0,5)(3,5)
\psbezier(6,2)(6,5)(6,5)(3,5)
\pscircle[linewidth=0.4pt,fillstyle=solid,fillcolor=white](1,1){0.2}
\pscircle[linewidth=0.4pt,fillstyle=solid,fillcolor=white](5,1){0.2}
\pscircle[linewidth=0.4pt,fillstyle=solid,fillcolor=black](3,3.5){0.2}
\pscircle[linewidth=0.4pt,fillstyle=solid,fillcolor=black](3,5){0.2}
\rput(3,4.25){$C$}
\rput(3,2.25){$C$}
\endpspicture
\,\oplus\,
\frac{1}{3!}\,
\pspicture[0.5](0,0)(4,3)
\psset{linewidth=\pstlw,xunit=0.5,yunit=0.5,runit=0.5}
\psset{arrowsize=2pt 2,arrowinset=0.2}
\psline(2,0)(2,1)
\psline(6,0)(6,1)
\psline{->}(2,0.5)(2,0)
\psline{->}(6,0.5)(6,0)
\psarc(2,2){1}{180}{360}
\psarc(1,3){1}{180}{360}
\psarc(6,2){1}{180}{360}
\psarc(7,3){1}{180}{360}
\psbezier(3,2)(3,3)(3,3)(4,3)
\psbezier(5,2)(5,3)(5,3)(4,3)
\psbezier(2,3)(2,4.5)(2,4.5)(4,4.5)
\psbezier(6,3)(6,4.5)(6,4.5)(4,4.5)
\psbezier(0,3)(0,6)(0,6)(4,6)
\psbezier(8,3)(8,6)(8,6)(4,6)
\pscircle[linewidth=0.4pt,fillstyle=solid,fillcolor=white](1,2){0.2}
\pscircle[linewidth=0.4pt,fillstyle=solid,fillcolor=white](2,1){0.2}
\pscircle[linewidth=0.4pt,fillstyle=solid,fillcolor=white](6,1){0.2}
\pscircle[linewidth=0.4pt,fillstyle=solid,fillcolor=white](7,2){0.2}
\pscircle[linewidth=0.4pt,fillstyle=solid,fillcolor=black](4,3){0.2}
\pscircle[linewidth=0.4pt,fillstyle=solid,fillcolor=black](4,4.5){0.2}
\pscircle[linewidth=0.4pt,fillstyle=solid,fillcolor=black](4,6){0.2}
\rput(4,5.25){$C$}
\rput(4,3.75){$C$}
\rput(4,2.25){$C$}
\endpspicture
\quad\ldots
\end{align}
\noindent
{\bf Note} that in the r.h.s inner lines $B$ and $C$ act on grade 1 
spaces only. Therefore $B^{\wedge}$ and $C^{\vee}$ are graded extensions
of $B$ and $C$.

The combinatorial factors $1/n!$ are not apparent if one 
has already taken account for the antisymmetry in terms like 
$B_{[i_1,i_2],[j_1j_2]}=B_{i_1j_1}B_{i_2j_2}-B_{i_1j_2}B_{i_2j_1}$. 
We introduce for the scalar and co-scalar product if extended to the whole 
space $\bigwedge V$ or $\bigvee V^*$ also the Sweedler notation
\begin{align}
B^{\wedge} &= B^{\wedge}_{(1)} \otimes B^{\wedge}_{(2)} \nn
C^{\vee} &= C^{\vee}_{(1)} \otimes C^{\vee}_{(2)} \,.
\end{align}
We note furthermore that for the Clifford co-product of $\Id$, based on 
the co-scalar product $C$ defined below, one finds
\begin{align}
\Delta_C(\Id) &= C^{\vee} = C^{\vee}_{(1)} \otimes C^{\vee}_{(2)}. 
\end{align}

\section{Product co-product duality}
\index{product co-product duality}

\subsection{By evaluation}

Having the evaluation established, we can explain the important concept 
of product co-product duality. Observe, that a co-vector might act on a 
{\it product\/} of 2 vectors $\omega(ab)$ and one can ask if the
co-vector can be 'distributed' on $a$ and $b$. Using tangles we obtain
\index{tangle!product co-product duality}
\begin{align}
\pspicture[0.5](0,0)(1.5,2.5)
\psset{linewidth=\pstlw,xunit=0.5,yunit=0.5,runit=0.5}
\psset{arrowsize=2pt 2,arrowinset=0.2}
\psline(0,5)(0,1)
\psline(1,5)(1,2)
\psline(3,5)(3,2)
\psline{<-}(0,4)(0,3.75)
\psline{->}(1,4)(1,3.75)
\psline{->}(3,4)(3,3.75)
\psarc(1,1){1}{180}{360}
\psarc(2,2){1}{180}{360}
\pscircle[linewidth=0.4pt,fillstyle=solid,fillcolor=white](2,1){0.2}
\pscircle[linewidth=0.4pt,fillstyle=solid,fillcolor=black](1,0){0.2}
\endpspicture
\quad=\quad
\pspicture[0.5](0,0)(3,2.5)
\psset{linewidth=\pstlw,xunit=0.5,yunit=0.5,runit=0.5}
\psset{arrowsize=2pt 2,arrowinset=0.2}
\psline(1,5)(1,4)
\psline(4,5)(4,3)
\psline(6,5)(6,2)
\psline(0,3)(0,2)
\psline{<-}(1,4.75)(1,4.5)
\psline{->}(4,4.5)(4,4.25)
\psline{->}(6,4.5)(6,4.25)
\psline(2,0)(4,0)
\psarc(1,3){1}{0}{180}
\psarc(3,3){1}{180}{360}
\psarc(2,2){2}{180}{270}
\psarc(4,2){2}{270}{360}
\pscircle[linewidth=0.4pt,fillstyle=solid,fillcolor=white](1,4){0.2}
\pscircle[linewidth=0.4pt,fillstyle=solid,fillcolor=black](3,0){0.2}
\pscircle[linewidth=0.4pt,fillstyle=solid,fillcolor=black](3,2){0.2}
\endpspicture
\end{align}
That is, one obtains $\omega(ab)=(\text{eval}\otimes\text{eval})
\Delta(\omega)(a\otimes b) = \omega_{(1)}(b)\omega_{(2)}(a)$, 
where the co-product $\Delta=m^*$ is the dualized product. Indeed, 
this can be done the other way around also
\begin{align}
\pspicture[0.5](0,0)(1.5,2.5)
\psset{linewidth=\pstlw,xunit=0.5,yunit=0.5,runit=0.5}
\psset{arrowsize=2pt 2,arrowinset=0.2}
\psline(3,5)(3,1)
\psline(2,5)(2,2)
\psline(0,5)(0,2)
\psline{<-}(3,4)(3,3.75)
\psline{->}(2,4)(2,3.75)
\psline{->}(0,4)(0,3.75)
\psarc(2,1){1}{180}{360}
\psarc(1,2){1}{180}{360}
\pscircle[linewidth=0.4pt,fillstyle=solid,fillcolor=white](1,1){0.2}
\pscircle[linewidth=0.4pt,fillstyle=solid,fillcolor=black](2,0){0.2}
\endpspicture
\quad=\quad
\pspicture[0.5](0,0)(3,2.5)
\psset{linewidth=\pstlw,xunit=0.5,yunit=0.5,runit=0.5}
\psset{arrowsize=2pt 2,arrowinset=0.2}
\psline(5,5)(5,4)
\psline(2,5)(2,3)
\psline(0,5)(0,2)
\psline(6,3)(6,2)
\psline(2,0)(4,0)
\psline{<-}(5,4.75)(5,4.5)
\psline{->}(2,4.5)(2,4.25)
\psline{->}(0,4.5)(0,4.25)
\psline(4,0)(4,0)
\psarc(5,3){1}{0}{180}
\psarc(3,3){1}{180}{360}
\psarc(2,2){2}{180}{270}
\psarc(4,2){2}{270}{360}
\pscircle[linewidth=0.4pt,fillstyle=solid,fillcolor=white](5,4){0.2}
\pscircle[linewidth=0.4pt,fillstyle=solid,fillcolor=black](3,0){0.2}
\pscircle[linewidth=0.4pt,fillstyle=solid,fillcolor=black](3,2){0.2}
\endpspicture
\end{align}
which shows how a product of co-vectors $\omega^\prime\omega^\pprime$
can be distributed over a vector $v$ as $(\omega^\prime\omega^\pprime)(v)
=(\text{eval}\otimes\text{eval})(\omega^\prime\otimes\omega^\pprime)
(v_{(1)}\otimes v_{(2)}) = \omega^\pprime(v_{(1)})\omega^\prime(v_{(2)})$.

In fact, the statement that an algebra over $A$ is dualized by
categorial duality into a co-algebra over $A$ and vice versa
is an equivalent assertion. The importance of this construction
cannot be overemphasized, since the whole theory of determinants,
permanents and their generalizations to supersymmetric spaces can be 
developed from this setting, \cite{grosshans:rota:stein:1987a}.
Furthermore, as we will demonstrate below this type of duality also
{\it yields\/} commutation relations.

\subsection{By scalar products}

Using the evaluation, we compose co-vector spaces and vector spaces,
which made it necessary to put arrows on the tangles. 
Since we have introduced cup and cap tangles for scalar and
co-scalar products, one can proceed to introduce co-products
from products, which are derived from these tangles and where the
entries are of the same type. However, this is no longer a duality in 
the above defined sense since it involves explicitely a scalar or co-scalar 
product. We will see, that an entirely new type of product will occur, 
the {\it contraction\/}. Due to our construction of the scalar product 
$B^{\wedge}$ as a graded morphism, $B^\wedge : \bigwedge V \rightarrow 
\bigvee V^*$, we have the important relation
\begin{align}
B^{\wedge}(a\wedge b) &= B^{\wedge}(a) \vee B^{\wedge}(b).
\end{align}
\index{outermorphism}
\index{scalar product!as graded morphism}
This was called outermorphism by Hestenes and Sobczyk 
\cite{hestenes:sobczyk:1992a}. The tangle equation is once more 
\begin{align}
\pspicture[0.5](0,0)(1.5,2.5)
\psset{linewidth=\pstlw,xunit=0.5,yunit=0.5,runit=0.5}
\psset{arrowsize=2pt 2,arrowinset=0.2}
\psline(0,5)(0,1)
\psline(1,5)(1,2)
\psline(3,5)(3,2)
\psline{->}(0,4)(0,3.75)
\psline{->}(1,4)(1,3.75)
\psline{->}(3,4)(3,3.75)
\psarc(1,1){1}{180}{360}
\psarc(2,2){1}{180}{360}
\pscircle[linewidth=0.4pt,fillstyle=solid,fillcolor=white](2,1){0.2}
\pscircle[linewidth=0.4pt,fillstyle=solid,fillcolor=black](1,0){0.2}
\rput(1,0.75){B}
\rput(2,1.75){$\wedge$}
\endpspicture
\quad=\quad
\pspicture[0.5](0,0)(3,2.5)
\psset{linewidth=\pstlw,xunit=0.5,yunit=0.5,runit=0.5}
\psset{arrowsize=2pt 2,arrowinset=0.2}
\psline(1,5)(1,4)
\psline(4,5)(4,3)
\psline(6,5)(6,2)
\psline(0,3)(0,2)
\psline{->}(1,4.5)(1,4.25)
\psline{->}(4,4.5)(4,4.25)
\psline{->}(6,4.5)(6,4.25)
\psline(2,0)(4,0)
\psarc(1,3){1}{0}{180}
\psarc(3,3){1}{180}{360}
\psarc(2,2){2}{180}{270}
\psarc(4,2){2}{270}{360}
\pscircle[linewidth=0.4pt,fillstyle=solid,fillcolor=white](1,4){0.2}
\pscircle[linewidth=0.4pt,fillstyle=solid,fillcolor=black](3,0){0.2}
\pscircle[linewidth=0.4pt,fillstyle=solid,fillcolor=black](3,2){0.2}
\rput(3,0.75){B}
\rput(3,2.75){B}
\endpspicture
\quad=:\quad
\pspicture[0.5](0,0)(1.5,2.5)
\psset{linewidth=\pstlw,xunit=0.5,yunit=0.5,runit=0.5}
\psset{arrowsize=2pt 2,arrowinset=0.2}
\psline(3,5)(3,1)
\psline(2,5)(2,2)
\psline(0,5)(0,2)
\psline{->}(3,4)(3,3.75)
\psline{->}(2,4)(2,3.75)
\psline{->}(0,4)(0,3.75)
\psarc(2,1){1}{180}{360}
\psarc(1,2){1}{180}{360}
\pscircle[linewidth=0.4pt,fillstyle=solid,fillcolor=white](1,1){0.2}
\pscircle[linewidth=0.4pt,fillstyle=solid,fillcolor=black](2,0){0.2}
\rput(2,0.75){B}
\rput(1,1.75){$\LL_B$}
\endpspicture
\end{align}
where we have defined the new product $\LL_B$, i.e. right contraction 
w.r.t. $B$. The defining tangle of the {\it right contraction\/} is thus:
\index{contraction!right}
\index{tangle!of right contraction} 
\begin{align}
\label{eqn:6-15}
\pspicture[0.5](0,0)(1,2)
\psset{linewidth=\pstlw,xunit=0.5,yunit=0.5,runit=0.5}
\psset{arrowsize=2pt 2,arrowinset=0.2}
\psline(0,4)(0,2)
\psline(2,4)(2,2)
\psline(1,1)(1,0)
\psarc(1,2){1}{180}{360}
\pscircle[linewidth=0.4pt,fillstyle=solid,fillcolor=white](1,1){0.2}
\rput(1,1.75){$\LL_B$}
\endpspicture
\quad:=\quad
\pspicture[0.5](0,0)(2,2)
\psset{linewidth=\pstlw,xunit=0.5,yunit=0.5,runit=0.5}
\psset{arrowsize=2pt 2,arrowinset=0.2}
\psline(1,4)(1,3)
\psline(0,2)(0,0)
\psline(4,4)(4,2)
\psarc(1,2){1}{0}{180}
\psarc(3,2){1}{180}{360}
\pscircle[linewidth=0.4pt,fillstyle=solid,fillcolor=white](1,3){0.2}
\pscircle[linewidth=0.4pt,fillstyle=solid,fillcolor=black](3,1){0.2}
\rput(3,0.25){$B$}
\rput(1,2){$\Delta^{\wedge}$}
\endpspicture
\end{align}
Of course, we can define in an analogous way the left contraction
\index{contraction!left}
\index{tangle!of left contraction}
\begin{align}
\pspicture[0.5](0,0)(1.5,2.5)
\psset{linewidth=\pstlw,xunit=0.5,yunit=0.5,runit=0.5}
\psset{arrowsize=2pt 2,arrowinset=0.2}
\psline(3,5)(3,1)
\psline(2,5)(2,2)
\psline(0,5)(0,2)
\psline{->}(3,4)(3,3.75)
\psline{->}(2,4)(2,3.75)
\psline{->}(0,4)(0,3.75)
\psarc(2,1){1}{180}{360}
\psarc(1,2){1}{180}{360}
\pscircle[linewidth=0.4pt,fillstyle=solid,fillcolor=white](1,1){0.2}
\pscircle[linewidth=0.4pt,fillstyle=solid,fillcolor=black](2,0){0.2}
\rput(1,2){$\wedge$}
\rput(2,0.75){B}
\endpspicture
\quad=\quad
\pspicture[0.5](0,0)(3,2.5)
\psset{linewidth=\pstlw,xunit=0.5,yunit=0.5,runit=0.5}
\psset{arrowsize=2pt 2,arrowinset=0.2}
\psline(5,5)(5,4)
\psline(2,5)(2,3)
\psline(0,5)(0,2)
\psline(6,3)(6,2)
\psline(2,0)(4,0)
\psline{->}(5,4.5)(5,4.25)
\psline{->}(2,4.5)(2,4.25)
\psline{->}(0,4.5)(0,4.25)
\psline(4,0)(4,0)
\psarc(5,3){1}{0}{180}
\psarc(3,3){1}{180}{360}
\psarc(2,2){2}{180}{270}
\psarc(4,2){2}{270}{360}
\pscircle[linewidth=0.4pt,fillstyle=solid,fillcolor=white](5,4){0.2}
\pscircle[linewidth=0.4pt,fillstyle=solid,fillcolor=black](3,0){0.2}
\pscircle[linewidth=0.4pt,fillstyle=solid,fillcolor=black](3,2){0.2}
\rput(3,0.75){B}
\rput(3,2.75){B}
\endpspicture
\quad=:\quad
\pspicture[0.5](0,0)(1.5,2.5)
\psset{linewidth=\pstlw,xunit=0.5,yunit=0.5,runit=0.5}
\psset{arrowsize=2pt 2,arrowinset=0.2}
\psline(0,5)(0,1)
\psline(1,5)(1,2)
\psline(3,5)(3,2)
\psline{->}(0,4)(0,3.75)
\psline{->}(1,4)(1,3.75)
\psline{->}(3,4)(3,3.75)
\psarc(1,1){1}{180}{360}
\psarc(2,2){1}{180}{360}
\pscircle[linewidth=0.4pt,fillstyle=solid,fillcolor=white](2,1){0.2}
\pscircle[linewidth=0.4pt,fillstyle=solid,fillcolor=black](1,0){0.2}
\rput(2,2){$\JJ_B$}
\rput(1,0.75){B}
\endpspicture
\end{align} 
which leads to the defining tangle for the {\it left contraction}:
\begin{align}
\label{eqn:6-17}
\pspicture[0.5](0,0)(1,2)
\psset{linewidth=\pstlw,xunit=0.5,yunit=0.5,runit=0.5}
\psset{arrowsize=2pt 2,arrowinset=0.2}
\psline(0,4)(0,2)
\psline(2,4)(2,2)
\psline(1,1)(1,0)
\psarc(1,2){1}{180}{360}
\pscircle[linewidth=0.4pt,fillstyle=solid,fillcolor=white](1,1){0.2}
\rput(1,1.75){$\JJ_B$}
\endpspicture
\quad:=\quad
\pspicture[0.5](0,0)(2,2)
\psset{linewidth=\pstlw,xunit=0.5,yunit=0.5,runit=0.5}
\psset{arrowsize=2pt 2,arrowinset=0.2}
\psline(3,4)(3,3)
\psline(4,2)(4,0)
\psline(0,4)(0,2)
\psarc(3,2){1}{0}{180}
\psarc(1,2){1}{180}{360}
\pscircle[linewidth=0.4pt,fillstyle=solid,fillcolor=white](3,3){0.2}
\pscircle[linewidth=0.4pt,fillstyle=solid,fillcolor=black](1,1){0.2}
\rput(1,0.25){$B$}
\rput(3,2){$\Delta^{\wedge}$}
\endpspicture
\end{align}

\noindent
{\bf Note} that these relations for the left and right contraction
are valid on the whole graded $A$-module, i.e. for {\it any grade\/}
and inhomogeneous element. It is to the best knowledge of the author 
the first time that such a formula is explicitely given. The same feature 
will be observed with the Rota-Stein cliffordization. We will identify below 
in the physics section contractions with respect to a 1-vector 
as {\it annihilation\/} operators and the wedge product with a 
1-vector as {\it creation\/} operators
\index{contraction!as annihilation}
\index{wedge product!as creation}
\begin{align}
&&v\JJ\,&\cong\, a_v &&&v\wedge\,&\cong\, a_v^\dagger &&\\
&&e_i\JJ\,&\cong\, a_i &&&e_i\wedge\,&\cong\, a_i^\dagger&& . 
\end{align} 
The above given generally valid relation will have an impact on 
calculations in quantum field theory. 

\index{dual isomorphism!from scalar products}
In fact, $B^{\wedge}$ and $C^{\vee}$ can be seen to be dual isomorphisms 
of a new kind. If one has a pairing 
${\langle} .\mid .{\rangle}_{\text{eval}}$, 
it is clear that vectors are moved to co-vectors by keeping the coefficients 
and altering the basis:
\begin{align}
\langle\omega\mid v{\rangle}_{\text{eval}} &= \omega(v) 
= \langle \epsilon^0\mid \omega^\star\JJ_{\delta}v{\rangle}_{\text{eval}}
= \omega^{\star i}v_i \langle \epsilon^0\mid e_0{\rangle}_{\text{eval}},
\end{align}
where $\omega=\omega_i\epsilon^i$ and $\omega^{\star}=\omega^{\star i}e_i
=\omega_j\delta^{ji}e_i$. Since there is no canonical dual isomorphism,
it is quite artificial to use $\delta$. If a scalar product and a co-scalar 
product are given, it is natural to use these maps $B$ and $C$ to move
vectors to co-vectors and vice versa in a pairing. We denote $B(u)$ the 
co-vector image of $u$ under the map $B$ and $C(\omega)$ the vector image
of $\omega$ under the map $C$.
\begin{align}
\langle\omega\mid u\wedge v{\rangle}_{\text{eval}} 
&=\langle \omega\LL_{\delta} B(u)\mid v{\rangle}_{\text{eval}}
= \langle \omega\LL_{B} u^{\star}\mid v{\rangle}_{\text{B,C}}
\end{align}
and
\begin{align}
\langle\omega\vee\rho\mid v{\rangle}_{\text{eval}} 
&= \langle \omega \mid C(\rho) \JJ_{\delta} v{\rangle}_{\text{eval}}
= \langle \omega \mid \rho^{\star}\JJ_{C}v {\rangle}_{\text{B,C}} .
\end{align}
This setting yields exactly the graded extension $B^{\wedge}$ and 
$C^{\vee}$ as introduced above. In physics, this will ensure that
the adjoint of a creation operator will be an annihilation operator
and not a polynomial of annihilation operators. Such a non-graded
extension of $B$ will lead to {\it polarization effects\/}. Note
that we have not assumed that $B\circ C = \Id_V$ and hence the
dualized dual is not in general identical with the original element. 

Let us explore the calculation rules of the tangle Eqns.
\ref{eqn:6-15} and \ref{eqn:6-17}. If we compute the left 
contraction on two 1-vectors $a$ and $b$, we get:
\begin{align}
\JJ_B(a\otimes b) &= (B\otimes \Id)\left(
(\Id\otimes \Delta)(a\otimes b)\right) \nn
&=(B\otimes \Id)(a\otimes b \otimes Id + a\otimes\Id \otimes b) \nn
&=B(a,b)\Id \,.
\end{align}
That is on 1-vectors the contraction product simply evaluates
to $\eta\circ B$. That is the first law of Chevalley deformation.
Two further relations are required which
describe how the contraction 'distributes' on the primary product.
We compute firstly the {\it left straightening law\/}:
\index{straightening law!for left multiplication}
\begin{align}
\pspicture[0.5](0,0)(1.5,3.5)
\psset{linewidth=\pstlw,xunit=0.5,yunit=0.5,runit=0.5}
\psset{arrowsize=2pt 2,arrowinset=0.2}
\psline(0,7)(0,5)
\psline(2,7)(2,5)
\psline(3,7)(3,3)
\psline(1,4)(1,3)
\psline(2,2)(2,0)
\psarc(2,3){1}{180}{360}
\psarc(1,5){1}{180}{360}
\pscircle[linewidth=0.4pt,fillstyle=solid,fillcolor=white](1,4){0.2}
\pscircle[linewidth=0.4pt,fillstyle=solid,fillcolor=white](2,2){0.2}
\rput(2,2.75){$\JJ_B$}
\endpspicture
\quad&=\quad
\pspicture[0.5](0,0)(2.5,3.5)
\psset{linewidth=\pstlw,xunit=0.5,yunit=0.5,runit=0.5}
\psset{arrowsize=2pt 2,arrowinset=0.2}
\psline(0,7)(0,5)
\psline(2,7)(2,5)
\psline(4,7)(4,5)
\psline(1,4)(1,3)
\psline(5,4)(5,0)
\psline(3,4)(3,3)
\psarc(4,4){1}{0}{180}
\psarc(2,3){1}{180}{360}
\psarc(1,5){1}{180}{360}
\pscircle[linewidth=0.4pt,fillstyle=solid,fillcolor=white](1,4){0.2}
\pscircle[linewidth=0.4pt,fillstyle=solid,fillcolor=white](4,5){0.2}
\pscircle[linewidth=0.4pt,fillstyle=solid,fillcolor=black](2,2){0.2}
\rput(2,2.75){$B$}
\endpspicture
\quad=\quad
\pspicture[0.5](0,0)(3.5,3.5)
\psset{linewidth=\pstlw,xunit=0.5,yunit=0.5,runit=0.5}
\psset{arrowsize=2pt 2,arrowinset=0.2}
\psline(0,7)(0,3)
\psline(2,7)(2,4)
\psline(6,7)(6,6)
\psline(7,5)(7,0)
\psline(6,4)(6,3)
\psarc(6,5){1}{0}{180}
\psarc(5,4){1}{0}{180}
\psarc(3,4){1}{180}{360}
\psarc(2,3){2}{180}{270}
\psarc(4,3){2}{270}{360}
\psline(2,1)(4,1)
\pscircle[linewidth=0.4pt,fillstyle=solid,fillcolor=white](6,6){0.2}
\pscircle[linewidth=0.4pt,fillstyle=solid,fillcolor=white](5,5){0.2}
\pscircle[linewidth=0.4pt,fillstyle=solid,fillcolor=black](3,3){0.2}
\pscircle[linewidth=0.4pt,fillstyle=solid,fillcolor=black](3,1){0.2}
\rput(3,1.75){$B$}
\rput(3,3.75){$B$}
\endpspicture
\nonumber \\[4ex]
&=\quad
\pspicture[0.5](0,0)(3.5,3.5)
\psset{linewidth=\pstlw,xunit=0.5,yunit=0.5,runit=0.5}
\psset{arrowsize=2pt 2,arrowinset=0.2}
\psline(0,7)(0,3)
\psline(2,7)(2,4)
\psline(5,7)(5,6)
\psline(7,4)(7,0)
\psline(4,5)(4,4)
\psline(5,4)(5,3)
\psarc(5,5){1}{0}{180}
\psarc(6,4){1}{0}{180}
\psarc(3,4){1}{180}{360}
\psarc(2,3){2}{180}{270}
\psarc(3,3){2}{270}{360}
\psline(2,1)(3,1)
\pscircle[linewidth=0.4pt,fillstyle=solid,fillcolor=white](5,6){0.2}
\pscircle[linewidth=0.4pt,fillstyle=solid,fillcolor=white](6,5){0.2}
\pscircle[linewidth=0.4pt,fillstyle=solid,fillcolor=black](3,3){0.2}
\pscircle[linewidth=0.4pt,fillstyle=solid,fillcolor=black](2.5,1){0.2}
\rput(2.5,1.75){$B$}
\rput(3,3.75){$B$}
\endpspicture
\quad=\quad
\pspicture[0.5](0,0)(1.5,3.5)
\psset{linewidth=\pstlw,xunit=0.5,yunit=0.5,runit=0.5}
\psset{arrowsize=2pt 2,arrowinset=0.2}
\psline(3,7)(3,5)
\psline(1,7)(1,5)
\psline(0,7)(0,3)
\psline(2,4)(2,3)
\psline(1,2)(1,0)
\psarc(1,3){1}{180}{360}
\psarc(2,5){1}{180}{360}
\pscircle[linewidth=0.4pt,fillstyle=solid,fillcolor=white](2,4){0.2}
\pscircle[linewidth=0.4pt,fillstyle=solid,fillcolor=white](1,2){0.2}
\rput(1,2.75){$\JJ_B$}
\rput(2,4.75){$\JJ_B$}
\endpspicture
\end{align}
That is, from the co-associativity of $\Delta$ and product co-product
duality w.r.t. $B$ we have {\it derived\/} the rule
\begin{align}
(u\wedge v) \JJ_B w &= u \JJ_B (v \JJ_B w),
\end{align}
where $u,v,w$ are arbitrary elements from $\bigwedge V$. This is the 
third law of Chevalley deformation.

To compute the {\it right straightening law\/} we have to compute
\index{straightening law!for right multiplication}
\begin{align}
\label{eqn:6-33}
\pspicture[0.5](0,0)(1.5,3.5)
\psset{linewidth=\pstlw,xunit=0.5,yunit=0.5,runit=0.5}
\psset{arrowsize=2pt 2,arrowinset=0.2}
\psline(3,7)(3,5)
\psline(1,7)(1,5)
\psline(0,7)(0,3)
\psline(2,4)(2,3)
\psline(1,2)(1,0)
\psarc(1,3){1}{180}{360}
\psarc(2,5){1}{180}{360}
\pscircle[linewidth=0.4pt,fillstyle=solid,fillcolor=white](2,4){0.2}
\pscircle[linewidth=0.4pt,fillstyle=solid,fillcolor=white](1,2){0.2}
\rput(1,2.75){$\JJ_B$}
\endpspicture
\quad&=\quad
\pspicture[0.5](0,0)(2,3.5)
\psset{linewidth=\pstlw,xunit=0.5,yunit=0.5,runit=0.5}
\psset{arrowsize=2pt 2,arrowinset=0.2}
\psline(0,7)(0,3)
\psline(2,7)(2,6)
\psline(4,7)(4,6)
\psline(3,5)(3,4)
\psline(4,3)(4,0)
\psarc(3,3){1}{0}{180}
\psarc(1,3){1}{180}{360}
\psarc(3,6){1}{180}{360}
\pscircle[linewidth=0.4pt,fillstyle=solid,fillcolor=white](3,5){0.2}
\pscircle[linewidth=0.4pt,fillstyle=solid,fillcolor=white](3,4){0.2}
\pscircle[linewidth=0.4pt,fillstyle=solid,fillcolor=black](1,2){0.2}
\rput(1,2.75){$B$}
\endpspicture
\quad=\quad
\pspicture[0.5](0,0)(3.5,3.5)
\psset{linewidth=\pstlw,xunit=0.5,yunit=0.5,runit=0.5}
\psset{arrowsize=2pt 2,arrowinset=0.2}
\psline(0,7)(0,2)
\psline(1,5)(1,3)
\psline(2,7)(2,6)
\psline(6,7)(6,6)
\psline(6,2)(6,0)
\psline(7,5)(7,3)
\psarc(2,5){1}{0}{180}
\psarc(6,5){1}{0}{180}
\psarc(1,2){1}{180}{360}
\psarc(2,3){1}{180}{360}
\psarc(6,3){1}{180}{360}
\psbezier(5,5)(5,4)(3,4)(3,3)
\psbezier[border=4pt,bordercolor=white](3,5)(3,4)(5,4)(5,3)
\pscircle[linewidth=0.4pt,fillstyle=solid,fillcolor=white](2,6){0.2}
\pscircle[linewidth=0.4pt,fillstyle=solid,fillcolor=white](2,2){0.2}
\pscircle[linewidth=0.4pt,fillstyle=solid,fillcolor=white](6,6){0.2}
\pscircle[linewidth=0.4pt,fillstyle=solid,fillcolor=white](6,2){0.2}
\pscircle[linewidth=0.4pt,fillstyle=solid,fillcolor=black](1,1){0.2}
\rput(1,1.75){$B$}
\endpspicture
\nonumber \\[4ex]
&=\quad
\pspicture[0.5](0,0)(5,3)
\psset{linewidth=\pstlw,xunit=0.5,yunit=0.5,runit=0.5}
\psset{arrowsize=2pt 2,arrowinset=0.2}
\psline(0,4)(0,3)
\psline(1,6)(1,5)
\psline(5,6)(5,5)
\psline(9,6)(9,5)
\psline(10,4)(10,2)
\psline(9,1)(9,0)
\psline(2,1)(4,1)
\psarc(1,4){1}{0}{180}
\psarc(5,4){1}{0}{180}
\psarc(9,4){1}{0}{180}
\psarc(9,2){1}{180}{360}
\psarc(2,3){2}{180}{270}
\psarc(3,4){1}{180}{360}
\psbezier(4,1)(5,1)(6,1)(6,2)
\psbezier(6,2)(6,3)(8,3)(8,4)
\psbezier[border=4pt,bordercolor=white](6,4)(6,3)(8,3)(8,2)
\pscircle[linewidth=0.4pt,fillstyle=solid,fillcolor=white](1,5){0.2}
\pscircle[linewidth=0.4pt,fillstyle=solid,fillcolor=white](5,5){0.2}
\pscircle[linewidth=0.4pt,fillstyle=solid,fillcolor=white](9,5){0.2}
\pscircle[linewidth=0.4pt,fillstyle=solid,fillcolor=white](9,1){0.2}
\pscircle[linewidth=0.4pt,fillstyle=solid,fillcolor=black](3,1){0.2}
\pscircle[linewidth=0.4pt,fillstyle=solid,fillcolor=black](3,3){0.2}
\rput(3,1.75){$B$}
\rput(3,3.75){$B$}
\endpspicture
\quad=\quad
\pspicture[0.5](0,0)(3,3)
\psset{linewidth=\pstlw,xunit=0.5,yunit=0.5,runit=0.5}
\psset{arrowsize=2pt 2,arrowinset=0.2}
\psline(1,7)(1,6)
\psline(4,7)(4,5)
\psline(6,7)(6,3)
\psline(4,1)(4,0)
\psline(0,5)(0,2)
\psbezier(0,2)(0,1)(3.5,4)(4,3)
\psline[border=4pt,bordercolor=white](3,4)(3,2)
\psarc(1,5){1}{0}{180}
\psarc(3,5){1}{180}{360}
\psarc(5,3){1}{180}{360}
\psarc(4,2){1}{180}{360}
\pscircle[linewidth=0.4pt,fillstyle=solid,fillcolor=white](1,6){0.2}
\pscircle[linewidth=0.4pt,fillstyle=solid,fillcolor=white](3,4){0.2}
\pscircle[linewidth=0.4pt,fillstyle=solid,fillcolor=white](5,2){0.2}
\pscircle[linewidth=0.4pt,fillstyle=solid,fillcolor=white](4,1){0.2}
\rput(3,4.75){$\JJ_B$}
\rput(5,2.75){$\JJ_B$}
\rput(0.5,1.5){eval}
\rput(3.85,3.6){eval}
\endpspicture
\end{align}
where we had to assume the identity
\begin{align}
\label{eqn:6-34}
\pspicture[0.5](0,0)(1,2)
\psset{linewidth=\pstlw,xunit=0.5,yunit=0.5,runit=0.5}
\psset{arrowsize=2pt 2,arrowinset=0.2}
\psbezier(2,4)(2,3)(0,3)(0,2)
\psbezier[border=4pt,bordercolor=white](0,4)(0,3)(2,3)(2,2)
\psline(2,2)(2,0)
\psline(0,2)(0,0)
\psline{->}(2,1.5)(2,1)
\psline{->}(0,2)(0,1.75)
\psline{->}(0,0)(0,0.25)
\pscircle[linewidth=0.4pt,fillstyle=solid,fillcolor=black](0,1){0.2}
\rput(0.75,1){$B$}
\endpspicture
\quad=\quad
\pspicture[0.5](0,0)(1,2)
\psset{linewidth=\pstlw,xunit=0.5,yunit=0.5,runit=0.5}
\psset{arrowsize=2pt 2,arrowinset=0.2}
\psbezier(2,2)(2,1)(0,1)(0,0)
\psbezier[border=4pt,bordercolor=white](0,2)(0,1)(2,1)(2,0)
\psline(2,4)(2,2)
\psline(0,4)(0,2)
\psline{->}(0,3.5)(0,3)
\psline{<-}(2,3.75)(2,4)
\psline{->}(2,2)(2,2.25)
\pscircle[linewidth=0.4pt,fillstyle=solid,fillcolor=black](2,3){0.2}
\rput(1.25,3){$B$}
\endpspicture
\end{align}
which is a requirement on the crossing. If this identity does not hold, 
we cannot move the $\JJ_B$ product 'under' the crossing but get not
a contraction but a different product. The above used identity reads 
in algebraic terms
\begin{align}
B_c^s\tau^{ab}_{sd} &= \tau^{ab}_{cs}B^s_d .
\end{align}
The tangle equation \ref{eqn:6-33} reads algebraically for arbitrary
multivectors $u,v,w$
\begin{align}
w\JJ_B (u\wedge v) &= (-1)^{(\vert w_{(1)}\vert\vert w_{(2)}\JJ_B u\vert)}
\, (w_{(2)}\JJ_B u)\wedge (w_{(1)}\JJ_B v)
\end{align}
If we assume the leftmost input $w$ to be a 1-vector $a$, we arrive 
with Gra{\ss}mann products and the Gra{\ss}mann graded switch as crossing, 
where Eqn. \ref{eqn:6-34} holds, at the following formula
\begin{align}
a\JJ_B(u\wedge v) &= (a\JJ_B u)\wedge v + \hat{u}\wedge(a \JJ_B v),
\end{align}
which is the second law of Chevalley deformation. However, our Hopf
algebraic result in Eqn. \ref{eqn:6-33} is valid for the input of any 
element of any grade, even inhomogeneous. The crossing has been 
replaced by $\hat{u}$, the grade involution. This is possible only if
the first factor is a 1-vector and shows that Chevalley deformation
is restricted by the fact that it does not properly deal with the 
crossing:
\index{grade involution!from crossing}
\begin{align}
\hat{\tau}(a \otimes u) &= (-1)^{\partial a\partial u}(u\otimes a) 
= ((-1)^{\partial u}u) \otimes a \nn
&= \hat{u} \otimes a.
\end{align}
Summarizing the formulas which we have just derived, we end up with the
rules of Chevalley deformation of a Gra{\ss}mann algebra, i.e. a
Clifford map which is given on 1-vectors $x\in V$ as
\index{map!Clifford}
\index{Clifford!map}
\index{Chevalley deformation}
\begin{align}
x &\rightarrow \gamma_x := x\JJ_B \,+\, x\wedge .
\end{align}
The operator $\gamma ~:~ V \otimes \bigwedge V\rightarrow \bigwedge V$ 
can be lifted to an action $\gamma ~:~\bigwedge V \otimes \bigwedge V
\rightarrow \bigwedge V$ by {\it recursive application\/} and linearity.
However, the Hopf algebraic counterparts are valid on the whole space and 
do not suffer any restriction on their input.
We can summarize the formulas which we have derived from Hopf gebraic 
considerations and compare them to the literature, e.g 
\cite{chevalley:1997a,bourbaki:1989a,budinich:trautmann:1988a,%
crumeyrolle:1990a,lounesto:2001a}.
Let $a,b \in V$, $u,v,w \in \bigwedge V$ it holds:
\index{contraction!rules}
\begin{align}
\label{eqn-DER-CON}
a\JJ_B b &= B(a,b) \nn
a\JJ_B(u\wedge v) &= (a\JJ_B u)\wedge v + \hat{u}\wedge(a \JJ_B v) \nn
(u\wedge v) \JJ_B w &= u \JJ_B (v \JJ_B w)
\end{align}

As a matter of fact, we could now develop the {\it co-contraction\/},
using the co-scalar product $C$ and derive analogous relations
as Eqns. \ref{eqn-DER-CON}.
\index{co-contraction}
\begin{align}
\pspicture[0.5](0,0)(1.5,3)
\psset{linewidth=\pstlw,xunit=0.5,yunit=0.5,runit=0.5}
\psset{arrowsize=2pt 2,arrowinset=0.2}
\psline(0,0)(0,4)
\psline(2,0)(2,4)
\psline(3,0)(3,5)
\psarc(1,4){1}{0}{180}
\psarc(2,5){1}{0}{180}
\pscircle[linewidth=0.4pt,fillstyle=solid,fillcolor=white](1,5){0.2}
\pscircle[linewidth=0.4pt,fillstyle=solid,fillcolor=black](2,6){0.2}
\rput(2,5.25){$C$}
\endpspicture
\quad=\quad
\pspicture[0.5](0,0)(3,3)
\psset{linewidth=\pstlw,xunit=0.5,yunit=0.5,runit=0.5}
\psset{arrowsize=2pt 2,arrowinset=0.2}
\psline(0,0)(0,4)
\psline(2,0)(2,3)
\psline(5,0)(5,2)
\psline(2,6)(4,6)
\psline(6,4)(6,3)
\psarc(3,3){1}{0}{180}
\psarc(5,3){1}{180}{360}
\psarc(2,4){2}{90}{180}
\psarc(4,4){2}{0}{90}
\pscircle[linewidth=0.4pt,fillstyle=solid,fillcolor=white](5,2){0.2}
\pscircle[linewidth=0.4pt,fillstyle=solid,fillcolor=black](3,4){0.2}
\pscircle[linewidth=0.4pt,fillstyle=solid,fillcolor=black](3,6){0.2}
\rput(3,5.25){$C$}
\rput(3,3.25){$C$}
\endpspicture
\quad&&\Rightarrow\quad
\pspicture[0.5](0,0)(1,3)
\psset{linewidth=\pstlw,xunit=0.5,yunit=0.5,runit=0.5}
\psset{arrowsize=2pt 2,arrowinset=0.2}
\psline(0,0)(0,3)
\psline(2,0)(2,3)
\psline(1,4)(1,6)
\psarc(1,3){1}{0}{180}
\pscircle[linewidth=0.4pt,fillstyle=solid,fillcolor=white](1,4){0.2}
\rput(1.1,3){$\Delta_{\JJ_C}$}
\endpspicture
\quad:=\quad
\pspicture[0.5](0,0)(2,3)
\psset{linewidth=\pstlw,xunit=0.5,yunit=0.5,runit=0.5}
\psset{arrowsize=2pt 2,arrowinset=0.2}
\psline(0,0)(0,3)
\psline(3,0)(3,2)
\psline(4,3)(4,6)
\psarc(1,3){1}{0}{180}
\psarc(3,3){1}{180}{360}
\pscircle[linewidth=0.4pt,fillstyle=solid,fillcolor=black](1,4){0.2}
\pscircle[linewidth=0.4pt,fillstyle=solid,fillcolor=white](3,2){0.2}
\rput(1,4.75){$C$}
\endpspicture
\end{align}
and
\begin{align}
\pspicture[0.5](0,0)(1.5,3)
\psset{linewidth=\pstlw,xunit=0.5,yunit=0.5,runit=0.5}
\psset{arrowsize=2pt 2,arrowinset=0.2}
\psline(3,0)(3,4)
\psline(1,0)(1,4)
\psline(0,0)(0,5)
\psarc(2,4){1}{0}{180}
\psarc(1,5){1}{0}{180}
\pscircle[linewidth=0.4pt,fillstyle=solid,fillcolor=white](2,5){0.2}
\pscircle[linewidth=0.4pt,fillstyle=solid,fillcolor=black](1,6){0.2}
\rput(1,5.25){$C$}
\endpspicture
\quad=\quad
\pspicture[0.5](0,0)(3,3)
\psset{linewidth=\pstlw,xunit=0.5,yunit=0.5,runit=0.5}
\psset{arrowsize=2pt 2,arrowinset=0.2}
\psline(6,0)(6,4)
\psline(4,0)(4,3)
\psline(1,0)(1,2)
\psline(4,6)(2,6)
\psline(0,4)(0,3)
\psarc(3,3){1}{0}{180}
\psarc(1,3){1}{180}{360}
\psarc(4,4){2}{0}{90}
\psarc(2,4){2}{90}{180}
\pscircle[linewidth=0.4pt,fillstyle=solid,fillcolor=white](1,2){0.2}
\pscircle[linewidth=0.4pt,fillstyle=solid,fillcolor=black](3,4){0.2}
\pscircle[linewidth=0.4pt,fillstyle=solid,fillcolor=black](3,6){0.2}
\rput(3,5.25){$C$}
\rput(3,3.25){$C$}
\endpspicture
\quad&&\Rightarrow\quad
\pspicture[0.5](0,0)(1,3)
\psset{linewidth=\pstlw,xunit=0.5,yunit=0.5,runit=0.5}
\psset{arrowsize=2pt 2,arrowinset=0.2}
\psline(2,0)(2,3)
\psline(0,0)(0,3)
\psline(1,4)(1,6)
\psarc(1,3){1}{0}{180}
\pscircle[linewidth=0.4pt,fillstyle=solid,fillcolor=white](1,4){0.2}
\rput(1.1,3){$\Delta_{\LL_C}$}
\endpspicture
\quad:=\quad
\pspicture[0.5](0,0)(2,3)
\psset{linewidth=\pstlw,xunit=0.5,yunit=0.5,runit=0.5}
\psset{arrowsize=2pt 2,arrowinset=0.2}
\psline(4,0)(4,3)
\psline(1,0)(1,2)
\psline(0,3)(0,6)
\psarc(3,3){1}{0}{180}
\psarc(1,3){1}{180}{360}
\pscircle[linewidth=0.4pt,fillstyle=solid,fillcolor=black](3,4){0.2}
\pscircle[linewidth=0.4pt,fillstyle=solid,fillcolor=white](1,2){0.2}
\rput(3,4.75){$C$}
\endpspicture
\end{align} 
In Sweedler notation these formulas are displayed as -- remember that 
$C^{\wedge}=\Delta_C(\Id)$:
\begin{align}
\Delta_{\JJ_C}(x) &= C^{\wedge}_{(1)} \otimes (C^{\wedge}_{(2)} \wedge x) \nn
\Delta_{\LL_C}(x) &= (x\wedge C^{\wedge}_{(1)})\otimes C^{\wedge}_{(2)}
\end{align}

It would be now possible to derive a co-Chevalley deformation based on a  
co-Clifford map.

\section{Cliffordization of Rota and Stein}

\index{cliffordization}
\index{product!quantized}
Cliffordization is a quite remarkable process. A product or 
co-product is {\it deformed\/} by cliffordization to yield a new 
{\it quantized\/} product, see \cite{oziewicz:1997b}. Deformation and 
quantization are therefore intimately related. In fact it turns out 
that imposing non-trivial commutation relations is equivalent to the choice of 
a bilinear form which gives rise to the deformation. It is remarkable to note 
that the process of cliffordization is quite ubiquitous in mathematics and not
restricted to quantum physics, \cite{rota:stein:1994a,rota:stein:1994b}.
Since we will use cliffordization mainly for Gra{\ss}mann exterior algebras
it should be emphasized that this method works for symmetric algebras
and even more general algebras also.

\subsection{Cliffordization of products}

Let $m_A~:~A\otimes A \rightarrow A$ be the product of a Hopf 
algebra $H(A,m_A,\eta,\Delta,\epsilon;S)$, or a convolution
$\Conv(A,A)$. As a prototype, the reader may think of a Gra{\ss}mann wedge
product. Now let a scalar product $B^\wedge$ be given on $A$ by 
exponentiation of $B$, which is represented by a cup tangle.

\index{tangle!cliffordization}
\index{tangle!sausage}
\index{cliffordization}
\begin{dfn}[Cliffordization]
A {\em Clifford product\/} (or circle product) $\&c$ on $A$ is defined via 
the tangle
\begin{align}
\pspicture[0.5](0,0)(1,3.5)
\psset{linewidth=\pstlw,xunit=0.5,yunit=0.5,runit=0.5}
\psset{arrowsize=2pt 2,arrowinset=0.2}
\psline(0,7)(0,4)
\psline(2,7)(2,4)
\psarc(1,4){1}{180}{360}
\psline(1,3)(1,0)
\pscircle[linewidth=0.4pt,fillstyle=solid,fillcolor=white](1,3){0.2}
\rput(1,3.75){$\&c$}
\endpspicture
\quad&:=\quad
\pspicture[0.5](0,0)(2.5,3.5)
\psset{linewidth=\pstlw,xunit=0.5,yunit=0.5,runit=0.5}
\psset{arrowsize=2pt 2,arrowinset=0.2}
\psline(1,7)(1,6)
\psline(5,7)(5,6)
\psarc(1,5){1}{0}{180}
\psarc(5,5){1}{0}{180}
\psarc(3,5){1}{180}{360}
\psarc(3,5){3}{180}{360}
\psline(3,2)(3,0)
\pscircle[linewidth=0.4pt,fillstyle=solid,fillcolor=black](3,4){0.2}
\pscircle[linewidth=0.4pt,fillstyle=solid,fillcolor=white](3,2){0.2}
\pscircle[linewidth=0.4pt,fillstyle=solid,fillcolor=white](1,6){0.2}
\pscircle[linewidth=0.4pt,fillstyle=solid,fillcolor=white](5,6){0.2}
\rput(1,5.25){$\Delta$}
\rput(5,5.25){$\Delta$}
\rput(3,4.75){$B^\wedge$}
\rput(3,2.75){$m_A$}
\endpspicture
\end{align}
where $B^\wedge$ is the bilinear form obtained from $B : V\otimes V
\rightarrow\openk$ by exponentiation.
\end{dfn}
This tangle can be easily remembered as 'sausage tangle', the
term was  coined by Oziewicz \cite{oziewicz:2001b}. Of course,
since $\&c$ is a $2\rightarrow 1$ map, it is a product. The new
product is totally defined by the structure tensors $\Delta,m_A$
and the scalar product $B$ of the primarily given Hopf algebra
or convolution. Cliffordization is a quite general process, 
see \cite{rota:stein:1994a,rota:stein:1994b} and it is by no means
restricted to Gra{\ss}mann Hopf algebras. It needs, in principle, only a 
convolution and a scalar product or even less restrictive two product
maps.

\begin{thrm}
\label{thrm:6-2}
The unit $\eta$ of $m_A$, if it exists, remains to be the unit
of the cliffordized product if (i) the unit $\eta$ is a 
co-algebra homomorphism of the primary co-product, and, (ii) the unit 
$\eta$ and counit $\epsilon$ are related via 
$B(A\otimes\eta)=\epsilon=B(\eta\otimes A)$.
\end{thrm}
\noindent
{\bf Proof:} The proof is given for the unit multiplied from the right,
the left multiplication case can be shown analogously.
\begin{align}
\pspicture[0.5](0,0)(1,3.5)
\psset{linewidth=\pstlw,xunit=0.5,yunit=0.5,runit=0.5}
\psset{arrowsize=2pt 2,arrowinset=0.2}
\psline(0,7)(0,4)
\psline(2,6)(2,4)
\psarc(1,4){1}{180}{360}
\psline(1,3)(1,0)
\pscircle[linewidth=0.4pt,fillstyle=solid,fillcolor=white](1,3){0.2}
\pscircle[linewidth=0.4pt,fillstyle=solid,fillcolor=white](2,6){0.2}
\rput(1,3.75){$\&c$}
\endpspicture
\quad=\quad
\pspicture[0.5](0,0)(2.5,3.5)
\psset{linewidth=\pstlw,xunit=0.5,yunit=0.5,runit=0.5}
\psset{arrowsize=2pt 2,arrowinset=0.2}
\psline(1,7)(1,6)
\psline(5,7)(5,6)
\psarc(1,5){1}{0}{180}
\psarc(5,5){1}{0}{180}
\psarc(3,5){1}{180}{360}
\psarc(3,5){3}{180}{360}
\psline(3,2)(3,0)
\pscircle[linewidth=0.4pt,fillstyle=solid,fillcolor=black](3,4){0.2}
\pscircle[linewidth=0.4pt,fillstyle=solid,fillcolor=white](3,2){0.2}
\pscircle[linewidth=0.4pt,fillstyle=solid,fillcolor=white](1,6){0.2}
\pscircle[linewidth=0.4pt,fillstyle=solid,fillcolor=white](5,6){0.2}
\pscircle[linewidth=0.4pt,fillstyle=solid,fillcolor=white](5,7){0.2}
\rput(3,3.25){$B^{\wedge}$}
\endpspicture
\qquad=\quad
\pspicture[0.5](0,0)(2.5,3.5)
\psset{linewidth=\pstlw,xunit=0.5,yunit=0.5,runit=0.5}
\psset{arrowsize=2pt 2,arrowinset=0.2}
\psline(1,7)(1,6)
\psarc(1,5){1}{0}{180}
\psarc(3,5){1}{180}{360}
\psarc(3,5){3}{180}{360}
\psline(3,2)(3,0)
\pscircle[linewidth=0.4pt,fillstyle=solid,fillcolor=black](3,4){0.2}
\pscircle[linewidth=0.4pt,fillstyle=solid,fillcolor=white](3,2){0.2}
\pscircle[linewidth=0.4pt,fillstyle=solid,fillcolor=white](1,6){0.2}
\pscircle[linewidth=0.4pt,fillstyle=solid,fillcolor=white](4,5){0.2}
\pscircle[linewidth=0.4pt,fillstyle=solid,fillcolor=white](6,5){0.2}
\rput(3,3.25){$B^{\wedge}$}
\endpspicture
\qquad=\quad
\pspicture[0.5](0,0)(1,3.5)
\psset{linewidth=\pstlw,xunit=0.5,yunit=0.5,runit=0.5}
\psset{arrowsize=2pt 2,arrowinset=0.2}
\psline(1,7)(1,5)
\psarc(1,4){1}{0}{180}
\psline(0,4)(0,0)
\pscircle[linewidth=0.4pt,fillstyle=solid,fillcolor=white](2,4){0.2}
\pscircle[linewidth=0.4pt,fillstyle=solid,fillcolor=white](1,5){0.2}
\endpspicture
\quad=
\pspicture[0.5](0,0)(1,3.5)
\psset{linewidth=\pstlw,xunit=0.5,yunit=0.5,runit=0.5}
\psset{arrowsize=2pt 2,arrowinset=0.2}
\psline(1,7)(1,0)
\endpspicture
\end{align}
In fact this means that we are dealing with an augmented algebra
and an augmented co-algebra as primary structure.

\begin{thrm}
\label{thrm:6-3}
The cliffordized product $\&c$ of a bi-associative Hopf gebra
or a bi-associative bigebra is associative under the condition
that the crossing fulfils the following symmetry requirement
\begin{align}
\tau^{ab}_{cd} = \tau^{bd}_{ac}
\quad&&
\pspicture[0.5](0,0)(1,1)
\psset{linewidth=\pstlw,xunit=0.5,yunit=0.5,runit=0.5}
\psset{arrowsize=2pt 2,arrowinset=0.2}
\psbezier(2,2)(2,1)(0,1)(0,0)
\psbezier[border=4pt,bordercolor=white](0,2)(0,1)(2,1)(2,0)
\endpspicture
\quad=\quad
\pspicture[0.5](0,0)(1,1)
\psset{linewidth=\pstlw,xunit=0.5,yunit=0.5,runit=0.5}
\psset{arrowsize=2pt 2,arrowinset=0.2}
\psbezier(0,2)(0,1)(2,1)(2,0)
\psbezier[border=4pt,bordercolor=white](2,2)(2,1)(0,1)(0,0)
\endpspicture
\end{align}
\end{thrm}
\noindent
{\bf Proof:} Looking at the tangles and using the fact that $m_A$ is an 
co-algebra homomorphism and $\Delta$ is an algebra homomorphism and 
product co-product duality between $m_A$ and $\Delta$ yields the result.

\index{Hopf algebra!cliffordized}
Rota and Stein showed, that the new structure 
$H(A,\&c,\eta,\Delta,\epsilon;S^{c})$ is a Hopf algebra. A point 
which was criticized at this level is that the product is deformed
by the sausage tangle, but the co-product remains to be the same. This
is quite unnatural. Moreover, since these products are no longer related 
by product co-product duality, a further cliffordization may lead to
non-associative products. This motivates the following definition:

\begin{dfn}[local and non-local products]
\label{dfn:6-4}
\index{product!local}
\index{co-product!local}
A product which possesses an augmentation $\epsilon$ such
that $\epsilon$ is an algebra homomorphism is called {\em local}. 
If the augmentation is not an algebra homomorphism the product is 
called {\em non-local}. The same notion is introduced by duality for 
co-products.
\end{dfn}

\begin{thrm}
\label{thrm:6-5}
\index{product!Clifford!non-local}
Products which arise from the process of cliffordization are
in general non-local.
\end{thrm}
\noindent
{\bf Proof:}
\begin{align}
\pspicture[0.5](0,0)(2.5,3.5)
\psset{linewidth=\pstlw,xunit=0.5,yunit=0.5,runit=0.5}
\psset{arrowsize=2pt 2,arrowinset=0.2}
\psline(1,7)(1,6)
\psline(5,7)(5,6)
\psarc(1,5){1}{0}{180}
\psarc(5,5){1}{0}{180}
\psarc(3,5){1}{180}{360}
\psarc(3,5){3}{180}{360}
\psline(3,2)(3,0)
\pscircle[linewidth=0.4pt,fillstyle=solid,fillcolor=black](3,4){0.2}
\pscircle[linewidth=0.4pt,fillstyle=solid,fillcolor=white](3,2){0.2}
\pscircle[linewidth=0.4pt,fillstyle=solid,fillcolor=white](1,6){0.2}
\pscircle[linewidth=0.4pt,fillstyle=solid,fillcolor=white](5,6){0.2}
\pscircle[linewidth=0.4pt,fillstyle=solid,fillcolor=white](3,0){0.2}
\endpspicture
\qquad=\quad
\pspicture[0.5](0,0)(2.5,3.5)
\psset{linewidth=\pstlw,xunit=0.5,yunit=0.5,runit=0.5}
\psset{arrowsize=2pt 2,arrowinset=0.2}
\psline(1,7)(1,6)
\psline(5,7)(5,6)
\psarc(1,5){1}{0}{180}
\psarc(5,5){1}{0}{180}
\psarc(3,5){1}{180}{360}
\pscircle[linewidth=0.4pt,fillstyle=solid,fillcolor=black](3,4){0.2}
\pscircle[linewidth=0.4pt,fillstyle=solid,fillcolor=white](0,5){0.2}
\pscircle[linewidth=0.4pt,fillstyle=solid,fillcolor=white](6,5){0.2}
\pscircle[linewidth=0.4pt,fillstyle=solid,fillcolor=white](1,6){0.2}
\pscircle[linewidth=0.4pt,fillstyle=solid,fillcolor=white](5,6){0.2}
\endpspicture
\qquad=\quad
\pspicture[0.5](0,0)(1,2)
\psset{linewidth=\pstlw,xunit=0.5,yunit=0.5,runit=0.5}
\psline(0,4)(0,2)
\psline(2,4)(2,2)
\psarc(1,2){1}{180}{360}
\pscircle[linewidth=0.4pt,fillstyle=solid,fillcolor=black](1,1){0.2}
\rput(1,-0.25){$B^{\wedge}$}
\endpspicture
\end{align}
where the last equality holds if the augmentation $\epsilon$ is
the counit of the primary co-product $\Delta$. Of course, locality is
preserved in the trivial case that $B^\wedge = \epsilon \otimes \epsilon$ 
holds true. But this does not lead to a new product and is not a
proper cliffordization.

\subsection{Cliffordization of co-products}

Having studied the cliffordization of products in some detail, we 
can shortly display the {\it co-cliffordization\/}. Let $\Delta ~:~A
\rightarrow A\otimes A$ be a co-product of a Hopf algebra $H(A,m_A,\eta,
\Delta,\epsilon;S)$ of a convolution $\Conv(A,A)$. Let a co-scalar product
$C^{\vee}$ be given, then we define

\index{co-cliffordization}
\begin{dfn}[Co-cliffordization]
A {\em co-Clifford product} $\Delta_c$ on $A$ is defined via the
following tangle. We employ the extension of the co-scalar product 
$C^\vee$, i.e. a cap tangle:
\begin{align}
\pspicture[0.5](0,0)(1,3.5)
\psset{linewidth=\pstlw,xunit=0.5,yunit=0.5,runit=0.5}
\psset{arrowsize=2pt 2,arrowinset=0.2}
\psline(0,0)(0,3)
\psline(2,0)(2,3)
\psarc(1,3){1}{0}{180}
\psline(1,4)(1,7)
\pscircle[linewidth=0.4pt,fillstyle=solid,fillcolor=white](1,4){0.2}
\rput(1,3.25){$\Delta_c$}
\endpspicture
\quad&:=\quad
\pspicture[0.5](0,0)(2.5,3.5)
\psset{linewidth=\pstlw,xunit=0.5,yunit=0.5,runit=0.5}
\psset{arrowsize=2pt 2,arrowinset=0.2}
\psline(1,0)(1,1)
\psline(5,0)(5,1)
\psarc(1,2){1}{180}{360}
\psarc(5,2){1}{180}{360}
\psarc(3,2){1}{0}{180}
\psarc(3,2){3}{0}{180}
\psline(3,5)(3,7)
\pscircle[linewidth=0.4pt,fillstyle=solid,fillcolor=black](3,3){0.2}
\pscircle[linewidth=0.4pt,fillstyle=solid,fillcolor=white](3,5){0.2}
\pscircle[linewidth=0.4pt,fillstyle=solid,fillcolor=white](1,1){0.2}
\pscircle[linewidth=0.4pt,fillstyle=solid,fillcolor=white](5,1){0.2}
\rput(1,1.75){$m_A$}
\rput(5,1.75){$m_A$}
\rput(3,2.25){$C$}
\rput(3,4.25){$\Delta$}
\endpspicture
\end{align}               
\end{dfn}

It would be worth to explore the structure $H(A,m_A,\eta,\Delta_c,
\epsilon;S^c)$ in the same way as Rota and Stein did for the 
cliffordization. However, we will concentrate on the case where both
products are deformed obtaining this asymmetric deformation as a
special case.

In an analogous sense, the definition \ref{dfn:6-4} and the assertions 
\ref{thrm:6-2}, \ref{thrm:6-3} and \ref{thrm:6-5} can be established 
for Clifford co-products too.

\subsection{Clifford maps for any grade}

A major drawback of Chevalley deformation of Gra{\ss}mann algebras 
is that the Clifford map $\gamma~:~ V\otimes \bigwedge V \rightarrow
\bigwedge V$, $x~:~ x\rightarrow \gamma_x$ is defined on 1-vectors only.
\index{Chevalley deformation!for any grade}
\index{Clifford map!for any grade}
\index{map!Clifford!of any grade}

However, since we found tangles for contractions of any grade, we can
now define a Clifford map for any grade using Hopf algebraic techniques.
Observe that the Rota and Stein 'sausage' tangle of cliffordization
can be rewritten as
\begin{align}
\pspicture[0.5](0,0)(2.5,3.5)
\psset{linewidth=\pstlw,xunit=0.5,yunit=0.5,runit=0.5}
\psset{arrowsize=2pt 2,arrowinset=0.2}
\psline(1,7)(1,6)
\psline(5,7)(5,6)
\psarc(1,5){1}{0}{180}
\psarc(5,5){1}{0}{180}
\psarc(3,5){1}{180}{360}
\psarc(3,5){3}{180}{360}
\psline(3,2)(3,0)
\pscircle[linewidth=0.4pt,fillstyle=solid,fillcolor=black](3,4){0.2}
\pscircle[linewidth=0.4pt,fillstyle=solid,fillcolor=white](3,2){0.2}
\pscircle[linewidth=0.4pt,fillstyle=solid,fillcolor=white](1,6){0.2}
\pscircle[linewidth=0.4pt,fillstyle=solid,fillcolor=white](5,6){0.2}
\endpspicture
\qquad&=\quad
\pspicture[0.5](0,0)(2,3.5)
\psset{linewidth=\pstlw,xunit=0.5,yunit=0.5,runit=0.5}
\psset{arrowsize=2pt 2,arrowinset=0.2}
\psline(1,7)(1,6)
\psline(4,7)(4,5)
\psarc(1,5){1}{0}{180}
\psarc(3,5){1}{180}{360}
\psarc(1.5,2.5){1.5}{180}{360}
\psline(0,5)(0,2.5)
\psline(3,4)(3,2.5)
\psline(1.5,1)(1.5,0)
\pscircle[linewidth=0.4pt,fillstyle=solid,fillcolor=white](1,6){0.2}
\pscircle[linewidth=0.4pt,fillstyle=solid,fillcolor=white](3,4){0.2}
\pscircle[linewidth=0.4pt,fillstyle=solid,fillcolor=white](1.5,1){0.2}
\rput(3,4.75){$\JJ_B$}
\endpspicture
\quad=\quad
\pspicture[0.5](0,0)(2,3.5)
\psset{linewidth=\pstlw,xunit=0.5,yunit=0.5,runit=0.5}
\psset{arrowsize=2pt 2,arrowinset=0.2}
\psline(3,7)(3,6)
\psline(0,7)(0,5)
\psarc(3,5){1}{0}{180}
\psarc(1,5){1}{180}{360}
\psarc(2.5,2.5){1.5}{180}{360}
\psline(4,5)(4,2.5)
\psline(1,4)(1,2.5)
\psline(2.5,1)(2.5,0)
\pscircle[linewidth=0.4pt,fillstyle=solid,fillcolor=white](3,6){0.2}
\pscircle[linewidth=0.4pt,fillstyle=solid,fillcolor=white](1,4){0.2}
\pscircle[linewidth=0.4pt,fillstyle=solid,fillcolor=white](2.5,1){0.2}
\rput(1,4.75){$\LL_B$}
\endpspicture
\end{align}
On 1-vectors this results in
\begin{align}
a \,\&c\, b &= a_{(1)}\wedge (a_{(2)} \JJ_B b) \nn
&= a \wedge (\Id \JJ_B b) + \Id \wedge (a\JJ_B b) \nn
&= a\wedge b + a\JJ_B b = \gamma_a \, b \\
a \,\&c\, b &= (a \LL_B b_{(1)}) \wedge  b_{(2)} \nn
&= (a \LL_B \Id) \wedge b + (a\LL b)\wedge \Id \nn
&=a\wedge b + b\LL_B a = \gamma_a \, b . 
\end{align}
Now, the above formula holds also in higher grades, and even when elements
are inhomogeneous. We compute an example where $a,b,x,y \in V$ and we 
Clifford multiply two step 2 extensors
\begin{align}
(a\wedge b) \,\&c\, (x\wedge y) &= 
(a\wedge b )\wedge (x\wedge y) + a\wedge (b \JJ_B ((x\wedge y)) \nn
&\phantom{=} -b\wedge (a \JJ_B ((x\wedge y)) + \Id\wedge((a\wedge b) \JJ_B (x\wedge y)) \nn
&=\phantom{+}
 a \wedge (b \wedge (x\wedge y) + b \JJ_B (x\wedge y)) \nn
&\phantom{=} +a \JJ_B  (b \wedge (x\wedge y) + b \JJ_B (x\wedge y)) 
-(a \JJ_B b) (x\wedge y) \nn
&=
\gamma_a (\gamma_b (x\wedge y)) -(a\JJ_B b) (x\wedge y) \nn
&=
(\gamma_a \wedge \gamma_b) (x \wedge y) \, . 
\end{align}
The case of co-cliffordization is handled along the same lines. We get
\index{co-Chevalley deformation}
\index{co-Clifford map}
\index{map!co-Clifford}
\begin{align}
\pspicture[0.5](0,0)(2.5,3.5)
\psset{linewidth=\pstlw,xunit=0.5,yunit=0.5,runit=0.5}
\psset{arrowsize=2pt 2,arrowinset=0.2}
\psline(1,0)(1,1)
\psline(5,0)(5,1)
\psarc(1,2){1}{180}{360}
\psarc(5,2){1}{180}{360}
\psarc(3,2){1}{0}{180}
\psarc(3,2){3}{0}{180}
\psline(3,5)(3,7)
\pscircle[linewidth=0.4pt,fillstyle=solid,fillcolor=black](3,3){0.2}
\pscircle[linewidth=0.4pt,fillstyle=solid,fillcolor=white](3,5){0.2}
\pscircle[linewidth=0.4pt,fillstyle=solid,fillcolor=white](1,1){0.2}
\pscircle[linewidth=0.4pt,fillstyle=solid,fillcolor=white](5,1){0.2}
\endpspicture
\qquad&=\quad
\pspicture[0.5](0,0)(2,3.5)
\psset{linewidth=\pstlw,xunit=0.5,yunit=0.5,runit=0.5}
\psset{arrowsize=2pt 2,arrowinset=0.2}
\psline(1,0)(1,1)
\psline(4,0)(4,2)
\psarc(1,2){1}{180}{360}
\psarc(3,2){1}{0}{180}
\psarc(1.5,4.5){1.5}{0}{180}
\psline(0,2)(0,4.5)
\psline(3,3)(3,4.5)
\psline(1.5,6)(1.5,7)
\pscircle[linewidth=0.4pt,fillstyle=solid,fillcolor=white](1,1){0.2}
\pscircle[linewidth=0.4pt,fillstyle=solid,fillcolor=white](3,3){0.2}
\pscircle[linewidth=0.4pt,fillstyle=solid,fillcolor=white](1.5,6){0.2}
\rput(3,1.25){$\Delta_{\JJ_C}$}
\endpspicture
\quad=\quad
\pspicture[0.5](0,0)(2,3.5)
\psset{linewidth=\pstlw,xunit=0.5,yunit=0.5,runit=0.5}
\psset{arrowsize=2pt 2,arrowinset=0.2}
\psline(3,0)(3,1)
\psline(0,0)(0,2)
\psarc(3,2){1}{180}{360}
\psarc(1,2){1}{0}{180}
\psarc(2.5,4.5){1.5}{0}{180}
\psline(4,2)(4,4.5)
\psline(1,3)(1,4.5)
\psline(2.5,6)(2.5,7)
\pscircle[linewidth=0.4pt,fillstyle=solid,fillcolor=white](3,1){0.2}
\pscircle[linewidth=0.4pt,fillstyle=solid,fillcolor=white](1,3){0.2}
\pscircle[linewidth=0.4pt,fillstyle=solid,fillcolor=white](2.5,6){0.2}
\rput(1,1.25){$\Delta_{\LL_C}$}
\endpspicture
\end{align}
and see that one can introduce a co-Clifford map. One finds obviously
for an element $a$ of {\it any\/} grade, even inhomogeneous,
\begin{align}
\Delta_C(a) &= (\wedge \otimes \Delta_{\JJ_C})(a_{(1)}\otimes a_{(2)}) \nn
&= (a_{(1)} \wedge C_{(1)}) \otimes  (C_{(2)}\wedge a_{(2)}) \, . 
\end{align}

\subsection{Inversion formulas}

A quite interesting point remains to be examined. Is it possible
to invert the cliffordization process. That is, given a deformed or
cliffordized product $\&c$ or co-cliffordized co-product $\Delta_c$, 
one can obtain back the undeformed product $m_A$ or co-product 
$\Delta$. This is done by the Rota and Stein inversion formulas 
\cite{rota:stein:1994a}, section 4, p.13059. In our notation using
Sweedlers convention about co-products, we find for Gra{\ss}mann-Clifford
products:
\index{inversion formulas}
\begin{align}
i) && B(u,v) &= 
\sum_{(u)(v)} S(u_{(1)})\wedge(u_{(2)} \,\&c\, v_{(1)})\wedge S(v_{(2)}) \nn
ii) && u\wedge v &=
\sum_{(u)(v)} \pm B(S(u_{(1)}),v_{(1)}) (u_{(2)}\,\&c\, v_{(2)}) \nn
iii) && u\wedge v &=
\sum_{(u)(v)} \pm B(u_{(1)},S(v_{(1)})) (u_{(2)}\,\&c\, v_{(2)})  
\end{align} 
where $S$ is the antipode of the undeformed Hopf algebra. The tangles
of this relations read:
\begin{align}
\pspicture[0.5](0,0)(1,3)
\psset{linewidth=\pstlw,xunit=0.5,yunit=0.5,runit=0.5}
\psset{arrowsize=2pt 2,arrowinset=0.2}
\psline(0,6)(0,4)
\psline(2,6)(2,4)
\psarc(1,4){1}{180}{360}
\psline(1,1)(1,0)
\pscircle[linewidth=0.4pt,fillstyle=solid,fillcolor=black](1,3){0.2}
\pscircle[linewidth=0.4pt,fillstyle=solid,fillcolor=white](1,1){0.2}
\rput(1,2.25){$B^{\wedge}$}
\endpspicture
\quad\overset{i)}{=}\quad
\pspicture[0.5](0,0)(3,3)
\psset{linewidth=\pstlw,xunit=0.5,yunit=0.5,runit=0.5}
\psline(1,6)(1,5)
\psline(5,6)(5,5)
\psarc(1,4){1.0}{0}{180}
\psarc(5,4){1.0}{0}{180}
\psarc(3,4){1.0}{180}{360}
\psarc(2,3){1.0}{270}{360}
\psarc(2,4){2.0}{180}{270}
\psarc(3,2){1.0}{180}{270}
\psarc(3,4){3.0}{270}{360}
\psline(3,1)(3,0)
\pscircle[linewidth=0.4pt,fillstyle=solid,fillcolor=white](5,5){0.2}
\pscircle[linewidth=0.4pt,fillstyle=solid,fillcolor=white](1,5){0.2}
\pscircle[linewidth=0.4pt,fillstyle=solid,fillcolor=white](2,2){0.2}
\pscircle[linewidth=0.4pt,fillstyle=solid,fillcolor=white](3,1){0.2}
\pscircle[linewidth=0.4pt,fillstyle=solid,fillcolor=white](3,3){0.2}
\rput*(0.25,3){$S$}
\rput*(5.75,3){$S$}
\rput(3,3.75){$\&c$}
\endpspicture 
\end{align}
and
\begin{align}
\pspicture[0.5](0,0)(1,3)
\psset{linewidth=\pstlw,xunit=0.5,yunit=0.5,runit=0.5}
\psset{arrowsize=2pt 2,arrowinset=0.2}
\psline(0,6)(0,4)
\psline(2,6)(2,4)
\psline(1,3)(1,0)
\psarc(1,4){1}{180}{360}
\pscircle[linewidth=0.4pt,fillstyle=solid,fillcolor=white](1,3){0.2}
\rput(1,3.75){$\wedge$}
\endpspicture     
\quad\overset{ii)}{=}\quad
\pspicture[0.5](0,0)(3,3)
\psset{linewidth=\pstlw,xunit=0.5,yunit=0.5,runit=0.5}
\psline(1,6)(1,5)
\psline(5,6)(5,5)
\psarc(1,4){1.0}{0}{180}
\psarc(5,4){1.0}{0}{180}
\psline(0,4)(0,2)
\psline(6,4)(6,2)
\psline{c-c}(4,4)(2,2)
\psline[border=4pt]{c-c}(2,4)(4,2)
\psarc(1,2){1.0}{180}{360}
\psarc(5,2){1.0}{180}{360}
\psline(5,1)(5,0)
\pscircle[linewidth=0.4pt,fillstyle=solid,fillcolor=black](1,1){0.2}
\pscircle[linewidth=0.4pt,fillstyle=solid,fillcolor=white](1,5){0.2}
\pscircle[linewidth=0.4pt,fillstyle=solid,fillcolor=white](5,1){0.2}
\pscircle[linewidth=0.4pt,fillstyle=solid,fillcolor=white](5,5){0.2}
\rput*(0,3){$S$}
\rput(1.0,1.75){$B$}
\rput(5.0,1.75){$\&c$}
\endpspicture
\quad\overset{iii)}{=}\quad
\pspicture[0.5](0,0)(3,3)
\psset{linewidth=\pstlw,xunit=0.5,yunit=0.5,runit=0.5}
\psline(1,6)(1,5)
\psline(5,6)(5,5)
\psarc(1,4){1.0}{0}{180}
\psarc(5,4){1.0}{0}{180}
\psline(0,4)(0,2)
\psline(6,4)(6,2)
\psline{c-c}(4,4)(2,2)
\psline[border=4pt]{c-c}(2,4)(4,2)
\psarc(1,2){1.0}{180}{360}
\psarc(5,2){1.0}{180}{360}
\psline(5,1)(5,0)
\pscircle[linewidth=0.4pt,fillstyle=solid,fillcolor=black](1,1){0.2}
\pscircle[linewidth=0.4pt,fillstyle=solid,fillcolor=white](1,5){0.2}
\pscircle[linewidth=0.4pt,fillstyle=solid,fillcolor=white](5,1){0.2}
\pscircle[linewidth=0.4pt,fillstyle=solid,fillcolor=white](5,5){0.2}
\rput(1.0,1.75){$B$}
\rput(5,1.75){$\&c$}
\rput*(2.1,2.2){$S$}
\endpspicture
\end{align}
It will be of great value for the later discussed applications to 
provide a few examples of these formulas. Let $x,y,z \in V$, 
one finds:
\begin{align}
B(x,y) &= x\,\&c\, y - x\wedge y \nn
x\wedge y\wedge z &= x \,\&c\, (y\wedge z) - B(x,y)z + B(x,z)y \nn
&= x \,\&c\, y \,\&c\, z -B(x,y)z + B(x,z)y - B(y,z)x.
\end{align}
These are the basic formulas which have been employed in 
\cite{fauser:1996b,fauser:1998a,fauser:2001b} to perform vertex 
normal-ordering, as it will be discussed below. 

The most remarkable fact is, that for the inversion formulas to hold
one needs to have an antipode $S$. It is also the antipode which is
needed in the Connes-Kreimer renormalization method, i.e the antipode
is hidden in the BPHZ formalism of perturbative renormalization. This 
gives strong evidence that quantum field theory should be formulated 
with Hopf gebras.

\section{Convolution algebra}

\index{algebra!convolution}
\index{convolution!algebra}
We have already defined the convolution using the structure tensors 
$\Delta$ and $m$ which are assumed to be associative here. We restrict
our discussion to the endomorphic case.
\begin{align}
\begin{array}{c@{\hskip 2truecm}c}
\Rnode{1}{A} & \Rnode{2}{A\otimes A} \\[8ex]
\Rnode{3}{A} & \Rnode{4}{A\otimes A}
\ncline{->}{1}{2}
\Aput{\Delta}
\ncline{->}{2}{4}
\Aput{f\otimes g}
\ncline{->}{4}{3}
\Bput{m}
\ncline{->}{1}{3}
\Bput{\star}
\end{array}
&&~
\pspicture[0.5](0,0)(1,3)
\psset{linewidth=\pstlw,xunit=0.5,yunit=0.5,runit=0.5}
\psset{arrowsize=2pt 2,arrowinset=0.2}
\psline{-}(1,6)(1,4)
\psline{-}(1,2)(1,0)
\psarc(1,3){1.0}{0}{360}
\pscircle[linewidth=0.4pt,fillstyle=solid,fillcolor=black](0,3){0.2}
\pscircle[linewidth=0.4pt,fillstyle=solid,fillcolor=black](2,3){0.2}
\pscircle[linewidth=0.4pt,fillstyle=solid,fillcolor=white](1,4){0.2}
\pscircle[linewidth=0.4pt,fillstyle=solid,fillcolor=white](1,2){0.2}
\rput(-0.75,3){$f$}
\rput(2.75,3){$g$}
\rput(1.5,5.5){$A$}
\rput(1.5,0.5){$A$}
\endpspicture            
\end{align}

This defines the convolution algebra $\Conv(A,m,\Delta)$ on the 
endomorphisms $f~:~ A\rightarrow A$. The convolution product is
denoted by $\star ~:~ \End A\otimes \End A \rightarrow \End A
\cong A\otimes A^*$. 

A convolution unit $u$ is defined as usual,
\index{unit!convolution}
\index{convolution!unit}
\begin{align}
f\star u ~&=~ f ~=~ u\star f
\end{align}
or in terms of tangles
\begin{align}
\label{eqn-CU}
\pspicture[0.5](0,0)(1,3)
\psset{linewidth=\pstlw,xunit=0.5,yunit=0.5,runit=0.5}
\psset{arrowsize=2pt 2,arrowinset=0.2}
\psline{-}(1,6)(1,5)
\psline{-}(1,1)(1,0)
\psarc(1,4){1.0}{0}{180}
\psarc(1,2){1.0}{180}{360}
\psline(0,2)(0,4)
\psline(2,2)(2,4)
\pscircle[linewidth=0.4pt,fillstyle=solid,fillcolor=white](1,5){0.2}
\pscircle[linewidth=0.4pt,fillstyle=solid,fillcolor=white](1,1){0.2}
\rput*(0,3){$f$}
\rput*(2,3){$u$}
\endpspicture
\quad=\quad
\pspicture[0.5](0,0)(0.5,3)
\psset{linewidth=\pstlw,xunit=0.5,yunit=0.5,runit=0.5}
\psset{arrowsize=2pt 2,arrowinset=0.2}
\psline{-}(1,6)(1,0)
\rput*(1,3){$f$}
\endpspicture
\quad=\quad
\pspicture[0.5](0,0)(1,3)
\psset{linewidth=\pstlw,xunit=0.5,yunit=0.5,runit=0.5}
\psset{arrowsize=2pt 2,arrowinset=0.2}
\psline{-}(1,6)(1,5)
\psline{-}(1,1)(1,0)
\psarc(1,4){1.0}{0}{180}
\psarc(1,2){1.0}{180}{360}
\psline(0,2)(0,4)
\psline(2,2)(2,4)
\pscircle[linewidth=0.4pt,fillstyle=solid,fillcolor=white](1,5){0.2}
\pscircle[linewidth=0.4pt,fillstyle=solid,fillcolor=white](1,1){0.2}
\rput*(0,3){$u$}
\rput*(2,3){$f$}
\endpspicture
\qquad\quad
\forall f.
\end{align}
The unit, if it exists, is unique, since we find for two units $u$ 
and $u^\prime$,
\begin{align}
\pspicture[0.5](0,0)(1,3)
\psset{linewidth=\pstlw,xunit=0.5,yunit=0.5,runit=0.5}
\psset{arrowsize=2pt 2,arrowinset=0.2}
\psline{-}(1,6)(1,5)
\psline{-}(1,1)(1,0)
\psarc(1,4){1.0}{0}{180}
\psarc(1,2){1.0}{180}{360}
\psline(0,2)(0,4)
\psline(2,2)(2,4)
\pscircle[linewidth=0.4pt,fillstyle=solid,fillcolor=white](1,5){0.2}
\pscircle[linewidth=0.4pt,fillstyle=solid,fillcolor=white](1,1){0.2}
\rput*(0,3){$u$}
\rput*(2,3){$u^\prime$}
\endpspicture
\quad=\quad
\pspicture[0.5](0,0)(0.5,3)
\psset{linewidth=\pstlw,xunit=0.5,yunit=0.5,runit=0.5}
\psset{arrowsize=2pt 2,arrowinset=0.2}
\psline{-}(1,6)(1,0)
\rput*(1,3){$u$}
\endpspicture
\quad=\quad
\pspicture[0.5](0,0)(1,3)
\psset{linewidth=\pstlw,xunit=0.5,yunit=0.5,runit=0.5}
\psset{arrowsize=2pt 2,arrowinset=0.2}
\psline{-}(1,6)(1,5)
\psline{-}(1,1)(1,0)
\psarc(1,4){1.0}{0}{180}
\psarc(1,2){1.0}{180}{360}
\psline(0,2)(0,4)
\psline(2,2)(2,4)
\pscircle[linewidth=0.4pt,fillstyle=solid,fillcolor=white](1,5){0.2}
\pscircle[linewidth=0.4pt,fillstyle=solid,fillcolor=white](1,1){0.2}
\rput*(0,3){$u$}
\rput*(2,3){$u$}
\endpspicture
\qquad\Rightarrow\quad
\pspicture[0.5](0,0)(0.5,3)
\psset{linewidth=\pstlw,xunit=0.5,yunit=0.5,runit=0.5}
\psset{arrowsize=2pt 2,arrowinset=0.2}
\psline{-}(1,6)(1,0)
\rput*(1,3){$u$}
\endpspicture
\quad=\quad
\pspicture[0.5](0,0)(0.5,3)
\psset{linewidth=\pstlw,xunit=0.5,yunit=0.5,runit=0.5}
\psset{arrowsize=2pt 2,arrowinset=0.2}
\psline{-}(1,6)(1,0)
\rput*(1,3){$u^\prime$}
\endpspicture
\end{align}
If the product $m$ possesses a unit $\eta$ and the co-product possesses 
a counit $\epsilon$, then the convolution unit $u$ is given by 
$u=\eta \circ \epsilon$, since we find
\begin{align}
\pspicture[0.5](0,0)(1,3)
\psset{linewidth=\pstlw,xunit=0.5,yunit=0.5,runit=0.5}
\psset{arrowsize=2pt 2,arrowinset=0.2}
\psline{-}(1,6)(1,5)
\psline{-}(1,1)(1,0)
\psarc(1,4){1.0}{0}{180}
\psarc(1,2){1.0}{180}{360}
\psline(0,2)(0,4)
\psline(2,2)(2,4)
\pscircle[linewidth=0.4pt,fillstyle=solid,fillcolor=white](1,5){0.2}
\pscircle[linewidth=0.4pt,fillstyle=solid,fillcolor=white](1,1){0.2}
\rput*(0,3){$f$}
\rput*(2,3){$u$}
\endpspicture
\quad=\quad
\pspicture[0.5](0,0)(1,3)
\psset{linewidth=\pstlw,xunit=0.5,yunit=0.5,runit=0.5}
\psset{arrowsize=2pt 2,arrowinset=0.2}
\psline{-}(1,6)(1,5)
\psline{-}(1,1)(1,0)
\psarc(1,4){1.0}{0}{180}
\psarc(1,2){1.0}{180}{360}
\psline(0,2)(0,4)
\pscircle[linewidth=0.4pt,fillstyle=solid,fillcolor=white](2,2){0.2}
\pscircle[linewidth=0.4pt,fillstyle=solid,fillcolor=white](2,4){0.2}
\pscircle[linewidth=0.4pt,fillstyle=solid,fillcolor=white](1,5){0.2}
\pscircle[linewidth=0.4pt,fillstyle=solid,fillcolor=white](1,1){0.2}
\rput*(0,3){$f$}
\endpspicture
\quad=\quad
\pspicture[0.5](0,0)(0.5,3)
\psset{linewidth=\pstlw,xunit=0.5,yunit=0.5,runit=0.5}
\psset{arrowsize=2pt 2,arrowinset=0.2}
\psline{-}(1,6)(1,0)
\rput*(1,3){$f$}
\endpspicture
\quad=\quad
\pspicture[0.5](0,0)(1,3)
\psset{linewidth=\pstlw,xunit=0.5,yunit=0.5,runit=0.5}
\psset{arrowsize=2pt 2,arrowinset=0.2}
\psline{-}(1,6)(1,5)
\psline{-}(1,1)(1,0)
\psarc(1,4){1.0}{0}{180}
\psarc(1,2){1.0}{180}{360}
\psline(2,2)(2,4)
\pscircle[linewidth=0.4pt,fillstyle=solid,fillcolor=white](0,2){0.2}
\pscircle[linewidth=0.4pt,fillstyle=solid,fillcolor=white](0,4){0.2}
\pscircle[linewidth=0.4pt,fillstyle=solid,fillcolor=white](1,5){0.2}
\pscircle[linewidth=0.4pt,fillstyle=solid,fillcolor=white](1,1){0.2}
\rput*(2,3){$f$}
\endpspicture
\quad=\quad
\pspicture[0.5](0,0)(1,3)
\psset{linewidth=\pstlw,xunit=0.5,yunit=0.5,runit=0.5}
\psset{arrowsize=2pt 2,arrowinset=0.2}
\psline{-}(1,6)(1,5)
\psline{-}(1,1)(1,0)
\psarc(1,4){1.0}{0}{180}
\psarc(1,2){1.0}{180}{360}
\psline(0,2)(0,4)
\psline(2,2)(2,4)
\pscircle[linewidth=0.4pt,fillstyle=solid,fillcolor=white](1,5){0.2}
\pscircle[linewidth=0.4pt,fillstyle=solid,fillcolor=white](1,1){0.2}
\rput*(0,3){$u$}
\rput*(2,3){$f$}
\endpspicture
\end{align}
where we used the unit and counit properties displayed in Eqns.
\ref{eqn:5-3} and \ref{eqn:5-21}.

Since we are interested in Gra{\ss}mann and Clifford Hopf gebras, 
we are dealing with unital associative algebras and co-gebras. From
product co-product duality it follows, that unital algebras are related
to counital co-gebras.

\begin{thrm}
If $m$ is a Gra{\ss}mann (Clifford) product and $\Delta$ is a 
Gra{\ss}mann (Clifford) co-product, then $\star$ is a unital 
convolution with unit $u=\eta\circ \epsilon$.
\end{thrm} 
\index{convolution!unital}
\index{convolution!Gra{\ss}mann!unital}
\index{convolution!Clifford!unital}

\noindent
{\bf Note} that we use product co-product duality to make statements
about the character of the mutual structure, we do {\it not\/} however
assert the relation $m^*=\Delta$ or equivalently $\Delta^*=m$.

Having established that Gra{\ss}mann or Clifford bi-convolutions are
unital, we can ask if an antipode exists. Recalling the axioms:
\index{tangle:antipode}
\begin{align}
\pspicture[0.5](0,0)(1,3)
\psset{linewidth=\pstlw,xunit=0.5,yunit=0.5,runit=0.5}
\psset{arrowsize=2pt 2,arrowinset=0.2}
\psline{-}(1,6)(1,5)
\psline{-}(1,1)(1,0)
\psarc(1,4){1.0}{0}{180}
\psarc(1,2){1.0}{180}{360}
\psline(0,2)(0,4)
\psline(2,2)(2,4)
\pscircle[linewidth=0.4pt,fillstyle=solid,fillcolor=white](1,5){0.2}
\pscircle[linewidth=0.4pt,fillstyle=solid,fillcolor=white](1,1){0.2}
\rput*(0,3){$S$}
\endpspicture 
\quad=\quad
\pspicture[0.5](0,0)(0.5,3)
\psset{linewidth=\pstlw,xunit=0.5,yunit=0.5,runit=0.5}
\psset{arrowsize=2pt 2,arrowinset=0.2}
\psline{-}(1,6)(1,0)
\rput*(1,3){$u$}
\endpspicture
\quad=\quad
\pspicture[0.5](0,0)(1,3)
\psset{linewidth=\pstlw,xunit=0.5,yunit=0.5,runit=0.5}
\psset{arrowsize=2pt 2,arrowinset=0.2}
\psline{-}(1,6)(1,5)
\psline{-}(1,1)(1,0)
\psarc(1,4){1.0}{0}{180}
\psarc(1,2){1.0}{180}{360}
\psline(0,2)(0,4)
\psline(2,2)(2,4)
\pscircle[linewidth=0.4pt,fillstyle=solid,fillcolor=white](1,5){0.2}
\pscircle[linewidth=0.4pt,fillstyle=solid,fillcolor=white](1,1){0.2}
\rput*(2,3){$S$}
\endpspicture
\quad=\quad
\pspicture[0.5](0,0)(0.5,3)
\psset{linewidth=\pstlw,xunit=0.5,yunit=0.5,runit=0.5}
\psset{arrowsize=2pt 2,arrowinset=0.2}
\psline{-}(1,6)(1,4)
\psline{-}(1,2)(1,0)
\pscircle[linewidth=0.4pt,fillstyle=solid,fillcolor=white](1,4){0.2}
\pscircle[linewidth=0.4pt,fillstyle=solid,fillcolor=white](1,2){0.2}
\rput*(1.75,4){$\epsilon$}
\rput*(1.75,2){$\eta$}
\endpspicture 
\end{align}
we can prove that the Gra{\ss}mann bi-convolution is antipodal. 
Therefore a Gra{\ss}mann algebra is graded and augmented,
by the counit, and that it is connected. Thus the recursive proof 
for the existence of teh antipode from page \pageref{proof-antipode} 
applies.

This argument does not apply to Clifford algebras. Clifford algebras 
are filtered algebras w.r.t. the $\openZ_n$-filtration inherited from the
Gra{\ss}mann algebra when viewed as endomorphism (sub)algebra
of $\bigwedge V$ or inherited from the tensor algebra via factorization.
The Clifford product is $\openZ_2$-graded only, and the counit if
taken as augmentation, yields a non-connected algebra.

\begin{thrm}[Oziewicz 1997 \cite{oziewicz:1997a}]
A Clifford bi-convolution with product $m^B$ based on the scalar product
$B$ and a co-product $\Delta^C$ based on the co-scalar product
$C$ is antipodal if and only if $C^{-1}\not=B$.
\end{thrm}
\index{convolution!antipodal}
\index{convolution!antipodeless}

In other words, if one uses the particular co-product which is gained by 
product co-product duality, the resulting bi-convolution 
$\Conv(A,m^B,\Delta^{B^{-1}})$ is antipodeless. This seemed to be a 
great drawback in the study of Clifford Hopf gebras, and led to the 
study of convolutions with independent product and co-product by 
Oziewicz and coworkers.

Regarding our analysis and recalling the idea of a Peano space
and the Gra{\ss}mann-Cayley algebra, it is quite natural to introduce
{\it independently\/} a wedge and a vee exterior product on $A$ and $A^*$.
Using now product co-product duality, one ends up with an independent product 
and co-product. This fact will be a major part of our analysis
of normal-ordering in quantum field theory. 

The basic point is that one can perform a cliffordization of a 
Gra{\ss}mann algebra w.r.t. a purely antisymmetric bilinear form. The 
resulting product is again an exterior product, but different from
the originally introduced wedge. Such an algebra will be called 
{\it quantum Gra{\ss}mann algebra\/}, for reasons given below.

\section{Crossing from the antipode}

We have up to now identified associative antipodal convolution
algebras with Hopf gebras. It remains to show, that the product is 
a co-gebra homomorphism and that the co-product is an algebra homomorphism. 
This condition is displayed by the following tangle, compare Eqn.
\ref{eqn:5-38}
\begin{align}
\pspicture[0.5](0,0)(1,3)
\psset{linewidth=\pstlw,xunit=0.5,yunit=0.5,runit=0.5}
\psline(0,6)(0,5)
\psline(2,6)(2,5)
\psarc(1,5){1.0}{180}{360}
\psline(1,4)(1,2)
\psarc(1,1){1.0}{0}{180}
\psline(0,1)(0,0)
\psline(2,1)(2,0)
\pscircle[linewidth=0.4pt,fillstyle=solid,fillcolor=white](1,4){0.2}
\pscircle[linewidth=0.4pt,fillstyle=solid,fillcolor=white](1,2){0.2}
\endpspicture
=\quad
\pspicture[0.5](0,0)(3,3)
\psset{linewidth=\pstlw,xunit=0.5,yunit=0.5,runit=0.5}
\psline(1,6)(1,5)
\psline(5,6)(5,5)
\psarc(1,4){1.0}{0}{180}
\psarc(5,4){1.0}{0}{180}
\psline(0,4)(0,2)
\psline(6,4)(6,2)
\psline{c-c}(4,4)(2,2)
\psline[border=4pt]{c-c}(2,4)(4,2)
\psarc(1,2){1.0}{180}{360}
\psarc(5,2){1.0}{180}{360}
\psline(1,1)(1,0)
\psline(5,1)(5,0)
\pscircle[linewidth=0.4pt,fillstyle=solid,fillcolor=white](1,1){0.2}
\pscircle[linewidth=0.4pt,fillstyle=solid,fillcolor=white](1,5){0.2}
\pscircle[linewidth=0.4pt,fillstyle=solid,fillcolor=white](5,1){0.2}
\pscircle[linewidth=0.4pt,fillstyle=solid,fillcolor=white](5,5){0.2}
\endpspicture 
\end{align}
where the crossing occurs in the r.h.s. because
of the crossed products involved in that calculation.

\index{tangle!crossing}
\index{crossing!definition from antipode}
\begin{thrm}[Oziewicz]
The crossing of a bi-associative antipodal bi-convolution is given as
\begin{align}
\pspicture[0.5](0,0)(1,1)
\psset{linewidth=\pstlw,xunit=0.5,yunit=0.5,runit=0.5}
\psbezier(2,2)(2,1)(0,1)(0,0)
\psbezier[border=4pt,bordercolor=white](0,2)(0,1)(2,1)(2,0)
\endpspicture
\quad&\equiv\quad
\pspicture[0.5](0,0)(3,4)
\psset{linewidth=\pstlw,xunit=0.5,yunit=0.5,runit=0.5}
\psline{-}(1,8)(1,7)
\psline{-}(5,8)(5,7)
\psline{-}(1,1)(1,0)
\psline{-}(5,1)(5,0)
\psline{-}(0,6)(0,2)
\psline{-}(6,6)(6,2)
\psline{-}(3,5)(3,3)
\psarc(1,6){1.0}{0}{180}
\psarc(5,6){1.0}{0}{180}
\psarc(3,6){1.0}{180}{360}
\psarc(3,2){1.0}{0}{180}
\psarc(1,2){1.0}{180}{360}
\psarc(5,2){1.0}{180}{360}
\pscircle[linewidth=0.4pt,fillstyle=solid,fillcolor=white](1,7){0.2}
\pscircle[linewidth=0.4pt,fillstyle=solid,fillcolor=white](5,7){0.2}
\pscircle[linewidth=0.4pt,fillstyle=solid,fillcolor=white](3,5){0.2}
\pscircle[linewidth=0.4pt,fillstyle=solid,fillcolor=white](3,3){0.2}
\pscircle[linewidth=0.4pt,fillstyle=solid,fillcolor=white](1,1){0.2}
\pscircle[linewidth=0.4pt,fillstyle=solid,fillcolor=white](5,1){0.2}
\rput*(0,4){S}
\rput*(6,4){S}
\endpspicture
\end{align}
\end{thrm}
\noindent
{\bf Proof:}
\begin{align}
\pspicture[0.5](0,0)(3,3)
\psset{linewidth=\pstlw,xunit=0.5,yunit=0.5,runit=0.5}
\psline(1,6)(1,5)
\psline(5,6)(5,5)
\psarc(1,4){1.0}{0}{180}
\psarc(5,4){1.0}{0}{180}
\psline(0,4)(0,2)
\psline(6,4)(6,2)
\psline{c-c}(4,4)(2,2)
\psline[border=4pt]{c-c}(2,4)(4,2)
\psarc(1,2){1.0}{180}{360}
\psarc(5,2){1.0}{180}{360}
\psline(1,1)(1,0)
\psline(5,1)(5,0)
\pscircle[linewidth=0.4pt,fillstyle=solid,fillcolor=white](1,1){0.2}
\pscircle[linewidth=0.4pt,fillstyle=solid,fillcolor=white](1,5){0.2}
\pscircle[linewidth=0.4pt,fillstyle=solid,fillcolor=white](5,1){0.2}
\pscircle[linewidth=0.4pt,fillstyle=solid,fillcolor=white](5,5){0.2}
\endpspicture 
\qquad&=\qquad
\pspicture[0.5](0,0)(4,5)
\psset{linewidth=\pstlw,xunit=0.5,yunit=0.5,runit=0.5}
\psline(1,10)(1,9)
\psline(7,10)(7,9)
\psline(0,8)(0,2)
\psline(1,7)(1,3)
\psline(4,4)(4,6)
\psline(7,7)(7,3)
\psline(8,8)(8,2)
\psline(1,1)(1,0)
\psline(7,1)(7,0)
\psarc(1,8){1.0}{0}{180}
\psarc(2,7){1.0}{0}{180}
\psarc(4,7){1.0}{180}{360}
\psarc(7,8){1.0}{0}{180}
\psarc(6,7){1.0}{0}{180}
\psarc(4,3){1.0}{0}{180}
\psarc(2,3){1.0}{180}{360}
\psarc(1,2){1.0}{180}{360}
\psarc(6,3){1.0}{180}{360}
\psarc(7,2){1.0}{180}{360}
\pscircle[linewidth=0.4pt,fillstyle=solid,fillcolor=white](1,9){0.2}
\pscircle[linewidth=0.4pt,fillstyle=solid,fillcolor=white](2,8){0.2}
\pscircle[linewidth=0.4pt,fillstyle=solid,fillcolor=white](4,6){0.2}
\pscircle[linewidth=0.4pt,fillstyle=solid,fillcolor=white](6,8){0.2}
\pscircle[linewidth=0.4pt,fillstyle=solid,fillcolor=white](7,9){0.2}
\pscircle[linewidth=0.4pt,fillstyle=solid,fillcolor=white](1,1){0.2}
\pscircle[linewidth=0.4pt,fillstyle=solid,fillcolor=white](2,2){0.2}
\pscircle[linewidth=0.4pt,fillstyle=solid,fillcolor=white](4,4){0.2}
\pscircle[linewidth=0.4pt,fillstyle=solid,fillcolor=white](6,2){0.2}
\pscircle[linewidth=0.4pt,fillstyle=solid,fillcolor=white](7,1){0.2}
\rput*(1,5){S}
\rput*(7,5){S}
\endpspicture
\qquad=
\nonumber
\end{align}
\begin{align}
=\qquad
\pspicture[0.5](0,0)(4,5)
\psset{linewidth=\pstlw,xunit=0.5,yunit=0.5,runit=0.5}
\psline(2,10)(2,9)
\psline(6,10)(6,9)
\psline(0,7)(0,3)
\psline(2,7)(2,3)
\psline(4,3)(4,7)
\psline(6,7)(6,3)
\psline(8,7)(8,3)
\psline(2,1)(2,0)
\psline(6,1)(6,0)
\psarc(2,8){1.0}{0}{180}
\psarc(1,7){1.0}{0}{180}
\psarc(4,8){1.0}{180}{360}
\psarc(6,8){1.0}{0}{180}
\psarc(7,7){1.0}{0}{180}
\psarc(4,2){1.0}{0}{180}
\psarc(1,3){1.0}{180}{360}
\psarc(2,2){1.0}{180}{360}
\psarc(7,3){1.0}{180}{360}
\psarc(6,2){1.0}{180}{360}
\pscircle[linewidth=0.4pt,fillstyle=solid,fillcolor=white](2,9){0.2}
\pscircle[linewidth=0.4pt,fillstyle=solid,fillcolor=white](1,8){0.2}
\pscircle[linewidth=0.4pt,fillstyle=solid,fillcolor=white](4,7){0.2}
\pscircle[linewidth=0.4pt,fillstyle=solid,fillcolor=white](6,9){0.2}
\pscircle[linewidth=0.4pt,fillstyle=solid,fillcolor=white](7,8){0.2}
\pscircle[linewidth=0.4pt,fillstyle=solid,fillcolor=white](1,2){0.2}
\pscircle[linewidth=0.4pt,fillstyle=solid,fillcolor=white](2,1){0.2}
\pscircle[linewidth=0.4pt,fillstyle=solid,fillcolor=white](4,3){0.2}
\pscircle[linewidth=0.4pt,fillstyle=solid,fillcolor=white](6,1){0.2}
\pscircle[linewidth=0.4pt,fillstyle=solid,fillcolor=white](7,2){0.2}
\rput*(2,5){S}
\rput*(6,5){S}
\endpspicture
\qquad&=\qquad
\pspicture[0.5](0,0)(1,4)
\psset{linewidth=\pstlw,xunit=0.5,yunit=0.5,runit=0.5}
\psline(0,8)(0,6)
\psline(2,8)(2,6)
\psline(1,5)(1,3)
\psline(0,2)(0,0)
\psline(2,2)(2,0)
\psarc(1,6){1.0}{180}{360}
\psarc(1,2){1.0}{0}{180}
\pscircle[linewidth=0.4pt,fillstyle=solid,fillcolor=white](1,5){0.2}
\pscircle[linewidth=0.4pt,fillstyle=solid,fillcolor=white](1,3){0.2}
\endpspicture
\end{align}
where we have used bi-associativity and the antipode axioms.

In fact, the 'crossing' is a planar graph, containing no 'over' 
or 'under' information, and has {\it a priori\/}
nothing to do with knots and links and their projections. This
tangle will also not fulfil in general the Reidemeister moves
of Eqn. \ref{eqn:4-55} or \ref{eqn:4-56}, even if the considered
structure had fulfilled Eqn. \ref{eqn:4-54}. Hence the name crossing
is quite misleading and should probably be replaced with
\index{scattering}
{\it scattering\/} or {\it transmutation}, since it is a generalized switch. 
However, we will stay with the term crossing since this tangle is employed
in building crossed products. It is of utmost importance to classify 
such crossings. Unfortunately not very much is known till now. For what 
type of structure tensors $m,\Delta$ is the crossing a pre-braid,
preserves a grading, filtration etc. ? 

Nevertheless, this theorem answers in part our main question about 
the independence of the structure tensors. We can reformulate it as 
follows: Given a bi-associative antipodal convolution, then the crossing 
is a function of the structure tensors $m,\Delta$ and $S$ which 
itself is a function of $m,\Delta$ -- in physicists notation:
\begin{align}
\tau &= \tau(m,\Delta) &&S=S(m,\Delta).
\end{align}

\noindent
{\bf Example:} We consider a 2-dimensional space $V$ and its dual space 
$V^*$. Let $\{ e_i\},~ i\in \{1,2\}$ be an arbitrary basis of $V$ and 
let $\{\epsilon^j\},~j\in \{1,2\}$ be the canonical dual basis of $V^*$ 
w.r.t. the basis of $V$ defined via $\epsilon^i(e_j)=\delta^i_j$. 
Introduce a  scalar product $B$ and a co-scalar product $C$ as
\begin{align}
B&\cong\left[
\begin{array}{cc}
a &b \\ c & d
\end{array}
\right]
&&
C\cong\left[
\begin{array}{cc}
u &z \\ w & v
\end{array}
\right].
\end{align}
The Clifford algebra $\Cl(V,B)$ has a Gra{\ss}mann basis $GB = \{
\Id, e_1, e_2, e_1\wedge e_2 \}$. The co-scalar product induces the
following co-product of $\Id$
\begin{align}
\Delta(\Id) &= C^{\wedge}_{(1)}\otimes C^{\wedge}_{(2)} \nn 
&=\Id\otimes \Id + 
u\,e_1\otimes e_1 + z\,e_1\otimes e_2 +
w\,e_2\otimes e_1 + v\,e_2\otimes e_2 
\nn
&+(zw-uv)\, e_1\wedge e_2\otimes e_1\wedge e_2.
\end{align}
Therefrom any co-product can be calculated by Rota-Stein
co-cliffordization.

The antipode is a linear operator on $\bigwedge V$ and can be represented
in the above defined basis and its dual basis on $\bigvee V^*$. One
finds using BIGEBRA \cite{BIGEBRA}
\index{antipode!bi-cliffordized in dim 2}
\begin{align}
S&= S^{a}_{\phantom{a}b}\, e_a\otimes \epsilon^b\,\cong\,
1/N \left[
\begin{array}{cccc}
1+(c-b)(w-z) & 0 & 0 & -c+b \\
0 & -1 & 0 & 0 \\
0 &  0 & -1 & 0 \\
z-w & 0 & 0 & 1
\end{array}\right] 
\nn
N&=(1-\tr(BC)+\det(BC)\,)
\end{align}

Only if {\it both\/} scalar products are symmetric, the antipode is 
grade preserving, hence a graded morphism. If either one or the other
scalar product is symmetric, $S$ has a triangular representation
as matrix and preserves the filtration of either $\bigwedge V$ or 
$\bigvee V^*$. In the case of a quantum Gra{\ss}mann Hopf algebra, 
i.e. $B\not=C^{-1}$, $B^T=-B$, $C^T=-C$, one arrives also at an 
antipode being a non-graded $\openk$-module morphism.
This will be of extreme importance in the theory of
perturbative renormalization according to Connes-Kreimer,
since it flaws the recursive formula for calculating the antipode,
see page \pageref{proof-antipode}. But this formula, enriched by the
renormalization scheme, is equivalent there to the Zimmermann
forest formula, \cite{connes:kreimer:1999a,kreimer:2000a,%
brouder:1999a,brouder:1999b,brouder:2000a}.

Regarding our example, current computer algebra can derive the crossing, 
which turns out to be a cumbersome expression, however, it is not able yet
to manage to calculate e.g. the minimal polynomial of the derived crossing, 
or to detect if it is a braid. A special case was, however, discussed 
in Fauser and Oziewicz \cite{fauser:oziewicz:2001a}.

\section{Local versus non-local products and co-products}

\subsection{Kuperberg Lemma 3.2. revisited}

\index{Kuperberg!Lemma 3.2}
In \cite{kuperberg:1991a} Kuperberg proved a lemma which establishes
some important relations widely used in the theory of 3-manifold 
invariants, e.g. \cite{kauffman:1999a}. Moreover, Kuperberg showed
that a certain invariant of Hopf algebras is connected to the 
Kauffman bracket \cite{kauffman:1991a}, i.e. the Jones polynomial.
This polynomial has however also an impact on quantum field theory 
as Witten showed \cite{witten:1989a}. It seems to be generally 
not well known that Hopf gebras can be defined without being connected.
Only recently such structures have been studied by Nill et al.
\cite{hausser:nill:1999a,nill:szlachanyi:wiesbrock:1998a,nill:1998a},
however, we reject the term weak Hopf algebra.

Kuperberg's lemma takes thus a central part in the theory of 3-manifold
invariants and quantum field theory as promoted by Witten. Regarding our
results, we will revisit this lemma. In our notation it reads

\begin{thrm}[Kuperberg \cite{kuperberg:1991a}, Lemma 3.2]
The following identities hold in {\em  any\/} Hopf algebra:
\begin{align}
&a)\quad
\pspicture[0.5](0,-0.5)(1,1.0)
\psset{linewidth=\pstlw,xunit=0.5,yunit=0.5,runit=0.5}
\psset{arrowsize=2pt 2,arrowinset=0.2}
\psline(0,3)(0,2)
\psline(2,3)(2,2)
\psarc(1,2){1}{180}{360}
\psline(1,1)(1,0)
\pscircle[linewidth=0.4pt,fillstyle=solid,fillcolor=white](1,0){0.2}
\pscircle[linewidth=0.4pt,fillstyle=solid,fillcolor=white](1,1){0.2}
\rput(1,1.75){m}
\rput(1,-0.75){$\epsilon_A$}
\endpspicture
\quad=\quad
\pspicture[0.5](0,-0.5)(1,1.0)
\psset{linewidth=\pstlw,xunit=0.5,yunit=0.5,runit=0.5}
\psset{arrowsize=2pt 2,arrowinset=0.2}
\psline(0,3)(0,0)
\psline(2,3)(2,0)
\pscircle[linewidth=0.4pt,fillstyle=solid,fillcolor=white](0,0){0.2}
\pscircle[linewidth=0.4pt,fillstyle=solid,fillcolor=white](2,0){0.2}
\rput(0,-0.75){$\epsilon_A$}
\rput(2,-0.75){$\epsilon_A$}
\endpspicture 
&&b)\quad
\pspicture[0.5](0,0)(1,3)
\psset{linewidth=\pstlw,xunit=0.5,yunit=0.5,runit=0.5}
\psset{arrowsize=2pt 2,arrowinset=0.2}
\psline(0,6)(0,4)
\psline(2,6)(2,4)
\psline(1,3)(1,0)
\psarc(1,4){1}{180}{360}
\pscircle[linewidth=0.4pt,fillstyle=solid,fillcolor=white](1,3){0.2}
\rput*(1,1){S}
\endpspicture
\quad=\quad
\pspicture[0.5](0,0)(1.5,3)
\psset{linewidth=\pstlw,xunit=0.5,yunit=0.5,runit=0.5}
\psset{arrowsize=2pt 2,arrowinset=0.2}
\psline(0,6)(0,4)
\psline(2,6)(2,4)
\psbezier(2,4)(2,3)(0,3)(0,2)
\psbezier[border=4pt,bordercolor=white](0,4)(0,3)(2,3)(2,2)
\psarc(1,2){1}{180}{360}
\psline(1,1)(1,0)
\pscircle[linewidth=0.4pt,fillstyle=solid,fillcolor=white](1,1){0.2}
\rput*(0,5){S}
\rput*(2,5){S}
\psbezier(2.25,4)(3,4)(2.25,2)(2.75,2)
\psbezier(2.25,0)(3,0)(2.25,2)(2.75,2)
\rput(4,2){$m^{op}$}
\endpspicture     
\nonumber \\[2ex]
&c)\quad
\pspicture[0.5](0,0)(0.5,3)
\psset{linewidth=\pstlw,xunit=0.5,yunit=0.5,runit=0.5}
\psset{arrowsize=2pt 2,arrowinset=0.2}
\psline(0.5,6)(0.5,0)
\pscircle[linewidth=0.4pt,fillstyle=solid,fillcolor=white](0.5,6){0.2}
\rput*(0.5,3){S}
\endpspicture
\quad=\quad
\pspicture[0.5](0,0)(0.5,3)
\psset{linewidth=\pstlw,xunit=0.5,yunit=0.5,runit=0.5}
\psset{arrowsize=2pt 2,arrowinset=0.2}
\psline(0.5,6)(0.5,0)
\pscircle[linewidth=0.4pt,fillstyle=solid,fillcolor=white](0.5,6){0.2}
\endpspicture
&&d)\quad
\pspicture[0.5](0,0)(1,3)
\psset{linewidth=\pstlw,xunit=0.5,yunit=0.5,runit=0.5}
\psset{arrowsize=2pt 2,arrowinset=0.2}
\psline{-}(1,6)(1,5)
\psline{-}(1,1)(1,0)
\psarc(1,4){1.0}{0}{180}
\psarc(1,2){1.0}{180}{360}
\psline(0,2)(0,4)
\psline(2,2)(2,4)
\pscircle[linewidth=0.4pt,fillstyle=solid,fillcolor=white](1,5){0.2}
\pscircle[linewidth=0.4pt,fillstyle=solid,fillcolor=white](1,1){0.2}
\rput*(2,3){$S^{-1}$}
\rput(1,1.75){$m^{op}$}
\endpspicture 
\quad=\quad
\pspicture[0.5](0,0)(0.5,3)
\psset{linewidth=\pstlw,xunit=0.5,yunit=0.5,runit=0.5}
\psset{arrowsize=2pt 2,arrowinset=0.2}
\psline{-}(0.5,6)(0.5,3.5)
\psline{-}(0.5,2.5)(0.5,0)
\pscircle[linewidth=0.4pt,fillstyle=solid,fillcolor=white](0.5,3.5){0.2}
\pscircle[linewidth=0.4pt,fillstyle=solid,fillcolor=white](0.5,2.5){0.2}
\rput(1,3.5){$\epsilon$}
\rput(1,2.5){$\eta$}
\endpspicture
\end{align}
\end{thrm}
\vfill\newpage
\noindent
Kuperberg's {\bf proof} of a) is as follows:
\begin{align}
\pspicture[0.5](0,-0.5)(1,1.0)
\psset{linewidth=\pstlw,xunit=0.5,yunit=0.5,runit=0.5}
\psset{arrowsize=2pt 2,arrowinset=0.2}
\psline(0,3)(0,2)
\psline(2,3)(2,2)
\psarc(1,2){1}{180}{360}
\psline(1,1)(1,0)
\pscircle[linewidth=0.4pt,fillstyle=solid,fillcolor=white](1,0){0.2}
\pscircle[linewidth=0.4pt,fillstyle=solid,fillcolor=white](1,1){0.2}
\endpspicture
\quad&\Rightarrow\quad
\pspicture[0.5](0,0)(3,3)
\psset{linewidth=\pstlw,xunit=0.5,yunit=0.5,runit=0.5}
\psline(1,6)(1,5)
\psline(5,6)(5,5)
\psarc(1,4){1.0}{0}{180}
\psarc(5,4){1.0}{0}{180}
\psline(0,4)(0,2)
\psline(6,4)(6,2)
\psline{c-c}(4,4)(2,2)
\psline[border=4pt]{c-c}(2,4)(4,2)
\psarc(1,2){1.0}{180}{360}
\psarc(5,2){1.0}{180}{360}
\psline(1,1)(1,0.25)
\psline(5,1)(5,0)
\pscircle[linewidth=0.4pt,fillstyle=solid,fillcolor=white](1,1){0.2}
\pscircle[linewidth=0.4pt,fillstyle=solid,fillcolor=white](1,5){0.2}
\pscircle[linewidth=0.4pt,fillstyle=solid,fillcolor=white](5,1){0.2}
\pscircle[linewidth=0.4pt,fillstyle=solid,fillcolor=white](5,5){0.2}
\pscircle[linewidth=0.4pt,fillstyle=solid,fillcolor=white](1,0.25){0.2}
\pspolygon[linestyle=dotted](-0.5,-0.5)(-0.5,2.5)(2.5,2.5)(2.5,-0.5)
\endpspicture      
\quad=\quad
\pspicture[0.5](0,0)(1,3)
\psset{linewidth=\pstlw,xunit=0.5,yunit=0.5,runit=0.5}
\psline(0,6)(0,5)
\psline(2,6)(2,5)
\psarc(1,5){1.0}{180}{360}
\psline(1,4)(1,2)
\psarc(1,1){1.0}{0}{180}
\psline(0,1)(0,0)
\psline(2,1)(2,0)
\pscircle[linewidth=0.4pt,fillstyle=solid,fillcolor=white](1,4){0.2}
\pscircle[linewidth=0.4pt,fillstyle=solid,fillcolor=white](1,2){0.2}
\pscircle[linewidth=0.4pt,fillstyle=solid,fillcolor=white](0,0){0.2}
\endpspicture         
\quad=\quad
\pspicture[0.5](0,0)(1,1.5)
\psset{linewidth=\pstlw,xunit=0.5,yunit=0.5,runit=0.5}
\psset{arrowsize=2pt 2,arrowinset=0.2}
\psline(0,3)(0,2)
\psline(2,3)(2,2)
\psarc(1,2){1}{180}{360}
\psline(1,1)(1,0)
\pscircle[linewidth=0.4pt,fillstyle=solid,fillcolor=white](1,1){0.2}
\endpspicture
\quad=
\nonumber \\[3ex]
&=\quad
\pspicture[0.5](0,0)(3,3)
\psset{linewidth=\pstlw,xunit=0.5,yunit=0.5,runit=0.5}
\psline(1,6)(1,5)
\psline(5,6)(5,5)
\psarc(1,4){1.0}{0}{180}
\psarc(5,4){1.0}{0}{180}
\psline(6,4)(6,2)
\psline[border=4pt]{c-c}(2,4)(4,2)
\psarc(5,2){1.0}{180}{360}
\psline(5,1)(5,0)
\pscircle[linewidth=0.4pt,fillstyle=solid,fillcolor=white](0,4){0.2}
\pscircle[linewidth=0.4pt,fillstyle=solid,fillcolor=white](4,4){0.2}
\pscircle[linewidth=0.4pt,fillstyle=solid,fillcolor=white](1,5){0.2}
\pscircle[linewidth=0.4pt,fillstyle=solid,fillcolor=white](5,1){0.2}
\pscircle[linewidth=0.4pt,fillstyle=solid,fillcolor=white](5,5){0.2}
\endpspicture
\quad=\quad
\pspicture[0.5](0,0)(3,3)
\psset{linewidth=\pstlw,xunit=0.5,yunit=0.5,runit=0.5}
\psline(1,6)(1,5)
\psline(5,6)(5,5)
\psarc(1,4){1.0}{0}{180}
\psarc(5,4){1.0}{0}{180}
\psline(0,4)(0,2)
\psline(6,4)(6,2)
\psline{c-c}(4,4)(2,2)
\psline[border=4pt]{c-c}(2,4)(4,2)
\psarc(5,2){1.0}{180}{360}
\psline(0,2)(0,0.5)
\psline(2,2)(2,0.5)
\psline(5,1)(5,0)
\pscircle[linewidth=0.4pt,fillstyle=solid,fillcolor=white](0,0.5){0.2}
\pscircle[linewidth=0.4pt,fillstyle=solid,fillcolor=white](2,0.5){0.2}
\pscircle[linewidth=0.4pt,fillstyle=solid,fillcolor=white](1,5){0.2}
\pscircle[linewidth=0.4pt,fillstyle=solid,fillcolor=white](5,1){0.2}
\pscircle[linewidth=0.4pt,fillstyle=solid,fillcolor=white](5,5){0.2}
\pspolygon[linestyle=dotted](-0.5,-0.5)(-0.5,2.5)(2.5,2.5)(2.5,-0.5)
\endpspicture
\quad&\Rightarrow\quad
\pspicture[0.5](0,-0.5)(1,1.0)
\psset{linewidth=\pstlw,xunit=0.5,yunit=0.5,runit=0.5}
\psset{arrowsize=2pt 2,arrowinset=0.2}
\psline(0,3)(0,0)
\psline(2,3)(2,0)
\pscircle[linewidth=0.4pt,fillstyle=solid,fillcolor=white](0,0){0.2}
\pscircle[linewidth=0.4pt,fillstyle=solid,fillcolor=white](2,0){0.2}
\endpspicture     
\end{align}
Furthermore, the proof of b) requires a), the proof of c) requires b) 
and the proof of d) requires c). Hence every assertion of the lemma is 
true if and only if a) holds.

But, we have already displayed a counterexample to property a)~, see
theorem \ref{thrm:6-5}! If a Hopf gebra has a non-local product, e.g. 
derived by cliffordization, we showed that the counit is not an 
algebra homomorphism, i.e. does not satisfy a). So, where is the error 
in the very suggestive proof?

There are two points to be criticized. First, we note that the proof 
uses a technique where to a tangle $x$ a 'helping tangle' $t$ or $b$
\index{tangle!helping}
is added and after some manipulations it is removed unaltered. That is 
this tangle acts like a catalyser in chemistry. However, if the
helping tangle $t$ is added on top it has to be monic and if the
helping tangle $b$ is added from below it has to be epi to guarantee
that a cancellation law holds, i.e. the tangle can be safely removed 
again. Let $x=y$ be the searched for tangle equation and $t$ from the top
and $b$ from the bottom added tangles, one computes in algebraic terms
\begin{align}
x &\Rightarrow\quad
xt=\ldots =yt \quad\Rightarrow\qquad x=y \nn
x &\Rightarrow\quad
bx=\ldots =by \quad\Rightarrow\qquad x=y
\end{align}
and a cancellation takes place. Now it is known \cite{oziewicz:2001a}
that the tangle $\Delta\circ m$, which occurs after the first
equality sign in the proof of a), is {\it not\/} invertible. Hence we
cannot assume a cancellation law to hold! Secondly, we noticed already
that the crossing is defined by the structure tensors $m,\Delta$ since 
an antipode exists. Hence one cannot assume that
\begin{align}
\pspicture[0.5](0,0)(1,1)
\psset{linewidth=\pstlw,xunit=0.5,yunit=0.5,runit=0.5}
\psline(0,2)(2,0)
\psline(2,2)(2,1)
\pscircle[linewidth=0.4pt,fillstyle=solid,fillcolor=white](2,1){0.2}
\endpspicture
\quad=\quad
\pspicture[0.5](0,0)(1,1)
\psset{linewidth=\pstlw,xunit=0.5,yunit=0.5,runit=0.5}
\psline(2,2)(0,0)
\psline[border=4pt,bordercolor=white](0,2)(2,0)
\pscircle[linewidth=0.4pt,fillstyle=solid,fillcolor=white](0,0){0.2}
\endpspicture
\end{align}
holds true, but this has to be proved. In fact, our counterexample,
obtained by cliffordization in theorem \ref{thrm:6-2} shows that
the assertion a) is not true for {\it any\/} Hopf (al)gebra. In fact,
if a Hopf algebra is {\it defined\/} to be connected, as in Sweedler
or Milnor Moore, Kuperberg's assertion is evidently true, but a) is
an axiom in this case and has not to be proved at all. We reformulate 
the Kuperberg lemma as follows:
\begin{thrm}
\index{Kuperberg!Lemma 3.2. revised}
In any Hopf algebra where the co-product is counital and the counit 
is an algebra homomorphism and the product is unital, the unit is 
a co-product homomorphism, and Eqn. \ref{eqn-CK} holds then the 
following identities hold:
\begin{align}
b)\qquad
\pspicture[0.5](0,0)(1,3)
\psset{linewidth=\pstlw,xunit=0.5,yunit=0.5,runit=0.5}
\psset{arrowsize=2pt 2,arrowinset=0.2}
\psline(0,6)(0,4)
\psline(2,6)(2,4)
\psline(1,3)(1,0)
\psarc(1,4){1}{180}{360}
\pscircle[linewidth=0.4pt,fillstyle=solid,fillcolor=white](1,3){0.2}
\rput*(1,1){S}
\endpspicture
\quad=\quad
\pspicture[0.5](0,0)(1.5,3)
\psset{linewidth=\pstlw,xunit=0.5,yunit=0.5,runit=0.5}
\psset{arrowsize=2pt 2,arrowinset=0.2}
\psline(0,6)(0,4)
\psline(2,6)(2,4)
\psbezier(2,4)(2,3)(0,3)(0,2)
\psbezier[border=4pt,bordercolor=white](0,4)(0,3)(2,3)(2,2)
\psarc(1,2){1}{180}{360}
\psline(1,1)(1,0)
\pscircle[linewidth=0.4pt,fillstyle=solid,fillcolor=white](1,1){0.2}
\rput*(0,5){S}
\rput*(2,5){S}
\psbezier(2.25,4)(3,4)(2.25,2)(2.75,2)
\psbezier(2.25,0)(3,0)(2.25,2)(2.75,2)
\rput(4,2){$m^{op}$}
\endpspicture     
\quad c)\qquad
\pspicture[0.5](0,0)(0.5,3)
\psset{linewidth=\pstlw,xunit=0.5,yunit=0.5,runit=0.5}
\psset{arrowsize=2pt 2,arrowinset=0.2}
\psline(0.5,6)(0.5,0)
\pscircle[linewidth=0.4pt,fillstyle=solid,fillcolor=white](0.5,6){0.2}
\rput*(0.5,3){S}
\endpspicture
\quad=\quad
\pspicture[0.5](0,0)(0.5,3)
\psset{linewidth=\pstlw,xunit=0.5,yunit=0.5,runit=0.5}
\psset{arrowsize=2pt 2,arrowinset=0.2}
\psline(0.5,6)(0.5,0)
\pscircle[linewidth=0.4pt,fillstyle=solid,fillcolor=white](0.5,6){0.2}
\endpspicture
\quad d)\quad
\pspicture[0.5](0,0)(1,3)
\psset{linewidth=\pstlw,xunit=0.5,yunit=0.5,runit=0.5}
\psset{arrowsize=2pt 2,arrowinset=0.2}
\psline{-}(1,6)(1,5)
\psline{-}(1,1)(1,0)
\psarc(1,4){1.0}{0}{180}
\psarc(1,2){1.0}{180}{360}
\psline(0,2)(0,4)
\psline(2,2)(2,4)
\pscircle[linewidth=0.4pt,fillstyle=solid,fillcolor=white](1,5){0.2}
\pscircle[linewidth=0.4pt,fillstyle=solid,fillcolor=white](1,1){0.2}
\rput*(2,3){$S^{-1}$}
\rput(1,1.75){$m^{op}$}
\endpspicture 
\quad=\quad
\pspicture[0.5](0,0)(0.5,3)
\psset{linewidth=\pstlw,xunit=0.5,yunit=0.5,runit=0.5}
\psset{arrowsize=2pt 2,arrowinset=0.2}
\psline{-}(0.5,6)(0.5,3.5)
\psline{-}(0.5,2.5)(0.5,0)
\pscircle[linewidth=0.4pt,fillstyle=solid,fillcolor=white](0.5,3.5){0.2}
\pscircle[linewidth=0.4pt,fillstyle=solid,fillcolor=white](0.5,2.5){0.2}
\rput(1,3.5){$\epsilon$}
\rput(1,2.5){$\eta$}
\endpspicture
\end{align}
\end{thrm}
\noindent
{\bf Note} that the product has thus to be {\it local\/} as defined
in definition \ref{dfn:6-4}. In terms of the previous chapter this
asserts that the algebra part of the Hopf algebra has to be an 
augmented connected algebra. We required the same for the coal-gebra. 

The proof of $b)$ requires furthermore that 
\begin{align}
\label{eqn-CK}
\pspicture[0.5](0,0)(1,1)
\psset{linewidth=\pstlw,xunit=0.5,yunit=0.5,runit=0.5}
\psset{arrowsize=2pt 2,arrowinset=0.2}
\psbezier(0,2)(0,1)(2,1)(2,0)
\psbezier[border=4pt,bordercolor=white](2,2)(2,1)(0,1)(0,0)
\rput*(0,0){$S$}
\endpspicture
\quad=\quad
\pspicture[0.5](0,0)(1,1)
\psset{linewidth=\pstlw,xunit=0.5,yunit=0.5,runit=0.5}
\psset{arrowsize=2pt 2,arrowinset=0.2}
\psbezier(0,2)(0,1)(2,1)(2,0)
\psbezier[border=4pt,bordercolor=white](2,2)(2,1)(0,1)(0,0)
\rput*(2,2){$S$}
\endpspicture
&&\text{and}
&&
\pspicture[0.5](0,0)(1,1)
\psset{linewidth=\pstlw,xunit=0.5,yunit=0.5,runit=0.5}
\psset{arrowsize=2pt 2,arrowinset=0.2}
\psbezier(0,2)(0,1)(2,1)(2,0)
\psbezier[border=4pt,bordercolor=white](2,2)(2,1)(0,1)(0,0)
\pscircle[linewidth=0.4pt,fillstyle=solid,fillcolor=white](2,0){0.2}
\endpspicture
\quad=\quad
\pspicture[0.5](0,0)(1,1)
\psset{linewidth=\pstlw,xunit=0.5,yunit=0.5,runit=0.5}
\psset{arrowsize=2pt 2,arrowinset=0.2}
\psline(0,2)(0,1)
\psline(2,2)(0,0)
\pscircle[linewidth=0.4pt,fillstyle=solid,fillcolor=white](0,1){0.2}
\endpspicture 
\, ,
\end{align}
which are further requirements on the crossing.

\subsection{Interacting and non-interacting Hopf gebras}

The observations of the previous section lead naturally to the following
questions:

{\bf Q1:} Is any crossing of an antipodal convolution, i.e. Hopf gebra,
which is derived from a local product and a local co-product 'knottish'
in that sense, that operators can be moved on the strings under and over
as done in Eqn. \ref{eqn-CK}?

If we remember the definition of the crossing as planar graph, this is
a non-trivial requirement. Looking at Gra{\ss}mann Hopf gebras, one finds
however that
\begin{align}
(S\otimes S)\circ\hat{\tau}(a\otimes b) &= 
(-)^{|a||b|}(S\otimes S)\circ {\sw}(a\otimes b) \nn
&=(-)^{|a||b|} {\sw} \circ (S\otimes S)(a\otimes b) =
\hat{\tau}\circ(S\otimes S)(a\otimes b) 
\end{align}
which fulfils also
\begin{align}
(f\otimes g)\circ \hat{\tau} &= 
\hat{\tau}\circ(g\otimes f)
\end{align}
for graded morphisms $f,g$.

{\bf Q2:} Does any Hopf gebra which possesses a crossing which fulfils the 
Reidemeister moves, and is thus 'knottish', has a local product and a local 
co-product? Or alternatively, are all Hopf gebras for 3-dimensional 
topology based on augmented connected algebras and connected augmented 
co-algebras? 
\index{Hopf gebra!knottish}

{\bf Q3:} Since cliffordized Hopf gebras possess non-local products,
and co-cliffordized Hopf gebras possess non-local co-products, are 
there 'knottish' such algebras?

{\bf Q4:} What is the topology behind cliffordized, i.e. quantized
Hopf gebras? Is this topology related to non-commutative geometry?

{\bf Q5:} Since Hopf algebras are generalizations of groups, which 
generalized groups are behind the non-local Hopf gebras? However,
some results in this direction are available 
\cite{fauser:2000e,fauser:1999a,fauser:1999b,%
ablamowicz:fauser:2000b,fauser:2000d}.

{\bf Q6:} Is there a reason from physics that cliffordized, i.e.
quantized, structures have to be used e.g. in quantum field theory?

We will not have occasion to answer these question in full detail in this 
treatise, but we will introduce a further notation which might suggest 
a direction towards the answers and which is motivated by demands of 
physics. Moreover we will show that normal-ordering is encoded this way, 
and a recent pre-print of Brouder \cite{brouder:2002a} shows that
cliffordization may be behind the curtain of Epstein-Glaser renormalization
of time-ordered products.

We will see below, that in quantum field theory Hopf gebras constitute
the structure of the generating functionals and of their algebraic
properties. It is quite suggestive, after examining this fact, to use
the counit as the vacuum expectation value 
\cite{fauser:1998a,fauser:2001e}. The below discussed topic of
normal-ordering deals with the connection of local and non-local
structures. This follows also from the fact, obtained in theorem 
\ref{thrm:6-5}, that the counit acting on the Clifford
product gives the cup tangle of the scalar product. Hence this tangle 
replaces the cup tangle in e.g. the Kauffman bracket. The crucial 
fact of the property a) of the Kuperberg Lemma 3.2, as discussed above 
is that the counit is an algebra homomorphism. If we look at the formula 
in terms of an expectation value, from 
\index{counit!as vacuum expectation value}
\begin{align}
\epsilon(ab) = \epsilon(a)\epsilon(b)
& &\langle ab\rangle = \langle a\rangle\langle b\rangle
\end{align}
it follows that one deals with a free theory as is well known that
the factorization of expectation values represents independence. A 
physically non-trivial theory has to have interactions between its 
constituting parts which renders the Hopf gebras with local
products to be a less interesting case. However, it is that case which
is employed, in the theory of knots and links, Kauffman bracket and 
Jones polynomial and therefore in Witten's approach to quantum field theory 
as described in \cite{witten:1989a}.

This motivates the following distinction:
\begin{dfn}
\index{Hopf gebra!interaction}
\index{Hopf gebra!non-interaction}
A Hopf gebra with local product and local co-product, i.e. a bilocal Hopf 
gebra, is called a {\em non-interacting\/} Hopf gebra, otherwise
the Hopf gebra is called {\em interacting\/} Hopf gebra.  
\end{dfn}
\noindent
{\bf Note} that a Hopf gebra is already called interacting if {\it one} 
of the involved products, product {\it or\/} co-product is non-local. 
If the co-product is non-local, then the co-product has no longer the form
\begin{align}
\Delta(x) &= x\otimes \Id + \Id \otimes x + 
{\sum_{(x)}}^\prime x_{(1)} \otimes x_{(2)} \, ,
\end{align}
which is used to derive the recursive form of the antipode. The primed
sum indicates proper sections of $x$, i.e. $x_{(i)}\not=\Id$. Hence, in 
that case the antipode formula also used by Connes and Kreimer in their
renormalization theory for perturbative quantum field theory cannot 
be established. Nevertheless an antipode and convolutive inverse
endomorphisms do exist in such cases too as we showed by direct calculation
using BIGEBRA.


%% file: antipode.tex
\chapter{Integrals, meet, join, unipotents, and `spinorial' antipode}

\section{Integrals}

We introduce a further structure in the convolution algebra, called
{\it integral\/}, see e.g. Sweedler \cite{sweedler:1969a}. 
\begin{dfn}
A {\em left/right integral} is an element $\mu_{L/R} \in A^*$, i.e.
a (multi) covector of the unital convolution $\Conv(A,A)$ which fulfils
\begin{align}
&\pspicture[0.5](0,0)(1,2.5)
\psset{linewidth=\pstlw,xunit=0.5,yunit=0.5,runit=0.5}
\psline(0,0)(0,2)
\psline(2,1)(2,2)
\psarc(1,2){1}{0}{180}
\psline(1,3)(1,5)
\pscircle[linewidth=0.4pt,fillstyle=solid,fillcolor=white](1,3){0.2}
\pscircle[linewidth=0.4pt,fillstyle=solid,fillcolor=black](2,1){0.2}
\rput(2,0.25){$\mu_R$}
\endpspicture
\quad=\quad
\pspicture[0.5](0,0)(0.5,2.5)
\psset{linewidth=\pstlw,xunit=0.5,yunit=0.5,runit=0.5}
\psline(1,0)(1,2)
\psline(1,3)(1,5)
\pscircle[linewidth=0.4pt,fillstyle=solid,fillcolor=white](1,2){0.2}
\pscircle[linewidth=0.4pt,fillstyle=solid,fillcolor=black](1,3){0.2}
\rput(1.75,2.75){$\mu_R$}
\endpspicture
&&
\pspicture[0.5](0,0)(1,2.5)
\psset{linewidth=\pstlw,xunit=0.5,yunit=0.5,runit=0.5}
\psline(0,1)(0,2)
\psline(2,0)(2,2)
\psarc(1,2){1}{0}{180}
\psline(1,3)(1,5)
\pscircle[linewidth=0.4pt,fillstyle=solid,fillcolor=white](1,3){0.2}
\pscircle[linewidth=0.4pt,fillstyle=solid,fillcolor=black](0,1){0.2}
\rput(0,0.25){$\mu_L$}
\endpspicture
\quad=\quad
\pspicture[0.5](0,0)(0.5,2.5)
\psset{linewidth=\pstlw,xunit=0.5,yunit=0.5,runit=0.5}
\psline(1,0)(1,2)
\psline(1,3)(1,5)
\pscircle[linewidth=0.4pt,fillstyle=solid,fillcolor=white](1,2){0.2}
\pscircle[linewidth=0.4pt,fillstyle=solid,fillcolor=black](1,3){0.2}
\rput(1.75,2.75){$\mu_L$}
\endpspicture
\end{align}
\end{dfn}
\noindent
In equation notation this reads for any $x$  
\begin{align}
&(\Id\otimes \mu_R)\Delta(x) = \mu_R(x)\Id
&&(\mu_L\otimes\Id)\Delta(x) = \mu_L(x)\Id \,.
\end{align}

To distinguish integrals from unit and counit, we use black bullets
in the graphical notation. Obviously zero is an integral, but a 
trivial one. Therefore we speak about proper or non-trivial integrals 
if $\mu_{R/L}\not=0$ is non-zero.

Using duality we define analogously cointegrals.
\begin{dfn}
A {\em left/right cointegral} is an element $e_{L/R} \in A$, i.e.
a (multi) vector of the counital convolution $\Conv(A,A)$ which fulfils
\begin{align}
&\pspicture[0.5](0,0)(1,2.5)
\psset{linewidth=\pstlw,xunit=0.5,yunit=0.5,runit=0.5}
\psline(0,5)(0,3)
\psline(2,4)(2,3)
\psarc(1,3){1}{180}{360}
\psline(1,2)(1,0)
\pscircle[linewidth=0.4pt,fillstyle=solid,fillcolor=white](1,2){0.2}
\pscircle[linewidth=0.4pt,fillstyle=solid,fillcolor=black](2,4){0.2}
\rput(2,4.75){$e_R$}
\endpspicture
\quad=\quad
\pspicture[0.5](0,0)(0.5,2.5)
\psset{linewidth=\pstlw,xunit=0.5,yunit=0.5,runit=0.5}
\psline(1,5)(1,3)
\psline(1,2)(1,0)
\pscircle[linewidth=0.4pt,fillstyle=solid,fillcolor=white](1,3){0.2}
\pscircle[linewidth=0.4pt,fillstyle=solid,fillcolor=black](1,2){0.2}
\rput(1.75,2.25){$e_R$}
\endpspicture
&&
\pspicture[0.5](0,0)(1,2.5)
\psset{linewidth=\pstlw,xunit=0.5,yunit=0.5,runit=0.5}
\psline(0,4)(0,3)
\psline(2,5)(2,3)
\psarc(1,3){1}{180}{360}
\psline(1,2)(1,0)
\pscircle[linewidth=0.4pt,fillstyle=solid,fillcolor=white](1,2){0.2}
\pscircle[linewidth=0.4pt,fillstyle=solid,fillcolor=black](0,4){0.2}
\rput(0,4.75){$e_L$}
\endpspicture
\quad=\quad
\pspicture[0.5](0,0)(0.5,2.5)
\psset{linewidth=\pstlw,xunit=0.5,yunit=0.5,runit=0.5}
\psline(1,5)(1,3)
\psline(1,2)(1,0)
\pscircle[linewidth=0.4pt,fillstyle=solid,fillcolor=white](1,3){0.2}
\pscircle[linewidth=0.4pt,fillstyle=solid,fillcolor=black](1,2){0.2}
\rput(1.75,2.25){$e_L$}
\endpspicture
\end{align}
\end{dfn}
\noindent
In equation notation this reads, for any $x\in \bigwedge V$
\begin{align}
&m(x\otimes e_R) = \epsilon(x)\, e_R
&&m(e_L\otimes x) = \epsilon(x)\, e_L \,.
\end{align}

Of course the term integral is taken since $(\mu_{L}\otimes Id) \circ \Delta$
($Id \otimes \mu_{R}) \circ \Delta$) is a linear form acting in such a way 
that the result is a scalar and the action is linear in the argument,
see Sweedler \cite{sweedler:1969a}.

\noindent
{\bf Example:} (continued)  We consider once more the Clifford biconvolution 
$\Cl(B,C)$ in $\dim V=2$, $\dim \bigwedge V = 2^2=4$ and $B,C$ as defined 
above in the chosen basis. We find the following
\begin{thrm}
For a Gra{\ss}mann biconvolution, i.e. a Gra{\ss}mann Hopf gebra,
$\bigwedge V = \Cl(0,0)$, i.e. $B\equiv 0 \equiv C$ identical zero,
there exists one and only one non-trivial left and right integral
$\mu = \epsilon^{12}$ where $\epsilon^{12}(e_{12})=1$ and 
$\epsilon^{12}(e_I)=0$, $I\not= (12)$, and there exists one and only
one non-trivial left and right cointegral $e=e_{12}$ where 
$e_{12}(\epsilon^{12})=1$ and $e_{12}(\epsilon^I)=0$, $I\not= (12)$.
\end{thrm}
\noindent
{\bf Proof:} by direct computation using CLIFFORD / BIGEBRA.

In general one finds, using physicists notation, $\gamma^5$ for
the elements of maximal grade, that $\mu=\gamma^5=\epsilon^n\vee
\ldots\vee \epsilon^1$ is the {\it unique integral\/} and 
$e=\gamma_5=e_1\wedge\ldots\wedge e_n$ is the {\it unique cointegral\/}
in the $n$-dimensional case of Gra{\ss}mann Hopf gebra. One should
thereby remember, that Gra{\ss}mann Hopf gebras are bi-augmented
and bi-connected and are thus well behaved. This situation changes 
drastically if we turn the products and co-products into non-local 
ones by cliffrodization.

\begin{thrm}
For a Clifford biconvolution $\Cl(B,C)$ ($\dim V=2$) as defined above 
one obtains no non-trivial integral unless $C\equiv 0$ and no non-trivial 
cointegral unless $B\equiv 0$, i.e. unless the cliffordization is 
trivial. 
\end{thrm}
\noindent
{\bf Proof:} by direct computation using CLIFFORD / BIGEBRA.

This result should be compared with various claims of the existence
of integrals, e.g. see \cite{kauffman:radford:1999a}. A theory of 
integrals for Hopf algebras having non-local products, called
quasi Hopf algebras, was developed in 
\cite{nill:szlachanyi:wiesbrock:1998a,nill:1998a,hausser:nill:1999a}.
Moreover, we have no doubts that these results can be generalized to
arbitray dimensions which needs an algebraic proof.

The result we found above for cliffordized and thereby non-local 
products and co-products agrees with the fact that Gra{\ss}mann 
algebras are faithfully represented on the module they are built over.
That is, the left/right regular representation $L_a b = ab$ ($R_a b = ba$)
is irreducible. This follows from the fact that $\Id$ is the only 
non-trivial idempotent element in $\bigwedge V$. Hence a minimal left/right 
ideal in $\bigwedge V$ is $\bigwedge V$ itself. This is a bi-ideal.

In the case of a cliffordized algebra one obtains new primitive 
idempotent elements $f_i^2=f_i$ with $\Id=\sum f_i$ and $f_if_j=f_jf_i$.
Such primitive idempotents generate left/right ideals which carry 
faithful representations, called {\it spinor representations}:
\begin{align}
{\cal S}_f &\cong \Cl\,f \nn
L_{\Cl} \, {\cal S}_f &\cong {\cal S}_f.
\end{align}
Let for the moment $\openk$ be the field of real number having no 
non-trivial involutive automorphisms. Any $f_i$ can take over the 
role of the unit $\eta~:~\openk\rightarrow A$ and 
$f_i\,\Cl\,f_i \cong \openk\, f_i$ can represent the base ring
the algebra is built over. In fact $f_i$ can take over the role of a 
cointegral on a left/right spinor space ${\cal S}_{f,L}, {\cal S}_{f,R}$,
which is a graded $\openk$-module:
\begin{align}
&\pspicture[0.5](0,0)(1,3)
\psset{linewidth=\pstlw,xunit=0.5,yunit=0.5,runit=0.5}
\psline(0,5)(0,3)
\psline(2,4)(2,3)
\psarc(1,3){1}{180}{360}
\psline(1,2)(1,0)
\pscircle[linewidth=0.4pt,fillstyle=solid,fillcolor=white](1,2){0.2}
\pscircle[linewidth=0.4pt,fillstyle=solid,fillcolor=black](2,4){0.2}
\rput(0,5.75){${\cal S}_{f,R}$}
\rput(2,4.75){$f$}
\endpspicture
\quad=\quad
\pspicture[0.5](0,0)(0.5,3)
\psset{linewidth=\pstlw,xunit=0.5,yunit=0.5,runit=0.5}
\psline(1,5)(1,3)
\psline(1,2)(1,0)
\pscircle[linewidth=0.4pt,fillstyle=solid,fillcolor=white](1,3){0.2}
\pscircle[linewidth=0.4pt,fillstyle=solid,fillcolor=black](1,2){0.2}
\rput(1.75,3){$\epsilon_f$}
\rput(1.75,2){$f$}
\endpspicture
&&
\pspicture[0.5](0,0)(1,3)
\psset{linewidth=\pstlw,xunit=0.5,yunit=0.5,runit=0.5}
\psline(0,4)(0,3)
\psline(2,5)(2,3)
\psarc(1,3){1}{180}{360}
\psline(1,2)(1,0)
\pscircle[linewidth=0.4pt,fillstyle=solid,fillcolor=white](1,2){0.2}
\pscircle[linewidth=0.4pt,fillstyle=solid,fillcolor=black](0,4){0.2}
\rput(2,5.75){${\cal S}_{f,L}$}
\rput(0,4.75){$f$}
\endpspicture
\quad=\quad
\pspicture[0.5](0,0)(0.5,3)
\psset{linewidth=\pstlw,xunit=0.5,yunit=0.5,runit=0.5}
\psline(1,5)(1,3)
\psline(1,2)(1,0)
\pscircle[linewidth=0.4pt,fillstyle=solid,fillcolor=white](1,3){0.2}
\pscircle[linewidth=0.4pt,fillstyle=solid,fillcolor=black](1,2){0.2}
\rput(1.75,3){$\epsilon_f$}
\rput(1.75,2){$f$}
\endpspicture
\end{align}
Here we have defined the counit $\epsilon_f$ to be $f\,\Cl\,f \mod f$
and the cointegral is given by $f$ itself. Integrals could be obtained
by categorial duality from this structure.

This is a tremendously important structure since it is directly related with
the representation of elementary particles in quantum field theory.
Moreover, the structure of the state space of a quantum system will
depend strongly on this fact. One finds a decomposition of the unit 
$\Id=\sum f_i$ which induces a direct sum decomposition of the 
representation space, i.e. left/right ideals. Also the counit will 
split along the same lines as $\epsilon=\sum \epsilon_f$. We
will use in this mathematical section only integrals and cointegrals
of Gra{\ss}mann Hopf gebras, see the discussion of meet and join below,
and we will not develop a theory of integrals and cointegrals
for Clifford biconvolution.

Moreover, in the physics sections below we will find that due to Wick
normal-ordering the situation there is much more involved. We will 
discuss these peculiarities there. 

\section{Meet and join}

In this section we will shortly explain in which way integrals and 
cointegrals are involved in Gra{\ss}mann-Cayley algebras.

Starting point is a Gra{\ss}mann Hopf gebra. The interpretation of
1-vectors is that of points of a projective space represented in a 
homogeneous way. That is, $a$ and $\alpha a$, $0\not=\alpha \in \openk$,
are the same point. In fact, the field $\openk$ does not play a major
role in what follows, but we will assume that $\openk$ is a field of 
characteristic 0.

The 'join' of two points $a,b$ is the line $l=\overline{ab}$ which
is represented by the exterior wedge product $\wedge$, i.e. $l=a\wedge b$.
Incidence of an arbitrary point $x$ with that line results in a 
linear dependence of the triple of vectors $a,b,x$ which results
in $x\wedge a\wedge b=0$. This can be found in Gra{\ss}mann 
\cite{grassmann:1878a,grassmann:1862a}. The incidence relation using 
the wedge or 'join' is a non-parametric relation. Such incidence
relations have been discussed recently e.g. in Conradt 
\cite{conradt:2000a}. It is obvious that three independent points 
constitute a plane etc. One notes therefore that the exterior wedge 
or 'join' product raises the grade and increases the dimension of
the geometric objects represented by them.

The dual question to the join is that of a 'meet'. Two lines 
$l=\overline{ab}$ and $m=\overline{cd}$, represented by arbitrary 
points $a,b$ and $c,d$ incident with them, meet eventually in a 
point $x$ or not. The meet is represented by the exterior vee-product
$\vee$, e.g. $l\vee m$. 

Using duality between points and planes in $\openP^3$, one finds that 
the meet constitutes an exterior vee-algebra $\bigvee V^*$ of planes.
That is, the meet of 2 planes is a line and the meet of 3 planes is 
a point. Hence the meet increases the grade in the space constituted 
from planes. It lowers the grade if planes and lines are seen to be 
represented by sets of points.

Gra{\ss}mann introduced the meet in the {\bf A2} by means of an 
Erg\"anzungs operator. This {\it Erg\"anzung\/} is related to an 
orthogonal complement and denoted by a vertical bar $\mid$. Gra{\ss}mann 
defined it by an implicit relation
\begin{align}
[a \wedge \vert a] &= 1
\end{align}
where $\vert a$ is the Erg\"anzung and $[\ldots]$ is a volume form as 
studied in the case of Peano space. In fact, that is Peano's source
\cite{peano:1888a}.

The meet, alias regressive product or `eingewandtes Produkt', as opposed 
to the exterior (progressive) product, was defined by Gra{\ss}mann 
\cite{grassmann:1862a} as
\begin{align}
\vert ( a\vee b) &:= (\vert a) \wedge (\vert b)\, . 
\end{align}
This relation is still projective and does not use a metric but depends 
on a symmetric correlation. One should compare this Gra{\ss}mannian 
definition with that of Hestenes, Sobczyk
\cite{hestenes:sobczyk:1992a} and Hestenes, Ziegler 
\cite{hestenes:ziegler:1991a} where inner products are used. This 
route was taken also by Conradt \cite{conradt:2000a,conradt:2000c}.

Recall that the bracket $[.,\ldots,.] ~:~\otimes^n V\rightarrow\openk$
was essentially identical to a determinant of the matrix of the 
vector components of its entries
\begin{align}
\det( a_1,\ldots,a_n ) &= [a_1,\ldots,a_n] \, .
\end{align}
But following Chevalley \cite{chevalley:1997a}, the determinant can
be calculated along the following line. Let $a_i = A(e_i)$ be the images
of some $e_i$, which constitute a basis of $V$, i.e. 
$[e_1,\ldots,e_n]\not=0$. Let $A^{\wedge}$ be the graded extension of $A$ 
on $\bigwedge V$, as we have extended the scalar and coscalar products 
above. One computes
\begin{align}
A^{\wedge}(e_1\wedge\ldots\wedge e_n) &= a_1\wedge\ldots\wedge a_n =
\alpha\, e_1\wedge\ldots\wedge e_n\nn
\det(A) &= \alpha = [a_1,\ldots, a_n] \, .
\end{align}
Not even an orientation is needed, the determinant respects only a relative
orientation between two sequences here. But this will change if a particular 
basis is selected and orientation is established relatively to such a 
{\it right handed basis\/}. Now we find that the cointegral of a 
Gra{\ss}mann Hopf gebra projects onto the highest grade element, hence 
on the determinant of that linear operator, which maps a certain basis 
fulfilling $[e_1,\ldots, e_n]=1$, i.e. linearily ordered and oriented,
onto the input of the bracket. We {\it define\/} the bracket using the 
unique Gra{\ss}mann Hopf integral

\begin{align}
[a_1,\ldots,a_n]_{\mu} 
\quad&\cong\quad
\pspicture[0.5](0,0)(1.5,2.25)
\psset{linewidth=\pstlw,xunit=0.5,yunit=0.5,runit=0.5}
\psline(1.5,0)(1.5,1)
\psline(0,2.5)(0,3.5)
\psline(3,2.5)(3,3.5)
\psarc(1.5,2.5){1.5}{180}{360}
\psline[linestyle=dotted](1,3)(2,3)
\pscircle[linewidth=0.4pt,fillstyle=solid,fillcolor=white](1.5,1){0.2}
\pscircle[linewidth=0.4pt,fillstyle=solid,fillcolor=black](1.5,0){0.2}
\rput(0,4){$a_1\wedge$}
\rput(1.5,4){$\ldots$}
\rput(3,4){$\wedge a_n$}
\rput(2,0){$\mu$}
\endpspicture
\end{align}
The Hopf algebraic definition of the meet as given by Doubilet, Rota
and Stein \cite{doubilet:rota:stein:1974a} reads
\begin{align}
a \vee b &= (a_1\wedge\ldots\wedge a_r) \vee (b_1\wedge\ldots\wedge b_s) \nn
&= [b_{(1)},a] \, b_{(2)} =  a_{(1)}\, [b,a_{(2)}] \nn
&= \pm [a,b_{(1)}] \, b_{(2)} =  \pm a_{(1)}\, [a_{(2)},b] 
\end{align}
and contains the bracket. That the bracket is not 
foreign to the Gra{\ss}mann Hopf gebra was discussed in the previous section. 
However, we can now see that the bracket involves a disguised integral. In 
the following tangles we use for clarity the last line for the meet and 
compute modulo signs, which is allowed in homogeneous coordinates of 
projective geometry. We define hence the meet as
\begin{align}
\pspicture[0.5](0,0)(1,2)
\psset{linewidth=\pstlw,xunit=0.5,yunit=0.5,runit=0.5}
\psline(1,0)(1,1)
\psarc(1,2){1}{180}{360}
\psline(0,2)(0,4)
\psline(2,2)(2,4)
\pscircle[linewidth=0.4pt,fillstyle=solid,fillcolor=white](1,1){0.2}
\rput(1,1.75){$\vee$}
\endpspicture
\quad:=\quad
\pspicture[0.5](0,0)(2,2)
\psset{linewidth=\pstlw,xunit=0.5,yunit=0.5,runit=0.5}
\psline(0,4)(0,2)
\psline(3,4)(3,3)
\psarc(3,2){1}{0}{180}
\psarc(1,2){1}{180}{360}
\psline(1,1)(1,0.5)
\psline(4,2)(4,0)
\pscircle[linewidth=0.4pt,fillstyle=solid,fillcolor=white](3,3){0.2}
\pscircle[linewidth=0.4pt,fillstyle=solid,fillcolor=white](1,1){0.2}
\pscircle[linewidth=0.4pt,fillstyle=solid,fillcolor=black](1,0.5){0.2}
\rput(1,1.75){$\wedge$}
\rput(3,1.75){$\Delta_\wedge$}
\endpspicture
\quad=\quad
\pspicture[0.5](0,0)(2,2)
\psset{linewidth=\pstlw,xunit=0.5,yunit=0.5,runit=0.5}
\psline(4,4)(4,2)
\psline(1,4)(1,3)
\psarc(1,2){1}{0}{180}
\psarc(3,2){1}{180}{360}
\psline(3,1)(3,0.5)
\psline(0,2)(0,0)
\pscircle[linewidth=0.4pt,fillstyle=solid,fillcolor=white](1,3){0.2}
\pscircle[linewidth=0.4pt,fillstyle=solid,fillcolor=white](3,1){0.2}
\pscircle[linewidth=0.4pt,fillstyle=solid,fillcolor=black](3,0.5){0.2}
\rput(3,1.75){$\wedge$}
\rput(1,1.75){$\Delta_\wedge$}
\endpspicture
\end{align}
A short calculation shows that the meet $\vee$ is associative. Note that
the r.h.s. consists of Kuperberg ladder diagrams truncated by the 
cointegral. The second equality was proved by Doubilet, Rota and Stein
\cite{doubilet:rota:stein:1974a}. 
Since the Kuperberg ladder tangles are invertible, we can derive the 
relations
\begin{align}
\pspicture[0.5](0,0)(2.5,3)
\psset{linewidth=\pstlw,xunit=0.5,yunit=0.5,runit=0.5}
\psline(0,6)(0,4)
\psline(3,6)(3,5)
\psline(1,3)(1,2)
\psline(4,3)(4,4)
\psline(2,1)(2,0)
\psline(5,2)(5,0)
\psarc(3,4){1}{0}{180}
\psarc(4,2){1}{0}{180}
\psarc(1,4){1}{180}{360}
\psarc(2,2){1}{180}{360}
\pscircle[linewidth=0.4pt,fillstyle=solid,fillcolor=white](1,3){0.2}
\pscircle[linewidth=0.4pt,fillstyle=solid,fillcolor=white](2,1){0.2}
\pscircle[linewidth=0.4pt,fillstyle=solid,fillcolor=white](3,5){0.2}
\pscircle[linewidth=0.4pt,fillstyle=solid,fillcolor=white](4,3){0.2}
\pscircle[linewidth=0.4pt,fillstyle=solid,fillcolor=black](2,0){0.2}
\rput*(2,4){$S$}
\endpspicture
\quad=\quad
\pspicture[0.5](0,0)(1,3)
\psset{linewidth=\pstlw,xunit=0.5,yunit=0.5,runit=0.5}
\psline(0,6)(0,2)
\psline(2,6)(2,0)
\pscircle[linewidth=0.4pt,fillstyle=solid,fillcolor=black](0,2){0.2}
\endpspicture
\quad\text{~and~}\quad
\pspicture[0.5](0,0)(2.5,3)
\psset{linewidth=\pstlw,xunit=0.5,yunit=0.5,runit=0.5}
\psline(5,6)(5,4)
\psline(2,6)(2,5)
\psline(4,3)(4,2)
\psline(1,3)(1,4)
\psline(3,1)(3,0)
\psline(0,2)(0,0)
\psarc(2,4){1}{0}{180}
\psarc(1,2){1}{0}{180}
\psarc(4,4){1}{180}{360}
\psarc(3,2){1}{180}{360}
\pscircle[linewidth=0.4pt,fillstyle=solid,fillcolor=white](4,3){0.2}
\pscircle[linewidth=0.4pt,fillstyle=solid,fillcolor=white](3,1){0.2}
\pscircle[linewidth=0.4pt,fillstyle=solid,fillcolor=white](2,5){0.2}
\pscircle[linewidth=0.4pt,fillstyle=solid,fillcolor=white](1,3){0.2}
\pscircle[linewidth=0.4pt,fillstyle=solid,fillcolor=black](3,0){0.2}
\rput*(3,4){$S$}
\endpspicture
\quad=\quad
\pspicture[0.5](0,0)(1,3)
\psset{linewidth=\pstlw,xunit=0.5,yunit=0.5,runit=0.5}
\psline(2,6)(2,2)
\psline(0,6)(0,0)
\pscircle[linewidth=0.4pt,fillstyle=solid,fillcolor=black](2,2){0.2}
\endpspicture
\end{align}
Having the $2\rightarrow 0$ tangle for the bracket, it is natural to 
introduce the product co-product duality w.r.t. this cup tangle:
\begin{align}
&\pspicture[0.5](0,-0.5)(1.5,2.5)
\psset{linewidth=\pstlw,xunit=0.5,yunit=0.5,runit=0.5}
\psset{arrowsize=2pt 2,arrowinset=0.2}
\psline(3,5)(3,1)
\psline(2,5)(2,2)
\psline(0,5)(0,2)
\psline(2,0)(2,-1)
\psarc(2,1){1}{180}{360}
\psarc(1,2){1}{180}{360}
\pscircle[linewidth=0.4pt,fillstyle=solid,fillcolor=white](1,1){0.2}
\pscircle[linewidth=0.4pt,fillstyle=solid,fillcolor=white](2,0){0.2}
\pscircle[linewidth=0.4pt,fillstyle=solid,fillcolor=black](2,-1){0.2}
\rput(1,1.75){$\vee$}
\endpspicture
\quad=\quad
\pspicture[0.5](0,-0.5)(3,2.5)
\psset{linewidth=\pstlw,xunit=0.5,yunit=0.5,runit=0.5}
\psset{arrowsize=2pt 2,arrowinset=0.2}
\psline(5,5)(5,4)
\psline(2,5)(2,3)
\psline(0,5)(0,2)
\psline(6,3)(6,2)
\psline(2,0)(4,0)
\psline(4,0)(4,0)
\psline(3,2)(3,1)
\psline(3,0)(3,-1)
\psarc(5,3){1}{0}{180}
\psarc(3,3){1}{180}{360}
\psarc(2,2){2}{180}{270}
\psarc(4,2){2}{270}{360}
\pscircle[linewidth=0.4pt,fillstyle=solid,fillcolor=white](5,4){0.2}
\pscircle[linewidth=0.4pt,fillstyle=solid,fillcolor=white](3,2){0.2}
\pscircle[linewidth=0.4pt,fillstyle=solid,fillcolor=white](3,0){0.2}
\pscircle[linewidth=0.4pt,fillstyle=solid,fillcolor=black](3,1){0.2}
\pscircle[linewidth=0.4pt,fillstyle=solid,fillcolor=black](3,-1){0.2}
\rput(3,2.75){$\wedge$}
\rput(5,2.25){$\Delta_{\wedge}$}
\endpspicture
\quad=\quad
\pspicture[0.5](0,-0.5)(1.5,2.5)
\psset{linewidth=\pstlw,xunit=0.5,yunit=0.5,runit=0.5}
\psset{arrowsize=2pt 2,arrowinset=0.2}
\psline(0,5)(0,1)
\psline(1,5)(1,2)
\psline(3,5)(3,2)
\psline(1,0)(1,-1)
\psarc(1,1){1}{180}{360}
\psarc(2,2){1}{180}{360}
\pscircle[linewidth=0.4pt,fillstyle=solid,fillcolor=white](2,1){0.2}
\pscircle[linewidth=0.4pt,fillstyle=solid,fillcolor=white](1,0){0.2}
\pscircle[linewidth=0.4pt,fillstyle=solid,fillcolor=black](1,-1){0.2}
\rput(2,2,75){$\vee$}
\endpspicture
\nonumber \\[2ex]
&\mu((a\vee b) \wedge c) \,=\, \mu(a \wedge ( b\vee c)) \, ,           
\end{align}
which is a straigthening formula also derived by Rota et al. In this 
particular sense, $\vee$ is a self dual product.

We know already from the Gra{\ss}mann Hopf gebra that by categorial
duality we can derive a Gra{\ss}mann co-product from the exterior wedge
product $\wedge \rightarrow \Delta_\wedge$. Along the same lines,
i.e. using eval, we can introduce a Gra{\ss}mann co-product for the meet 
$\vee \rightarrow \Delta_\vee$. This notion depends on the unique
integral of the Gra{\ss}mann Hopf gebra
\begin{align}
\pspicture[0.5](0,0)(1,2)
\psset{linewidth=\pstlw,xunit=0.5,yunit=0.5,runit=0.5}
\psline(1,4)(1,3)
\psarc(1,2){1}{0}{180}
\psline(0,0)(0,2)
\psline(2,0)(2,2)
\pscircle[linewidth=0.4pt,fillstyle=solid,fillcolor=white](1,3){0.2}
\rput(1,2.25){$\Delta_\vee$}
\endpspicture
\quad:=\quad
\pspicture[0.5](0,0)(2,2)
\psset{linewidth=\pstlw,xunit=0.5,yunit=0.5,runit=0.5}
\psline(0,2)(0,0)
\psline(3,1)(3,0)
\psarc(3,2){1}{180}{360}
\psarc(1,2){1}{0}{180}
\psline(1,3)(1,3.5)
\psline(4,2)(4,4)
\pscircle[linewidth=0.4pt,fillstyle=solid,fillcolor=white](3,1){0.2}
\pscircle[linewidth=0.4pt,fillstyle=solid,fillcolor=white](1,3){0.2}
\pscircle[linewidth=0.4pt,fillstyle=solid,fillcolor=black](1,3.5){0.2}
\rput(3,2.25){$\wedge$}
\rput(1,2.25){$\Delta_\wedge$}
\endpspicture
\quad=\quad
\pspicture[0.5](0,0)(2,2)
\psset{linewidth=\pstlw,xunit=0.5,yunit=0.5,runit=0.5}
\psline(4,0)(4,2)
\psline(1,0)(1,1)
\psarc(1,2){1}{180}{360}
\psarc(3,2){1}{0}{180}
\psline(3,3)(3,3.5)
\psline(0,2)(0,4)
\pscircle[linewidth=0.4pt,fillstyle=solid,fillcolor=white](1,1){0.2}
\pscircle[linewidth=0.4pt,fillstyle=solid,fillcolor=white](3,3){0.2}
\pscircle[linewidth=0.4pt,fillstyle=solid,fillcolor=black](3,3.5){0.2}
\rput(1,2.25){$\wedge$}
\rput(3,2.25){$\Delta_\wedge$}
\endpspicture
\end{align}

The full symmetry of this structure was already noted by A. Lotze 1955 
\cite{lotze:1955a}, using Hopf gebras only implicitly in the combinatorics
of indices. Lotze showed that the exterior product derived from the meet
along the same lines as the meet itself was obtained, is again the join!
\begin{align}
\pspicture[0.5](0,0)(1,2)
\psset{linewidth=\pstlw,xunit=0.5,yunit=0.5,runit=0.5}
\psline(1,0)(1,1)
\psarc(1,2){1}{180}{360}
\psline(0,2)(0,4)
\psline(2,2)(2,4)
\pscircle[linewidth=0.4pt,fillstyle=solid,fillcolor=white](1,1){0.2}
\rput(1,1.75){$\wedge$}
\endpspicture
\quad:=\quad
\pspicture[0.5](0,0)(2,2)
\psset{linewidth=\pstlw,xunit=0.5,yunit=0.5,runit=0.5}
\psline(0,4)(0,2)
\psline(3,4)(3,3)
\psarc(3,2){1}{0}{180}
\psarc(1,2){1}{180}{360}
\psline(1,1)(1,0.5)
\psline(4,2)(4,0)
\pscircle[linewidth=0.4pt,fillstyle=solid,fillcolor=white](3,3){0.2}
\pscircle[linewidth=0.4pt,fillstyle=solid,fillcolor=white](1,1){0.2}
\pscircle[linewidth=0.4pt,fillstyle=solid,fillcolor=black](1,0.5){0.2}
\rput(1,1.75){$\vee$}
\rput(3,1.75){$\Delta_\vee$}
\endpspicture
\quad=\quad
\pspicture[0.5](0,0)(2,2)
\psset{linewidth=\pstlw,xunit=0.5,yunit=0.5,runit=0.5}
\psline(4,4)(4,2)
\psline(1,4)(1,3)
\psarc(1,2){1}{0}{180}
\psarc(3,2){1}{180}{360}
\psline(3,1)(3,0.5)
\psline(0,2)(0,0)
\pscircle[linewidth=0.4pt,fillstyle=solid,fillcolor=white](1,3){0.2}
\pscircle[linewidth=0.4pt,fillstyle=solid,fillcolor=white](3,1){0.2}
\pscircle[linewidth=0.4pt,fillstyle=solid,fillcolor=black](3,0.5){0.2}
\rput(3,1.75){$\vee$}
\rput(1,1.75){$\Delta_\vee$}
\endpspicture
\end{align}
Denoting the Erg\"anzung in modern notation by a star $\star$, the full
mathematical structure turns out to be the {\it Gra{\ss}mann-Cayley 
double Hopf gebra\/} or {\it fourfold algebra}
\begin{align}
GC(\bigwedge V, \bigvee V^*, \wedge, \Delta_\wedge, \vee,\Delta_\vee,\star)
\end{align}
which possesses also units, counits, antipodes $S^\wedge$, $S^\vee$, 
integrals and cointegrals. Furthermore this structure is the vector 
space analog of the Boolean algebra of sets, the algebra of 
logical inference.

One can check by easy computations that the cointegral is the unit
of the meet while the integral is the unit of the meet co-product. The 
integrals and cointegrals obtained from the meet play the same role
for the wedge again.

Unfortunately we have no further opportunity to discuss the geometry
behind this interesting highly symmetric algebraic structure in this 
work.

\section{Crossings}

We examine some properties of the crossing derived from the antipode.
This will be done for our $\dim V=2$ example. This is not a theory of 
the crossing, but it gives valuable hints how the crossing behaves.

\noindent
{\bf Example:} (continued)
Let $\dim V=2$ and $B,C$ be the arbitrary scalar and co-scalar products
of the Clifford biconvolution as in the previously discussed cases.

\begin{thrm}
The Clifford biconvolution $\Cl(B,C)$ is commutative as an algebra
\begin{align}
m\circ \hat{\tau} &= m 
\end{align}
if the scalar product is identically zero, $B\equiv 0$, and the 
co-scalar product is symmetric, $C^T=C$.\\
The Clifford bi-convolution $\Cl(B,C)$ is commutative as a co-gebra
\begin{align}
\hat{\tau}\circ \Delta &= \Delta 
\end{align}
if the co-scalar product is identically zero, $C\equiv 0$, and the 
scalarproduct is symmetric, $B^T=B$.
\end{thrm}
\noindent
{\bf Proof:} by direct computation using CLIFFORD / BIGEBRA.

We will furthermore check if we can derive a non-knottish skein relation
for the crossing, as various such relations have been suggested by 
Oziewicz. Let $\dim V=2$ and $\Cl(B,C)$ be an arbitrary Clifford 
biconvolution. Are there solutions to the following skein relation, 
where the Kuperberg ladders are involved?
\begin{align}
\pspicture[0.5](0,0)(1,2.5)
\psset{linewidth=\pstlw,xunit=0.5,yunit=0.5,runit=0.5}
\psline(0,5)(0,4)
\psline(2,5)(2,4)
\psline(0,1)(0,0)
\psline(2,1)(2,0)
\psarc(1,1){1}{0}{180}
\psarc(1,4){1}{180}{360}
\psline(1,3)(1,2)
\pscircle[linewidth=0.4pt,fillstyle=solid,fillcolor=white](1,3){0.2}
\pscircle[linewidth=0.4pt,fillstyle=solid,fillcolor=white](1,2){0.2}
\endpspicture
\quad=\quad
q\,\,
\pspicture[0.5](0,0)(2,2)
\psset{linewidth=\pstlw,xunit=0.5,yunit=0.5,runit=0.5}
\psline(0,2)(0,0)
\psline(3,1)(3,0)
\psarc(3,2){1}{180}{360}
\psarc(1,2){1}{0}{180}
\psline(1,3)(1,4)
\psline(4,2)(4,4)
\pscircle[linewidth=0.4pt,fillstyle=solid,fillcolor=white](3,1){0.2}
\pscircle[linewidth=0.4pt,fillstyle=solid,fillcolor=white](1,3){0.2}
\endpspicture
\quad+ t\,\,
\pspicture[0.5](0,0)(2,2)
\psset{linewidth=\pstlw,xunit=0.5,yunit=0.5,runit=0.5}
\psline(4,0)(4,2)
\psline(1,0)(1,1)
\psarc(1,2){1}{180}{360}
\psarc(3,2){1}{0}{180}
\psline(3,3)(3,4)
\psline(0,2)(0,4)
\pscircle[linewidth=0.4pt,fillstyle=solid,fillcolor=white](1,1){0.2}
\pscircle[linewidth=0.4pt,fillstyle=solid,fillcolor=white](3,3){0.2}
\endpspicture
\end{align}
As a result we obtain pairs of scalar and coscalar products $(B,C)$
such that the skein relation holds true for $q=0$ or $t=0$. However,
we find also solutions of the form
\begin{align}
q&=1-t 
&& B=\left[\begin{array}{cc}
\frac{1}{u} & -\frac{z}{uv} \\
0 & \frac{1}{v}
\end{array}\right]
&& C=\left[\begin{array}{cc}
u & z \\ 
0 & v
\end{array}\right]
\nn
q&=1-t 
&& B=\left[\begin{array}{cc}
-\frac{vc}{w} & \frac{vuc+w}{w^2}\\
c &  -\frac{uc}{w}
\end{array}\right]
&& C=\left[\begin{array}{cc}
u & \frac{vuc+w}{wc} \\
w & v
\end{array}\right] \, .
\end{align}
It would be an interesting task to investigate which sort of skein
relations can appear in this way and in which way such skein relations
can be used in physics. However we will follow another route here.

\section{Convolutive unipotents}

An unipotent element $x$ of an unital algebra fulfils $x^2=\Id$.
Idempotent elements which square to themselves $f^2=f$, are related to 
unipotents. That is, from every non-trivial unipotent element $x$ one 
can construct two idempotent elements which are mutually orthogonal
\begin{align}
f_{\pm} &= \frac{1}{2}(\Id\pm x) &&\Rightarrow && f^2_{\pm}=f_{\pm},
\qquad f_+ f_- = f_- f_+ \, .
\end{align}
The knowledge of commuting unipotents is thus closely related to
that of (primitive) idempotents. But primitive idempotents generate
minimal left and right ideals which carry irreducible faithful
representations of algebras, if this algebra is fully reducible.

In a Gra{\ss}mann algebra all elements but $\Id$ are nilpotent, i.e. 
$x^n=0$ for some $n\in \openZ$. Hence $\Id$ and $0$ are the only 
unipotent elements which are at the same time idempotents. That is,
$\bigwedge V$ is faithfully and irreducibly represented on $\bigwedge V$
seen as $\openk$-module.

Cliffordization changes this fact drastically. The non-local Clifford 
product generates, depending on the base field and the involved bilinear 
form, a certain number $r$ of commuting unipotents. From these $r$ commuting
unipotents, $2^r$ primitive idempotents are constructed. The number $r$
is the Radon-Hurwitz number found by Hurwitz and Radon during their 
studies of the composition of quadratic forms 
\cite{hurwitz:1923a,radon:1923a}. Benno Eckmann showed that there is
a group theoretical root of this number \cite{eckmann:1943a} and
provided a short proof. Later, Hasler Whitney showed, that 
the Radon-Hurwitz number is related to the number of independent 
vector fields on spheres \cite{whitney:1935a}. That is, the Radon-Hurwitz
number is directly related to topological properties of the underlying 
group manifolds. This provides a relation to the topological relation
of Hopf algebras and the process of cliffordization of Hopf gebras.

Since we saw that the convolution establishes a group like structure
via the convolution and finally the Hopf gebra. We can now start to
study the endormorphsims of $\bigwedge V$ forming a convolution Hopf 
algebra and ask if there are convolutive unipotents in Hopf gebras.
\begin{thrm}
In a Gra{\ss}mann Hopf gebra the convolutive unit $u$ is up to the 
sign the only unipotent element.
\end{thrm}

Hence the unit and zero $(u,0)$ are the only idempotent endomorphisms
under the Gra{\ss}mann convolution product. That is, we expect the totality 
of all graded endomorphisms $\End \bigwedge V$ to form an irreducible
representation space for the convolution as was the case for the 
Gra{\ss}mann algebra itself. Once more cliffordization changes the situation.

\noindent
{\bf Example:} (continued) We consider the $\dim V=2$ Clifford 
bi-convolution $\Cl(B,C)$.
\begin{thrm}
There are more than 90 non-trivial $T\not=0$ convolution unipotent
solutions of the equation $T\star T = u$ in the Clifford biconvolution 
$\Cl(B,C)$ of $\dim V=2$. [Not all of them being independent.] 
In particular, among the solutions there are singular endomorphisms 
$\det T=0$ and non-singular endomorphisms $\det T\not=0$.
\end{thrm}
\noindent
{\bf Proof:} by direct computation using CLIFFORD / BIGEBRA.

We expect thus a non-trivial and highly interesting representation theory 
of Clifford biconvolutions and Clifford Hopf gebras. However, we cannot
enter this topic here.

\subsection{Convolutive 'adjoint'}

The notion of an adjoint operator belongs to the theory of inner product
spaces. One defines the adjoint operator $A^*$ to be the operator shifted 
into the left slot of the inner product
\begin{align}
\langle x A^* \mid y\rangle &=
\langle x \mid A y\rangle \, .
\end{align}
A certain type of 'adjointness' is the inverse which is defined w.r.t. 
a product, not a bilinear form. In the case of a group one has
\begin{align}
\label{eqn:grp-ID}
g \star h &= \Id \star g^{-1}h = gh^{-1}\star \Id \, .
\end{align}
This gives an idea to introduce a {\it convolutive adjoint\/}, 
resp. inverse, along the lines of the group inversion. However, in the
case of groups the identity of Eqn. \ref{eqn:grp-ID} is trivial,
since $\star$ is the group multiplication identical to juxtaposition.
But in the case of a Hopf gebra or a convolution algebra, the convolution 
product is {\it not\/} identical with the repeated application or 
composition of endomorphisms:
\begin{align}
\text{group}&&& g\star g^{-1} = \Id = gg^{-1} \nn
\text{conv. alg.}&&& g\star g^{-1} = \Id \star g^\star g^{-1}
=  g^{-1} g^\star\star\Id \, .
\end{align}

\subsection{A square root of the antipode}

The convolutive inverse is mediated by the antipode. However, we can use
the above found unipotents which are related to the antipode itself
by convolutive adjointness
\begin{align}
 T \star T &= u \nn
\Rightarrow\qquad \Id \star T^\star T &= T\, T^\star \star \Id = u \,. 
\end{align}
The second line is the defining relation for the antipode. Hence we 
find
\begin{align}
T^\star T &= S = TT^\star
\end{align}
or equivalently
\begin{align}
S(T)T &= S = TS(T) \, .
\end{align}
This mechanism is related to M\"obius inversion of polynoms, which we 
unfortunately cannot examine here. However, we can report that using
BIGEBRA we have been able to find a great variety of operators $T$,
which proves their existence in some special cases. Moreover,
we find invertible and singular $T$s which will induce a rich
representation theory.

In a certain sense, the operator $T$ is like a square root of $S$ and 
could be called {\it spinorial antipode\/} since it parallels by analogy
the ${\bf\rm spin} \rightarrow SO$ or ${\bf\rm pin}\rightarrow O$ covering.
It is not yet clear, if such operators are connected to coverings
of the topological spaces behind the Hopf algebra.

\subsection{Symmetrized product co-procduct tangle} 

As an application, we can symmetrize the Kuperberg ladder tangles.
Using the convolution unipotent $T$, we define
\begin{align}
\pspicture[0.5](0,0)(2,2)
\psset{linewidth=\pstlw,xunit=0.5,yunit=0.5,runit=0.5}
\psline(0,4)(0,2)
\psline(3,4)(3,3)
\psarc(3,2){1}{0}{180}
\psarc(1,2){1}{180}{360}
\psline(1,1)(1,0)
\psline(4,2)(4,0)
\pscircle[linewidth=0.4pt,fillstyle=solid,fillcolor=white](3,3){0.2}
\pscircle[linewidth=0.4pt,fillstyle=solid,fillcolor=white](1,1){0.2}
\rput*(2,2){$T$}
\endpspicture
\quad\Rightarrow\qquad
\pspicture[0.5](0,0)(2.5,3)
\psset{linewidth=\pstlw,xunit=0.5,yunit=0.5,runit=0.5}
\psline(0,6)(0,4)
\psline(3,6)(3,5)
\psline(1,3)(1,2)
\psline(4,3)(4,4)
\psline(2,1)(2,0)
\psline(5,2)(5,0)
\psarc(3,4){1}{0}{180}
\psarc(4,2){1}{0}{180}
\psarc(1,4){1}{180}{360}
\psarc(2,2){1}{180}{360}
\pscircle[linewidth=0.4pt,fillstyle=solid,fillcolor=white](1,3){0.2}
\pscircle[linewidth=0.4pt,fillstyle=solid,fillcolor=white](2,1){0.2}
\pscircle[linewidth=0.4pt,fillstyle=solid,fillcolor=white](3,5){0.2}
\pscircle[linewidth=0.4pt,fillstyle=solid,fillcolor=white](4,3){0.2}
\rput*(2,4){$T$}
\rput*(3,2){$T$}
\endpspicture
\quad=\quad
\pspicture[0.5](0,0)(1,3)
\psset{linewidth=\pstlw,xunit=0.5,yunit=0.5,runit=0.5}
\psline(0,6)(0,0)
\psline(2,6)(2,0)
\endpspicture
\nonumber \\[2ex]
\pspicture[0.5](0,0)(2,2)
\psset{linewidth=\pstlw,xunit=0.5,yunit=0.5,runit=0.5}
\psline(4,4)(4,2)
\psline(1,4)(1,3)
\psarc(1,2){1}{0}{180}
\psarc(3,2){1}{180}{360}
\psline(3,1)(3,0)
\psline(0,2)(0,0)
\pscircle[linewidth=0.4pt,fillstyle=solid,fillcolor=white](1,3){0.2}
\pscircle[linewidth=0.4pt,fillstyle=solid,fillcolor=white](3,1){0.2}
\rput*(2,2){$T$}
\endpspicture
\quad\Rightarrow\qquad
\pspicture[0.5](0,0)(2.5,3)
\psset{linewidth=\pstlw,xunit=0.5,yunit=0.5,runit=0.5}
\psline(5,6)(5,4)
\psline(2,6)(2,5)
\psline(4,3)(4,2)
\psline(1,3)(1,4)
\psline(3,1)(3,0)
\psline(0,2)(0,0)
\psarc(2,4){1}{0}{180}
\psarc(1,2){1}{0}{180}
\psarc(4,4){1}{180}{360}
\psarc(3,2){1}{180}{360}
\pscircle[linewidth=0.4pt,fillstyle=solid,fillcolor=white](4,3){0.2}
\pscircle[linewidth=0.4pt,fillstyle=solid,fillcolor=white](3,1){0.2}
\pscircle[linewidth=0.4pt,fillstyle=solid,fillcolor=white](2,5){0.2}
\pscircle[linewidth=0.4pt,fillstyle=solid,fillcolor=white](1,3){0.2}
\rput*(2,2){$T$}
\rput*(3,4){$T$}
\endpspicture
\quad=\quad
\pspicture[0.5](0,0)(1,3)
\psset{linewidth=\pstlw,xunit=0.5,yunit=0.5,runit=0.5}
\psline(2,6)(2,0)
\psline(0,6)(0,0)
\endpspicture
\end{align} 
That is, these tangles are also unipotents. However, our choice implies 
further alterations in the theory. If we stay with the crossing of
Oziewicz, but change $S\rightarrow T$, we get
\begin{align}
\pspicture[0.5](0,0)(1,1)
\psset{linewidth=\pstlw,xunit=0.5,yunit=0.5,runit=0.5}
\psbezier(2,2)(2,1)(0,1)(0,0)
\psbezier[border=4pt,bordercolor=white](0,2)(0,1)(2,1)(2,0)
\endpspicture
\quad&:=\quad
\pspicture[0.5](0,0)(3,4)
\psset{linewidth=\pstlw,xunit=0.5,yunit=0.5,runit=0.5}
\psline{-}(1,8)(1,7)
\psline{-}(5,8)(5,7)
\psline{-}(1,1)(1,0)
\psline{-}(5,1)(5,0)
\psline{-}(0,6)(0,2)
\psline{-}(6,6)(6,2)
\psline{-}(3,5)(3,3)
\psarc(1,6){1.0}{0}{180}
\psarc(5,6){1.0}{0}{180}
\psarc(3,6){1.0}{180}{360}
\psarc(3,2){1.0}{0}{180}
\psarc(1,2){1.0}{180}{360}
\psarc(5,2){1.0}{180}{360}
\pscircle[linewidth=0.4pt,fillstyle=solid,fillcolor=white](1,7){0.2}
\pscircle[linewidth=0.4pt,fillstyle=solid,fillcolor=white](5,7){0.2}
\pscircle[linewidth=0.4pt,fillstyle=solid,fillcolor=white](3,5){0.2}
\pscircle[linewidth=0.4pt,fillstyle=solid,fillcolor=white](3,3){0.2}
\pscircle[linewidth=0.4pt,fillstyle=solid,fillcolor=white](1,1){0.2}
\pscircle[linewidth=0.4pt,fillstyle=solid,fillcolor=white](5,1){0.2}
\rput*(0,4){T}
\rput*(6,4){T}
\endpspicture
\end{align} 
which leads to a deformation of the crossed products
\begin{align}
\pspicture[0.5](0,0)(1,3)
\psset{linewidth=\pstlw,xunit=0.5,yunit=0.5,runit=0.5}
\psline(0,6)(0,5)
\psline(2,6)(2,5)
\psarc(1,5){1.0}{180}{360}
\psline(1,4)(1,2)
\psarc(1,1){1.0}{0}{180}
\psline(0,1)(0,0)
\psline(2,1)(2,0)
\pscircle[linewidth=0.4pt,fillstyle=solid,fillcolor=white](1,4){0.2}
\pscircle[linewidth=0.4pt,fillstyle=solid,fillcolor=white](1,2){0.2}
\endpspicture
\quad=\quad
\pspicture[0.5](0,0)(3,3)
\psset{linewidth=\pstlw,xunit=0.5,yunit=0.5,runit=0.5}
\psline(1,6)(1,5)
\psline(5,6)(5,5)
\psarc(1,4){1.0}{0}{180}
\psarc(5,4){1.0}{0}{180}
\psline(0,4)(0,2)
\psline(6,4)(6,2)
\psline{c-c}(4,4)(2,2)
\psline[border=4pt]{c-c}(2,4)(4,2)
\psarc(1,2){1.0}{180}{360}
\psarc(5,2){1.0}{180}{360}
\psline(1,1)(1,0)
\psline(5,1)(5,0)
\pscircle[linewidth=0.4pt,fillstyle=solid,fillcolor=white](1,1){0.2}
\pscircle[linewidth=0.4pt,fillstyle=solid,fillcolor=white](1,5){0.2}
\pscircle[linewidth=0.4pt,fillstyle=solid,fillcolor=white](5,1){0.2}
\pscircle[linewidth=0.4pt,fillstyle=solid,fillcolor=white](5,5){0.2}
\rput*(0,3){$T$}
\rput*(6,3){$T$}
\endpspicture
\end{align} 
If this is not asserted, we can no longer prove from the crossing
that $\Delta$ and $m$ are algebra and cogebra morphisms, which is however
a main feature of the Hopf gebra.

A detailed study of the representation theory of convolution algebras
will provide further information in which way such a generalization
has meaningfully to be developed.

%% file: gen_cli_prod.tex
\chapter{Generalized cliffordization}

Rota and Stein developed indeed a much more general concept of 
cliffordization \cite{rota:stein:1994a} as we have till now used. 
They used a product like mapping 
$\,\&r\, : \bigwedge V \times \bigwedge V \rightarrow \bigwedge V$ with
a non-scalar target as the `cup' tangle in the cliffordization. Such a 
tangle is no longer a `sausage', but has a third line which is internal 
and downwards. However, it is not clear why a product should be deformed 
by another {\it a priori} given product. However, quantum mechanics uses 
complex valued (anti)commutators which are maps 
$[.,.] : V \times V \rightarrow \openk$. If the target is extended to 
$\bigwedge V$ then the (anti)commutator products are no longer complex 
valued but operator valued. Since the (anti)commutator algebra
emerges itself from a cliffordization, as demonstrated above, this would
imply a cliffordization of a cliffordized product, hence a cliffordization 
of second order. In terms of tangles a sausage tangle inside a second 
more general cliffordization. While this is interesting in its own right,
we will generalize cliffordization to the most general case where the
cup tangle is scalar valued, i.e. where the tangle remains to be a 
sausage. It will turn out, using a result of Brouder \cite{brouder:2002a},
who used Hopf algebraic formulas of Pinter \cite{pinter:2000a,pinter:2001a}, 
that a scalar valued bilinear form may be already sufficient to be able to 
introduce renormalized time-ordered operator products and correlation 
functions as required by Epstein-Glaser renormalization 
\cite{epstein:glaser:1973a}. This mechanism introduces the renormalization
parameters via a scalar valued $Z$-pairing. While this is a special bilinear
form, i.e. cup tangle, we consider in this chapter the general case. 
This will allow us to derive some of the defining relations of the 
Z-pairing of renormalization by reasonable assumptions about the
resulting cliffordized product.

\section{Linear forms on $\bigwedge V \times \bigwedge V$}

In the preceeding section we were interested in generalizing Clifford 
algebras of a quadratic form $\Cl(V,g)$ to quantum Clifford algebras of 
an arbitrary bilinear form $\Cl(V,B)$. The form 
\begin{align}
B : V \times V \rightarrow \openk
\end{align}
is defined on $V$ and has $n^2$ independent parameters. This form has to be
extended to a bilinear form
\begin{align}
B^{\wedge} : \bigwedge V \times \bigwedge V \rightarrow \openk 
\end{align}
which obviously has also only $n^2$ parameters. However, a general 
bilinear form $\BF : \bigwedge V\times\bigwedge V\rightarrow\openk$ 
has $2^n\times 2^n= 4^n$ parameters. The aim of this chapter is
to investigate what kind of restrictions follows for $\BF$ if we
assert some properties to the cliffordized product based on this
bilinear form.

First of all, we recall what kind of bilinear form we have used
till now. The bilinear form $B^{\wedge}$, which we defined by exponentiation
or as a Laplace pairing, is graded in the following sense (on two extensors 
and extended by bilinearity):
\begin{align}
B^{\wedge}(x_1\wedge\ldots\wedge x_r,
y_1\wedge\ldots y_s) &=\left\{
\begin{array}{cl}
B(\Id,\Id) & \text{if $r=0=s$} \\
B(x,y)     & \text{if $r=1=s$}\\
(-)^{\frac{r(r-1)}{2}} \det([B(x_i,y_j)]) & \text{if $r=s$} \\
0          & \text{otherwise.}
\end{array}\right.
\end{align}
Note that we have used a different indexing in the first factor as
previously, which results in the prefactor in front of the determinants
due to reversion of factors. In a standard basis this yields a matrix 
representation where $B^{\wedge}$ decomposes into a block structure
\begin{align}
[B^{\wedge}] &=
\left[\begin{array}{ccccc}
1 & 0 & \ldots & \ldots & \ldots\\
0 & B & 0 & \ldots & \ldots\\
\vdots & 0 & B^2 & 0 & \ldots\\
\vdots & \vdots & 0 & \ddots & \ldots\\
\vdots & \vdots & \vdots & \vdots & \ddots
\end{array}\right]
\end{align}
where $B^r : \bigwedge^r V \times \bigwedge^r V \rightarrow \openk$ is
an ${r \choose n}\times{r\choose n}$-matrix and a function of $B$. 
The grading enforces the off-block entries to be zero. One sees easily 
that the $B^r$'s are built from minors of the original bilinear form $B$.

\noindent
{\bf Example:} Let $\dim V = 3$ and $B$ be represented in an arbitrary 
basis $\{ e_i \}$ as a $3\times 3$-matrix. $B^0$ is defined to have the 
value $1$. This requirement will allow the unit of an augmented
connected Hopf gebra to stay to be the unit w.r.t. the cliffordized
product, i.e. $\Id \,\&\!c\, \Id = B(\Id,\Id)\Id = \Id$. $B^1$ is 
obviously identical to $B$. An entry of $B^2$ is given as
\begin{align}
B^2(e_i\wedge e_j , e_k\wedge e_l) &= (-)^{2(2-1)/2}
(B(e_i,e_k)B(e_i,e_l)-B(e_j,e_l)B(e_i,e_k))
\end{align} 
and all of them make up a $3\times 3$-matrix, while 
$B^3= -det([B(e_i,e_j)])$ is the determinant of $B$. Taking the trace
of $B^\wedge$, one obtains the invariants of the bilinear form $B$,
which could be reformulated using the eigenvalues $\lambda_i$ of $B$
as $\{1,\sum \lambda_i, \sum_{i<j} \lambda_i\lambda_j, 
\sum_{i<j<k} \lambda_i\lambda_j\lambda_k = \prod \lambda_i\}$.

A first generalization would break up this correlation of the ${
n \choose r}\times{r\choose n}$-block-matrices arising from $B$. While $B$ 
has $n^2$ parameters, such a general graded bilinear form would come up 
with $\sum_r {n \choose r} = 2^n$ parameters. If we would introduce a 
grade operator (particle number operator) such maps would commute.

However, it is at this stage more convenient to introduce an ungraded 
bilinear form $\BF : \bigwedge V \times \bigwedge V \rightarrow \openk$
without any restrictions having $4^n$ free parameters. Alternatively this
map can be regarded as a {\it linear form\/} on the space $\bigwedge V \times 
\bigwedge V$. The tangle does not change at all and remains to be a cup 
tangle. However, the properties of the product obtained by cliffordization,
and hence the (Hopf) algebra built from it change {\it dramatically\/}!
In fact this generalization will be sufficient to allow us to incorporate
renormalization directly into a cliffordization process. While we showed 
in \cite{fauser:2001b} that time- and normal-ordered correlation
functions and operator products are related by a cliffordization process,
C. Brouder \cite{brouder:2002a} noted that Epstein-Glaser renormalization 
\cite{epstein:glaser:1973a} may be incorporated into this process too. It 
was the achievement of Pinter \cite{pinter:2000a,pinter:2001a} to prove the 
equivalence of Epstein-Glaser renormalization, which resides in position
space, to the Connes-Kreimer renormalization 
\cite{kreimer:1998a,kreimer:2000a,connes:kreimer:2000a,connes:kreimer:2001a}
which is equivalent to the BPHZ renormalization in momentum space and
the forest formulas there. 

We define the {\it generalized Clifford product\/} 
\begin{align}
\pspicture[0.5](0,0)(1,3.5)
\psset{linewidth=\pstlw,xunit=0.5,yunit=0.5,runit=0.5}
\psset{arrowsize=2pt 2,arrowinset=0.2}
\psline(0,7)(0,4)
\psline(2,7)(2,4)
\psarc(1,4){1}{180}{360}
\psline(1,3)(1,0)
\pscircle[linewidth=0.4pt,fillstyle=solid,fillcolor=white](1,3){0.2}
\rput(1,3.75){$\&r$}
\endpspicture
\quad&:=\quad
\pspicture[0.5](0,0)(2.5,3.5)
\psset{linewidth=\pstlw,xunit=0.5,yunit=0.5,runit=0.5}
\psset{arrowsize=2pt 2,arrowinset=0.2}
\psline(1,7)(1,6)
\psline(5,7)(5,6)
\psarc(1,5){1}{0}{180}
\psarc(5,5){1}{0}{180}
\psarc(3,5){1}{180}{360}
\psarc(3,5){3}{180}{360}
\psline(3,2)(3,0)
\pscircle[linewidth=0.4pt,fillstyle=solid,fillcolor=black](3,4){0.2}
\pscircle[linewidth=0.4pt,fillstyle=solid,fillcolor=white](3,2){0.2}
\pscircle[linewidth=0.4pt,fillstyle=solid,fillcolor=white](1,6){0.2}
\pscircle[linewidth=0.4pt,fillstyle=solid,fillcolor=white](5,6){0.2}
\rput(1,5.25){$\Delta$}
\rput(5,5.25){$\Delta$}
\rput(3,4.75){$\BF$}
\rput(3,2.75){$\wedge$}
\endpspicture
\end{align}
as in Rota-Stein \cite{rota:stein:1994a}, but now for a general cup
tangle. Algebraically this is equivalent to
\begin{align}
x \,\&r\, y &= \BF(x_{(2)},y_{(1)}) x_{(1)} \wedge y_{(2)}.
\end{align}
This product is called a `generalized' Clifford product, since it leads
to algebras which are Clifford like, but different to classical
Clifford algebras. We agree to call this Clifford algebras also quantum 
Clifford algebras (QCA).

\section{Properties of generalized Clifford products}

We have unfortunately not the opportunity to develop a theory of
generalized Clifford products, so we concentrate on some essential
properties which we want to assert on the product emerging from 
cliffordization to be able to utilize it in quantum field theory. 
Such properties of the product will have a direct impact on the possible 
form of the bilinear form $\BF$. Our original wedge products were 
associative and unital. We will use these wedges to model operator 
products in QFT. Since we are interested especially in renormalized 
time-ordered operator products inside of Wick monomials and correlation 
functions, we will study associative unital generalized Clifford products. 
These restrictions enforced on the generalized Clifford product $\,\&r\,$ 
will allow us to derive the assertions made on the renormalization 
parameters and the $Z$-grading by Brouder as a direct consequence.

\subsection{Units for generalized Clifford products}

We can firstly ask under which conditions the unit of the Gra{\ss}mann 
algebra remains to be the unit of the cliffordized product $\,\&r\,$. 
Hence we have the condition
\begin{align}
\pspicture[0.5](0,0)(2.5,3.5)
\psset{linewidth=\pstlw,xunit=0.5,yunit=0.5,runit=0.5}
\psset{arrowsize=2pt 2,arrowinset=0.2}
\psline(1,6.5)(1,6)
\psline(5,7)(5,6)
\psarc(1,5){1}{0}{180}
\psarc(5,5){1}{0}{180}
\psarc(3,5){1}{180}{360}
\psarc(3,5){3}{180}{360}
\psline(3,2)(3,0)
\pscircle[linewidth=0.4pt,fillstyle=solid,fillcolor=black](3,4){0.2}
\pscircle[linewidth=0.4pt,fillstyle=solid,fillcolor=white](3,2){0.2}
\pscircle[linewidth=0.4pt,fillstyle=solid,fillcolor=white](1,6){0.2}
\pscircle[linewidth=0.4pt,fillstyle=solid,fillcolor=white](5,6){0.2}
\pscircle[linewidth=0.4pt,fillstyle=solid,fillcolor=white](1,6.5){0.2}
\rput(3,4.75){$\BF$}
\endpspicture
\quad\quad=\quad
\pspicture[0.5](0,0)(1.0,3.5)
\psset{linewidth=\pstlw,xunit=0.5,yunit=0.5,runit=0.5}
\psline(1,7)(1,0)
\endpspicture
\quad=\quad
\pspicture[0.5](0,0)(2.5,3.5)
\psset{linewidth=\pstlw,xunit=0.5,yunit=0.5,runit=0.5}
\psset{arrowsize=2pt 2,arrowinset=0.2}
\psline(1,7)(1,6)
\psline(5,6.5)(5,6)
\psarc(1,5){1}{0}{180}
\psarc(5,5){1}{0}{180}
\psarc(3,5){1}{180}{360}
\psarc(3,5){3}{180}{360}
\psline(3,2)(3,0)
\pscircle[linewidth=0.4pt,fillstyle=solid,fillcolor=black](3,4){0.2}
\pscircle[linewidth=0.4pt,fillstyle=solid,fillcolor=white](3,2){0.2}
\pscircle[linewidth=0.4pt,fillstyle=solid,fillcolor=white](1,6){0.2}
\pscircle[linewidth=0.4pt,fillstyle=solid,fillcolor=white](5,6){0.2}
\pscircle[linewidth=0.4pt,fillstyle=solid,fillcolor=white](5,6.5){0.2}
\rput(3,4.75){$\BF$}
\endpspicture
\end{align}
Because the co-product of the Gra{\ss}mann Hopf gebra is connected and 
augmented, the unit is a co-algebra homomorphism as discussed above.
In formulas: $\Delta\circ \eta = \eta \otimes\eta$. The same holds true
for the counit. The condition on $\BF$ which asserts that $\eta$ becomes 
the unit w.r.t. $\,\&r\,$ reads then
\begin{align}
\pspicture[0.5](0,0)(1,1)
\psset{linewidth=\pstlw,xunit=0.5,yunit=0.5,runit=0.5}
\psarc(1,2){1}{180}{360}
\pscircle[linewidth=0.4pt,fillstyle=solid,fillcolor=white](0,2){0.2}
\pscircle[linewidth=0.4pt,fillstyle=solid,fillcolor=black](1,1){0.2}
\rput(1,1.75){$\BF$}
\endpspicture
\quad=\quad
\pspicture[0.5](0,0)(1,1)
\psset{linewidth=\pstlw,xunit=0.5,yunit=0.5,runit=0.5}
\psline(1,2)(1,1)
\pscircle[linewidth=0.4pt,fillstyle=solid,fillcolor=white](1,1){0.2}
\endpspicture
\quad=\quad
\pspicture[0.5](0,0)(1,1)
\psset{linewidth=\pstlw,xunit=0.5,yunit=0.5,runit=0.5}
\psarc(1,2){1}{180}{360}
\pscircle[linewidth=0.4pt,fillstyle=solid,fillcolor=white](2,2){0.2}
\pscircle[linewidth=0.4pt,fillstyle=solid,fillcolor=black](1,1){0.2}
\rput(1,1.75){$\BF$}
\endpspicture
\end{align}
or in formulas
\begin{align}
\BF(\eta,X) &= \epsilon(X) = \BF(X,\eta)\quad\quad \forall X\in \bigwedge V.
\end{align} 
Using a matrix representation, this implies the following block
structure for $\BF$
\begin{align}
[\BF] &=
\left[\begin{array}{ccccc}
1 & 0 & \ldots & 0 & \ldots\\
0 & B_{1,1} & \ldots & B_{1,n} & \ldots\\
\vdots & \vdots & \ddots & \ldots & \ldots\\
0 & B_{n,1} & \vdots & B_{n,n} & \ldots\\
\vdots & \vdots & \vdots & \vdots & \ddots
\end{array}\right]
\end{align}
which is a modest restriction. The renormalization scheme of 
Brouder-Epstein-Glaser (BEG) fulfils exactly this requirement.

However, it is an incidence and not an automatism that the Gra{\ss}mann 
unit $\Id$ remains to be the algebra unit under the generalized Clifford 
product. An algebra homomorphism maps only unit onto unit and has not
to have the unit as an invariant element. In fact very general elements 
of $\bigwedge V$ can be made to be a left/right unit w.r.t. the new 
Clifford product $\,\&r\,$ by a suitable choice of the bilinear form 
$\BF$. 

Let $GB=\{\Id, e_1,\ldots, e_n, e_1\wedge e_2, \ldots \}$ be a canonical
Gra{\ss}mann basis of $\bigwedge V$. An element
\begin{align}
X &= X_0\Id + \sum_i X_i e_i + \sum_{i<j} X_{ij} e_i \wedge e_j + 
\ldots X_{1\ldots n}e_1 \wedge \ldots \wedge e_n
\end{align}
where some or all $X_{i_r,\ldots,i_s}$ are non-zero is a left
unit if the corresponding rows of $\BF$ fulfil some constraints.

\noindent
{\bf Example:} We consider a 3-dimensional space $V$, $\dim \bigwedge V 
= 2^3=8$ and obtain using BIGEBRA the following matrix representation of 
$\BF$ if we assert that
\begin{align}
X &= X_0 \Id + \sum X_{i} e_i
\end{align}
is the left unit of $\,\&r\,$
\begin{align}
[\BF] &=
\left[\begin{array}{cccc}
\frac{1}{X_0} & 0 & \ldots & 0\\
0 & \ldots & \ldots & 0 \\
0 & \ldots & \ldots & 0 \\
0 & \ldots & \ldots & 0 \\
B_{4,1} & \ldots & \ldots & B_{4,8}\\
\vdots & \ldots & \ldots & \vdots \\
B_{8,1} & \ldots & \ldots & B_{8,8}
\end{array}\right]
\end{align}
This shows that one can have non-trivial such units.

The most interesting case is however given by an element $X$ which is a
primitive idempotent element. Of course, this requires the basic product
and/or co-product to be non-local and we deal no longer with an augmented 
connected Hopf gebra, which makes this construction peculiar.
But, such an element can easily be turned into a left or right unit by 
the above mechanism. Idempotent elements are connected with minimal left 
ideals and `spinor' representations. Probably more important is the fact, 
that the same process can be established for counits. However, counits can 
be related to the vacuum structure of a QFT. Generalized Clifford co-products 
allow hence a great variety of possible candidates for vacua. However, 
we will see that other demands restrict $\BF$ further.

\subsection{Associativity of generalized Clifford products}

Since renormalized time-ordered products are required to remain to be 
associative, we ask next which conditions are necessary for $\BF$ to 
ensure associativity of $\,\&r\,$. It will turn out, that this is a quite 
strong restriction and that one might be forced to drop associativity. 
Recall that the definition of $\,\&r\,$ was
\begin{align}
u\;\&r\;v &= \BF(u_{(2)}, v_{(1)}) u_{(1)} \wedge v_{(2)}
\end{align}
with an arbitrary bilinear form $\BF$, $\,\&r\,$ may not even be unital.
We need further more that the co-product $\Delta$ is an algebra homomorphism
\begin{align}
\Delta(a\wedge b) &= (a\wedge b)_{(1)}\otimes (a\wedge b)_{(2)} \nn
&= (-1)^{\vert a_{(2)}\vert\vert b_{(1)}\vert}
(a_{(1)}\wedge b_{(1)})\otimes (a_{(2)}\wedge b_{(2)}) .
\end{align}
Compute
\begin{align}
(u \,\&r\, v) \,\&r\, w
&= \BF(u_{(2)}, v_{(1)}) (u_{(1)} \wedge v_{(2)}) \,\&r\, w \nn
&= (-1)^{\vert u_{(12)}\vert\vert v_{(21)}\vert}
   \BF(u_{(2)}, v_{(1)})\BF(u_{(12)} \wedge v_{(22)}, w_{(1)})
    (u_{(11)} \wedge v_{(21)} \wedge w_{(2)}) \nn
&= (-1)^{\vert u_{(21)}\vert\vert v_{(2)}\vert}
   \BF(u_{(22)}, v_{(11)})\BF(u_{(21)} \wedge v_{(12)}, w_{(1)})
    (u_{(1)} \wedge v_{(2)} \wedge w_{(2)}) \nn
\end{align}
where we have used co-associativity and (graded) co-commutativity of 
the Gra{\ss}mann co-product, which results in the replacements 
$(u_{(11)},u_{(12)},u_{(2)}) \rightarrow (u_{(1)},u_{(21)},u_{(22)})$, 
$(v_{(1)},v_{(21)},v_{(22)}) \rightarrow (v_{(11)},v_{(12)},v_{(2)})$
and $v_{(22)}\otimes v_{(21)} = (-1)^{\vert v_{(22)}\vert\vert v_{(21)}\vert}
v_{(21)}\otimes v_{(22)}$. In the same manner we compute
\begin{align}
u \,\&r\, (v \,\&r\, w)
&= u \,\&r\, (\BF(v_{(2)}, w_{(1)}) (v_{(1)} \wedge w_{(2)})) \nn
&= (-1)^{\vert v_{(12)}\vert\vert w_{(21)}\vert}
   \BF(u_{(2)}, v_{(11)} \wedge w_{(21)})\BF(v_{(2)}, w_{(1)})
    (u_{(1)} \wedge v_{(12)} \wedge w_{(22)}) \nn
&= (-1)^{\vert v_{(2)}\vert\vert w_{(22)}\vert}
   \BF(u_{(2)}, v_{(11)} \wedge w_{(12)})\BF(v_{(12)}, w_{(11)})
    (u_{(1)} \wedge v_{(2)} \wedge w_{(2)})
\end{align}
where once more co-associativity and co-commutativity have been used. 
The requirement of associativity
\begin{align}
(u\,\&r\,v) \,\&r\, w &=
u\,\&r\,(v \,\&r\, w)
\end{align} 
implies that the coefficients of the above equations have, after a
renaming, to fulfil
\begin{align}
(-1)^{\vert u_{(1)}\vert}
\BF(u_{(2)}, v_{(1)})\,\BF(u_{(1)} \wedge v_{(2)}, w)
&=
(-1)^{\vert w_{(2)}\vert}
\BF(u, v_{(1)} \wedge w_{(2)})\,\BF(v_{(2)}, w_{(1)})
\end{align}
which can be rewritten using the antipode as
\begin{align}
\BF(u_{(2)}, v_{(1)})\,\BF(S(u_{(1)}) \wedge v_{(2)}, w)
&=
\BF(u, v_{(1)} \wedge S(w_{(2)}))\,\BF(v_{(2)}, w_{(1)}) .
\end{align}
This requirement should be compared with Brouder's coupling identities
\cite{brouder:2002a}. The coupling identity of Brouder is close but not 
identical to the conditions given by Rota-Stein for a Laplace Hopf 
algebra \cite{rota:stein:1994a}.

We can try to simplify further this condition by employing 
product co-product duality. That is, we can Laplace expand the
bilinear forms having a wedge product in one of its arguments
\begin{align}
\BF(u,v\wedge w) &= \BF(u_{(1)},w)\BF(u_{(2)},v).
\end{align}
This yields
\begin{align}
(-1)^{\vert u_{(1)}\vert}
\BF(u_{(2)}, v_{(1)})\,\BF(u_{(1)},w_{(2)})\,\BF(v_{(2)},w_{(1)})
=\nn
(-1)^{\vert w_{(2)}\vert}
\BF(u, v_{(1)} \wedge w_{(2)})\,\BF(v_{(2)}, w_{(1)}) .
\end{align}
Cancelling out the common factor and renaming yield the product
co-product duality up to a sign. We get
\begin{align}
\BF(u,v\wedge w) &= (-1)^{\vert w\vert\vert u_{(1)}\vert}
\BF(u_{(1)},w)\,\BF(u_{(2)},v) .
\end{align}
If $\BF$ is a graded bilinear form, we have shown
\begin{thrm}
A graded pairing which obeys product co-product duality results in 
an associative cliffordization.
\end{thrm}

On the other hand, if we assume that $\Id$ remains to be a unit
we find additionally the following set of relations
\begin{align}
\BF(\Id,\Id)&=1 \nn
\BF(\Id,X)  &= \epsilon(X) \nn
\BF(X,\Id)  &= \epsilon(X) \nn
\end{align}
which can be completed by product co-product duality and 
the requirement that $\BF^{-1}$ is the convolutional inverse
of $\BF$, in formulas
\begin{align}
\epsilon(u)\epsilon(v) &=
\BF(u_{(1)},v_{(2)})\,\BF^{-1}(u_{(2)},v_{(1)})
\end{align} 
which results in 
\begin{align}
\BF^{-1}(u,v)&=\BF(S(u),v) \nn
\BF^{-1}(u,v)&=\BF(u,S(v)) \nn
\BF(u,v)&=\BF(S(u),S(v))
\end{align}
where we have assumed that $S^2=\Id$ is involutive. Such a structure
is called a co-quasitriangular structure \cite{majid:1995a}. This
structure will be investigated elsewhere. The above condition
derived from associativity may be addressed as a $2$-cocycle,
even if product co-product duality does not hold.

\subsection{Commutation relations and generalized Clifford products}

Quantum field theory needs not only Wick monomials, normal- and 
time-ordered products and correlation functions, but employs also
a canonical field quantization. Since the quantization encodes 
properties of the system under consideration they should not be
altered by renormalization. Therefore, since we want to go for
renormalization, one has to assert that (anti)commutation relations 
are not altered by the generalized cliffordization process and thus 
by the renormalization process. But, this is true only for the 
formulation of the commutation relations on $V \times V$, higher terms 
might be altered by renormalization effects!

We investigate what kind of assertion is required on $\BF$ to guarantee
that the basic (anti)\-commutation relations of the generators 
(field operators) remain unaltered. We have to demand
\begin{align}
\{\psi_1,\psi_2\}_+ &= 2\,g_{12} \quad\quad \text{or equivalently}\nn
\{   e_1,e_2   \}_+ &= 2\,g_{12}
\end{align} 
where $g$ is extended as usual to $\bigwedge V$ by exponentiation which 
yields $g^\wedge$. However, replacing the $\wedge$ in the anticommutator
by the generalized Clifford product $\&r$ one ends up with additional 
terms. Note that
\begin{align}
e_i \,\&r\, e_j 
&= \BF(e_{i(2)},e_{j(1)}) e_{i(1)} \wedge e_{j(2)} \nn
&= \BF(\Id,\Id) e_i \wedge e_j +
 \BF(\Id,e_j) e_i + \BF(e_i,\Id) e_j + \BF(e_i,e_j) \Id
\end{align}
which comes up with the unwanted second and third terms. Hence our 
demand is that
\begin{align}
\BF(\Id,a) &= 0 \quad\quad \forall a \in V \nn
\BF(a,\Id) &= 0 \quad\quad \forall a \in V.
\end{align}
It turns out, that this is the case in the definitions of the pairings
for renormalization. We will see later, that the $Z$-pairing has to be 
a $\openZ_2$-graded morphism, which gives a second argument for this 
relation.

\subsection{Laplace expansion i.e. product co-product duality implies 
exponentially generated bilinear forms}

This section comes up with a peculiarity about the Laplace expansion.
Rota and Stein introduced a so called Laplace Hopf algebra, which is 
an augmented connected Hopf algebra where the product may be deformed by
cliffordization and the bilinear form permits Laplace expansion 
\cite{rota:stein:1994a}. In fact they gave two more relations which we will 
not consider here. Also Brouder used Laplace Hopf algebras in his work
on renormalization \cite{brouder:2002a} which is the main interest of our 
study too.

We made in this work good use of product co-product duality, and Laplace 
expansions. This motivated to examine in which way the condition that
the wedge product and the Gra{\ss}mann co-product are related.
In form of a tangles product co-product duality for tow given such
structures imposes restrictions on $\BF$. This reads in tangles
\begin{align}
\pspicture[0.5](0,0)(1.5,2.5)
\psset{linewidth=\pstlw,xunit=0.5,yunit=0.5,runit=0.5}
\psset{arrowsize=2pt 2,arrowinset=0.2}
\psline(3,5)(3,1)
\psline(2,5)(2,2)
\psline(0,5)(0,2)
\psline{->}(3,4)(3,3.75)
\psline{->}(2,4)(2,3.75)
\psline{->}(0,4)(0,3.75)
\psarc(2,1){1}{180}{360}
\psarc(1,2){1}{180}{360}
\pscircle[linewidth=0.4pt,fillstyle=solid,fillcolor=white](1,1){0.2}
\pscircle[linewidth=0.4pt,fillstyle=solid,fillcolor=black](2,0){0.2}
\rput(1,2){$\wedge$}
\rput(2,0.75){$\BF$}
\endpspicture
\quad=\quad
\pspicture[0.5](0,0)(3,2.5)
\psset{linewidth=\pstlw,xunit=0.5,yunit=0.5,runit=0.5}
\psset{arrowsize=2pt 2,arrowinset=0.2}
\psline(5,5)(5,4)
\psline(2,5)(2,3)
\psline(0,5)(0,2)
\psline(6,3)(6,2)
\psline(2,0)(4,0)
\psline{->}(5,4.5)(5,4.25)
\psline{->}(2,4.5)(2,4.25)
\psline{->}(0,4.5)(0,4.25)
\psline(4,0)(4,0)
\psarc(5,3){1}{0}{180}
\psarc(3,3){1}{180}{360}
\psarc(2,2){2}{180}{270}
\psarc(4,2){2}{270}{360}
\pscircle[linewidth=0.4pt,fillstyle=solid,fillcolor=white](5,4){0.2}
\pscircle[linewidth=0.4pt,fillstyle=solid,fillcolor=black](3,0){0.2}
\pscircle[linewidth=0.4pt,fillstyle=solid,fillcolor=black](3,2){0.2}
\rput(3,0.75){$\BF$}
\rput(3,2.75){$\BF$}
\endpspicture
\end{align}
and
\begin{align}
\pspicture[0.5](0,0)(1.5,2.5)
\psset{linewidth=\pstlw,xunit=0.5,yunit=0.5,runit=0.5}
\psset{arrowsize=2pt 2,arrowinset=0.2}
\psline(0,5)(0,1)
\psline(1,5)(1,2)
\psline(3,5)(3,2)
\psline{->}(0,4)(0,3.75)
\psline{->}(1,4)(1,3.75)
\psline{->}(3,4)(3,3.75)
\psarc(1,1){1}{180}{360}
\psarc(2,2){1}{180}{360}
\pscircle[linewidth=0.4pt,fillstyle=solid,fillcolor=white](2,1){0.2}
\pscircle[linewidth=0.4pt,fillstyle=solid,fillcolor=black](1,0){0.2}
\rput(1,0.75){$\BF$}
\rput(2,1.75){$\wedge$}
\endpspicture
\quad=\quad
\pspicture[0.5](0,0)(3,2.5)
\psset{linewidth=\pstlw,xunit=0.5,yunit=0.5,runit=0.5}
\psset{arrowsize=2pt 2,arrowinset=0.2}
\psline(1,5)(1,4)
\psline(4,5)(4,3)
\psline(6,5)(6,2)
\psline(0,3)(0,2)
\psline{->}(1,4.5)(1,4.25)
\psline{->}(4,4.5)(4,4.25)
\psline{->}(6,4.5)(6,4.25)
\psline(2,0)(4,0)
\psarc(1,3){1}{0}{180}
\psarc(3,3){1}{180}{360}
\psarc(2,2){2}{180}{270}
\psarc(4,2){2}{270}{360}
\pscircle[linewidth=0.4pt,fillstyle=solid,fillcolor=white](1,4){0.2}
\pscircle[linewidth=0.4pt,fillstyle=solid,fillcolor=black](3,0){0.2}
\pscircle[linewidth=0.4pt,fillstyle=solid,fillcolor=black](3,2){0.2}
\rput(3,0.75){$\BF$}
\rput(3,2.75){$\BF$}
\endpspicture
\end{align}
If we write this in algebraic terms using Sweedler notation, we obtain 
for arbitrary elements $u,v,w \in \bigwedge V$ 
\begin{align}
\BF(u\wedge v, w) &= \BF(u,w_{(2)})\BF(v,w_{(1)}) 
\hskip 1.5truecm \text{and} \nn
\BF(u, v \wedge w) &= \BF(u_{(1)},w)\BF(u_{(2)},v) .
\end{align}
This is in fact the Laplace expansion into rows or columns, where the
Hopf gebraic method allows to expand in a single step into a couple
of rows or columns. It was shown in \cite{grosshans:rota:stein:1987a}
that the Laplace pairing implies an exponential
extension of $B$ to $B^\wedge$. This was used also in our consideration
about associativity. It is obvious from the exponential representation,
that bilinear forms can be added in the following way
\begin{align}
B &= g + F \nn
B^\wedge &= e^{B}\,=\, e^{g+F} \,=\, e^{g}\,e^{F}\,e^{\frac{1}{2}[g,F]}
\times\ldots
\end{align}
which is a Hausdorff like formula. Pairings which are obtained from
co-boundaries, as will be done below, result in exponentially 
generated bilinear forms. However, since they have to be calculated
via Hausdorff like formulas, also here the Hopf gebra approach is
indispensable to manage the complexity of the calculations. On the other 
hand, it will be interesting to study Hausdorff like formulas on their
own right using Hopf gebraic methods. We have to add, that a Clifford
product and its dualized co-product do {\it not\/} form a Hopf gebra
but only a biconvolution, since no antipode exists. This can be 
circumvented if one introduces different orderings for vectors and
co-vectors \cite{fauser:oziewicz:2001a}.

\section{Renormalization group and $Z$-pairing}

\subsection{Renormalization group}

To be able to formulate renormalized time-ordered products Brouder had to 
introduce renormalization parameters. This is done by using scalar valued
linear forms $Z : \bigwedge V\rightarrow \openk$. This introduces $2^n$
independent renormalization parameters and even infinitely many such
parameters if an uncountable or contiguous index set is used. We require
that the linear forms $Z$ do form a group under the convolution product
\begin{align}
Z \star Z^\prime &= Z^\pprime \nn
Z \star Z^{-1} &= u \,=\,\eta\circ\epsilon . 
\end{align}
We assume that the convolution is w.r.t. a Gra{\ss}mann Hopf gebra or 
a symmetric such Hopf gebra, which might be called Weyl Hopf gebra. In 
this case it is possible to deduce the convolutive unit to be 
$u=\eta\circ\epsilon$. Such Hopf gebras are augmented and connected.

In such cases one is able to define the inverse $Z^{-1}$ of $Z$ by the 
well established recursive way. We require that $Z$ is normalized and 
acts trivially on the generators
\begin{align}
&&Z(1)=1&&Z(e_i)=0.&&
\end{align}
Furthermore, we use the notion of proper cuts indicated by a prime at
the sum over the terms of the co-product
\begin{align}
\Delta^\prime(x) &= \Delta(X) - \Id\otimes x - x\otimes\Id
\,=\, {\sum}^\prime x_{(1)}\otimes x_{(2)}. 
\end{align}
We stress once more that this is possible for non-interacting Hopf 
gebras only. From
\begin{align}
Z\star Z^{-1}(x) &= \eta\circ\epsilon(x) \,=\, x_0 \nn
&= Z(x_{(1)})Z^{-1}(x_{(2)}) \nn
&=Z^{-1}(x)+Z(x)+{\sum}^\prime x_{(1)}\otimes x_{(2)}
\end{align}
one finds
\begin{align}
Z^{-1} &=(x_0-Z(x))-{\sum}^\prime x_{(1)}\otimes x_{(2)}
\end{align}
for the inverse. Especially one finds
\begin{align}
&&Z^{-1}(x) =1 && Z^{-1}(x) =0\quad\forall x\in V&&
\end{align}
showing that $Z^{-1}(x)$ belongs to the class of linear forms
which we consider. 

These linear forms constitute a group under the convolution,
which is called {\it renormalization group\/}. We ask under
which conditions this group is abelian. That is, we want to
have $A\star B = B\star A$ for our linear forms. We compute 
this convolution as follows:
\begin{align}
A\star B(x) &= A(x_{(1)}) B(x_{(2)}) \nn
&=(-1)^{\vert x_{(1)}\vert\vert x_{(2)}\vert}
  A(x_{(2)}) B(x_{(1)}) \nn
&=(-1)^{(\vert x_{(1)}\vert\vert x_{(2)}\vert 
       +\vert A(x_{(2)})\vert\vert B(x_{(1)})\vert)}
  B(x_{(1)}) A(x_{(2)}) \nn
&=(-1)^{(\vert x_{(1)}\vert\vert x_{(2)}\vert 
       +\vert A(x_{(2)})\vert\vert B(x_{(1)})\vert)}
  B\star A(x) 
\end{align}
We define $\openZ_2$-graded maps as follows
\begin{align}
C : \bigwedge V^\pm \rightarrow \bigwedge V^\pm
\end{align}
and one obtains in this case that
\begin{align}
\vert x\vert &= \vert C(x)\vert \quad \mod 2
\end{align}
holds. If the above discussed maps $A$ and $B$ are $\openZ_2$-graded,
we find for the prefactor
\begin{align}
(-1)^{\vert x_{(1)}\vert\vert x_{(2)}\vert 
       +\vert A(x_{(2)})\vert\vert B(x_{(1)})\vert}
&= (-)^{2\vert x_{(1)}\vert\vert x_{(2)}\vert  } \,=\, +1
\end{align}
and the corresponding convolution product is commutative, hence the
renormalization group is abelian. Since our maps $Z$ shall be normalized,
we have to choose them to be even, i.e. they are trivial on the odd parts
\begin{align}
Z &: \bigwedge^{2n+1} V \rightarrow 0.
\end{align}

Any linear map can be written as a linear operator $Z^\prime$ followed
by the counit
\begin{align}
Z &= \epsilon\circ Z^\prime&& Z^\prime : \bigwedge V \rightarrow \bigwedge V
\end{align}
which introduces however additional spurious or `gauge' parameters due to 
the projection. $Z^\prime$ can now be written as a bilinear form, that is 
as a cup tangle by `bending up' one leg using evaluation, acting then in
the following way
\begin{align}
&&\partial Z^\prime &: \bigwedge V \times \bigvee V^* \rightarrow \openk ,&& 
\end{align} 
i.e. a cup tangle on this `quantum double' of spaces. This allows to
apply a boundary operator which results in the following co-chain,
see \cite{majid:1995a}
\begin{align}
\partial Z^\prime(u,v^*) &=
Z(u_{(1)})Z(v^*_{(2)}) Z^{-1}(u_{(2)}\wedge v^*_{(1)}).
\end{align}
In terms of tangles this reads
\begin{align}
\pspicture[0.5](0,0)(1,2.5)
\psset{linewidth=\pstlw,xunit=0.5,yunit=0.5,runit=0.5}
\psline(0,4)(0,3)
\psline(2,4)(2,3)
\psline{->}(0,4)(0,3.5)
\psline{<-}(2,3.75)(2,3)
\psarc(1,3){1}{180}{360}
\pscircle[linewidth=0.4pt,fillstyle=solid,fillcolor=black](1,2){0.2}
\rput(1,0.75){$\partial Z^\prime$}
\endpspicture
\quad&=
\quad
\pspicture[0.5](0,0)(3,2.5)
\psset{linewidth=\pstlw,xunit=0.5,yunit=0.5,runit=0.5}
\psline(1,5)(1,4)
\psline(5,5)(5,4)
\psline{->}(1,5)(1,4.5)
\psline{<-}(5,4.75)(5,4)
\psarc(1,3){1}{0}{180}
\psarc(5,3){1}{0}{180}
\psarc(3,3){1}{180}{360}
\psline(3,2)(3,1)
\pscircle[linewidth=0.4pt,fillstyle=solid,fillcolor=white](1,4){0.2}
\pscircle[linewidth=0.4pt,fillstyle=solid,fillcolor=white](5,4){0.2}
\pscircle[linewidth=0.4pt,fillstyle=solid,fillcolor=white](3,2){0.2}
\rput(0,2.25){$Z$}
\rput(6,2.25){$Z$}
\rput(3,0.25){$Z^{-1}$}
\endpspicture
\end{align}
and one can use commutativity of the convolution product to remove
ambiguties at which place, left or right, the linear forms have to be 
applied. If we introduce a self dual space ${\bf V}=V\oplus V^*$, we 
can neglect the duality, i.e. the arrows in the tangles and arrive at 
the $Z$-pairing employed by Brouder.

\subsection{Renormalized time-ordered products as generalized 
Clifford products}

It was proved in \cite{fauser:2001b} that normal-ordered and time-ordered
operator products can be defined as wedge and cliffordized wedge products,
where the latter was called dotted wedge product and emerged from a 
cliffordization w.r.t. a fully antisymmetric bilinear form $F^{\wedge}$.
While the wedge can be identified with the time-ordered products the dotted
wedge corresponds to the normal-ordered case. In our previously
published works, we had restricted this mechanism to Hamilton formalism,
which relays on a one-time formulation 
\cite{stumpf:fauser:pfister:1993a,fauser:1996c,fauser:1996b,%
fauser:stumpf:1997a,fauser:1998a}
for convenience and to make contact to the well established formalism
of Stumpf and coll., see \cite{stumpf:borne:1994a,borne:lochak:stumpf:2001a}.
The problem of renormalization is avoided in a different manner by Stumpf,
since the main interest was in non-linear spinor field theory, which is 
non-renormalizable. The transition from multi-time to one-time formulations
is intimately connected to renormalization. But, the Hopf gebraic method is 
purely combinatorial and does not make any assumption about the nature
of the coordinates of the fields and applies for multi-time correlation
functions too.

In the following we introduce with Brouder \cite{brouder:2002a} some 
additional pairings which can be used to define cliffordizations
and generalized Clifford products. This is done with hindsight to
come up with a structure suitable for renormalization. These pairings
and the resulting cliffordization are easily established in BIGEBRA
and we used this device to check the results given here.

Since we want to have normal- and time-ordered operator products
and correlation functions, we have to have a bilinear form $B=g+F$,
$g^T=g$, $F^T=-F$ in the case of fermions. As long as algebraic 
relations are considered, this mechanism works for bosons in a similar 
manner. In fact the original work of Brouder deals with bosons.

To come up with a generalized product, we have to add to the bilinear form 
$B^{\wedge}$ additional parameters, the parameters of renormalization. We
use the renormalization group.

Having $Z$ and $Z^{-1}$ established we are able to define the Z-pairing
as a product deformation of the time-ordered product, i.e. the wedge.
This is the above defined co-chain and reads as follows
\begin{align}
&\partial Z : \bigwedge V \times \bigwedge V \rightarrow \openk\circ\Id \nn
&\partial Z(x,y) = \sum 
Z(x_{(1)})Z(y_{(1)})Z^{-1}(x_{(2)}\wedge y_{(2)}) .
\end{align}
In \cite{fauser:1996b} we had shown that in reordering processes 
which connect bilinear forms like $g^\wedge$ and $B^\wedge$ no additional 
divergencies occur! A grade preserving bilinear form $F^\wedge$ which 
arrises from exponentiating seemed not be able to mediate renormalization.
We have
\begin{align}
{\cal N}(\psi_1,\ldots,\psi_n) &= \psi_1 \dwedge \ldots \dwedge \psi_n&&
\text{normal-ordered product} \nn
{\cal T}(\psi_1,\ldots,\psi_n) &= \psi_1 \wedge  \ldots \wedge  \psi_n&&
\text{time-ordered product.}
\end{align}
We obtained a `transition formula' from the cliffordization 
\cite{fauser:2001b}
\begin{align}
u \dwedge v &= F^\wedge (u_{(2)},v_{(1)}) (u_{(1)} \wedge v_{(2)}).
\end{align}
Which holds in both directions. Choosing the counit of the 
normal-ordered algebra, i.e. w.r.t. the dotted wedge, to be the vacuum,
we note that $\dwedge$ is a local product.
However already the time-ordered product $\wedge$ is then non-local with 
great implications for the vacuum structure. Remember that as {\it algebras}
both structures are isomorphic, which probably rendered their distinction
to be difficult for a long time. Both such structures are {\it quantized}
due to the introduction of a symmetric bilinear form $g$ extended by 
exponentiation to $g^\wedge$ in the same manner.

Brouder imposed on the $Z$-mappings the following restrictions, which
we derived above from certain assertions about the properties of the
product.

$Z(\Id)=\Id$, $Z(a)=0$ $\forall a\in V$ and arbitrary else. This results
using the recursive formula for the inverse in the following conditions
for $Z^{-1}$. $Z^{-1}(\Id)=\Id$, $Z^{-1}(a)=0$ $\forall a \in V$. These 
conditions necessarily establish that the unit $\Id$ remains to be the 
unit of the renormalized time-ordered product algebra and that the
anticommutation relations are not altered {\it on the generators}!
Furthermore, Brouder proved, in the case of bosons, that the $Z$-pairing
fulfils the coupling identity which ensures associativity of the product 
as we have shown above.

Properties of the $Z$-pairing if used as cup tangle in a cliffordization are:
\begin{itemize}
\item $Z$-cliffordization yields an unital algebra with $\Id$ as unit.
\item $Z$-cliffordization preserves associativity, since it fulfils the co-chain or coupling condition.
\item The $Z$-pairing is even, the renormalization group is therefore abelian.
\item $Z$-cliffordization respects the quantization. That is 
$Z$-cliffordization does not alter the (anti)commutation relations of the 
generating space $V\times V$.
\end{itemize}

Before we can proceed to renormalized time-ordered products, we
have to combine the bilinear form obtained from quantization and 
time-ordering and the $Z$-pairing. This is done using Hopf gebra 
methods and yields the total bilinear form $\BF$
\begin{align}
\BF(u,v) &= \partial Z(u_{(1)},v_{(2)}) B^\wedge(u_{(2)},v_{(1)})
\end{align}
where we arranged the order of the entries in the r.h.s. to avoid 
crossings. This differs from Brouder, since he used bosonic algebras
there were no sign problem for him. This could be called the Hausdorff 
bilinear form of $B^\wedge$ and $\partial Z$ for reasons given above.
The proof that the cliffordization w.r.t. the bilinear form $\BF$
yields the renormalized time-ordered product can be found in 
\cite{brouder:2002a}.

The remaining problem is to fix the multiplicative renormalization
constants by some arguments from physics. In fact they have to be chosen
to subtract the divergencies emerging in perturbative QFT calculations.
A mathematical basis of axioms would however provide a finite theory from 
the beginning and not order by order. Furthermore, if one is interested in 
non-Fock vacua, e.g. in QCD for studying confinement etc. or for the
calculation of composites, there is no way out of an {\it a priori}
renormalization which should be based on mathematical arguments.

Brouder \cite{brouder:2002a}, using results of Pinter 
\cite{pinter:2000a,pinter:2001a}, has shown that the operator
products and correlation function obtained from $\BF$-cliffordization, 
which includes time-ordering, normal-ordering for the counit, and 
renormalization mediated by the $Z$-pairing is equivalent to the 
Epstein-Glaser renormalization of time-ordered products 
\cite{epstein:glaser:1973a}. This allows us to identify
the product $\&r$ with a renormalized time-ordered operator product
and the coefficients derived therefrom, see Brouder's paper. 
The cliffordization results in a tremendous simple formula for the 
renormalized Green functions. Let ${\cal S}$ be the $S$-matrix, one has
\begin{align}
(\psi_1 \wedge \psi_2 \wedge {\cal T}({\cal S}^\wedge(\psi_n)) &= 
G_{12}{\cal T}({\cal S}^\wedge(\psi_{n+2})) \nn
(\psi_1 \circ \psi_2 \circ {\cal T}^{\text{ren}}({\cal S}^\wedge(\psi_n)) 
&= 
G^{\text{ren}}_{12}{\cal T}^{\text{ren}}({\cal S}^\wedge(\psi_{n+2})).
\end{align}
As a result we see that renormalization is still within the cliffordization 
scheme but needs a generalized Clifford product based on a bilinear form. 
This has a tremendous impact e.g. on the structure of higher commutators, 
the generic antipode, the generic crossing of the renormalized product etc. 
This structure will be considered in more detail elsewhere.

The Brouder formulation of Epstein-Glaser renormalization of time-ordered
products turns out to be {\it one} but probably not the most general or 
unique generalization of a Clifford product arising from cliffordization.
An axiomatization of renormalized cliffordizations would be most desirable.
In this sense our investigations of this chapter are preliminary.

%% file: qft.tex
\chapter{(Fermionic) quantum field theory and Clifford Hopf gebra}

In this chapter we develop a formulation of fermionic quantum field
theory (QFT) based on Hopf gebraic methods. We concentrate on fermions,
however, the bosonic case runs along the same lines. Bosons will
occur in our treatment of spinor quantum electrodynamic. This approach 
to QFT would have not been possible without the versatile development of 
functional quantum field theory by Stumpf and coll. 
\cite{pfister:1987a,pfister:1990a,stumpf:borne:1994a,kerschner:1994a,%
grimm:1994a,fauser:1996c,borne:lochak:stumpf:2001a}. 
Since the method is readily available in 
two monographs we will not develop functional QFT here, but deliver 
only those parts which are necessary for our transition to a Hopf 
gebraic treatment. It will turn out, that this is not a mere translation.
The Hopf gebraic formulation will clearly separate concepts which have
not properly been distinguished in conventional QFT. Only this
mathematically sound tool allows to deliver powerful and efficient
formulas which will be useful in calculations too.

The aim of this chapter is to provide an algebraic skeleton for QFT 
which is dealing with all peculiarities of QFT as far as the algebraic 
parts are concerned. Renormalization is treated along the lines of the 
last chapter following Brouder \cite{brouder:2002a}. However, the present 
treatment is still not ready-to-use, say for $n$-th order perturbation 
calculations in QED. But it is possible already now to formalize a great 
variety of QF theoretical calculations to such an extend that computer 
algebra systems as CLIFFORD/BIGEBRA for Maple can be used to evaluate 
expressions. This includes e.g. the derivation of functional equations 
or renormalization.  

Path integrals are, at this stage of the development, not appropriate, 
since they are to compact and formal to provide access e.g. to the peculiar 
combinatorics of renormalization. Nevertheless, since path integrals are 
formal solutions of the Schwinger-Dyson hierarchy equations it may be 
possible to come up with a mathematical background for their foundation. 
Reviewing our techniques, as developed so far in this treatise, it may 
be conjectured that the path measure will be related to an umbral operator 
\cite{kisil:1998a}. 
Such an operator is a linear form on our Hopf gebra having very peculiar 
combinatorial and algebraic properties. Umbrae cannot be transformed e.g. 
without taking peculiarities into account etc. Therefore it seems to be convenient 
to develop firstly the more basic functional differential point  of view, 
but using generating functionals.

\section{Field equations}

To be able to define a QFT we need two informations: (i) a field equation 
and (ii) (anti)commu\-ta\-tion rules. The Lagrangian point of view is more 
sophisticated and mostly chosen to incorporate symmetries. However it 
is not essential to the theory. Field equations are formulated in terms
of field operators $\psi_I$, where the multiindex $I$ may contain any
sort of indices, including uncountable indices or continuous variables
$I=(\Lambda,r,t,\alpha,A,\ldots)$. For a compact presentation, and
anticipating the mathematical structure which we are going to use, we
include at the same time the index $\Lambda$ which distinguishes field
operators and their conjugated or dual field operators. An analogous notation
is chosen for bosonic field operators $B_K$.

As the name suggests, field operators are elements of an operator 
algebra. The algebra structure is encoded in the quantization which is
given by (anti)commutation relations
\begin{align}
\{\psi_{I_1},\psi_{I_2}\}_+ = A_{I_1I_2} &&
[B_{K_1},B_{K_2}]_-&= C_{K_1K_2} .
\end{align}
Note that the above relations contain, due to our index doubling, all
commutators including the zero ones. It is obvious, that $A$ is symmetric
and $C$ is antisymmetric.

An example of a field equation would be a non-linear spinor field,
which comes up with the following field equation. This equation
is written in a Schr\"odinger like form which will allow us later
to pass easily to generating functionals and functional equations.
\begin{align}
i\partial_0 \psi_I &= D_{II_1}\psi_{I_1}
        +gV_I^{\{I_1I_2I_2\}_{as}}
        \psi_{I_1}\psi_{I_2}\psi_{I_3} .
\end{align}
$D_{I_1I_2}$ is (almost) the Dirac operator and $V_I^{I_1I_2I_3}$
is a local interaction vertex. Obviously this equation is non-renormalizable,
but very oftenly used in phenomenological models. A slightly more
general spinor equation was used by Stumpf to recover the
standard model of elementary particle theory via composite calculations.
Functional QFT was developed to manage the problems arising from such 
a task since it requires non-perturbative methods.

A second example is spinor QED, i.e. a Dirac field coupled to a
vector boson, the photon. The field equations are derived from the 
classical equations
\begin{align}
i)  &&&\partial_\nu F^{\mu\nu}+\frac{ie_0}{2}
      \Psi C\gamma^{\mu}\sigma^2\Psi=0 &&& \nn
ii) &&& (i\gamma^\mu\partial_\mu-m_0)\Psi
       +e_0 A_\mu\gamma^\mu\sigma^3\Psi=0 &&& \nn
iii)&&&F^{\mu\nu}:=\partial^\mu A^\nu-\partial^\nu A^\mu, &&&
\end{align}
where we have suppressed the indices. $C$ is the charge conjugation matrix
of Dirac theory and $\sigma$ encodes the double index information, taking
care if a field is a adjoint or not. Details may be found in
\cite{fauser:1996c,fauser:stumpf:1997a,fauser:1998a}. 
Introducing Coulomb gauge, eliminating the longitudinal part 
of the vector potential $A$ and introducing an index doubled field
$B_K$ for the canonical pair ($A_K,E_K$) of bosonic fields allows
us to write the field equations as
\begin{align}
i\partial_0\Psi_{I_1}&=D_{I_1I_2}\Psi_{I_2}+W^k_{I_1I_2}B_k\Psi_{I_2}
                 +U_{I_1}^{I_2I_3I_4}\Psi_{I_2}\Psi_{I_3}\Psi_{I_4}\nn
i\partial_0 B_{K_1}&=L_{K_1K_2}B_{K_2}+J_{K_1}^{I_1I_2}\Psi_{I_1}\Psi_{I_2}
\end{align}
where we have used the following abbreviations
\begin{align}
I                  &:= \{\alpha,\Lambda,{\bf r}\}\qquad
K                  := \{k,\eta,{\bf z}\}\qquad
P^{tr}             := 1-\Delta^{-1}\nabla\otimes\nabla\nn
D_{I_1I_2}         &:=-(i\gamma_0\gamma^k\partial_k-\gamma_0
                       m)_{\alpha_1\alpha_2}\delta_{\Lambda_1\Lambda_2}
                       \delta ({\bf r}_1-{\bf r}_2)\nn
W^K_{I_1I_2}       &:=e_0(\gamma_0\gamma^k)_{\alpha_1\alpha_2}
                       \delta({\bf r}_1-{\bf r}_2)
                       \delta({\bf r}_1-{\bf z})\delta_{1\eta}
                       \sigma^3_{\Lambda_1\Lambda_2}\nn
U_{I_1}^{I_2I_3I_4}&:=-\frac{i}{8\pi}e_0^2\left[
                       (C\gamma_0)_{\alpha_2\alpha_3}
                       \delta_{\alpha_1\alpha_4}
                       \sigma^2_{\Lambda_2\Lambda_3}
                       \sigma^2_{\Lambda_1\Lambda_4}
                       \frac{\delta({\bf r}_2-{\bf r}_3)
                       \delta({\bf r}_1-{\bf r}_4)}{
                       \vert{\bf r}_1-{\bf r}_2\vert}\right]_{
                       as\{I_2I_3I_4\}} \nn
L_{K_1K_2}         &:=i\delta({\bf z}_2-{\bf z}_1)\delta_{k_1k_2}
                       \delta_{\eta_11}\delta_{\eta_22}
                       +i\Delta({\bf z}_1)\delta({\bf z}_1-{\bf z}_2)
                       \delta_{k_1k_2}\delta_{\eta_12}
                       \delta_{\eta_21}\nn
J_K^{I_1I_2}       &:=-\frac{1}{2}e_0 P^{tr}({\bf z}-{\bf r}_1)
                       \delta({\bf r}_1-{\bf r}_2)(C\gamma^k)_{
                       \alpha_1\alpha_2}\delta_{2\eta}
                       \sigma^2_{\Lambda_1\Lambda_2}
\end{align}
and impose canonical quantization
\begin{align}
i)  &&& \{\Psi_{I_1},\Psi_{I_2}\}^t_+:=A_{I_1I_2}=C\gamma_0\sigma^1
       \delta({\bf r}_1-{\bf r}_2) &&&\nn
ii) &&& [B_K,\Psi_{I}]^t_-:=0 &&&\nn
iii)&&& [B_{K_1},B_{K_2}]_-^t=:C_{K_1K_2}. &&&
\end{align}
In the above field equation we find the Laplacian $L_{K_1K_2}$,
a non-linear self interaction term for the spinor field 
$U_{I_1}^{I_2I_3I_4}$ which emerges from the longitudinal vector
potential, i.e. the Coulomb interaction and {\it two\/} coupling
terms. One is a boson-fermion coupling $W^K_{I_1I_2}$ and the other
is a fermion-boson coupling $J_{K_1}^{I_1I_2}$. Our presentation
is not covariant. It is well known that consistency implies that 
quantization has to be done
in Hamilton formulation, which is also not covariant. 
However, we do not loose any information as long as we perform exact 
manipulations with this system of equations since we could transform back 
to the covariant picture. However, in the chosen form it is much easier 
to distinguish the various terms and to appreciate their physical 
content.

The task is now to provide a Schwinger-Dyson hierarchy for these
field equations from which in principle results can be computed.
To do this, one has to pass from a single field to a hierarchy
of field equations which will be encoded by a functional.

\section{Functionals}

We will define fermionic functionals in terms of quantum Clifford 
algebras. The bosonic case runs along the same lines by analogy,
as far as no topological assertions have to be made. 

A {\it generating function\/} is used to encode a sequence
of numbers ($a_1,a_2,a_3,\ldots$) by a function, such that they appear 
as the coefficients in a polynomial expansion. 
\begin{align}
f(t) &= \sum \frac{a_n}{n!}t^n \nn
\frac{t}{e^t-1} &= \sum \frac{B_n}{n!} t^n \nn
B_n &= \frac{\partial^n}{\partial t^n} \frac{t}{e^t-1}\vert_{t=0}
\end{align}
where the $B_n$ are the Bernoulli numbers. To get a certain coefficient,
we can use the $n$th-derivative and evaluate it at the point ($t=0$). The idea 
is to generalize this technique in that way, that the coefficients are 
functions, or even distributions. 

Having an operator algebra, we need a semi ordering as we have introduced
for bases in Clifford algebras in chapter 2. Such ordered monomials 
span linearly the space on which the operator algebra is built over. If 
we choose a (semi) ordering ${\cal P}$, we obtain the reduced words  
\begin{align}
e_{I_1,\ldots,I_n} &:= {\cal P}(\psi_{I_1},\ldots,\psi_{I_n}),
\end{align}
which constitute a basis. Note that such a basis is usually assumed 
to have the symmetry of the fields, i.e. in our case antisymmetric
for fermions and symmetric for bosons. This is not necessary, as we saw 
in the case of a Clifford basis versus a Gra{\ss}mann basis, 
but solely used in QFT. However, see the $U(2)$-model discussed below. 
The generators of this non-commutative
polynomial ring are the {\it Schwinger sources\/} $j_{I}$ for fermions
and $b_K$ for bosons. They have to reflect the (anti)symmetry 
of the ordered fields, since we are interested in Gra{\ss}mann or Weyl 
bases. Therefore they have to span Gra{\ss}mann or Weyl algebras.
The derivation operators are written as $\partial_I$ for fermions
and $\delta_K$ for bosons. The commutation relations are ($\delta_{I_1I_2}$,
$\delta_{K_1K_2}$ are Kronecker symbols)
\begin{align}
   \{ j_{I_1}, j_{I_2} \}_+ \,=\, 0 \,=
\, \{ \partial_{I_1}, \partial_{I_2} \}_+
&& \{ j_{I_1}, \partial_{I_2} \}_+ \,=\, \delta_{I_1I_2}
\end{align}
for the fermions, while the bosonic sources fulfil
\begin{align}
   [ b_{K_1}, b_{K_2} ]_- \,=\, 0 \,=
\, [\delta_{K_1}, \delta_{K_2} ]_-
&& [\delta_{K_1}, b_{K_2} ]_- \,=\, \delta_{K_1K_2} .
\end{align}
We have adopted the convention to incorporate the factor $2$ into the
definition of the bilinear forms. This will lead to a factor $1/2$
in contractions, but follows the QFT standard. It is clear, that
the $j_I$ sources generate a Gra{\ss}mann algebra, while the $b_K$ 
sources do the same for a Weyl algebra or symmetric algebra. We denote by 
$V$ the space generated by the $j_I$s and by $\bigwedge V$ the whole space,
similar we use $W$ and $\Sym(W)$ for the Weyl algebra. The product
between the sources is defined to be the wedge product for the $j_I$s,
the vee-product for the duals $\partial_I$ and juxtaposition for the
bosons. The last setting is slightly to narrow, but we will deal with 
bosons only occasionally. 

From the previously obtained results, we know, that these algebras
are Hopf gebras. We can even adjoin a unit map $\eta$ and a counit
map $\epsilon$ and know that these algebras are bi-augmented bi-connected
Hopf gebras. The crossing is the graded switch $\tau$ for fermions 
and the non-graded switch $\sw$ for bosons. 

The physicist's notation for the counit uses functional `vacuum' states. 
One defines
\begin{align}
& j_I\mid0\rangle_F = 0 \qquad\qquad {}_F\!\langle 0\mid \partial_I \nn 
\epsilon^\wedge(U) &=  {}_F\!\langle 0\mid 
\;\sum_{\text{grades}} U^{I_1,\ldots,I_n}
j_{I_1}\wedge\ldots\wedge j_{I_n} 
\mid0\rangle_F \nn
&=U^0 .
\end{align}
In fact, $\epsilon^\wedge$ is the projection onto the coefficient
of the identity element $\Id$, since usually one assumes normalization
so that ${}_F\!\langle 0\mid \Id \mid0\rangle_F =1$.

To expand a reduced word $e_{I_1,\ldots,I_n}$ into the $j$-sources,
we need a mapping ${\rm\sf p} : \bigwedge V \rightarrow \openk$ which
establishes the ordering ${\cal P}$ of the field operators. Before we 
can proceed, we have to give the definition of the field operators
in terms of the sources $j,\partial$. 

We will make use of the Chevalley deformation, knowing  that this can be
generalized by Hopf gebraic means, and define the field 
operator to be a Clifford map. Let $e_{I_1,\ldots,I_n}$ be a reduced
monomial w.r.t. the ordering ${\cal P}$. We have to define two field
operators which allow to add a single field operator to this basis monom,
one adds from the left and one adds from the right. These field operators
read
\begin{align}
\psi_I &= \partial_i + \frac{1}{2}P_{IL}\,j_L\;\wedge \nn
\psi^{op}_I &= \partial_i - \frac{1}{2}P_{LI}\,j_L\;\wedge  .
\end{align}
Note that the indices of the field operator acting by opposite 
multiplication are reversed, due to the fact that $\Cl^{op}(V,P)$ and 
$\Cl(V,-P^T)$ are isomorphic. This could be called a Pieri formula for 
fermionic quantum field theory. The bosonic case is treated similarly.
With our pre-knowledge, we identify the action of such an operator
as a Clifford product and write a circle for this product, where we
leave the ${~}^{op}$ with the field operator abusing its meaning as 
opposite product. We obtain
\begin{align}
\psi_I\circ e_{I_1,\ldots,I_n} &= e_{I,I_1,\ldots,I_n}
\,=\, {\cal P}(\psi_I,\psi_{I_1},\ldots,\psi_{I_n})  \nn
\psi^{op}_I\circ e_{I_1,\ldots,I_n} &= e_{I_1,\ldots,I_n,I}
\,=\, {\cal P}(\psi_{I_1},\ldots,\psi_{I_n},\psi_I) .
\end{align}
Now, given an element $U \in \bigwedge V$, we can expand it into
our $j$-basis using the projection $ \pi^a_p$ as follows
\begin{align}
\pi^a_p(e_{I_1,\ldots,I_n})
&=\pi^a_p\big({\cal P}(\psi_{I_1},\ldots,\psi_{I_n})\big) \nn
&= \langle \partial_{I_n}\vee\ldots\vee\partial_{I_1}
{\cal P}(\psi_{I_1},\ldots,\psi_{I_n})\,\rangle_{\pi^a_p}\;
j_{I_1}\wedge\ldots\wedge j_{I_n} \mid 0 \rangle_F \nn
&=\langle 0\mid \partial_{I_n}\vee\ldots\vee\partial_{I_1}
{\cal P}(\psi_{I_1},\ldots,\psi_{I_n})\mid a\rangle
\;j_{I_1}\wedge\ldots\wedge j_{I_n} \mid 0\rangle_F \nn
&= \rho({I_1},\ldots,{I_n}\mid a)\;
j_{I_1}\wedge\ldots\wedge j_{I_n} \mid 0\rangle_F .
\end{align}
$\pi^a_p$ is a sort of grade projection operator. We are able
to learn several things from this calculation:
\begin{itemize}
\item[i)] The functional vacuum is related to the counit of the 
Gra{\ss}mann algebra of the Schwinger sources.
\item[ii)] The counit of the Schwinger sources is {\it not\/} directly 
related to 
the physical vacuum, which we have denoted by a bra-ket notation also. 
Moreover, this is a second and independent linear form $\langle\ldots
\rangle_{\pi^a_p} : \bigwedge V  \rightarrow \openk$. This
linear form may be parametrized by a physical state $\vert a\rangle$.
If this state is the vacuum state, we deal with vacuum expectation values.
\item[iii)] The physical vacuum depends also on the ordering ${\cal P}$.
The mapping ${\rm\sf p} :\bigwedge V\rightarrow \openk$ is identical 
to the {\it correlation functions\/} $\rho({I_1},\ldots,{I_n}\mid a)$.
\item[iv)] This projection encodes all of the combinatorics of QFT,
as we will see below. 
\end{itemize}
The immediate question is, if there is a second Gra{\ss}mann exterior 
or Clifford 
product, which turns the above projection $\pi^a_p$ into the counit
of this algebra. Moreover we have to ask, what kind of product is induced
by this type of projection inside the reduced words. In the case of
quantum Clifford algebras we know already that we can choose e.g. a 
Clifford basis, a Gra{\ss}mann wedge basis or a dotted Gra{\ss}mann 
wedge basis.

Note that this projection is exactly the same as we have introduced
for the renormalization parameters in the previous chapter. In fact,
we can show the following result. Given a mapping ${\rm\sf p} :\bigwedge V
\rightarrow \openk$ we can define an algebra homomorphism as follows.
Let ${\cal P}$ be the operator acting on $x\in \bigwedge V$ by convolution
in the following way (${\rm\sf p}^{-1}$ defined as ${\rm\sf p}^{-1}
\star {\rm\sf p} \,=\, \eta\circ\epsilon$)
\begin{align}
&&{\cal P}(x) = {\rm\sf p}(x_{(1)})x_{(2)} &&
{\cal P}^{-1}(x) = {\rm\sf p}^{-1}(x_{(1)})x_{(2)}&&
\end{align}
which is assumed to be commutative, i.e.
\begin{align}
{\cal P} &= {\rm\sf p}\star \Id \,=\, \Id \star {\rm\sf p} \nn
{\cal P}^{-1} &= {\rm\sf p}^{-1}\star \Id \,=\, \Id \star {\rm\sf p}^{-1},
\end{align}
then the product inside the ordering is obtained by the following 
homomorphism ($x,y\in\bigwedge V$)
\begin{align}
&&{\cal P}(x\circ^{\rm\sf p} y) = {\cal P}(x) \wedge {\cal P}(y),&&
{\cal P}^{-1}(x\wedge y) = {\cal P}^{-1}(x) \circ^{\rm\sf p}{\cal P}^{-1}(y),&&
\end{align}
moreover the circle product $\circ^p$ is the cliffordization of the 
undeformed wedge product w.r.t. the bilinear form (up to a sign)
\begin{align}
\partial {\cal P}(u,v) &= 
{\rm\sf p}(u_{(1)}){\rm\sf p}(v_{(2)})
{\rm\sf p}^{-1}(u_{(2)}\wedge v_{(1)}) \nn
\partial {\cal P}^{-1}(u,v) &= 
{\rm\sf p}^{-1}(u_{(1)}){\rm\sf p}^{-1}(v_{(2)})
{\rm\sf p}(u_{(2)}\wedge v_{(1)})
\end{align}
This bilinear form is a co-chain, as we have discussed in the previous chapter.
Indeed, it looks as if the `renormalization group' as discussed there
is much more an `ordering group'. We will see, that renormalization might
be addressed as a sort of residual re-ordering if normal-ordering
is done w.r.t. the free propagator. Reorderings are invertible since
we demanded that the endomorphisms ${\rm\sf p}$ do form a group under 
convolution.

To prove the above given statement, one has to use the fact that 
${\cal P}(x)={\rm\sf p}(x_{(1)})x_{(2)} \,=\, x_{(1)}{\rm\sf p}(x_{(2)})$ 
and that the co-product is an algebra homomorphism and the cocommutativity.
In terms of tangles this reads (up to signs):
\begin{align}
\pspicture[0.5](0,0)(1,3.5)
\psset{linewidth=\pstlw,xunit=0.5,yunit=0.5,runit=0.5}
\psset{arrowsize=2pt 2,arrowinset=0.2}
\psline(0,7)(0,4)
\psline(2,7)(2,4)
\psarc(1,4){1}{180}{360}
\psline(1,3)(1,0)
\pscircle[linewidth=0.4pt,fillstyle=solid,fillcolor=white](1,3){0.2}
\rput(1,3.75){$\circ^{\rm\sf p}$}
\endpspicture
\quad:=\quad
\pspicture[0.5](0,0)(2.5,3.5)
\psset{linewidth=\pstlw,xunit=0.5,yunit=0.5,runit=0.5}
\psset{arrowsize=2pt 2,arrowinset=0.2}
\psline(1,7)(1,6)
\psline(4,7)(4,6)
\psarc(1,5){1}{0}{180}
\psarc(4,5){1}{0}{180}
\psarc(2,5){2}{180}{270}
\psarc(3,5){2}{270}{360}
\psline(2.5,3)(2.5,2)
\psline(2,3)(3,3)
\psarc(2.5,1){1}{0}{180}
\psline(3.5,1)(3.5,0)
\pscircle[linewidth=0.4pt,fillstyle=solid,fillcolor=white](1,6){0.2}
\pscircle[linewidth=0.4pt,fillstyle=solid,fillcolor=white](4,6){0.2}
\pscircle[linewidth=0.4pt,fillstyle=solid,fillcolor=white](2.5,3){0.2}
\pscircle[linewidth=0.4pt,fillstyle=solid,fillcolor=white](2.5,2){0.2}
\rput(2,4.25){${\rm\sf p}$}
\rput(3,4.25){${\rm\sf p}$}
\rput(1.5,0.25){${\rm\sf p^{-1}}$}
\endpspicture
\quad=\quad
\pspicture[0.5](0,0)(3.5,3.5)
\psset{linewidth=\pstlw,xunit=0.5,yunit=0.5,runit=0.5}
\psset{arrowsize=2pt 2,arrowinset=0.2}
\psline(1,7)(1,6)
\psline(6,7)(6,6)
\psarc(1,5){1}{0}{180}
\psarc(2,4){1}{0}{180}
\psarc(6,5){1}{0}{180}
\psarc(5,4){1}{0}{180}
\psarc(2,4){1}{180}{270}
\psarc(5,4){1}{270}{360}
\psarc(2,3){2}{180}{270}
\psarc(5,3){2}{270}{360}
\psline(0,5)(0,3)
\psline(7,5)(7,3)
\psline(2,3)(5,3)
\psline(2,1)(5,1)
\psline(3.5,1)(3.5,0)
\pscircle[linewidth=0.4pt,fillstyle=solid,fillcolor=white](1,6){0.2}
\pscircle[linewidth=0.4pt,fillstyle=solid,fillcolor=white](2,5){0.2}
\pscircle[linewidth=0.4pt,fillstyle=solid,fillcolor=white](6,6){0.2}
\pscircle[linewidth=0.4pt,fillstyle=solid,fillcolor=white](5,5){0.2}
\pscircle[linewidth=0.4pt,fillstyle=solid,fillcolor=white](3.5,1){0.2}
\rput(3,3.5){${\rm\sf p}$}
\rput(4,3.5){${\rm\sf p}$}
\rput(3.5,2){${\rm\sf p^{-1}}$}
\psline(3.5,3)(3.5,2.5)
\pscircle[linewidth=0.4pt,fillstyle=solid,fillcolor=white](3.5,3){0.2}
\endpspicture
\end{align}
Hence we can state, under the above given assumptions, that if an 
algebra morphism is constructed from a linear form by convolution 
then it induces a cliffordization where the bilinear form is a cochain. 

After this preliminary consideration we can write down a generating functional.
For this purpose, we add up the reduced monoms including the prefactor, that
is the $n$-point correlation function, to build up a general element of
$\bigwedge V$. This reads, if we indicate also the state $\vert a>$ 
w.r.t. the transition matrix elements, as
\begin{align}
\vert {\cal P}(j,a)\rangle_F &= \sum_{i=0}^{n} \frac{i^n}{n!} 
\rho_n(I_1,\ldots,I_n\mid a)\;j_{I_1}\wedge\ldots\wedge j_{I_n} 
\vert 0\rangle_F\nn
\rho_n(I_1,\ldots,I_n\mid a)
&= <0\mid {\cal P}(\psi_{I_1},\ldots,\psi_{I_n}) \vert a> . 
\end{align}
With respect to the endomorphisms of $\bigwedge V$, this is a `state' and has
thus also transformation properties. The implementation of the Poincar\'e 
group symmetries on functional spaces e.g. 
is discussed in great detail in \cite{stumpf:borne:1994a}. 

\section{Functional equations}

We have now generating functionals at our disposal. The next step is to
implement the dynamics on such functionals. They code directly the
Schwinger-Dyson hierarchy, which is the hierarchy of the coupled $n$-point
correlation functions. Our goal is to derive a Schr\"odinger like
equation for such functional states. It should however be noted that
our basis elements, the reduced words $e_{I_1,\ldots,I_n}$ are neither
normalized nor orthogonal and cannot directly be interpreted in physical
terms. If one assumes a Fock representation, the ordinary perturbative
treatment does apply.

Our starting point is the Heisenberg equation
\begin{align}
i\dot{\psi}_I &= [\psi, H]_- \nn
H &= H[\psi].
\end{align}
where $H[\psi]$ is assumed to generate a one parameter family of
automorphisms by integrating the above equation
\begin{align}
\psi_I(t) &= e^{iHt} \psi_I(0) e^{-iHt}.
\end{align}
This equation translates into the following functional equation
\begin{align}
i\partial_0 \vert {\cal P}(j,p)\rangle_F &= H[j,\partial]^{\cal P}
\vert {\cal P}(j,p)\rangle_F
\end{align}
and our task is to calculate the {\it functional Hamiltonian\/} 
$H[j,\partial]^{\cal P}$,
which depends on the chosen ordering. We can use the above defined 
Clifford maps to perform this task. This results in
\begin{align}
i\dot{e}_{I_1,\ldots,I_n} 
&= [e_{I_1,\ldots,I_n}, H[\psi] ]_- 
\,=\, e_{I_1,\ldots,I_n} H[\psi]-H[\psi] e_{I_1,\ldots,I_n} \nn
&= H[\psi^{op}] e_{I_1,\ldots,I_n} -H[\psi] e_{I_1,\ldots,I_n}.
\end{align}
The opposite product allows to write the functional Hamiltonian as an endomorphism 
acting from the left alone. Having an equation for a reduced word, i.e. a basis
monom, we can add up the hierarchy and obtain the functional equation
\begin{align}
i\partial_0 \mid {\cal P}(j,a)\rangle_F &=  ( H[\psi^{op}] - H[\psi] ) \mid
{\cal P}(j,a)\rangle_F \nn
&= H[j,\partial]^{\cal P} \mid {\cal P}(j,a)\rangle_F.
\end{align}
Since we can read of the functional form of $H[\psi]$ directly from the 
equations of motion, we are immediately ready to calculate the functional
Hamiltonian by replacing the field operators with the appropriate 
Clifford map. This is done as follows
\begin{align}
H[j,\partial]^{{\cal P}} &= H[\psi^{op}] - H[\psi] 
\,=\, H[\partial-1/2\, P^Tj] - H[\partial +1/2\, Pj].
\end{align}
Since our field equations are in general polynomial in the interaction terms, 
the calculations can be performed very quickly. In fact, this could be recast
in Hopf gebraic form. For the interaction term of the spinor field theory this reads
\begin{align}
{\cal P}^{-1}(V_{I}^{I_1I_2I_3}\,\psi_{I_1}\wedge \psi_{I_1}\wedge \psi_{I_1})
&=V_{I}^{I_1I_2I_3}\, {\cal P}^{-1}(\psi_{I_1})\,\circ^{{\rm\sf p^{-1}}} 
{\cal P}^{-1}(\psi_{I_2})\,\circ^{{\rm\sf p^{-1}}} {\cal P}^{-1}(\psi_{I_3}) .
\end{align}
If we demand, in difference to the requirements for renormalization discussed
in the previous chapter, that ${\cal P}^{-1}(\psi_{I_1})=\psi_{I_1}$ we end up
with the same term, but w.r.t. the new circle product. The opposite field
operators imply a reversion of the circle products, i.e. a right action.

The crucial point in this consideration is, that we assume that the original 
classical field equations have to be formulated with the same wedge product
which we used for forming the generating functionals. In other words, the 
ordering is an algebra homomorphism which induces in our case quantization 
and ordering in a single step. The quantization stems for fermions from the 
symmetric parts of $P$ while the ordering depends on the antisymmetric
part. In the case of bosons these symmetries are interchanged. The antisymmetric 
part will show up in the next section to be related to the propagator of the theory.

\section{Vertex renormalization}

Since we deal with non-linear terms like the above discussed $\psi^3$ term, 
one has to ask if there occur ordering problems. Indeed, it is well known
from standard treatments, that one has to `remove' vertex singularities.
However, such singularities emerge also in the re-ordering from, say
time- to normal-ordering and vice versa. It is hence only possible to remove 
such singular contributions in {\it one\/} ordering, say the normal-ordering. 
We have shown in Ref. \cite{fauser:1996c,fauser:1996b}, that in the present 
formalism {\it no additional singularities\/} emerge from a re-ordering.

The point is, that in the standard treatment one does not write down the
product. There, the transition is done simply by adding contraction terms,
some of them are diagonal $P_{I_1I_1}\cong P(x,x)$ and diverge. This fact is
usually discussed verbally and as a solution one comes up with a vertex
renormalization denoted by colons
\begin{align}
:\,\psi_{I_1}\psi_{I_2}\psi_{I_3}\,: \;&=\;
\psi_{I_1}\psi_{I_2}\psi_{I_3} 
+ \text{contraction terms} - \text{singularities.}
\end{align} 
Nevertheless, this moves around the singularity only from one picture
into another and does not yield a solution of the problem, 
since if the singularities are substracted in the time-ordered picture, they
reappear e.g. in the normal-ordered formulation. 

Introducing the proper products allows to get rid of these singularities in
{\it all\/} orderings. That is, a transition from one picture into another
does not introduce new spurious singularities. We calculate
\begin{align}
\psi_{I_1}\circ^{{\rm\sf p^{-1}}}\psi_{I_2}\circ^{{\rm\sf p^{-1}}}\psi_{I_3} &=
\psi_{I_1}\wedge\psi_{I_2}\wedge\psi_{I_3}
-P_{I_1I_2}\psi_{i_3}
-P_{I_2I_3}\psi_{i_1}
-P_{I_3I_1}\psi_{i_2}
\end{align}
and no diagonal singular terms occur. This outcome motivates to study
if this mechanism applies directly to functional equations.

\section{Time- and normal-ordering}

In this section we report our findings from Refs. 
\cite{pfister:1987a,fauser:1996c,fauser:1996b,%
fauser:stumpf:1997a,fauser:1998a}, therefore
we do not provide the calculational details, but try to exhibit the newly
established Hopf gebraic aspects. The former calculations are already 
considerable more efficient than the derivation of functional equations by means 
of e.g. Hausdorff formulas. Our `replacement' formalism is tied to the 
exponential representation of the functionals and goes back to Anderson
\cite{anderson:1954a}\footnote{I thank C. Brouder for pointing out this 
reference to me.}.

Before we go into the details of examples, we examine the transition from 
time- to normal-ordered functionals. The time-ordered functional is defined
w.r.t. the wedge product and has $\tau$-functions as coefficients. We 
abbreviate the indices by numbers for convenience.
\begin{align}
\left\vert {\cal T}(j,a)\right \rangle^\wedge
&=
\sum_{n=0}^\infty \frac{i^n}{n!} \tau_n(1,\ldots,n\mid a)
j_1\wedge\ldots\wedge j_n \mid 0 \rangle_F .
\end{align}
The normal-ordered functional is also {\it expanded w.r.t. the wedge
product\/} and has the $\phi$-functions as coefficients.
\begin{align}
\left\vert {\cal N}(j,a)\right \rangle^\wedge
&=
\sum_{n=0}^\infty \frac{i^n}{n!} \phi_n(1,\ldots,n\mid a)
j_1 \wedge\ldots\wedge j_n \mid 0 \rangle_F\,.
\end{align}
Since both fucntionals are expanded in the {\it same basis\/}, it is clear
that they are different in their content. In fact, the normal-ordered
functional corresponds to one-particle irreducible correlation functions
in perturbative QFT. Time- and normal ordered functionals can be related
by the expenentiated propagator $F_{I_1I_2}$ in the following manner
\begin{align}
&&
\left\vert {\cal T}(j,a)\right\rangle^\wedge\,=\,
e^{-\frac{1}{2} F_{I_1I_2}j_{I_1}\wedge j_{I_2}}
\left\vert {\cal N}(j,a)\right\rangle^\wedge 
&&
\left\vert {\cal N}(j,a)\right\rangle^\wedge\,=\,
e^{\frac{1}{2} F_{I_1I_2}j_{I_1}\wedge j_{I_2}}
\left\vert {\cal T}(j,a)\right\rangle^\wedge \,
\end{align}
where we have indicated that the functionals are expanded in the
wedge basis. In Hopf gebraic terms, this transiton is mediated by 
an algebra homomorphism in the following way. Let ${\cal N}$ be the 
algebra morphism which transformes from time- to normal-ordering by
changing the basis from the wedge to the dotted wedge basis
\begin{align}
{\cal N}(u \wedge v) &= {\cal N}(u)\dwedge{\cal N}(v).
\end{align}
Applying this operation to a time ordered functional yields
the normal-ordered such functional
\begin{align}
{\cal N}\big(\left\vert {\cal T}(j,a)\right\rangle^\wedge \big)
&=
\left\vert {\cal N}(j,a)\right\rangle^{\dwedge}.
\end{align}
But this functional has the {\it same\/} expansion coefficients, i.e. it stays
with the $\tau$-functions! Only after we have re-expanded the 
normal-ordered functional into the wegde basis, we end up with the above
given result. This reads explicitely as
\begin{align}
\mid {\cal N}(j,a) \rangle^{\dwedge} &=
\sum \frac{i^n}{n!}\, \tau_n(1,\ldots,n\mid a)\; 
j_{I_1} \dwedge \ldots \dwedge j_{I_n} \mid 0 \rangle_F.\nn
\mid {\cal N}(j,a) \rangle^{\wedge} &=
\sum \frac{i^n}{n!}\, \phi_n(1,\ldots,n\mid a)\; 
j_{I_1} \wedge \ldots \wedge j_{I_n} \mid 0 \rangle_F\, .
\end{align} 
The two expansions shows, that the connenction between the 
$\tau$-functions and the 
$\phi$-functions is not directly mediated by the reordering, but by
the re-expansion of the re-ordered functionals in the wedge basis.
This is not an artifact of the theory, but is closely related to the
fact the we have to chose a unique counit which acts as a projection
onto the identity element. This counit depends on the chosen basis. 
If we use $\epsilon^\wedge$ we have to expand {\it all\/} functionals
in this particular wedge basis. This does on the other hand imply, that 
we have to use different ordering mappings (linear forms on $\bigwedge V$) e.g. 
${\rm\sf t} : \bigwedge V \rightarrow\openk$ and 
${\rm\sf n} : \bigwedge V\rightarrow\openk$
for obtaining the time- or normal-ordered correlation functions. In 
this functions one encodes thereby the information about the physical 
vaccum structure of the theory. This should not be confused with
the functional `vacuum' $\vert 0\rangle_F$ which does not contain
physical informations.

\subsection{Spinor field theory}

We start considering the non-linear spinor field theory. Its Hamiltonian 
is displayed as
\begin{align}
H[\psi]^\wedge
&=\frac{1}{2}A_{I_1I_3}D_{I_3I_2}\psi_{I_1}{\wedge}\psi_{I_2} 
+\frac{g}{4}A_{I_1I_5}V_{I_1}^{I_2I_3I_4}
\psi_{I_1}{\wedge}\psi_{I_2}{\wedge}\psi_{I_3}{\wedge}\psi_{I_4}
\end{align}
where we have introduced explicitly the wedge product. Remember that 
$A_{I_1I_2}$ is the anticommutator of the fields in index doubled formulation.
Since we know, that a Clifford map $\partial\pm 1/2\,Aj\wedge$ is a map
into the Clifford algebra $\Cl(V,A)$ and $\Cl^{op}(V,A)=\Cl(V,-A^T)$,
we are able to identify the corresponding cliffordization with the time-ordered
functional equation. The field operators have to be chosen as 
\begin{align}
\psi &=\frac{1}{i} \partial_i -\frac{i}{2}A_{II_1}j_{I_1}\,\wedge \nn
\psi^{op} &=\frac{1}{i} \partial_i +\frac{i}{2}A_{II_1}j_{I_1}\,\wedge 
\end{align}
to stay compatible with QF theoretic conventions. 
Note that this differs from our 
previously chosen conventions by a swap between product and opposite product
resulting just in a renaming of the fields.
If we assume the state $\vert a\rangle$ to be an eigenstate of the Hamiltonian, 
we end up with the renormalized energy eigen-functional equation
\begin{align}
E_{a0}\vert{\cal T}(j,a)\rangle^\wedge_F&=
\Big\{D_{I_1I_2}j_{I_1}\partial_{I_2}\nn
&
+
gj_{I_1}V_{I_1}^{I_2I_3I_4}\Big(
\partial_{I_3}\partial_{I_2}\partial_{I_1}
+\frac{1}{4}
A_{I_4I_4^\prime}A_{I_3I_3^\prime}
j_{I_4^\prime}j_{I_3^\prime}\partial_{I_2}\Big)\Big\}
\vert{\cal T}(j,a)\rangle^\wedge_F 
\end{align} 
where the energy value $E_{a0}=E_a-E_0$ is the difference of the energy
of the state $\vert a\rangle$ w.r.t. the vacuum energy $E_0$ and thus 
renormalized. One identifies the terms as follows: The $D$ term is the 
kinetic part of the Dirac operator, the $V$ term has two parts, the
interaction part $j\partial^3$ and a quantization part $j^2\partial^2$.
This functional equation is time-ordered. 

However, for composite calculation one needs the normal-ordered 
functional equation. Usually this is obtained by the deviation over the 
intermediate step of the time-ordered equation. But we have another opportunity,
we can simply introduce a different cliffordization based on another ordering.
Two orderings which are based on antisymmetric operator products or
correlation functions can differ only by an antisymmetric part. We can use
the {\it propagator\/} of the theory to perform this transition. The field
operators translate as
\begin{align}
\psi^{op} &= \frac{1}{i}\partial + \frac{i}{2} Aj\dwedge =
\frac{1}{i}\partial+\frac{i}{2}Aj\wedge+iFj\wedge \nn
\psi      &= \frac{1}{i}\partial - \frac{i}{2} Aj\dwedge =
\frac{1}{i}\partial-\frac{i}{2}Aj\wedge+iFj\wedge
\end{align}
and we can easily calculate the normal-ordered energy functional equation
in a {\it single step\/} as 
\begin{align}
E_{a0}\vert {\cal N}(j,a)\rangle^\wedge_F
&=
H[j,\partial]^{\dwedge} \vert {\cal N}(j,a)\rangle^\wedge_F\nn
&=\Big\{
D_{I_1I_2}j_{I_1}\partial_{I_2}
-D_{I_1I_3}F_{I_3I_2}j_{I_1}j_{I_2}\nn
&
+gV_{I_1}^{I_2I_3I_4}\Big[
 j_{I_1^\prime}\partial_{I_4}\partial_{I_3}\partial_{I_2}
-3F_{I_2I_2^\prime}j_{I_1}j_{I_2^\prime}
\partial_{I_3}\partial_{I_4}\nn
&
+(3F_{I_3I_3^\prime}F_{I_2I_2^\prime}
  +\frac{1}{4}A_{I_3I_3^\prime}A_{I_2I_2^\prime})
  j_{I_1}j_{I_3^\prime}j_{I_2^\prime}
  \partial_{I_4}\nn
&
+(3F_{I_3I_3^\prime}F_{I_2I_2^\prime}
  +\frac{1}{4}A_{I_3I_3^\prime}A_{I_2I_2^\prime})
  F_{I_1I_1^\prime}
  j_{I_1}j_{I_3^\prime}j_{I_2^\prime}j_{I_1^\prime}\Big]
\vert {\cal N}(j,a)\rangle^\wedge_F. 
\end{align} 
This equation is of greater complexity, but has to be taken as starting point 
for e.g. composite calculations. One finds the same terms as in the time-ordered
case, but also new terms constituting exchange and quantization terms. 

\subsection{Spinor quantum electrodynamics}

Dealing with a coupling theory is slightly more involved. While we can 
immediately apply our method also to the bosonic field operators $B_K$,
we have to reconsider the commutation relations for bosons and fermions
\begin{align}
i)  &&& \{\Psi_{I_1},\Psi_{I_2}\}^t_+:=A_{I_1I_2}=C\gamma_0\sigma^1
       \delta({\bf r}_1-{\bf r}_2)\nn
ii) &&& [B_K,\Psi_{I}]^t_-:=0\nn
iii)&&& [B_{K_1},B_{K_2}]_-^t=:C_{K_1K_2}.
\end{align}
The second equation states that the bosons are considered to be elementary and
are not functions of the fermionic fields. If we try to derive the Hamiltonian, 
we find two coupling terms. It was shown in \cite{stumpf:fauser:pfister:1993a}
that a consistency condition is required to ensure that they occur. It is a 
remarkable fact, that this condition arises from the fact that one has to
demand that the functional equations are independent on the ordering of
the bosons and fermions vice versa
\begin{align}
\vert {\cal P}(a,j,b)\rangle^\wedge &= \sum_{n,m=0}^{\infty}
\frac{i^n}{n!m!} \rho^s(I_1,\ldots,I_n,K_1,\ldots,K_m\mid a)\,\times\nn
&\hskip 0.75truecm \times\;j_{I_1}\wedge\ldots\wedge j_{I_n} \;
b_{K-1}\ldots b_{K_m} \mid0\rangle_F \nonumber \\[2ex]
\rho^1(I_1,\ldots,I_n,K_1,\ldots,K_m\mid a)&=
\langle 0\mid {\cal P}^f(\psi_{I_1},\ldots,\psi_{I_n})
{\cal P}^b(B_{K_1},\ldots,B_{K_m}) \mid a\rangle \nn
\rho^2(I_1,\ldots,I_n,K_1,\ldots,K_m\mid a)&=
\langle 0\mid {\cal P}^b(B_{K_1},\ldots,B_{K_m})
{\cal P}^f(\psi_{I_1},\ldots,\psi_{I_n}) \mid a\rangle 
\end{align}
We denoted the ordering by ${\cal P}$ which specialized to ${\cal P}^f$
for fermions and ${\cal P}^b$ for bosons.
The requirement that the functional equations for the $\rho^1$ hierarchy
is equivalent to that of the $\rho^2$ correlation function results in 
the following {\it reaction relation} which ensures that action and reaction
between fermions and bosons are mutually equal 
\begin{align}
C_{K_1K}W^K_{I_1I_2}\Psi_{I_2}&=2A_{I_1I}J^{II_2}_K\Psi_{I_2}.
\end{align}
It is remarkable, that from this equation one is able to {\it compute\/} 
the anticommutator $C_{K_1K_2}$ if the commutator of the fermions $A_{I_1I_2}$ 
is given and vice versa, see \cite{pfister:1987a,stumpf:fauser:pfister:1993a}.
Finally this relation can be used to eliminate one of the interaction terms
in the Hamiltonian in favour of the other, we choose
\begin{align}
H[\Psi ,B]^\wedge&=\frac{1}{2}A_{I_1I_3}D_{I_3I_2}\Psi_{I_1}\wedge\Psi_{I_2}
            +\frac{1}{2}A_{I_1I_3}W^K_{I_3I_2}B_K\Psi_{I_1}\wedge\Psi_{I_2}\nn
          &+\frac{1}{4}A_{I_1I_5}U_{I_5}^{I_2I_3I_4}
             \Psi_{I_1}\wedge\Psi_{I_2}\wedge\Psi_{I_3}\wedge\Psi_{I_4}\nn
          &+\frac{1}{2}C_{K_1K_3}L_{K_3K_2}B_{K_1}B_{K_2}.
\end{align}
The analogous transition to the functional equation yields for the 
time-ordered energy equation
\begin{align}
E_{0a}\vert{\cal T}(a,j,b)\rangle^\wedge
&=\big\{D_{I_1I_2}j_{I_1}\partial_{I_2}
+W^K_{I_1I_2}j_{I_1}\partial_{I_2}\partial^b_K
+L_{K_1K_2}b_{K_1}\partial^b_{K_2}\nn
&+U_{I_1}^{I_2I_3I_4}j_{I_1}(\partial_{I_2}\partial_{I_3}\partial_{I_4}
  -\frac{1}{4}A_{I_3I_3^\prime}A_{I_2I_2^\prime}j_{I_2^\prime}j_{I_3^\prime}
  j_{I_4})\nn
&+J_K^{I_1I_2}b_K(\partial_{I_1}\partial_{I_2}+\frac{1}{4}A_{I_1I_1^\prime}
  A_{I_2I_2^\prime}j_{I_1^\prime}j_{I_2^\prime})\big\}\vert{\cal T}
  (a,j,b)\rangle^\wedge \, .
\end{align}
If we are interested in the normal-ordered energy equation, where only the
fermions are normal-ordered, we have to add the propagator term in the 
Clifford map for the fermionic field operators and get as result
\begin{align}
E_{0a}\vert{\cal N}(a,j,b)\rangle^\wedge
=&\Big\{
 D_{I_1I_2}j_{I_1}\partial_{I_2}-D_{I_1I_3}F_{I_3I_2}j_{I_1}j_{I_2}\nn
&
+W^K_{I_1I_2}\big[
j_{I_1}\partial_{I_2}-F_{I_2I_2^\prime}j_{I_1}j_{I_2^\prime}\big]
\partial^b_K\nn
&
+J_k^{I_1I_2}b_K\big[
 \partial_{I_1}\partial_{I_2}
 -2F_{I_1I_1^\prime}j_{I_1^\prime}\partial_{I_2}\nn
&+(F_{I_1I_1^\prime}F_{I_2I_2^\prime}+\frac{1}{4}
    A_{I_1I_1^\prime}A_{I_2I_2^\prime})j_{I_1^\prime}j_{I_2^\prime}
    \big]\nn
&
+U_{J}^{I_1I_2I_3}j_{J}
\big[
\partial_{I_1}\partial_{I_2}\partial_{I_3}
-3F_{I_3I_4}j_{I_4}\partial_{I_2}\partial_{I_1}\nn
&
+(3F_{I_3I_4}F_{I_2I_5}+\frac{1}{4}
   A_{I_3I_4}A_{I_2I_5})j_{I_4}j_{I_5}\partial_{I_1}\nn
&
-(F_{I_3I_4}F_{I_2I_5}F_{I_1I_6}+\frac{1}{4}
  A_{I_3I_4}A_{I_2I_5}A_{I_1I_6})j_{I_4}j_{I_5}j_{I_6}\big]\nn
&
+L_{K_1K_2}b_{K_1}\partial^b_{K_2}
\Big\}\vert{\cal N}(a,j,b)\rangle^\wedge \,.
\end{align}
This equation can be used as a starting point to calculate positronium
bound states, see Ref. \cite{fauser:stumpf:1997a}.

As a general rule, we see from this calculations, that one can 
perform the following cliffordization process in Hopf gebraic terms
\begin{align}
{\cal P}(H[\psi,B]^\wedge) &= H[j,\partial,b,\delta]^{
\circ^{{\rm\sf p^f}}\circ^{{\rm\sf p^b}}}
\end{align}
where ${\rm\sf p}^f$ and ${\rm\sf p}^b$ are scalar valued ordering and
quantization maps inducing the cliffordization in the fermion and boson
sectors. This structure will be investigated elsewhere.

\subsection{Renormalized time-ordered products}

We have discussed already in the previous chapter the method, introduced by
Brouder \cite{brouder:2002a}, which allows to rewrite Epstein-Glaser
renormalization in Hopf algebraic terms. The Epstein-Glaser formalism
comes up with a renormalized time-ordered product in position
space, while the BPHZ renormalization, also employed by Connes and Kreimer 
\cite{kreimer:1998a,kreimer:2000a,connes:kreimer:2000a,connes:kreimer:2001a},
resides in momentum space. It was Pinter who established a clear and 
to Hopf algebras related formulation of the Epstein-Glaser theory 
\cite{pinter:2000a,pinter:2001a}. Finally Brouder realized that
this mechanism is a disguised cliffordization. In our formalism,
we have to add simply a new bilinear form $Z$ which introduces the
renormalization parameters. Since the reorderings, including the
renormalization, form a group under convolution, we can introduce
an operator ${\cal Z}$ and a linear form ${\rm\sf z}$ as done above
to introduce the renormalization. The renormalized time ordered functional
is then achieved by
\begin{align}
{\cal Z}\big(\left\vert {\cal T}(j,a)\right\rangle^\wedge \big)
&=
\left\vert {\cal Z}(j,a)\right\rangle^{\circ^{\rm\sf z}}
\,=\,
\left\vert {\cal Z}(j,a)\right\rangle^{\wedge}
\end{align}
where the last step introduces renormalization in the correlation 
functions by re-expanding the functional. The whole combinatorics of
this process is encoded in this singe and harmless looking equation!

The crucial point is to investigate what kind of physical reason is 
behind this additional reordering. There are two possibilities: 

(i) Since usually one reorders by the free propagator and not w.r.t. the exact 
propagator, there is a 
deficit in the ordering process which leads to singularities and has 
to be removed. In other words, one would expect to obtain no singularities at 
all if one would use the exact propagator of the theory.

(ii) We have studied generalized cliffordizations in the previous 
chapter. It might be possible, that renormalization is needed due to 
the fact, that the quantization and reordering process which is usually 
performed can come up only with exponentially generated bilinear forms,
e.g. we had the extension of $B$ into $B^\wedge$. Such bilinear forms
are related to theories which possess {\it only\/} two-particle 
interactions. If physics needs {\it inevitably\/} non-exponentially 
generated bilinear forms, such a contribution can be introduced by 
renormalization and the $Z$-pairing. 

These possibilities will be studied elsewhere.

\section{On the vacuum structure}

While the preceding sections dealt with realistic QF theories, we will
discuss the peculiarities occuring from the vacuum structure in a $U(1)$-
and $U(2)$-model. This will allow to be very explicite while being not
bothered with complications of a realistic theory. But, already the $U(2)$- 
model, if it is considered as describing a fiber on the space of modes,
is a realistic model of BCS superconductivity and provides even 
generalizations. A detailed exposition including the relation to 
an analogous $C^*$-algebraic treatment can be found in Ref. 
\cite{fauser:2001e}. It was in fact this work which initiated 
the study of time- and normal-ordering and generally QFT in 
Clifford Hopf gebraic terms. 

\subsection{One particle Fermi oscillator, $U(1)$}

In this section we study the simplest possible model, which consists
of a single fermionic particle. We are interested in the degrees
of freedom of the fiber only, so we suppress a momentum index, which 
could however be added without altering our consideration. The CAR
algebra of a single fermion is created by two generators $\{a,a^\dagger\}$
which we denote also by $\{e_1,e_2\}$ in the index doubled formulation, 
i.e. the index describes the adjointness of the operator. The CAR relations
read
\begin{align}
&&\{ a,a^\dagger\}_+ = \Id&&\text{others zero.}&&
\end{align}
The adjoint map is the algebra antihomomorphism which interchanges
$a$ and $a^\dagger$. This algebra can be turned into a $C^*$-algebra.

Reformulating the CAR relation in Clifford Hopf gebraic terms
does not allow to fix the bilinear form $B$, but only its symmetric
part, which encodes thereby the quantization. We introduce therefore
a parameter $\nu$, which represents the antisymmetric part, in the 
following way
\begin{align}
[B_\nu] &=
   \begin{array}{|cc|} 0 & \nu \\  1-\nu & 0 \end{array}
=  \begin{array}{|cc|} 0 & \frac{1}{2} \\  \frac{1}{2} & 0 \end{array}
  +\begin{array}{|cc|} 0 & -\frac{1}{2}+\nu \\ 
                      \frac{1}{2}-\nu & 0 \end{array} \nn
&= [g]+[F_{\nu}].
\end{align}
This form is chosen for convenience to be able to make contact to 
$C^*$-algebraic calculations done by Kerschner \cite{kerschner:1994a}
and ourselves \cite{fauser:2000e}. The cliffordization is performed w.r.t. 
this bilinear form. We denote the Clifford product by juxtaposition and the
contraction is given as $\JJ_v$ to indicate its dependence on the parameter. 
The underlying space of the Gra{\ss}mann algebra is based on the wedge 
$\wedge$ product.

It is an easy task to recompute the CAR relations in Clifford terms
\begin{align}
e_i e_j + e_j e_i &= B_\nu(e_i,e_j) + B_\nu(e_j,e_i) \nn
&= 2\,g(e_i,e_j) = \delta_{i,n+1-j}
\end{align}
showing that only the symmetric part of $B_\nu$ enters the 
quantization. 

We want to investigate the meaning of the parameter $\nu$. For this 
purpose, we introduce the counit $\epsilon^\wedge$ and compute the 
{\it `vacuum' expectation values\/} of the elements in a Clifford basis. 
Remember, that the operator product is given by the Clifford product 
and that physicists do commonly write down only expressions using this 
product. We get
\begin{align}
\epsilon^\wedge(\Id) &= 1\qquad\text{normalization} \nn
\epsilon^\wedge(a\,a^\dagger) &= 
\epsilon^\wedge(e_1\circ e_2) = 
\epsilon^\wedge(e_1\wedge e_2 + B_{12}) = \nu \nn
\epsilon^\wedge(a) &= 0\nn
\epsilon^\wedge(a^\dagger) &= 0 .
\end{align}
If we require that our state is positive, we get from the above result
and form $\epsilon^\wedge(a^\dagger\,a)=1-\nu$ the condition 
$0\le \nu\le 1$ for $\nu$ or equivalently $\det(B_\nu)<0$. The convex set of 
positive, normalized, linear functionals on the CAR algebra is thus
parameterized by $\nu\in [0,1]$. 

The reader should note the difference in the description, while the
usual treatment comes up with a variety of states acting on a fixed
operator algebra, we have a unique state, the {\it counit\/} and 
parameterize the operator product by adding the antisymmetric part
$F_\nu$.

In physics, one introduces a Fock vacuum by the following requirement
\begin{align}
a\vert0\rangle_{\cal F} &= 0 .
\end{align} 
This, however does fix the value $\nu$ immediately! One finds
$0=\epsilon^\wedge(a^\dagger\, a)=\epsilon^\wedge(1-a\,a^\dagger)
= 1-\nu$ and hence $\nu=1$. A basis of the algebra under this condition
is given by the Fock space basis
\begin{align}
\{\,\vert0\rangle_{\cal F}\,,\, a^\dagger\,\vert0\rangle_{\cal F}\,\} 
\end{align}
which is two dimensional, and in fact a spinor representation. However,
our treatment is totally arbitrary w.r.t. the name of the operators, and
we could have introduced a dual Fock space demanding that
\begin{align}
a^\dagger\vert0\rangle_{{\cal F}^*} &= 0
\end{align} 
which would have resulted in the basis
\begin{align}
\{\,\vert0\rangle_{{\cal F}^*}\,,\, a\,\vert0\rangle_{{\cal F}^*}\,\} \,,
\end{align}
the span of which we call dual Fock space. It can be shown that this setting 
corresponds to the parameter $\nu=0$. What happens for $\nu\in ]0,1[$ ?

While we found two dimensional representations for $\nu=1$ and $\nu=0$,
we get a 4-dimensional representation in the general case, rendering the
algebra to be indecomposable. In other words, 
\begin{align}
a\,a^\dagger \quad\text{and}\quad 
a^\dagger\,a
\end{align}
are almost idempotents if and only if $\nu=1$ and $\nu=0$. States with
$\nu\in[0,1[$ can be described as linear combinations of this two states
and come up to be mixed states. It can be shown, that the time-ordered 
case is obtained if $\nu=1/2$, in which case the antisymmetric part
$F_v$ of our bilinear form is not present. Renormalization does not 
make any sense in this almost to trivial example. Since we can come up
with a particle number operator which acts on the operators
\begin{align}
{}[N,a]_- &= -a \nn
{}[N,a^\dagger]_- &= \phantom{-}a^\dagger
\end{align}
we call this a $U(1)$-model. It turns out, that $N$ depends on $\nu$
in the following way
\begin{align}
N &= (v-\frac{1}{2})-e_1\wedge e_2 \,=\, \nu\,\Id +a\,a^\dagger .
\end{align}
This is a Lie group generator only if $\nu = 1/2$, otherwise one deals with a
central extended Lie group.

A detailed study of families of idempotents parameterized by a parameter
$\nu$ will be given elsewhere.

\subsection{Two particle Fermi oscillator, $U(2)$}

While the one particle case is not very interesting, we gain
a resonable interesting model already in the next dimension, having
two particles and hence four creation and annihilation operators.
We have the CAR relations ($\alpha,\beta \in (1,2)$)
\begin{align}
&\{a_\alpha, a_\beta \} = 0 = \{a_\alpha^\dagger,
a_\beta^\dagger \} \nn
&\{a_\alpha, a_\beta^\dagger \}  = \delta_{\alpha,\beta} \Id
\end{align}
which we will encode once more by index doubling as $\{a_\alpha, a_\beta,
a_\alpha^\dagger, a_\beta^\dagger \} = \{e_1,e_2,e_3,e_4\}$. While in 
the $U(1)$-model we had only a single operator at our disposal, we can 
implement in the 2-dimensional case a $U(2)$ action. Let $N,S_k$ be the 
generators of the Lie group $U(2)$, we define 
\begin{align}
& [N,S_k] = 0, \quad [S_k,S_l] = i\epsilon_{klm}S_m &\nn
& S^\dagger = S,\quad N^\dagger = N & \nn
& [S_k, a_\alpha] = \sigma_k^{\alpha \beta} a_\beta , \quad
 [S_k, a_\alpha^\dagger] = \hat{\sigma}_k^{\alpha \beta}
a_\beta^\dagger & \nn
& [N,a_\alpha] = + a_\alpha,\quad
  [N,a_\alpha^\dagger] = - a_\alpha^\dagger \,. &
\end{align}
This are the defining relation of the $U(2)$ generators, two reality
conditions and finally their action on the CAR generators. 
The relations are not independent. A basis
using operator products(!) can be given by looking for the eigen 
states of $S_3,\sum S_kS^k$ and $N$, which we denote as $(s_3,s(s+1),q)$.

\begin{table}[t]
{\hfill
\begin{tabular}[t]{|c|c||c|c|c|}
\hline\hline
\multicolumn{5}{c}{\bf Table 1.} \\
\hline\hline
No: & $A \in CAR$ & $s_3$ & $s(s+1)$ & $q$ \\
\hline\hline
$g_1$   & $\Id$ & 0 & 0 & 0 \\
$g_2$   & $\frac{1}{2}(a_1a_1^\dagger+a_2 a_2^\dagger)$ & 0 & 0 & 0 \\
$g_3$   & $a_1 a_2 a_2^\dagger a_1^\dagger$ & 0 & 0 & 0 \\
\hline
$g_4$   & $a_1 a_2^\dagger$ & 1 & 2 & 0 \\
$g_5$   & $\frac{1}{2}(a_1a_1^\dagger-a_2 a_2^\dagger)$ & 0 & 2 & 0 \\
$g_6$   & $a_2 a_1^\dagger$ & -1 & 2 & 0 \\
\hline
$g_7$   & $a_1$ & $\frac{1}{2}$ & $\frac{3}{4}$ & 1 \\
$g_8$   & $a_1 a_2 a_2^\dagger$ & $\frac{1}{2}$ & $\frac{3}{4}$ & 1 \\
$g_9$   & $a_2$ & $-\frac{1}{2}$ & $\frac{3}{4}$ & 1 \\
$g_{10}$& $a_2 a_1 a_1^\dagger$ & $-\frac{1}{2}$ & $\frac{3}{4}$ & 1 \\
\hline
$g_{11}$& $a_2^\dagger$ & $\frac{1}{2}$ & $\frac{3}{4}$ & -1 \\
$g_{12}$& $a_1 a_1^\dagger a_2^\dagger$ & $\frac{1}{2}$ & $\frac{3}{4}$ & -1 \\
$g_{13}$& $a_1^\dagger$ & $-\frac{1}{2}$ & $\frac{3}{4}$ & -1 \\
$g_{14}$& $a_2 a_2^\dagger a_1^\dagger$ & $-\frac{1}{2}$ & $\frac{3}{4}$ & -1 \\
\hline
$g_{15}$& $a_1 a_2$ & 0 & 0 & 2 \\
$g_{16}$& $a_1^\dagger a_2^\dagger$ & 0 & 0 & -2 \\
\hline
\end{tabular}\hfill}
\caption{\label{Tabelle 1 in Kap. 8}
Eigenvectors of the $U(2)$ and their $U(2)$ quantum numbers. 
Operator products are Clifford products.}
\end{table} 

We ask for a linear form $\omega_{\nu w}$ which is positive and normalized.
Such a linear form is characterized by its action on a basis and we have 
to compute the expectation values w.r.t. this form for all 16 states as 
given in table \ref{Tabelle 1 in Kap. 8}. However, we are interested in 
such states only which are invariant under the $U(2)$ action. Under this 
requirement we find for the non-zero expectation values
\begin{align}
\omega_{\nu w}(\Id) &= 1 \nn
\omega_{\nu w}(a_1\,a_1^\dagger) &= \nu \nn
\omega_{\nu w}(a_2\,a_2^\dagger) &= \nu \nn
\omega_{\nu w}(a_1\,a_2\,a_2^\dagger\,a_1^\dagger) &= w .
\end{align}
This is the result of Kerschner \cite{kerschner:1994a}, which he obtained 
by $C^*$-algebraic considerations. However, note that the above basis is
not antisymmetric and the operator products cannot be seen as Wick 
monomials or correlation functions. This fact will lead below to a 
{\it renormalization\/} of the above displayed expectation values.

Let us introduce a bilinear form which allows us to cliffordize the
Gra{\ss}mann algebra over the 4 generators $\{e_i\}$ in such a manner,
that the CAR relations hold. The most general form is
\begin{align}
{}[B(e_i,e_j)] &= 
\begin{array}{|cccc|}
  0 & u & q & r   \\
 -u & 0 & s & t   \\
 -q & 1-s & 0 & m \\
1-r & -t & -m & 0
\end{array} \, .
\end{align}
This yields a quantum Clifford algebra and our `vacuum' state is the
counit w.r.t. the wedge product $\epsilon^\wedge$. We want to express
the $U(2)$ generators in terms of the generators. From the CAR 
relations one obtains that
\begin{align}
Q^\prime &:= a_1 a_1^\dagger + a_2 a_2^\dagger.
\end{align}
However, this operator has a non-vanishing expectation value and we have 
to {\it renormalize\/} it. This reads
\begin{align}
Q&= a_1 a_1^\dagger + a_2 a_2^\dagger -(r+s)\,\Id 
  = a_1\wedge a_1^\dagger + a_2 \wedge a_2^\dagger.
\end{align}
From this display we see that the operator $Q$ has to be defined in 
the wedge basis. The same applies for the basis vectors in our above
given table and the other $U(2)$ generators. We find
\begin{align}
&g_4^\prime = g_4 -q,\hskip 2truecm g_5^\prime = g_5 -\frac{1}{2}(r-s)& \nn
&g_6^\prime = g_6 -t,\hskip 2truecm g_{15}^\prime = g_{15} -u& \nn
&g_{16}^\prime = g_{16} +m\,.&
\end{align}
After this renormalization we can {\it derive\/} the parameters 
$(\nu, w)$ of 
the `vacuum' state $\omega_{\nu w}$ from the data given by the bilinear 
form that is from quantization and from the `propagator' which enter 
the cliffordization process. While the symmetric 
part is obtained due to canonical quantization, the antisymmetric part 
was in QFT related to the propagator. In our $U(2)$ example, the 
propagator is given as
\begin{align}
{} [ F_{ij} ] &= 
   [ \langle \frac{1}{2}[a_i , a_j^\dagger] \rangle_0^\wedge ] 
=\begin{array}{|cc|} -r +1/2 & q \\ 
                      t & -s+1/2 
\end{array} \, .
\end{align}
This shows, that $\nu$ and $w$ are functions of the parameters $r,s,q,t$. 
This model was discussed in great detail, in Ref. \cite{fauser:2001e}, 
including its direct link to BCS theory, the gap equation etc. The reader
is invited to consult the original source for this details. 

\begin{figure}
\centering
\includegraphics[height=0.75\textwidth]{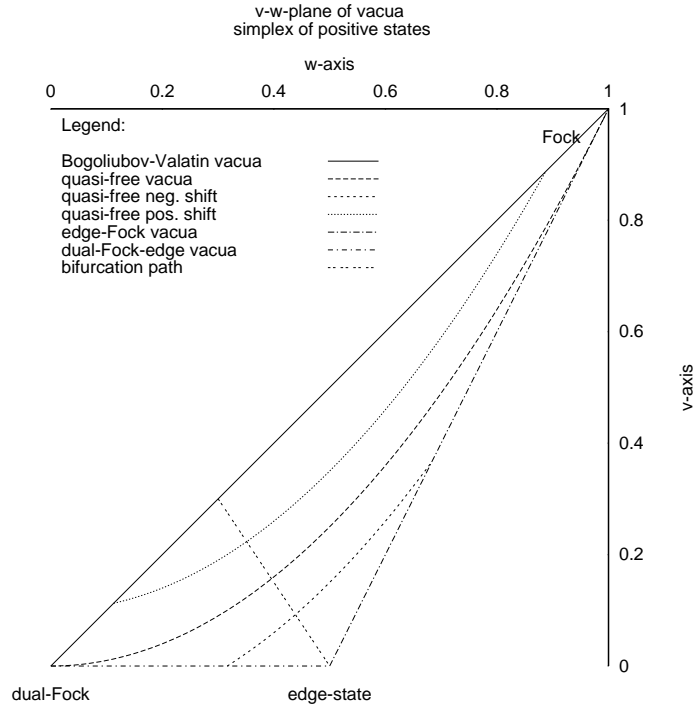}
\caption{Plane of vacuum states for an $U(2)$-model.}
\label{fig1}
\end{figure} 

To complete our discussion we want to study the vacuum structure of
the present $U(2)$-model a little further. First of all, we consider 
under which condition the state $\omega_{\nu w}$ is positive. An analogous
consideration as above yields
\begin{align}
&0 < w < \nu < 1& \nn
&2\nu-1 < w \, . &
\end{align}
We draw a diagram, see Figure \ref{fig1}, 
where every point in the affine Euclidean plane 
corresponds to a state $\omega_{\nu w}$. The positive states form
a triangle. We want to discuss the states in and on the borders of the
triangle.

Let us impose the Fock and dual Fock space conditions and
see which point in the plane corresponds to it 
\begin{align}
&&\prod_{i\in I} a_i \mid 0 >_{\cal F} \,=\, 0 &&& \text{Fock space}&&\nn
&&\prod_{i\in I} a_i^\dagger \mid 0 >_{{\cal F}^*} \,=\, 0  &&&
\text{dual Fock space.} &&
\end{align}
We obtain e.g. for the Fock space condition the following expectation values
\begin{align}
0 &= < a_i^\dagger a_i >_{\cal F} \, =\, 
<\Id - a_i a_i^\dagger >_{\cal F} \, =\,
1- \left\{
\begin{array}{l} 
r\quad i=1 \\
s\quad i=2\, .
\end{array} \right.
\end{align}
from which we deduce $\nu=1$. Furthermore we find
\begin{align}
0 &= < a_2^\dagger a_1^\dagger a_1 a_2>_{\cal F} \,=\,
      <a_2^\dagger a_2 - \Id + a_1^\dagger a_1 + 
       a_1 a_2 a_2^\dagger a_1^\dagger >_{\cal F} \nn
  &= 0-1+0+w,
\end{align}
which yields $w=1$. An analogous computation can be done for the 
dual Fock case, and we get
\begin{align}
\omega_{\cal F} &=
 \omega_{\nu w}\Big\vert_{\nu=1\, w=1} \,=\, \omega_{11},
\hskip 2truecm
\omega_{{\cal F}^*} = 
   \omega_{\nu w}\Big\vert_{\nu=0\, w=0} \,=\, \omega_{00}.
\end{align}
Hence we can identify the up-right and down-left edges of the trinagle
of sets to be the Fock and dual Fock state. The representation space 
is in both cases 4-dimensional and reads as follows
\begin{align}
{\cal V}_{\cal F} &= \{ \Id \mid 0 >_{\cal F},
                  a_1^\dagger \mid 0 >_{\cal F},
                  a_2^\dagger \mid 0 >_{\cal F},
                a_1^\dagger a_2^\dagger \mid 0 >_{\cal F} \} \nn
{\cal V}_{{\cal F}^*} &= \{ \Id \mid 0 >_{{\cal F}^*},
                  a_1 \mid 0 >_{{\cal F}^*},
                  a_2 \mid 0 >_{{\cal F}^*},
                a_1 a_2 \mid 0 >_{{\cal F}^*} \}.
\end{align}
The line which connects these two edges can be reached by 
Bogoliubov-Valatin transformations. These states are usually employed
in BCS theory for condensates.

Quasi free states are defined to have no higher correlations, i.e. there
exists a transformation into a free theory, see 
\cite{bratelli:robinson:1979a,bratelli:robinson:1981a}. We can ask, 
which states in our plane do not possess higher correlations 
($\kappa_n=0$, for all $n>1$, $\kappa_n$ is defined below). Hence we 
have to assert that 
\begin{align}
\kappa_1(a_\alpha a^\dagger_\beta) 
 &= \omega_{\nu w}(a_\alpha a^\dagger_\beta)\,=\, \nu \nn
0 \,=\,\kappa_2(a_{\alpha_1}a_{\alpha_2} a^\dagger_{\beta_1}a^\dagger_{\beta_2})
 &= \omega_{\nu w}(a_{\alpha_1}a_{\alpha_2} 
                   a^\dagger_{\beta_1}a^\dagger_{\beta_2}) \nn
    &\phantom{=}
     +\omega_{\nu w}(a_{\alpha_1} a^\dagger_{\beta_1}) 
      \omega_{\nu w}(a_{\alpha_2} a^\dagger_{\beta_2}) \nn
    &\phantom{=}
     -\omega_{\nu w}(a_{\alpha_1} a^\dagger_{\beta_2}) 
      \omega_{\nu w}(a_{\alpha_2} a^\dagger_{\beta_1}) \nn
 &= w-\nu^2
\end{align}
holds. From $\kappa_2=0$ we find a parabola in our diagram, which connects
the Fock and dual Fock states and shows that these states are quasi free too.
Having no higher correlations means that there is no interaction, hence
these states build a border between regions having interactions of possibly
different type, e.g. attracting or rejecting. Since we know that the line
which connects Fock and dual Fock states is related to BCS theory, and since
one has a condensate due to an attractive interaction, we may address
the area between the parabola of quasi free states and the line of 
Bogoliubov-Valatin states as the condensate area.

Since every positive state can be written as a convex combination
of extremal states, it remains to discuss the third edge of the triangle, 
which we
call with Kerschner `edge`-state and denote it by $\omega_{\cal E}$. 
We know that this state is at the position $\nu=1/2$, $w=0$. It is 
easy to see that this condition leads to the following 8-dimensional
space
\begin{align}
\omega_{\cal E} &= \omega_{\nu w}\Big\vert_{\nu=1/2\, w=0} 
\,=\, \omega_{1/2\, 0} \nn
{\cal V}_{\cal E} &= \{
\mid 0 >_{\cal E}\,, 
a_1\,\mid 0 >_{\cal E}\,, 
a_2\mid 0 >_{\cal E}\,, 
a_1^\dagger\,\mid 0 >_{\cal E}, 
a_2^\dagger\,\mid 0 >_{\cal E}\,, \nn
&\phantom{=}\; 
\frac{1}{2}(a_1 a_1^\dagger - a_2 a_2^\dagger )\mid 0 >_{\cal E}\,,
a_1 a_2^\dagger\,\mid 0 >_{\cal E}, 
a_2 a_1^\dagger\, \mid 0 >_{\cal E} \,\}.
\end{align}
It is remarkable that in this set a spin triplet occurs which is not 
present in the Fock or dual Fock space. Moreover, we find spin up and
down particles and antiparticles (annihilators w.r.t. the Fock vacuum!).

If one derives a gap-equation, see Ref. \cite{fauser:2001e}, one
notes that the discriminant is negative for states in the area between 
the parabola of quasi free states and the edge-state, which disallows
two solutions. On the other hand, if one looks at states between the
quasi free parabola and the Bogoliubov-Valatin states (left border line)
one has two solutions and a gap. This gap can be related to the
common energy gab of BCS theory. 

Having discussed roughly the vacuum states which arise from 
$\epsilon^\wedge$ by Hopf algebraic means, especially by cliffordization, 
we close this comprehensive treatise.

However we want to remark that this is only the starting point into 
a new and exciting field, which we await to be fruitful for studies
in various directions. Hopf gebras will help us to understand
what quantization means geometrically, a new approach to renormalization
is opened, the vacuum/state space structures of a theory can be explored,
dynamics is related to states directly, which will have interesting 
consequences, and many more. We await to enter hopfish times and quantum 
Clifford algebras will play a major role.

%% file: cliff_big.tex
\chapter{CLIFFORD and BIGEBRA packages for Maple}

\section{Computer algebra and Mathematical physics}

Computer algebra was a major tool to investigate the topics which have 
been presented in this work. We had the opportunity to state even some
theorems which we proved in low dimensions by direct calculations. Of 
course, the strength of a Computer Algebra System (CAS) is not to prove
general theorems, but to provide a general area to explore mathematics 
and physics in an {\it experimental\/} way. Moreover, a CAS can help to
surmount difficulties which would not be tractable at the moment by 
analytical, algebraical or arithmetical methods. E.g. when we computed 
the antipode of a two dimensional Clifford bi-convolution algebra this took
some hours of computing time on a present day state of the art computer
with lots of RAM. Only after the solution is found, it is an easy task to 
check by hand, so not relying on the computer any more, that this is 
indeed the searched antipode. A much wider area is opened by the 
possibility to check own and other people's assertions and claims simply 
by evaluating them in special cases. While this cannot lead to a proof, 
many such assertions can be disproved. This leads at the end to a 
refinement of their formulations and eventually to an idea how to prove
such mathematical assertions by generalizing the generic case. Also in
this work, we had the opportunity to find out many shortcomings of 
statements found in the literature. As a prominent example may be recalled
the distinction between interacting and non-interacting, i.e. connected and 
non-connected, Hopf gebras. A simple re-calculation of standard material
led to the fact that a Clifford Hopf gebra cannot be connected which 
stems from the non-locality of the cliffordization. Seeing the problem
was essential to come up with a solution.

We want to summarize the cutting edge points which were valuable to
the present research and which will become for sure a common tool 
in research in future times.

\begin{itemize}
\item	{\bf Check Assertions:}
If one has a prejudice that some assertion should be true in an algebraic
setting, randomly chosen special cases can give confidence into such
a belief. More boldly, a single counterexample can put down the whole 
business immediately. This might look distracting but saves a tremendous 
amount of work, since only such assertions remain for being proved which
are already tested to some amount and have a particular chance to generalize
to a theorem.
\item	{\bf Computations:}
It should not be underestimated how time consuming it is to evaluate
lengthy computations. While the CAS cannot substitute a sever knowledge
of the mathematics behind and a  sound physical concept to work on, it can
help to compute with much fewer errors than any calculation by hand can
provide. Moreover, using a CAS one can reach areas which are un-tractable
by hand-written calculations simply by its mere length.
\item 	{\bf Develop new Mathematics:}
Since new mathematical tools are not shipped with a CAS, one has to
develop ones own functionality as an add to the common features of 
such a system. E.g. Maple \cite{MAPLE}\footnote{Maple is a registered 
trademark of Maple Waterloo Software, 
see {\tt\small http://www.maplesoft.com/}}
comes already with a tremendous ability to deal with many parts of algebra,
but it was not able to deal with Gra{\ss}mann and Clifford algebras and 
Hopf gebras. The development of such a device was a major impulse to 
investigate the mathematical structure in great depth. In fact, if you can
teach the mathematics to a computer you have really understood the case.
\item	{\bf Experimental Mathematics:}
Having the opportunity to deal with a CAS opens the field of experimental 
mathematics. This includes partly the other topics of this list, but should
not underestimated in its own dynamics. Exploring mathematics by doing 
particular experiments justifying or deceasing {\it own\/} assumptions
is of extreme value to be able to enter a field fast and in a secure 
and solid way. This leads immediately to the next item.
\item	{\bf Teaching:}
{\it Experimental Mathematics\/} may be regarded as an additional tool in
teaching complicated mathematics. Students can see what type of behaviour 
some algebraic or physical structures have before the try to understand 
or perform on their own a proof to master finally the topic. The CAS
enables dealing in a concrete way with mathematical structures. Visualisation,
erasing of miss-conceptions, and allowing a neat approach to complicated
technicalities have already boosted up the field of non-linear dynamics.
This field enjoyed a renaissance after the advent of sufficiently fast 
computers to handle the numerics. However, CAS is much more valuable 
since it really develops the algebraic understanding of the mathematical
subject. 
\end{itemize} 

The particular CAS we use here is Maple V rel 5.1. Perhaps any reasonable 
general such tool could be employed. However, the already existing package 
CLIFFORD, developed by Rafa{\l} Ab{\l}amowicz \cite{CLIFFORD}, which I had
enjoyed to use for now a couple of years, was reason enough for this choice.

In the next section we will give some hints how CLIFFORD can be used for
computations in Clifford algebra. However, since there is a valuable and 
well developed online help consisting of approx. 150 help-pages, we stay with 
those features which were actually used in this work and which were essential
for the development and design of the BIGEBRA package. The latest version
of CLIFFORD is Cliff5 (i.e. version 5). CLIFFORD will be developed jointly
in future with Rafa{\l} Ab{\l}amowicz. 

The section on the BIGEBRA package will describe in a very cursory way 
the essential features which have been used to establish the assertions
and theorems stated in this works. Some proofs have been by ``direct 
computation using CLIFFORD/BIGEBRA'' and we feel responsible to exemplify
the abilities of CLIFFORD/BIGEBRA to give some hints how this was 
established. Full confidence can however be obtained only by looking at 
the particular, sometimes long-winding, Maple worksheets containing the actual
computations. BIGEBRA was developed in close cooperation jointly with Rafa{\l}
Ab{\l}amowicz.

\section{The CLIFFORD Package -- rudiments of version 5}

The CLIFFORD package was developed by Rafa{\l} Ab{\l}amowicz since 1996. It
is available from his web-server at {\tt http://math.tntech.edu/rafal/}.
From version 5 onwards the package comes together with the additional BIGEBRA
package and is developed jointly with the author. Since there is an extensive
online documentation, included into the Maple online help system, with 
help-page for every function we give only a look-and-feel description of
those functions which are needed later in the BIGEBRA examples.

To load the CLIFFORD package we simply type in the following command: 
\begin{Verbatim}
restart:with(Cliff5):
\end{Verbatim}
This has loaded the package and offers now to perform calculations in
Gra{\ss}mann and Clifford algebras. First of all, let us show how to
select a Clifford algebra and how to assemble a basis, particular, and 
general elements. Such elements will be called Clifford or Gra{\ss}mann 
polynoms, Clifford or Gra{\ss}mann monoms with or without a scalar pre-factor.
We compute over general algebraic expressions dealing thus with Clifford or
Gra{\ss}mann modules. A colon suppresses the output of the command, while
a semi-colon ends a statement and returns its output. The generators of the
algebras are denoted as {\tt e1,e2,e3,$\ldots$,ea,eb,$\ldots$}.
\begin{Verbatim}
dim_V:=2:                 ## set dim. of generating space
B:=linalg[diag](1$dim_V); ## diagonal Euclidean metric; $ short for seq.
\end{Verbatim}
\relax
\[
B :=  \left[ 
{\begin{array}{rr}
1 & 0 \\
0 & 1
\end{array}}
 \right] 
\]
\relax
\begin{Verbatim}
bas:=cbasis(dim_V);       ## get a basis spanning the Algebra
\end{Verbatim}
\relax
\[
\mathit{bas} := [\mathit{Id}, \,\mathit{e1}, \,\mathit{e2}, \,
\mathit{e1we2}]
\]
\relax
\begin{Verbatim}
p1:=e1we2;                ## notion for 'e1 wedge e2'
\end{Verbatim}
\relax
\[
\mathit{p1} := \mathit{e1we2}
\]
\relax
\begin{Verbatim}
p2:=a*e1+b*e1we2-4*Id;    ## a Grassmann polynom, Id is the unit
\end{Verbatim}
\relax
\[
\mathit{p2} := a\,\mathit{e1} + b\,\mathit{e1we2} - 4\,\mathit{Id
}
\]
\relax
\begin{Verbatim}
p3:=x*eaweb+ec;           ## a Grassmann polynom with symbolic indices	
\end{Verbatim}
\relax
\[
\mathit{p3} := x\,\mathit{eaweb} + \mathit{ec}
\]
\relax
\begin{Verbatim}
X:=add(_X[i]*bas[i],i=1..2^dim_V); ## a general element 
\end{Verbatim}
\relax
\[
X := {\mathit{\_X}_{1}}\,\mathit{Id} + {\mathit{\_X}_{2}}\,
\mathit{e1} + {\mathit{\_X}_{3}}\,\mathit{e2} + {\mathit{\_X}_{4}
}\,\mathit{e1we2}
\]
\relax
Since the wedge product $\wedge$ was already used internally for building 
the Gra{\ss}mann basis, we start by exemplifying the usage of the wedge 
product.
\begin{Verbatim}
wedge(e1,e2);             ## wedge of e1 and e2
\end{Verbatim}
\relax
\[
\mathit{e1we2}
\]
\relax
\begin{Verbatim}
&w(e1,e2);                ## short form for wedge
\end{Verbatim}
\relax
\[
\mathit{e1we2}
\]
\relax
\begin{Verbatim}
e1 &w e2;                 ## infix form for wedge
\end{Verbatim}
\relax
\[
\mathit{e1we2}
\]
\relax
\begin{Verbatim}
&w(p1,p2);                ## wedge on particular elements
\end{Verbatim}
\relax
\[
 - 4\,\mathit{e1we2}
\]
\relax
\begin{Verbatim}
&w(X,X);                  ## square of a general element 	
\end{Verbatim}
\relax
\[
{\mathit{\_X}_{1}}^{2}\,\mathit{Id} + 2\,{\mathit{\_X}_{2}}\,{
\mathit{\_X}_{1}}\,\mathit{e1} + 2\,{\mathit{\_X}_{3}}\,{\mathit{
\_X}_{1}}\,\mathit{e2} + 2\,{\mathit{\_X}_{4}}\,{\mathit{\_X}_{1}
}\,\mathit{e1we2}
\]
\relax
Given the Gra{\ss}mann algebra as above, we have also contractions at 
our disposal. The contractions act w.r.t. the chosen bilinear form $B$,
which could also be symbolic or unassigned at all. The (left) contraction
acts as a graded derivation on the module generated by the above given 
basis. It also established the bilinear form. To manipulate Gra{\ss}mann
basis elements we need also a device to put them into a standard order, i.e. 
the function `{\tt reorder}' and a function which constitutes the grading, 
i.e. `{\tt gradeinv}'. The eigenspace of {\tt gradeinv} are exactly
the even and odd elements.  
\begin{Verbatim} 
map(gradeinv,bas);               ## map means `apply to the list`
\end{Verbatim}
\relax 
\[
[\mathit{Id}, \, - \mathit{e1}, \, - \mathit{e2}, \,\mathit{e1we2
}]
\]
\relax
\begin{Verbatim} 
map(i->1/2*(i+gradeinv(i)),bas); ## even elements
\end{Verbatim}
\relax 
\[
[\mathit{Id}, \,0, \,0, \,\mathit{e1we2}]
\]
\relax
\begin{Verbatim} 
map(i->1/2*(i-gradeinv(i)),bas); ##  odd elements
\end{Verbatim}
\relax 
\[
[0, \,\mathit{e1}, \,\mathit{e2}, \,0]
\]
\relax
\begin{Verbatim} 
linalg[matrix](dim_V,dim_V,(i,j)->LC(e.i,e.j));## contraction on vectors
\end{Verbatim}
\relax 
\[
 \left[ 
{\begin{array}{cc}
\mathit{Id} & 0 \\
0 & \mathit{Id}
\end{array}}
 \right] 
\]
\relax
\begin{Verbatim} 
## derivation property, LC taken w.r.t. the unassigned bilinear form `K'
LC(e1,e1we2,K)=LC(e1,e1,K) &w e2+gradeinv(e1)*LC(e1,e2,K); 
\end{Verbatim}
\relax 
\[
{K_{1, \,1}}\,\mathit{e2} - {K_{1, \,2}}\,\mathit{e1}={K_{1, \,1}
}\,\mathit{e2} - \mathit{e1}\,{K_{1, \,2}}\,\mathit{Id}
\]
\relax
It is well know that the Clifford product of a 1-vector can be established
as an endomorphism on the Gra{\ss}mann basis underlying the Clifford 
algebra. Such a particular endomorphism is called a {\em Clifford map}. 
The Clifford product in CLIFFORD ver. 5 is however based on the Hopf 
algebraic process of Cliffordization.
\begin{Verbatim}
CliMap:=proc(x,u,B) LC(x,u,B)+wedge(x,u) end: ## the Clifford map
CliMap(e1,Id,B);    ## contraction part is zero
\end{Verbatim}
\relax
\[
\mathit{e1}
\]
\relax
\begin{Verbatim}
CliMap(e1,e1,K);    ## wedge part is zero
\end{Verbatim}
\relax
\[
{K_{1, \,1}}\,\mathit{Id}
\]
\relax
\begin{Verbatim}
CliMap(e1,e2,K);    ## Clifford product w.r.t. the bilinear form `K'
\end{Verbatim}
\relax
\[
{K_{1, \,2}}\,\mathit{Id} + \mathit{e1we2}
\]
\relax
\begin{Verbatim}
CliMap(e2,e1we2,B); ## action on a bi-vector
\end{Verbatim}
\relax
\[
 - \mathit{e1}
\]
\relax
\begin{Verbatim}
cmul(e2,e1we2);     ## compare with the builtin Clifford product
\end{Verbatim}
\relax
\[
 - \mathit{e1}
\]
\relax
Of course, the Clifford product has to be extended to a general first 
argument. This can be done by using the rules given in the main text.
Since more features of CLIFFORD are explained in the following section 
which describes the BIGEBRA package, we end by exemplifying the 
{\tt clisolve} facility. This function allows to solve equations
in Gra{\ss}mann and Clifford algebras either for particular elements
and their coefficients or for arbitrary elements. We will show how to 
find idempotents. Remember that we had defined an arbitrary element $X$.
\begin{Verbatim}
X; ## general element; dim_V = 2
\end{Verbatim}
\relax
\[
{\mathit{\_X}_{1}}\,\mathit{Id} + {\mathit{\_X}_{2}}\,\mathit{e1}
 + {\mathit{\_X}_{3}}\,\mathit{e2} + {\mathit{\_X}_{4}}\,\mathit{
e1we2}
\]
\relax
\begin{Verbatim}
sol:=map(allvalues,clisolve(cmul(X,X)-X,X));
\end{Verbatim}
\relax
\begin{eqnarray*}
\lefteqn{\mathit{sol} := [0, \,\mathit{Id}, \,{\displaystyle 
\frac {1}{2}} \,\mathit{Id} + {\displaystyle \frac {1}{2}} \,
\sqrt{1 + 4\,{\mathit{\_X}_{4}}^{2}}\,\mathit{e1} + {\mathit{\_X}
_{4}}\,\mathit{e1we2}, \,{\displaystyle \frac {1}{2}} \,\mathit{
Id} - {\displaystyle \frac {1}{2}} \,\sqrt{1 + 4\,{\mathit{\_X}_{
4}}^{2}}\,\mathit{e1} + {\mathit{\_X}_{4}}\,\mathit{e1we2}, } \\
 & & {\displaystyle \frac {1}{2}} \,\mathit{Id} + {\mathit{\_X}_{
2}}\,\mathit{e1} + {\displaystyle \frac {1}{2}} \,\sqrt{ - 4\,{
\mathit{\_X}_{2}}^{2} + 4\,{\mathit{\_X}_{4}}^{2} + 1}\,\mathit{
e2} + {\mathit{\_X}_{4}}\,\mathit{e1we2},  \\
 & & {\displaystyle \frac {1}{2}} \,\mathit{Id} + {\mathit{\_X}_{
2}}\,\mathit{e1} - {\displaystyle \frac {1}{2}} \,\sqrt{ - 4\,{
\mathit{\_X}_{2}}^{2} + 4\,{\mathit{\_X}_{4}}^{2} + 1}\,\mathit{
e2} + {\mathit{\_X}_{4}}\,\mathit{e1we2}]\mbox{\hspace{148pt}}
\end{eqnarray*}
\relax
\begin{Verbatim}
## re-compute the equation to check for correctness 
sol_square:=map(i->clicollect(simplify(cmul(i,i)-i)),sol);
\end{Verbatim}
\relax
\[
\mathit{sol\_square} := [0, \,0, \,0, \,0, \,0, \,0]
\]
\relax
All functions come with well developed help-pages. They can be reached by 
typing {\tt ?function} at the Maple commandline or searching the help
of Maple. A general help-page for the entire package and its sub-packages 
is available by typing {\tt ?Clifford[intro]}. A general introduction to 
Maple and its programming facilities may be found e.g. in \cite{wright:2002a}.

\section{The BIGEBRA Package}

This appendix provides only a very basic look-and-feel explanation of
the BIGEBRA package. The online documentation of BIGEBRA comes with over 
100 printed pages and should be consulted as reference. However, we felt 
it necessary to exhibit BIGEBRA's abilities here, since it was used to 
prove some statements in the text.

The BIGEBRA package (version 0.16) loads automatically the CLIFFORD package 
since the latter package is internally needed. We suppress the startup messages
by setting \verb|_SILENT| to $true$.
\begin{Verbatim}
restart:_CLIENV[_SILENT]:=true:with(Bigebra):
\end{Verbatim}
\begin{maplettyout}
Warning, new definition for drop_t
Warning, new definition for gco_d_monom
Warning, new definition for gco_monom
Warning, new definition for init
\end{maplettyout}
The particular functions of BIGEBRA are described below very shortly
to give an overview. For detailed help-pages and much more detailed
examples use the Maple online help by typing {\tt ?Bigebra,<function>}.

\subsection{{\tt \&cco} -- Clifford co-product}

The internal computation of the Clifford co-product is done by Rota-Stein
co-cliffordization as explained in the main text. The Clifford co-product has 
therefore to be initialized before the first usage, since it needs internally 
the Clifford co-product of the unit element, i.e. the `cap' tangle. 
Furthermore one needs also a co-scalar product which is stored in the 
matrix {\tt BI} (or left undefined), the dimension of the base space, defined
in {\tt dim\_V}, can range between 1 and 9. We have to set:
\begin{Verbatim}
dim_V:=2:
BI:=linalg[matrix](dim_V,dim_V,[a,b,c,d]);
\end{Verbatim}
\relax
\[
\mathit{BI} :=  \left[ 
{\begin{array}{cc}
a & b \\
c & d
\end{array}}
 \right] 
\]
\relax
\begin{Verbatim}
make_BI_Id():
&cco(e1);
\end{Verbatim}
\relax
\[
(\mathit{Id}\,\mathrm{\&t}\,\mathit{e1}) - b\,(\mathit{e1}\,
\mathrm{\&t}\,\mathit{e1we2}) - d\,(\mathit{e2}\,\mathrm{\&t}\,
\mathit{e1we2}) + (\mathit{e1}\,\mathrm{\&t}\,\mathit{Id}) + c\,(
\mathit{e1we2}\,\mathrm{\&t}\,\mathit{e1}) + d\,(\mathit{e1we2}\,
\mathrm{\&t}\,\mathit{e2})
\]
\relax
The most remarkable fact is that the Clifford co-product of the unit element
{\tt Id} is not \verb+&t(Id,Id)+ but
\begin{Verbatim}
&cco(Id);
\end{Verbatim}
\relax
\[
(\mathit{Id}\,\mathrm{\&t}\,\mathit{Id}) + a\,(\mathit{e1}\,
\mathrm{\&t}\,\mathit{e1}) + c\,(\mathit{e2}\,\mathrm{\&t}\,
\mathit{e1}) + b\,(\mathit{e1}\,\mathrm{\&t}\,\mathit{e2}) + d\,(
\mathit{e2}\,\mathrm{\&t}\,\mathit{e2}) + (c\,b - d\,a)\,(
\mathit{e1we2}\,\mathrm{\&t}\,\mathit{e1we2})
\]
\relax
The Clifford co-product is however co-associative.

\subsection{{\tt \&gco} -- Gra{\ss}mann co-product}

The Gra{\ss}mann co-product is the basic function of the BIGEBRA package, 
since the Clifford co-product is derived by the process of co-cliffordization.
It turns out that the Gra{\ss}mann co-product is a combinatorial function 
on the index set of Gra{\ss}mann multi-vectors, this is used in the package 
to get a fast evaluation of this function. The Gra{\ss}mann co-product is 
that of a connected and augmented co-algebra, which we called non-interacting 
Hopf  gebra in the main text.
\begin{Verbatim}
&gco(Id);    ## this is as expected
\end{Verbatim}
\relax
\[
\mathit{Id}\,\mathrm{\&t}\,\mathit{Id}
\]
\relax
\begin{Verbatim}
&gco(e1);
\end{Verbatim}
\relax
\[
(\mathit{Id}\,\mathrm{\&t}\,\mathit{e1}) + (\mathit{e1}\,\mathrm{
\&t}\,\mathit{Id})
\]
\relax
\begin{Verbatim}
&gco(e1we2); ## sum over splits
\end{Verbatim}
\relax
\[
(\mathit{Id}\,\mathrm{\&t}\,\mathit{e1we2}) + (\mathit{e1}\,
\mathrm{\&t}\,\mathit{e2}) - (\mathit{e2}\,\mathrm{\&t}\,\mathit{
e1}) + (\mathit{e1we2}\,\mathrm{\&t}\,\mathit{Id})
\]
\relax
Note that in the last case the sum is over all splits which are compatible
with the permutation symmetry of the factors. The signs are such that 
multiplying back gives {\em for each term\/} the original input. Hence we 
get two to the power of the grade of the element as a prefactor:
\begin{Verbatim}
eval(subs(`&t`=wedge,[op(%)]));
\end{Verbatim}
\relax
\[
[\mathit{e1we2}, \,\mathit{e1we2}, \,\mathit{e1we2}, \,\mathit{
e1we2}]
\]
\relax
\begin{Verbatim}
eval(`+`(op(%)));
\end{Verbatim}
\relax
\[
4\,\mathit{e1we2}
\]
\relax

\subsection{{\tt \&gco\_d} -- dotted Gra{\ss}mann co-product}

The dotted Gra{\ss}mann co-product is taken with respect to a different 
filtration of the Gra{\ss}mann algebra under consideration. This different
filtration is represented by the dotted wedge basis built w.r.t the dotted 
wedge product $\dot\wedge$. The dotted Gra{\ss}mann co-product is a wrapper 
function which translates the wedge basis elements into the dotted wedge 
basis ones, computes there the regular Gra{\ss}mann co-product and transforms
back the tensor product into the undotted basis. For examples see the online 
help of BIGEBRA.

\subsection{{\tt \&gpl\_co} -- Gra{\ss}mann Pl\"ucker co-product}

The Gra{\ss}mann-Pl\"ucker co-product evaluates the co-product w.r.t. 
the {\tt meet} (resp. {\tt \&v}) product of hyperplanes since it can be 
shown that the meet is an exterior product for hyperplanes. If we represent 
hyperplanes using Pl\"ucker coordinates, we can ask for a co-product on 
these Pl\"ucker coordinatized hyperplanes, which is in fact related to the 
wedge product of the points. For examples see the online help of BIGEBRA.

\subsection{{\tt \&map} -- maps products onto tensor slots}

The {\tt \&map} function extends product to be able to act on tensors. 
For instance one wants to wedge or Clifford multiply a tensor, say
{\tt \&t(e1,e2we3,e1we2)}, in two adjacent slots of the tensor. 
This is achieved as
\begin{Verbatim}
dim_V:=4:
&map(&t(e1,e2we3,e1we4),2,wedge);
\end{Verbatim}
\relax
\[
\mathit{e1}\,\mathrm{\&t}\,\mathit{e1we2we3we4}
\]
\relax
\begin{Verbatim}
&map(&t(e1,e2we3,e1we4),1,cmul);
\end{Verbatim}
\relax
\[
(\mathit{e1we2we3}\,\mathrm{\&t}\,\mathit{e1we4}) + {B_{1, \,2}}
\,(\mathit{e3}\,\mathrm{\&t}\,\mathit{e1we4}) - {B_{1, \,3}}\,(
\mathit{e2}\,\mathrm{\&t}\,\mathit{e1we4})
\]
\relax
Any $2\rightarrow 1$ mapping can be applied to tensors by this device. As
most of the BIGEBRA and CLIFFORD functions this is a multilinear mapping.

\subsection{{\tt \&t} -- tensor product}

The tensor product is a basic feature of the BIGEBRA package. The tensor 
product is an unevaluated product which is multilinear over any Maple 
expression which is not a CLIFFORD basis element. That is we are able to 
compute over Clifford modules. However, re-defining the Clifford type
{\tt type/cliscalar} one can change the behaviour. A few examples are
\begin{Verbatim}
&t(a*e1,3*e2+5*e3);
\end{Verbatim}
\relax
\[
3\,a\,(\mathit{e1}\,\mathrm{\&t}\,\mathit{e2}) + 5\,a\,(\mathit{
e1}\,\mathrm{\&t}\,\mathit{e3})
\]
\relax
\begin{Verbatim}
&t(e1,sin(x)*e2,e1we2);
\end{Verbatim}
\relax
\[
\mathrm{sin}(x)\,\mathrm{\&t}(\mathit{e1}, \,\mathit{e2}, \,
\mathit{e1we2})
\]
\relax
\begin{Verbatim}
&t(e1we2*z,-e3/z+t*e4,e2/t);
\end{Verbatim}
\relax
\[
 - {\displaystyle \frac {\mathrm{\&t}(\mathit{e1we2}, \,\mathit{
e3}, \,\mathit{e2})}{t}}  + z\,\mathrm{\&t}(\mathit{e1we2}, \,
\mathit{e4}, \,\mathit{e2})
\]
\relax
The tensor product allows studying decomposition and periodicity theorems.
One can handle multi-particle Clifford algebra, compute in different 
Clifford algebras, e.g. different bilinear or quadratic forms, and is 
able to investigate tangles of Gra{\ss}mann Hopf gebras and Clifford 
convolution algebras. A computation of a Gra{\ss}mann or Clifford antipode
would be impossible without this device. Moreover, also more geometric
notions as the {\tt meet} or {\tt \&v} (vee) product benefit from this 
structure.

\subsection{{\tt \&v} -- vee-product, i.e. meet}

The meet or vee-product computes the join of two extensors. It constitutes 
an exterior product on its own right, but on hyperplanes, not on points.
If hyperplanes are identified which the duals of points, which needs a
correlation and introduces a bilinear form, a complete dual approach to
the Gra{\ss}mann-Cayley algebra and its deformed structure the Clifford
convolution algebra is obtained. A few examples are:
\begin{Verbatim}
dim_V:=3:B:='B': ## unassign B
meet(e1we2,e2we3),&v(e1we2,e2we3); ## meet and &v are the same
\end{Verbatim}
\relax
\[
 - \mathit{e2}, \, - \mathit{e2}
\]
\relax
\begin{Verbatim}
&v(e1we2+e2we3,e2we3+e1we3);       ## acts on polynoms too
\end{Verbatim}
\relax
\[
 - \mathit{e1} - \mathit{e2} + \mathit{e3}
\]
\relax
Note that the meet introduces signs and it is the oriented meet of the
support of the extensor which describes the linear subspace. Of course 
a geometrical meaning of polynomial such objects is not obvious, but 
the meet nevertheless inherits linearity from its construction. The meet 
is calculated using the Peano {\tt bracket} and the co-product as
\[
{\rm meet}(x,y) = x_{(1)} [y,x_{(2)}] = [y_{(1)},x] y_{(2)}
\]
where the order of factors is important. The bracket can be understood 
in hopfish terms too.

\subsection{{\tt bracket} -- the Peano bracket}

The Peano bracket and Peano algebra was introduced by Rota et al. 
\cite{doubilet:rota:stein:1974a,barnabei:brini:rota:1985a} and called in 
the first paper Cayley algebra. However, Peano introduced the bracket as
a device to define Gra{\ss}mann's regressive product in dimension three,
see \cite{peano:1888a}. We showed in the main text that the Peano bracket
can be derived using a non-trivial integral of the Gra{\ss}mann Hopf gebra.
\[
[x,y] = h( x \wedge y) 
\]
where $h(x) : \bigwedge V \rightarrow \openk$ is a non-trivial integral.
In the case of the Gra{\ss}mann Hopf gebra this is the projection onto 
the highest grade element. BIGEBRA needs thus no bilinear form to define 
the bracket but only a maximal dimension. The {\tt bracket} function takes 
any number of arguments, wedges them together and projects onto the highest 
grade, e.g.
\begin{Verbatim}
dim_V:=3:
bracket(e1we2we3),bracket(e1,e2,e3);
\end{Verbatim}
\relax
\[
1, \,1
\]
\relax
\begin{Verbatim}
dim_V:=4:
bracket(e1we2,e2we3we4); ##   0 expected
\end{Verbatim}
\relax
\[
0
\]
\relax
\begin{Verbatim}
bracket(a*e1we2,b*e3we4);## a*b expected  
\end{Verbatim}
\relax
\[
a\,b
\]
\relax

\subsection{{\tt contract} -- contraction of tensor slots}

Given a tensor with at least two slots, contract allows to map a 
$2 \rightarrow 0$ mapping onto adjacent such slots. The tensor elements
can be seen as vectors or co-vectors, so we have in fact 4 types of 
contractions.
\begin{Verbatim}
contract(&t(e1,e1,e2),1,EV);         ## evaluation on slots 1,2
\end{Verbatim}
\relax
\[
\mathrm{\&t}(\mathit{e2})
\]
\relax
\begin{Verbatim}
contract(&t(e1,e1we2,e3we4),2,bracket); ## bracket on slots 2,3
\end{Verbatim}
\relax
\[
\mathrm{\&t}(\mathit{e1})
\]
\relax

\subsection{{\tt define} -- Maple define, patched}

The {\tt define} facility of Maple turned out to be not very useful for 
defining multilinear associative functions. It showed up to compute wrong
results and was not designed to handle an arbitrary base ring. BIGEBRA
patches define so that {\tt type/cliscalar} is used for scalars and that 
any function defined with define like {\tt define(`\&r`,flat,multilinear)}
to be associative, i.e. flat and multilinear. For further information
see the online help-page of BIGEBRA.

\subsection{{\tt drop\_t} -- drops tensor signs}

This is a helper function to drop the tensor sign {\tt \&t} from Clifford
expressions, i.e. tensors of rank one. For technical reasons the tensor 
sign is not automatically dropped.
\begin{Verbatim}
drop_t(&t(a*e1+b*e1we2));
\end{Verbatim}
\relax
\[
a\,\mathit{e1} + b\,\mathit{e1we2}
\]
\relax

\subsection{{\tt EV} -- evaluation map}

The evaluation map is given by the action of co-vectors on vectors acting 
in the natural way. If a canonical co-basis $\theta^a$ is defined, one
finds $\theta^a(e_b) = \delta^a_b$ where $\delta$ is the Kronecker symbol.
The user has to take care in which tensor slot the co-vectors reside, since
they are, unfortunately, displayed by the same basis symbols $eiwej$ etc.
The evaluation map acts on any multivector polynom in $\bigwedge V$.
\begin{Verbatim}
EV(e1,a*Id+b*e1+c*e2+d*e1we2);   ## b expected
\end{Verbatim}
\relax
\[
b
\]
\relax
\begin{Verbatim}
EV(e1we2,e1we2),EV(e1we2,e2we3); ## 1,0 expected
\end{Verbatim}
\relax
\[
1, \,0
\]
\relax

\subsection{{\tt gantipode} -- Gra{\ss}mann antipode}

The Gra{\ss}mann antipode is the antipode of the Gra{\ss}mann Hopf gebra.
The most remarkable fact is that this antipode map is equivalent to the 
main involution of a Clifford algebra of the same space or the main 
involution of the Gra{\ss}mann algebra
\begin{Verbatim}
dim_V:=3:
bas:=cbasis(dim_V);
\end{Verbatim}
\relax
\[
\mathit{bas} := [\mathit{Id}, \,\mathit{e1}, \,\mathit{e2}, \,
\mathit{e3}, \,\mathit{e1we2}, \,\mathit{e1we3}, \,\mathit{e2we3}
, \,\mathit{e1we2we3}]
\]
\relax
\begin{Verbatim}
map(gantipode,bas);
\end{Verbatim}
\relax
\[
[\mathit{Id}, \, - \mathit{e1}, \, - \mathit{e2}, \, - \mathit{e3
}, \,\mathit{e1we2}, \,\mathit{e1we3}, \,\mathit{e2we3}, \, - 
\mathit{e1we2we3}]
\]
\relax
\begin{Verbatim}
map(gradeinv,bas);
\end{Verbatim}
\relax
\[
[\mathit{Id}, \, - \mathit{e1}, \, - \mathit{e2}, \, - \mathit{e3
}, \,\mathit{e1we2}, \,\mathit{e1we3}, \,\mathit{e2we3}, \, - 
\mathit{e1we2we3}]
\]
\relax

\subsection{{\tt gco\_unit} -- Gra{\ss}mann co-unit}

Since the co-gebra structure is obtained by categorical duality, the 
Gra{\ss}mann co-gebra possesses a co-unit. This can be exemplified as 
follows:
\begin{Verbatim}
_X:=add(X[i]*bas[i],i=1..2^dim_V);    ## arbitrary element
\end{Verbatim}
\relax
\[
\mathit{\_X} := {X_{1}}\,\mathit{Id} + {X_{2}}\,\mathit{e1} + {X
_{3}}\,\mathit{e2} + {X_{4}}\,\mathit{e3} + {X_{5}}\,\mathit{
e1we2} + {X_{6}}\,\mathit{e1we3} + {X_{7}}\,\mathit{e2we3} + {X_{
8}}\,\mathit{e1we2we3}
\]
\relax
\begin{Verbatim}
simplify(drop_t(gco_unit(&gco(_X),1)) - _X); ## 0 expected
\end{Verbatim}
\relax
\[
0
\]
\relax
\begin{Verbatim}
simplify(drop_t(gco_unit(&gco(_X),2)) - _X); ## 0 expected
\end{Verbatim}
\relax
\[
0
\]
\relax

\subsection{{\tt gswitch} -- graded (i.e. Gra{\ss}mann) switch}

The graded switch is the natural switch of the Gra{\ss}mann Hopf gebra. 
It is not the generic switch of a Clifford algebra if the bilinear form 
is not identical zero. The graded switch swaps two adjacent factors
of a tensor and counts the minus signs arising from the reordering 
of the factors.
\begin{Verbatim}
gswitch(&t(e1,e2,e3we4),1); ## - expected
\end{Verbatim}
\relax
\[
 - \mathrm{\&t}(\mathit{e2}, \,\mathit{e1}, \,\mathit{e3we4})
\]
\relax
\begin{Verbatim}
gswitch(&t(e1,e2,e3we4),2); ## + expected
\end{Verbatim}
\relax
\[
\mathrm{\&t}(\mathit{e1}, \,\mathit{e3we4}, \,\mathit{e2})
\]
\relax

\subsection{{\tt help} -- main help-page of BIGEBRA package}

This is not a function of the package, but the main help-page of the 
BIGEBRA package. It can be accessed in a Maple session by typing 
{\tt ?Bigebra,help}. The main help-page gives an alphabetic listing
of BIGEBRA functions, links it to CLIFFORD and provides some literature
from which place some of the algorithms and mathematics have been taken.
The reader is urged to look up this page.

\subsection{{\tt init} -- init procedure}

BIGEBRA needs a tricky init procedure to patch load the package and patch
the Maple {\tt define} function. Init loads BIGEBRA, then the tensor product
{\tt \&t} is defined which loads the define code into the session. Then 
BIGEBRA is loaded a second time to overwrite in the memory the unsuited
parts of define. {\tt Init} loads CLIFFORD, i.e. Cliff5, if it was not 
already loaded.

\subsection{{\tt linop/linop2} -- action of a linear operator on a Clifford polynom}

Since we have been interested in tangle equations like the definition of 
the antipode. The action of certain operators on a tensor slot is therefore
necessary. Sometimes it is useful to have matrix representations of
such operators and {\tt linop} provides this facility. {\tt linop2}
is the same function which acts however on two adjacent tensor slots, 
hence we have
\[
{\rm linop} \in \End \bigwedge V
\]
\[
{\rm linop2} \in \End \bigwedge V \otimes \bigwedge V  
\]

\subsection{{\tt make\_BI\_Id} -- cup tangle need for {\tt \&cco}}

This function computes the cap tangle for a certain co-scalar product
either unassigned or defined as a matrix named {\tt BI}. See either {\tt
\&cco} above or the online help-page of BIGEBRA.

\subsection{{\tt mapop/mapop2} -- action of an operator on a tensor slot}

While {\tt linop(2)} defines a linear operator as an endomorphism on 
$\bigwedge V$ seen as linear space. The function {\tt mapop(2)}
allows to apply these operators to any tensor slot of a tensor or to
any two adjacent tensor slots. For some example and the usage see
the help-page of BIGEBRA.

\subsection{{\tt meet} -- same as {\tt \&v} (vee-product)}

The meet is a synonym for the {\tt \&v} (vee-) product. However, in the
BIGEBRA package the meet and vee-products are computed differently, we
have
\[
meet(x,y) = x_{(1)} [y , x_{(2)}]
\]
while 
\[
\&v(x,y) = [y_{(1)} , x ] y_{(2)}]
\]
This allows to check that both definitions are equivalent. This 
computation can be found, together with many geometric applications
and some benchmarks in the online help-page for the {\tt meet}
in the BIGEBRA package.

\subsection{{\tt pairing} -- A pairing w.r.t. a bilinear form}

The pairing is a decorated cup tangle, where the decoration describes
the bilinear form used to convert one element into a co-vector, i.e. 
a scalar product. The pairing s graded and can be defined as follows
\[
\langle x \mid y \rangle = \left\{
\begin{array}{cl}
\pm{\rm det}(\langle x_i \mid y_j \rangle) & \text{if grade $x$ = grade $y$}\\
0 & \text{otherwise}
\end{array}\right.
\]
where $x,y$ are extensors of $\bigwedge V$ and the pairing is extended
by bilinearity. For explicite examples see the online help-page of the
BIGEBRA package.

\subsection{{\tt peek} -- extract a tensor slot}

This is a technical function used mostly internally to be able to access
certain tensor slots. For explicite examples and the correct syntax see 
the online help-page of the BIGEBRA package.

\subsection{{\tt poke} -- insert a tensor slot}

This is a technical function used mostly internally to be able to insert
Clifford elements as new tensor slots in an arbitrary tensor polynomial. 
For explicite examples and the correct syntax see the online help-page of 
the BIGEBRA package.

\subsection{{\tt remove\_eq} -- removes tautological equations}

This is a technical function used mostly internally. It drops tautological 
equations in a set of equations. For explicite examples see the online 
help-page of the BIGEBRA package. 

\subsection{{\tt switch} -- ungraded switch}

The switch simply swaps adjacent tensor slots, no sign is computed.
\begin{Verbatim}
switch(&t(e1,e2,e3we4),1);
\end{Verbatim}
\relax
\[
\mathrm{\&t}(\mathit{e2}, \,\mathit{e1}, \,\mathit{e3we4})
\]
\relax
\begin{Verbatim}
switch(&t(e1,e2,e3we4),2);
\end{Verbatim}
\relax
\[
\mathrm{\&t}(\mathit{e1}, \,\mathit{e3we4}, \,\mathit{e2})
\]
\relax

\subsection{{\tt tcollect} -- collects w.r.t. the tensor basis}

This is a function which is needed to customise the output of some BIGEBRA 
functions for inputting it into other such functions. Furthermore it allows
a better comparison of tensor polynomials. For explicite examples see the 
online help-page of the BIGEBRA package. 

\subsection{{\tt tsolve1} -- tangle solver}

The tangle solver is an extension of the CLIFFORD function {\tt clisolve}.
It allows to solve for {\em endomorphisms} acting in $n \rightarrow 1$
tangles, therefore the name. Most of the axioms and definitions of 
Gra{\ss}mann Hopf gebras and Clifford bi-convolution algebras are of this 
type. The online help-page for {\tt tsolve1} comes up with explicite
computations of the unit for Gra{\ss}mann convolution, the Gra{\ss}mann
antipode and some facts about integrals in Gra{\ss}mann and Clifford 
bi-convolutions. For explicite examples see the online help-page of the 
BIGEBRA package. 

\subsection{{\tt VERSION} -- shows the version of the package}

This command is issued as {\tt VERSION();} and returns some information
about the release of the BIGEBRA package.

\subsection{{\tt type/tensorbasmonom} -- new Maple type}

To be able to facilitate symbolic computations Maple provides a type
checking system. BIGEBRA as CLIFFORD use this device and define some 
new types extending this mechanism. A tensorbasmonom is any expression
which is an extensor without any prefactor, e.g.
\begin{Verbatim}
type(&t(e1,e2,e3),tensorbasmonom);   ## true  expected
\end{Verbatim}
\relax
\[
\mathit{true}
\]
\relax
\begin{Verbatim}
type(a*&t(e1,e2),tensorbasmonom);    ## false expected
\end{Verbatim}
\relax
\[
\mathit{false}
\]
\relax
\begin{Verbatim}
type(&t(e1)+&t(e2),tensorbasmonom);  ## false expected
\end{Verbatim}
\relax
\[
\mathit{false}
\]
\relax
\begin{Verbatim}
type(a*sin(x)*e1we3,tensorbasmonom); ## false expected
\end{Verbatim}
\relax
\[
\mathit{false}
\]
\relax

\subsection{{\tt type/tensormonom} -- new Maple type}

A tensormonom is a tensorbasmonom possibly having a prefactor from the
ring the tensor product is built over. This type is inclusive in that
way that a tensorbasmonom is also considered to be a tensormonom.
\begin{Verbatim}
type(&t(e1,e2,e3),tensormonom);   ## true  expected
\end{Verbatim}
\relax
\[
\mathit{true}
\]
\relax
\begin{Verbatim}
type(a*&t(e1,e2),tensormonom);    ## true  expected
\end{Verbatim}
\relax
\[
\mathit{true}
\]
\relax
\begin{Verbatim}
type(&t(e1)+&t(e2),tensormonom);  ## false expected
\end{Verbatim}
\relax
\[
\mathit{false}
\]
\relax
\begin{Verbatim}
type(a*sin(x)*e1we3,tensormonom); ## false expected
\end{Verbatim}
\relax
\[
\mathit{false}
\]
\relax

\subsection{{\tt type/tensorpolynom} -- new Maple type}

A tensor polynom is a sum of tensormonoms. This type is also inclusive. 
\begin{Verbatim}
type(&t(e1,e2,e3),tensorpolynom);   ## true  expected
\end{Verbatim}
\relax
\[
\mathit{true}
\]
\relax
\begin{Verbatim}
type(a*&t(e1,e2),tensorpolynom);    ## true  expected
\end{Verbatim}
\relax
\[
\mathit{true}
\]
\relax
\begin{Verbatim}
type(&t(e1)+&t(e2),tensorpolynom);  ## true  expected
\end{Verbatim}
\relax
\[
\mathit{true}
\]
\relax
\begin{Verbatim}
type(a*sin(x)*e1we3,tensorpolynom); ## false expected
\end{Verbatim}
\relax
\[
\mathit{false}
\]
\relax